\newcount\mgnf\newcount\tipi\newcount\tipoformule\newcount\greco

\tipi=2          
\tipoformule=0   


\global\newcount\numsec
\global\newcount\numfor
\global\newcount\numtheo
\global\advance\numtheo by 1

\def\senondefinito#1{\expandafter\ifx\csname#1\endcsname\relax}

\def\SIA #1,#2,#3 {\senondefinito{#1#2}%
\expandafter\xdef\csname #1#2\endcsname{#3}\else
\write16{???? ma #1,#2 e' gia' stato definito !!!!} \fi}

\def\etichetta(#1){(\veroparagrafo.\veraformula)%
\SIA e,#1,(\veroparagrafo.\veraformula) %
\global\advance\numfor by 1%
\write15{\string\FU (#1){\equ(#1)}}%
\write16{ EQ #1 ==> \equ(#1) }}

\def\letichetta(#1){\veroparagrafo.\verotheo
\SIA e,#1,{\veroparagrafo.\verotheo}
\global\advance\numtheo by 1
\write15{\string\FU (#1){\equ(#1)}}
\write16{ Sta \equ(#1) == #1 }}

\def\tetichetta(#1){\veroparagrafo.\veraformula 
\SIA e,#1,{(\veroparagrafo.\veraformula)}
\global\advance\numfor by 1
\write15{\string\FU (#1){\equ(#1)}}
\write16{ tag #1 ==> \equ(#1)}}

\def\FU(#1)#2{\SIA fu,#1,#2 }

\def\etichettaa(#1){(A\veroparagrafo.\veraformula)%
\SIA e,#1,(A\veroparagrafo.\veraformula) %
\global\advance\numfor by 1%
\write15{\string\FU (#1){\equ(#1)}}%
\write16{ EQ #1 ==> \equ(#1) }}

\def\BOZZA{
\def\alato(##1){%
 {\rlap{\kern-\hsize\kern-1.4truecm{$\scriptstyle##1$}}}}%
\def\aolado(##1){%
 {
{
 \rlap{\kern-1.4truecm{$\scriptstyle##1$}}}}}
}

\def\alato(#1){}
\def\aolado(#1){}

\def\veroparagrafo{\number\numsec}
\def\veraformula{\number\numfor}
\def\verotheo{\number\numtheo}

\def\Eq(#1){\eqno{\etichetta(#1)\alato(#1)}}
\def\eq(#1){\etichetta(#1)\alato(#1)}
\def\leq(#1){\leqno{\aolado(#1)\etichetta(#1)}}
\def\teq(#1){\tag{\aolado(#1)\tetichetta(#1)\alato(#1)}}
\def\Eqa(#1){\eqno{\etichettaa(#1)\alato(#1)}}
\def\eqa(#1){\etichettaa(#1)\alato(#1)}
\def\eqv(#1){\senondefinito{fu#1}$\clubsuit$#1
\write16{#1 non e' (ancora) definito}%
\else\csname fu#1\endcsname\fi}
\def\equ(#1){\senondefinito{e#1}\eqv(#1)\else\csname e#1\endcsname\fi}

\def\Lemma(#1){\aolado(#1)Lemma \letichetta(#1).$\;$}%
\def\Theorem(#1){{\aolado(#1)Theorem \letichetta(#1).$\;$}}%
\def\Proposition(#1){\aolado(#1){Proposition \letichetta(#1).$\;$}}%
\def\Corollary(#1){{\aolado(#1)Corollary \letichetta(#1).$\;$}}%
\def\Remark(#1){{\noindent\aolado(#1){\bf Remark \letichetta(#1).$\;$}}}%
\def\Definition(#1){{\noindent\aolado(#1){\bf Definition
\letichetta(#1).$\;$}}}
\def\Example(#1){\aolado(#1) Example \letichetta(#1)$\!\!$\hskip-1.6truemm}

\def\include#1{
\openin13=#1.aux \ifeof13 \relax \else
\input #1.aux \closein13 \fi}

\openin14=\jobname.aux \ifeof14 \relax \else
\input \jobname.aux \closein14 \fi
\openout15=\jobname.aux


{\count255=\time\divide\count255 by 60 \xdef\hourmin{\number\count255}
        \multiply\count255 by-60\advance\count255 by\time
   \xdef\hourmin{\hourmin:\ifnum\count255<10 0\fi\the\count255}}

\def\oramin{\hourmin }

\def\data{\number\day/\ifcase\month\or january \or february \or march \or
april \or may \or june \or july \or august \or september
\or october \or november \or december \fi/\number\year;\ \oramin}

\newcount\pgn \pgn=1
\def\foglio{\number\numsec:\number\pgn
\global\advance\pgn by 1}
\def\foglioa{A\number\numsec:\number\pgn
\global\advance\pgn by 1}

\footline={\rlap{\hbox{\copy200}}\hss\tenrm\folio\hss}

\def\TIPIO{
\font\setterm=amr7 
\def \settepunti{\def\rm{\fam0\setterm}
\textfont0=\setterm   
\normalbaselineskip=9pt\normalbaselines\rm
}\let\nota=\settepunti}

\def\TIPITOT{
\font\twelverm=cmr12
\font\twelvei=cmmi12
\font\twelvesy=cmsy10 scaled\magstep1
\font\twelveex=cmex10 scaled\magstep1
\font\twelveit=cmti12
\font\twelvett=cmtt12
\font\twelvebf=cmbx12
\font\twelvesl=cmsl12
\font\ninerm=cmr9
\font\ninesy=cmsy9
\font\eightrm=cmr8
\font\eighti=cmmi8
\font\eightsy=cmsy8
\font\eightbf=cmbx8
\font\eighttt=cmtt8
\font\eightsl=cmsl8
\font\eightit=cmti8
\font\sixrm=cmr6
\font\sixbf=cmbx6
\font\sixi=cmmi6
\font\sixsy=cmsy6
\font\twelvetruecmr=cmr10 scaled\magstep1
\font\twelvetruecmsy=cmsy10 scaled\magstep1
\font\tentruecmr=cmr10
\font\tentruecmsy=cmsy10
\font\eighttruecmr=cmr8
\font\eighttruecmsy=cmsy8
\font\seventruecmr=cmr7
\font\seventruecmsy=cmsy7
\font\sixtruecmr=cmr6
\font\sixtruecmsy=cmsy6
\font\fivetruecmr=cmr5
\font\fivetruecmsy=cmsy5
\textfont\truecmr=\tentruecmr
\scriptfont\truecmr=\seventruecmr
\scriptscriptfont\truecmr=\fivetruecmr
\textfont\truecmsy=\tentruecmsy
\scriptfont\truecmsy=\seventruecmsy
\scriptscriptfont\truecmr=\fivetruecmr
\scriptscriptfont\truecmsy=\fivetruecmsy
\def \eightpoint{\def\rm{\fam0\eightrm}
\textfont0=\eightrm \scriptfont0=\sixrm \scriptscriptfont0=\fiverm
\textfont1=\eighti \scriptfont1=\sixi   \scriptscriptfont1=\fivei
\textfont2=\eightsy \scriptfont2=\sixsy   \scriptscriptfont2=\fivesy
\textfont3=\tenex \scriptfont3=\tenex   \scriptscriptfont3=\tenex
\textfont\itfam=\eightit  \def\it{\fam\itfam\eightit}%
\textfont\slfam=\eightsl  \def\sl{\fam\slfam\eightsl}%
\textfont\ttfam=\eighttt  \def\tt{\fam\ttfam\eighttt}%
\textfont\bffam=\eightbf  \scriptfont\bffam=\sixbf
\scriptscriptfont\bffam=\fivebf  \def\bf{\fam\bffam\eightbf}%
\tt \ttglue=.5em plus.25em minus.15em
\setbox\strutbox=\hbox{\vrule height7pt depth2pt width0pt}%
\normalbaselineskip=9pt
\let\sc=\sixrm  \let\big=\eightbig  \normalbaselines\rm
\textfont\truecmr=\eighttruecmr
\scriptfont\truecmr=\sixtruecmr
\scriptscriptfont\truecmr=\fivetruecmr
\textfont\truecmsy=\eighttruecmsy
\scriptfont\truecmsy=\sixtruecmsy
}\let\nota=\eightpoint}

\newfam\msbfam   
\newfam\truecmr  
\newfam\truecmsy 
\newskip\ttglue
\ifnum\tipi=0\TIPIO \else\ifnum\tipi=1 \TIPI\else \TIPITOT\fi\fi

\def\a{\alpha}
\def\b{\beta}
\def\d{\delta}
\def\e{\epsilon}

\def\vf{\varphi}
\def\g{\gamma}
\def\k{\kappa}
\def\l{\lambda}
\def\r{\rho}
\def\s{\sigma}
\def\t{\tau}
\def\th{\theta}

\def\z{\zeta}
\def\o{\omega}
\def\D{\Delta}
\def\L{\Lambda}
\def\G{\Gamma}
\def\O{\Omega}
\def\S{\Sigma}

\def\del #1{\frac{\partial^{#1}}{\partial\l^{#1}}}

\def\E{{I\kern-.25em{E}}}
\def\N{{I\kern-.25em{N}}}
\def\M{{I\kern-.25em{M}}}
\def\R{{I\kern-.25em{R}}}
\def\Z{{Z\kern-.425em{Z}}}
\def\eu{{1\kern-.25em\hbox{\sm I}}}
\def \1{{1\kern-.25em{\rm I}}}
\def\C{{I\kern-.64em{C}}}
\def\P{{I\kern-.25em{P}}}
\def\eop{{ \vrule height7pt width7pt depth0pt}\par\bigskip}

\def\del{\partial}


\def\AA{{\cal A}}
\def\BB{{\cal B}}
\def\CC{{\cal C}}
\def\DD{{\cal D}}
\def\EE{{\cal E}}
\def\FF{{\cal F}}
\def\GG{{\cal G}}

\def\II{{\cal I}}
\def\JJ{{\cal J}}

\def\LL{{\cal L}}
\def\OO{{\cal O}}
\def\PP{{\cal P}}
\def\SS{{\cal S}}
\def\TT{{\cal T}}
\def\NN{{\cal N}}
\def\MM{{\cal M}}
\def\WW{{\cal W}}
\def\VV{{\cal V}}

\def\LL{{\cal L}}

\def\RR{{\cal R}}
\def\QQ{{\cal Q}}
\def\YY{{\cal Y}}

\def\chap #1#2{\line{\ch #1\hfill}\numsec=#2\numfor=1}

\def\wh{\widehat}

\def\sqr#1#2{{\vcenter{\vbox{\hrule height.#2pt
     \hbox{\vrule width.#2pt height#1pt \kern#1pt
   \vrule width.#2pt}\hrule height.#2pt}}}}
\def\qed{ $\mathchoice\sqr64\sqr64\sqr{2.1}3\sqr{1.5}3$}


\newcount\foot
\foot=1
\def\note#1{\footnote{${}^{\number\foot}$}{\ftn #1}\advance\foot by 1}
\def\tag #1{\eqno{\hbox{\rm(#1)}}}
\def\frac#1#2{{#1\over #2}}
\def\sfrac#1#2{{\textstyle{#1\over #2}}}
\def\text#1{\quad{\hbox{#1}}\quad}
\def\newpage{\vfill\eject}

\def\proof{{\noindent\pr Proof: }}

\def\remark{\noindent{\bf Remark: }}
\def\thanks{\noindent{\bf Aknowledgements: }}
\font\pr=cmbxsl10


\font\ch=cmbx12
\font\ftn=cmr8

\font\it=cmti10
\font\bf=cmbx10
\font\sm=cmr7

%
\catcode`\X=12\catcode`\@=11
\def\n@wcount{\alloc@0\count\countdef\insc@unt}
\def\n@wwrite{\alloc@7\write\chardef\sixt@@n}
\def\n@wread{\alloc@6\read\chardef\sixt@@n}
\def\crossrefs#1{\ifx\alltgs#1\let\tr@ce=\alltgs\else\def\tr@ce{#1,}\fi
   \n@wwrite\cit@tionsout\openout\cit@tionsout=\jobname.cit
   \write\cit@tionsout{\tr@ce}\expandafter\setfl@gs\tr@ce,}
\def\setfl@gs#1,{\def\@{#1}\ifx\@\empty\let\next=\relax
   \else\let\next=\setfl@gs\expandafter\xdef
   \csname#1tr@cetrue\endcsname{}\fi\next}
\newcount\sectno\sectno=0\newcount\subsectno\subsectno=0\def\r@s@t{\relax}
\def\resetall{\global\advance\sectno by 1\subsectno=0
  \gdef\firstpart{\number\sectno}\r@s@t}
\def\resetsub{\global\advance\subsectno by 1
   \gdef\firstpart{\number\sectno.\number\subsectno}\r@s@t}
\def\v@idline{\par}\def\firstpart{\number\sectno}
\def\l@c@l#1X{\firstpart.#1}\def\gl@b@l#1X{#1}\def\t@d@l#1X{{}}
\def\m@ketag#1#2{\expandafter\n@wcount\csname#2tagno\endcsname
     \csname#2tagno\endcsname=0\let\tail=\alltgs\xdef\alltgs{\tail#2,}%
  \ifx#1\l@c@l\let\tail=\r@s@t\xdef\r@s@t{\csname#2tagno\endcsname=0\tail}\fi
   \expandafter\gdef\csname#2cite\endcsname##1{\expandafter
     \ifx\csname#2tag##1\endcsname\relax?\else{\rm\csname#2tag##1\endcsname}\fi
    \expandafter\ifx\csname#2tr@cetrue\endcsname\relax\else
     \write\cit@tionsout{#2tag ##1 cited on page \folio.}\fi}%
   \expandafter\gdef\csname#2page\endcsname##1{\expandafter
     \ifx\csname#2page##1\endcsname\relax?\else\csname#2page##1\endcsname\fi
     \expandafter\ifx\csname#2tr@cetrue\endcsname\relax\else
     \write\cit@tionsout{#2tag ##1 cited on page \folio.}\fi}%
   \expandafter\gdef\csname#2tag\endcsname##1{\global\advance
     \csname#2tagno\endcsname by 1%
   \expandafter\ifx\csname#2check##1\endcsname\relax\else%
\fi
   \expandafter\xdef\csname#2check##1\endcsname{}%
   \expandafter\xdef\csname#2tag##1\endcsname
     {#1\number\csname#2tagno\endcsnameX}%
   \write\t@gsout{#2tag ##1 assigned number \csname#2tag##1\endcsname\space
      on page \number\count0.}%
   \csname#2tag##1\endcsname}}%
\def\m@kecs #1tag #2 assigned number #3 on page #4.%
   {\expandafter\gdef\csname#1tag#2\endcsname{#3}
   \expandafter\gdef\csname#1page#2\endcsname{#4}}
\def\re@der{\ifeof\t@gsin\let\next=\relax\else
    \read\t@gsin to\t@gline\ifx\t@gline\v@idline\else
    \expandafter\m@kecs \t@gline\fi\let \next=\re@der\fi\next}
\def\t@gs#1{\def\alltgs{}\m@ketag#1e\m@ketag#1s\m@ketag\t@d@l p
    \m@ketag\gl@b@l r \n@wread\t@gsin\openin\t@gsin=\jobname.tgs \re@der
    \closein\t@gsin\n@wwrite\t@gsout\openout\t@gsout=\jobname.tgs }
\outer\def\localtags{\t@gs\l@c@l}
\outer\def\globaltags{\t@gs\gl@b@l}
\outer\def\newlocaltag#1{\m@ketag\l@c@l{#1}}
\outer\def\newglobaltag#1{\m@ketag\gl@b@l{#1}}

\def\t@gsoff#1,{\def\@{#1}\ifx\@\empty\let\next=\relax\else\let\next=\t@gsoff
   \expandafter\gdef\csname#1cite\endcsname{\relax}
   \expandafter\gdef\csname#1page\endcsname##1{?}
   \expandafter\gdef\csname#1tag\endcsname{\relax}\fi\next}
\def\verbatimtags{\let\ift@gs=\iffalse\ifx\alltgs\relax\else
   \expandafter\t@gsoff\alltgs,\fi}
\catcode`\X=11 \catcode`\@=\active
\localtags
%
\setbox200\hbox{$\scriptscriptstyle \data $}
\global\newcount\numpunt
\magnification=1000
\hoffset=0.cm
\baselineskip=14pt
\parindent=12pt
\lineskip=4pt\lineskiplimit=0.1pt
\parskip=0.1pt plus1pt

\hyphenation{small}

\catcode`\@=11

 \centerline {\bf { One-dimensional  random field Kac's model:} }
\vskip.3truecm
\centerline {\bf  localization of the phases \footnote{$^*$}
   {\eightrm Supported by: CNR-CNRS-Project 8.005,
INFM-Roma;  MURST/Cofin 01-02/03-04; FAPERJ Projects E-26/150.940-99 and
E-26/151.905/00; CNPq-CNR Project: 91.0069/00-0}}
\vskip1truecm
 \centerline{\bf Marzio Cassandro} 
\centerline {Dipartimento di Fisica,
Universit\`a di Roma ``La Sapienza'',} 
\centerline {INFM-Sez. di Roma. P.le A.
Moro, 00185 Roma Italy. }
\centerline {cassandro@roma1.infn.it}
\bigskip
\centerline{ \bf Enza Orlandi}
\centerline{ Dipartimento di Matematica,
Universit\`a  di Roma Tre,}
\centerline{ L.go S.Murialdo 1, 00156 Roma, Italy.}
\centerline{ orlandi@mat.uniroma3.it} 
\bigskip
\centerline{\bf Pierre Picco }
\centerline{ CPT--CNRS, Luminy, Case 907, F-13288
Marseille Cedex 9, France}
\centerline{ and CMM-CNRS Blanco Encalada 2120, Santiago,
Chile.}
\centerline{  picco@cpt.univ-mrs.fr}
\bigskip
\centerline{ \bf Maria Eulalia Vares\footnote{$^{**}$}
{\eightrm Partially supported by CNPq.}
}
\centerline{ CBPF, Rua Dr. Xavier Sigaud, 150.}
\centerline{ 22290-180, Rio de Janeiro, RJ,Brasil.}
\centerline{ eulalia@cbpf.br.}

\vskip1truecm

  { \bf Abstract}
We study the typical profiles of a one dimensional random field
Kac model, 
for values of the temperature and magnitude  of the field 
in the region of the two  absolute minima for the free
energy of the corresponding random field  Curie Weiss model.
We show that, for a set of realizations of the random field of 
overwhelming  probability,  the localization of the two phases
corresponding to the previous minima is completely determined.
Namely, we are able to construct random intervals  tagged with a 
sign, where typically,
with respect to the infinite volume Gibbs measure, the profile is
rigid and takes, according to the sign,  one of the two values 
corresponding to the  previous  minima.
Moreover,  we characterize the transition from one phase to the other.  
The analysis  extends the one done by Cassandro, Orlandi and Picco in
[\rcite{COP1}].
\vskip.5cm
\noindent {\bf Key Words and Phrases}: phase transition,
random walk, random environment, Kac potential.\hfill\break
\vskip.5cm

\newpage
\bigskip \bigskip
\chap {1 Introduction}1
\numsec= 1 \numfor= 1 \numtheo=1

\medskip

We consider a  one-dimensional spin system interacting via a ferromagnetic 
two-body Kac
potential  and  external random field given by independent
Bernoulli variables. Problems where  a stochastic contribution is added to 
the energy of the
system arise naturally in condensed matter physics where the presence of the 
impurities
causes the microscopic structure to vary from point to point. 
Some of the vast literature on these
topics may been seen consulting  [\rcite {Ah}-\rcite {BRZ}], [\rcite {BK}], [\rcite
{FFS}-\rcite {IM}], [\rcite {Ku}], [\rcite {SW}].

 Kac's potentials
is a short way to  denote two-body ferromagnetic interactions with
range $\frac 1\g $, where $\g$ is a dimensionless parameter such
that when $\g \to 0$, i.e. very long range, the strength of the
interaction becomes very weak, but in such a way that the total
interaction between one spin and all the others is finite.  They
were  introduced in [\rcite {KUH}], and then generalized in
[\rcite {LP}], to provide  a rigorous proof of the 
validity of the van der Waals
theory of a {liquid--vapor} phase transition.  Performing first the
thermodynamic limit of the spin system interacting via Kac's
potential, and then the limit of infinite range, $\g \to 0$, they
rigorously derived the Maxwell rule. This implies that   the free
energy  of the system is the convex envelope of the corresponding
free energy for the Curie-Weiss model. 
This leads to two spatially homogeneous phases, corresponding to the
two points of minima of  the free energy of the  
Curie-Weiss model.
Often we will call $+$ phase the one associated to the positive
minimizer, and $-$ phase the one associated to the negative
minimizer.  For $\g$ fixed and different from zero, there are
several papers trying to understand qualitatively and
quantitatively the features of systems with long, but finite range
interaction. (See for instance [\rcite {CP}], [\rcite {LMP}], 
[\rcite {BZ}], [\rcite {FP}].)
In the one dimensional case, the analysis [\rcite {COP}] for Ising
spin and [\rcite {Bo}]  for more general spin, gives  a
satisfactory description of the typical profiles.

Similar type of analysis holds for Ising spin systems interacting
via a Kac potential and external random field. In this paper,
extending the analysis done in [\rcite {COP1}], 
we   study,  for $\g$ small but different from
zero, in one dimension,  the typical profiles of the system
for all the values of the temperature and magnitude  of the field in the
region of two  absolute minima for the free energy of the
corresponding random field Curie Weiss model, whose behavior is closely
connected with the local behavior of the random field Kac model.   
Through a
block-spin transformation, the microscopic system is mapped into a
system on $L^{\infty}(\R) \times L^{\infty}(\R)$, for which the
length of interaction becomes of order one (the macroscopic
system). It has been proven in [\rcite {COP1}] that if  the system is 
considered on an interval of
length $  \frac 1 \gamma (\log \frac 1 \gamma)^p$, $p\ge 2$, then 
for intervals
whose length in macroscopic scale is of order $ \frac {1} {\gamma
\log\log \frac 1\gamma}$, the typical block spin profile is rigid, 
taking one of the two values corresponding to the minima of the
free energy for the random field Curie Weiss model, or makes at
most one transition from one of the minima to the other. { This holds for}
almost all realizations of the field. It was
also proven that   the
typical profiles are not rigid over any interval of length at
least $L_1(\gamma)=\frac 1 \gamma (\log \frac 1 \gamma)(\log\log
\frac 1\gamma)^{2+\rho}$, for any $\rho>0$. In  [\rcite {COP1}] 
the results  are shown for values of
the temperature and magnitude of the field in a subset of the 
region of two  absolute minima for the
free energy of the corresponding random {field}  Curie Weiss model.

In the present work we show  that,  on a set of realizations 
of the random field of probability that goes to 1 when $\g \downarrow 0$
,  we can construct random intervals of length of order
$\frac1 \g$ to 
which we associate a sign in such 
a way that the magnetization profile is rigid on these
intervals and, according to the sign, they belong to the $+$ or { to the}
$-$ phase.
A description of the transition from one phase to the other is 
also discussed.

The main problem in the proof of the previous results is the ``non
locality'' of the system, due to the presence of the random field.
Within  a run of positively magnetized blocks of length 1 in macro scale,
the ferromagnetic interaction will favor the persistence of blocks 
positively magnetized. The effect of the random 
magnetic  fields is related to  the  sum over these blocks
of the random magnetic fields.
It is relatively easy to see that  the fluctuations of the sum of the
random field over  intervals of order in macro scale  
$\frac 1 \g$ are the relevant ones.
But this is not enough.  To
determine the beginning, the end of the random interval, and the
sign attributed to it, it is essential to verify  other local
requirements for the random field. We need a detailed
analysis of the  sum  
of the random field in all subintervals of the large interval of
order  $\frac 1 \g$. In fact it could happen that even though at large
the random field undergoes to a positive (for example) fluctuation,
locally there  are negative fluctuations which make not
convenient (in terms of the total free energy) for the system to have a
magnetization profile close to the + phase in that interval.

Another problem in our analysis is due to the fact that  the
previously mentioned  block-spin transformation  
gives rise to a random  multibody  potential. 
Using a deviation inequality [\rcite{LT}], it turns out that 
for our analysis it is enough to compute the Lipschitz
norm of this  multibody potential. This is done by using 
cluster expansion tools to  represent this multibody potential  
as an absolute convergent series. 

The plan of the paper is the following. In Section 2 we give the
description of the model and present the main results. In Section
3 we   prove  probability estimates on functions of the random
field  which will allow us to construct the random intervals
  together with the corresponding sign.
In
Section 4  we  show   that,  typically, the magnetization profiles
are rigid over the macroscopic scale $\frac \e \g$, for any
$\e>0$, provided $\g$ is small enough.  This is an important
intermediate result. In Section 5 we finally prove the theorems
stated in Section 2. In Section 6 we prove some technical 
results needed in Section 5. 
In Section 7, we present a rather short, self contained and complete
proof of the convergence of  
the cluster expansion  for our  model. This is a standard tool
in Statistical Mechanics, but the application to this model is new.
\smallskip

{\bf Acknowledgments}  

This paper was motivated by a discussion with Errico Presutti, where the main
question was raised and an heuristic answer suggested.
We further thank him  for useful suggestions and  comments.   
Two of us (MC  and EO) acknowledge the hospitality of the CPT-CNRS
(Luminy) and of the IMPA and CBPF  (Rio de Janeiro).  
 PP acknowledges  the  hospitality of the Department of
Mathematics of Universit\`a di Roma--Tre, the Department of Physics 
of Universit\`a di Roma ``La Sapienza", the IMPA and  CBPF  (Rio de Janeiro), 
and the CMM (Santiago, Chili) where a part of this work was done. 
PP thanks also Roberto Cominetti (CMM) for an  enlightening  discussion. 
MEV acknowledges the hospitality of 
the Department of Mathematics of Universit\`a di Roma 
``Tor Vergata" and of the 
Department of Physics of Universit\`a di Roma ``La Sapienza".


\vskip 1.cm
\chap{2 Description of the model and  main results}2
\numsec= 2
\numfor= 1
\numtheo=1

 Let $(\O,\AA,\P)$ be a probability space on which we have defined
$h \equiv \{h_i\}_{i\in \Z}$, a family of
independent, identically distributed Bernoulli random variables with
$ \P[h_i=+1]=\P[h_i=-1]=1/2$. They represent random signs of external
magnetic fields
acting on a spin system on $\Z$, and whose magnitude  is denoted by
$\th>0$.
The configuration space is $\SS\equiv  \{-1,+1\}^\Z$.
If $\s \in \SS$ and $i\in \Z$,
$\s_i$ represents the value of the spin at site $i$. 
The pair interaction among spins is given by a Kac potential of the form
$J_\g(i-j)\equiv \g J(\g (i-j)),\, \g>0$,
on which one requires, for $r \in \R$: (i) $J(r)\geq 0$ (ferromagnetism);
 (ii) $J(r)=J(-r)$ (symmetry); (iii) $J(r)\le ce^{-c'|r|}$ for $c,c'$
positive constants
(exponential decay); (iv) $\int J(r)dr=1$ (normalization).
For simplicity we fix $J(r)= \1_{ [|r|\le 1/2]}$,
 though the behavior is the same under the above conditions.

For $\L \subseteq \Z$ we set $\SS_\L=\{-1,+1\}^\L$; its elements
are usually denoted by $\s_\L$; also,
if $\s \in \SS$, $\s_\L$ denotes its restriction to $\L$. Given
$\L\subset \Z$  finite and a realization of the magnetic fields,
the free boundary condition hamiltonian in the  volume $\L$ is given
by
$$
H_\g(\s_\L)[\o]= -\frac 12 \sum_{(i,j) \in \L \times \L}
 J_\g(i-j) \s_i \s_j -\th\sum_{i\in \L} h_i[\o]\s_i,
\Eq(2.1)
$$
which is then a random variable on $(\O,\AA,\P)$. In the following we drop
the $\o$ from the notation. 

The corresponding {\sl Gibbs measure} on the finite volume $\L$,
at inverse temperature $\b>0$ and free
boundary condition  is then a random variable with values
on the space of probability measures on $\SS_\L$. We denote it by 
$\mu_{\b,\th,\g,\L}$  and it is defined by
$$
\mu_{\b,\th,\g,\L}(\s_\L)
= \frac 1{Z_{\b,\th,\g,\L}} \exp\{-\b H_\g(\s_\L)\} \quad \quad
\s_\L \in \SS_\L,
\Eq(2.3)
$$
where  $Z_{\b,\th,\g,\L}$ is the normalization factor usually
called  partition function.

To take into account the interaction between the spins in $\L$ and
those outside $\L$ we set
$$
W_{\g}(\s_\L,\s_{\L^c}) = -\sum_{i\in \L} \sum_{j\in \L^c} J_\g(i-j)
\s_i \s_j.
\Eq(2.2)
$$
If $\tilde \s \in \SS$, the Gibbs
measure on the finite volume $\L$ and boundary condition
$\tilde \s_{\L^c}$ is  the  random probability measure on $\SS_\L$,
denoted by $\mu_{\b,\th,\g,\L}^{\tilde \s_{\L^c}}$  and
defined by
$$
\mu_{\b,\th,\g,\L}^{\tilde \s_{\L^c}}(\s_{\L})=
\frac 1{Z_{\b,\th,\g,\L}^{\tilde \s_{\L^c}}}
\exp\left\{- \b (H_{\g}(\s_{\L})+ W_{\g}(\s_{\L},\tilde
\s_{\L^c}))\right\},
\Eq(2.300)
$$
where again the partition function $Z_{\b,\th,\g,\L}^{\tilde
\s_{\L^c}}$ is  the normalization factor.

Given a  realization of $h$ and   $\g>0$, there is a unique
weak-limit of
$ \mu_{\b,\th,\g,\L}$ along a family of volumes $\L_{L}= [-L,L] \cap \Z$,  $  L
\in \N $; such limit
is called the infinite volume Gibbs measure
$\mu_{\b,\th,\g}$. The limit does not depend on the boundary
conditions, which may be taken
$h$-dependent, but it is a random element, i.e.,
different realizations of $h$ give a priori
different infinite
volume Gibbs measures.

As in [\rcite {COP}] and [\rcite {COP1}], our analysis of the large
scale profiles under $\mu_{\b,\th,\g}$
in the limit of $\g \downarrow  0$
involves a block spin transformation, which
transforms our microscopic system on $\Z$
into a {\it macroscopic} system on $\R$. Since the interaction length
is $\g^{-1}$,
one starts by a suitable scale transformation such that  on the
 new  scale, which we call {\it  the macroscopic scale}, the
interaction length becomes one.  Therefore,  a macroscopic
 volume, always taken as an interval  $I \subseteq \R$,
corresponds to the microscopic volume
$  \L(I)=\g^{-1}I\cap \Z$.
 The results will always 
be expressed in the macroscopic scale. 
The block spin transformation involves a ``coarse graining".
 Before making this precise let us set
some notations and  basic definitions,
mostly from [\rcite {COP1}].
\smallskip
\noindent  Given a rational positive number  $\d$,  $\DD_\d$ denotes the partition
of $\R$ into (macroscopic) intervals $\tilde A_\d(x)=((x-1)\d, x\d ]$
where $x \in \Z$. 
 If
$I\subset \R$ denotes a macroscopic interval
we let ${\cal C}_{\d}(I)=\{x \in\Z; \tilde A_\d(x) \subseteq I\}$.  In the following 
we will consider, if not explicitly  written,
intervals  always  in  macroscopic scale and 
$\DD_\d$--measurable, i.e.,     $I=\cup\{\tilde A_\d(x);
x \in{\cal C}_{\d}(I)\}$. 
\smallskip
The coarse graining will involve a scale $0<\d^*(\g)<1$ satisfying
certain conditions of smallness and will be the smallest scale.
The elements of $\DD_{\d^*}$ will be denoted by
$\tilde A(x)$, with $x \in \Z$.
The blocks $\tilde A(x)$ correspond to intervals of length $\d^*$ in  the
macroscopic scale  and induce a
partition of $\Z$ into blocks (in microscopic scale) of  order $\d^*
\g^{-1}$, hereby denoted by $A(x)=\{ i \in \Z ; i \g \in \tilde
A(x)\}=\{a(x)+1,\dots,a(x+1)\}$;  for notational
simplicity, if no confusion arises, we omit  to write the explicit dependence on 
  $\g,\d^*$.
We assume for convenience, that $\g=2^{-n}$ for some integer $n$, with
$\d^*$
such that $ \d^* \g^{-1}$ is an integer, so that $a(x)= x \d^*
\g^{-1}$, with $x \in \Z$.  We assume that $\d^*\g^{-1}\uparrow \infty$.

Given a realization $h[\o]\equiv (h_i[\o])_{i\in \Z}$, we set
$ A^+(x)= \big\{ i \in A(x); h_i[\o]=+1 \big\}$
and  $A^-(x)= \big\{ i \in A(x); h_i[\o]=-1 \big \}$. 
Let $ \l(x)\equiv {\rm sgn}( |A^+(x)|-(2\g)^{-1}\d^*)$, where ${\rm sgn}$
is the sign
function, with the convention that ${\rm  sgn}(0)=0$. For convenience
we assume  $\d^* \g^{-1}$ to be even, in which  case:
$$
\P\left[ \l(x)=0 \right]= 2^{-\d^* \g^{-1}}
{{\d^* \g^{-1}}\choose {\d^* \g^{-1}/2}}.
\Eq(2.7)
$$
Of course $\l(x)$ is a symmetric random variable.
When $\l(x)=\pm 1$ we set
$$
l(x)\equiv
\inf \{l > a(x)\, : \sum_{j=a(x)+1}^l \1_{\{ A^{\l(x)}(x)\}}(j)\geq
\d^* \g^{-1}/2\}
\Eq(2.8)
$$
and consider the following decomposition of $A(x)$:
$B^{\l(x)}(x)= \left\{ i\in A^{\l(x)}(x); i \le l(x) \right\}$
and $ B^{-\l(x)}(x)=A(x)\setminus B^{\l(x)}(x)$.
When $\l(x)=0$ we set $B^+(x)=A^+(x)$ and $B^-(x)=A^-(x)$.
We set $ D(x)\equiv A^{\l(x)}(x)\setminus B^{\l(x)}(x)$. In this
way,   the set $B^{\pm}(x)$ depend on  the 
  realizations of $\o$, but  the cardinality 
$|B^{\pm}(x)|=\d^* \g^{-1}/2$  is the same for all realizations.  
We  define
$$
  m^{\d^*}(\pm,x, \s)=\frac {2\g}{\d^*} \sum_{i\in B^{\pm}(x)} \s_i.
\Eq(2.10)
$$
We have 
$$
\frac \g{\d^*} \sum_{i\in A(x)}\s_i=
 \frac 12 (m^{\d^*}(+,x,\s)+m^{\d^*}(-,x,\s))
\Eq(2.12)
$$
and
$$
\frac \g{\d^*} \sum_{i\in A(x)} h_i\s_i =\frac 12
(m^{\d^*}(+,x,\s)-m^{\d^*}(-,x,\s)) +
\l(x) \frac {2\g}{\d^*} \sum_{i\in D(x)} \s_i .
\Eq(2.13)
$$
Given a  volume $\L\subseteq \Z$ in the original microscopic spin
system, it corresponds to the macroscopic volume $I=\g \L=\{\g i;
i \in \L\}$, assumed to be $\DD_{\d^*}$--measurable to avoid rounding problems.
 The block spin transformation, as considered in [\rcite {COP1}], is the
random map which associates to the spin configuration $\s_\L$
the vector $(m^{\d^*}(x, \s))_{x \in {\cal C}_{\d^*}(I)}$, where
$m^{\d^*}(x, \s)=(m^{\d^*}(+,x, \s),m^{\d^*}(-,x, \s))$, with
values in the set
$$
 \MM_{\d^*} (I)\equiv \prod_{x\in {\cal C}_{\d^*} (I)}
\left  \{-1, -1+ \frac {4\g}{\d^*}, -1+ \frac {8\g}{\d^*},\dots,
1- \frac {4 \g}{\d^*},1 \right\}^2.
\Eq (2.14)
$$

As in [\rcite {COP1}], we  use the same notation $\mu_{\b,\th,\g,\L}$
to denote
both, the Gibbs measure on
$\SS_\L$, and
the probability measure it induces on $ \MM_{\d^*}(I)$, through the
block spin transformation, i.e.,
a coarse grained version of the original measure.
Analogously, the infinite volume limit (as $\L\uparrow\Z$) of the
laws of the block
spin 
$(m^{\d^*}(x, \s))_{x \in {\cal C}_{\d^*}(I)}$
under the Gibbs measure
will also be denoted by  $\mu_{\b,\th,\g}$. If
$\;\lim_{\g \downarrow 0} \d^*(\g) =0$, this limiting measure
 will be supported by
$${\cal T}= \{ m \equiv (m_1,m_2) \in L^{\infty}(\R) \times L^{\infty}(\R);
\|m_1\|_\infty\vee\|m_2\|_\infty\le 1\}.  \Eq (space) $$ 
 
\medskip
\noindent  To denote a generic element in
$\MM_{\d^*}(I)$ we write
 $$
 m^{\d^*}_I\equiv(m^{\d^*}(x))_{x \in {\cal C}_{\d^*}(I)}
\equiv (m^{\d^*}_1(x),m^{\d^*}_2(x))_{x \in {\cal C}_{\d^*}(I)}. \Eq (2.14a)
$$
Since  $I$ is $\DD_{\d^*}$-- measurable, we can identify $m^{\d^*}_I$ with
the element of ${\cal T}$ which equals $m^{\d^*}(x)$ on each
$\tilde A(x)= ((x-1)\d^*, x\d^* ] $ for $x \in {\cal C}_{\d^*}(I)$, and 
vanishes outside
$I$. 
We denote by $T$, the
linear bijection on ${\cal T}$ 
defined by
$$
(Tm)(x)= (-m_2(x),-m_1(x))
\quad \forall x\in \R. \Eq(4.5P)
$$

\medskip
While analysing  some specific  block spin configurations, as in section 4, 
one encounters  a relevant functional that can be 
expressed as $\FF+ \g \GG$, where $\FF$ is deterministic and
$\GG$ is stochastic.

For the definition of $\FF$ we recall the relation of the local
behavior of the Random Field Kac model with the corresponding 
Random Field Curie-Weiss model.
The last one is  obtained when the
volume
$|\L|=\g^{-1}$  and its  canonical free energy $f_{\b,\th}(m_1,m_2)$
is given by
$$
f_{\b,\th}(m_1,m_2)=
-\frac{(m_1+m_2)^2}{8}-\frac{\th}{2}(m_1-m_2)+\frac{1}{2\b}(\II(m_1)+\II(m_2)),
\Eq (2.14000)
$$
for $(m_1,m_2) \in [-1,+1]^2$  and  
$\II(m)=\frac{(1+m)}{2}\log\big(\frac{1+m}{2}\big)+\frac{(1-m)}{2}\log\big(\frac{1-m}{2}\big)$.

Let us state  
some properties of $f_{\b,\th}(m_1,m_2)$. This will
point out the proper range of $\b,\th$ to be considered.
Differentiating \eqv (2.14000) we see that $(m_1,m_2) \in [-1,1]^2$ is
a  critical point of
$f_{\b,\th}(\cdot,\cdot)$ if and only if
$$
\eqalign{m_1&= \tanh (\b(m_1+m_2)/2 +\b\th)\cr
m_2&= \tanh (\b(m_1+m_2)/2 -\b\th).} \Eq (2.16)
$$
 The  sum of  the two   equations  in \eqv (2.16) is closed with
respect to  $\tilde m=(m_1+m_2)/2$ 
$$
\tilde m=\frac12 \tanh\b(\tilde m+\th) + \frac12 
\tanh\b(\tilde m-\th)\equiv g_{\b}(\tilde m,\th).
 \Eq(2.1800)
$$
It can be proved that
$$
\eqalign{&
1<\b<3/2  \;, \; 0<\th < \th_{1,c}(\b)\equiv\frac 1\b {\rm artanh}
(1-\b^{-1})^{1/2}; \text{ or }\cr &
3/2\le \b <+\infty  \;, \; 0<\th \le \th_{1,c}(\b)
}
\Eq(2.17)
$$
are necessary and sufficient for the existence of exactly three
solutions, $\tilde m=-\tilde m_{\b},0,\tilde m_{\b}$,
(with $\tilde m_{\b}>0$) to   equation
 \eqv(2.1800), verifying
$$
\frac{\del g_{\b}}{\del m}(\tilde m_{\b},\th) <1.
\Eq(2.17000)
$$
To simplify notations we do not write explicitly the dependence on $\th$ of $\tilde m_{\b}$. 
The result  on the solutions of  \eqv (2.1800)   implies that, setting   
$$  m_{\b,1}=\tanh \b(\tilde m_{\b}+\th); \qquad \qquad 
 m_{\b,2}=\tanh \b(\tilde m_{\b}-\th),  \Eq (phases) $$ 
  $m_\b=(m_{\b,1},m_{\b,2})$ and $Tm_\b=(-m_{\b,2},-m_{\b,1})$
are  solutions of \eqv  (2.16) corresponding to 
  the two global minima  of $f_{\b,\th}(\,\cdot\,)$,
$f_{\b,\th}(m_\b)=f_{\b,\th}(Tm_\b)$.  We  denote  $m_\b$  
  the $+$ phase and $Tm_\b$   the $-$ phase.
\smallskip 
\noindent {\bf Remark.}
Concerning equation \eqv(2.1800)
the following can also be proven: $\tilde m=0$ is the unique
solution, if $0<\b\le 1$. For $1<\b<3/2$, $\th\ge \th_{1,c}(\b)$,
again the unique solution is $\tilde m=0$ and  $\lim_{\th\uparrow
\th_{1,c}(\b)} \tilde m_{\b,\th}=0$. For $\b\ge 3/2$, there exists
$\th_{3,c}(\b)>\th_{1,c}(\b)$ such that for
$\th_{1,c}(\b)<\th<\th_{3,c}(\b)$ there exist  five solutions,
$\tilde m=-\tilde m_{2,\b,\th}, -\tilde m_{1,\b,\th},0, \tilde
m_{1,\b,\th},\tilde m_{2,\b,\th}$, with $0<\tilde
m_{1,\b,\th}<\tilde m_{2,\b,\th}$; when $\th\uparrow
\th_{3,c}(\b)$, $\tilde m_{1,\b,\th}\uparrow \tilde m_{3,\b}>0,
\tilde m_{2,\b,\th}\downarrow \tilde m_{3,\b}$, where $\tilde
m_{3,\b}=g_{\b}(\tilde m_{3,\b},\th)$ but $\frac{\del g_{\b}}{\del
m}(\tilde m_{3,\b},\th_{3,c}(\b))=1$; at last when
$\th>\th_{3,c}(\b)$, $\tilde m=0$ is the only solution.
 Property  \eqv(2.17000) will be constantly used in this work. In particular
we will not treat the case $\th=\th_{3,c}$.
\smallskip 
Throughout the work we assume \eqv(2.17) to be  satisfied, so that
$f_{\b,\th}(m_1,m_2)$ has  exactly three critical points, two
points of minima around which  $f_{\b,\th}(\,\cdot\,)$ is quadratic
and a local maximum. 
Moreover there exists a strictly positive constant
$\k (\b,\th)$ so that for each $m \in [-1,+1]^2$
$$
f_{\b,\th}(m)-f_{\b,\th}(m_{\b}) \geq \k (\b,\th)
\min\{\|m-m_\b\|^2_1,\|m-Tm_\b\|^2_1\}, \Eq(2.19)
$$
where 
$ \| \cdot \|_1$  the $\ell^1$ norm in $\R^2$ and 
$m_\b=(m_{\b,1},m_{\b,2})$, see \eqv (phases).

\remark 
Note that for $1<\b<3/2$, as $\th \uparrow \th_{1,c}$ we have
$\k (\b,\th) \downarrow 0$, but under \eqv(2.17) we have
always $\k (\b,\th)>0$. Since we want to work in the whole region
of $\b,\th$ where \eqv(2.17) is satisfied a little  care 
of $\k (\b,\th)$ will be taken.

We   introduce the so called ``excess  free energy
functional" $\FF(m)$,  $m \in {\cal T}$:
 $$
 \eqalign {& \FF(m)= \FF(m_1,m_2) \cr & = \frac 1 4 \int \int
J(r-r') \left [\tilde m(r)- \tilde m(r')\right ]^2 dr dr'+ \int
\left [ f_{\b,\th}(m_1(r),m_2(r)) - f_{\b,\th}(m_{\b,1},m_{\b,2})
\right] dr }\Eq (AP.1) $$   with
 $ f_{\b,\th}(m_1,m_2) $  given by \eqv (2.14000)
  and $\tilde
m(r)=(m_1(r)+m_2(r))/2$.
The functional $\FF$ is  well defined and non-negative,
although it may  take the value $+\infty$.
Clearly, the absolute minimum of $\FF$ is attained at the functions
constantly equal to the minimizers  of
$f_{\b,\th}$. 
$\FF$ represents the continuum approximation
of the deterministic contribution to the free 
energy of the system (cf. \equ(3.11P)) subtracted  by 
$f_{\b,\th}(m_{\b})$,
the free  energy of the homogeneous phases. 
Notice that  $\FF$  is invariant  under
the  $T$-transformation, defined in  \eqv (4.5P).
It has been proven in [\rcite {COP2}] that under 
the condition $m_1(0)+m_2(0)=0$,
there exists a unique  minimizer 
$\bar m=(\bar m_1,\bar m_2)$, 
of $\FF$ over the set
$$
\MM_{\infty}=\{(m_1,m_2)
\in {\cal T}; \liminf_{r\to -\infty} m_i(r) < 0<\liminf_{r\to +\infty} m_i(r),
i=1,2 \}.
 \Eq(AP.5)
$$
Without the condition $m_1(0)+m_2(0)=0$, there is a continuum of minimizers,
all other minimizers are translates of $\bar m$.
The minimizer 
$\bar m(\cdot )$ is infinitely differentiable.
Furthermore, there exists positive constant $c$ depending only on
$\b$ and $\th$  such  that
$$
\eqalign{
   \| \bar m(r)-m_\b\|_1 & \le ce^{-\a |r|}, \, \text{if} \,  r>0;\cr
  \|\bar  m (r)-Tm_\b\|_1 & \le ce^{-\a |r|}, \, \text{if} \ r< 0, \cr
 }
\Eq(V.1)
$$
where  $\a=\a(\b,\th)>0$ is given by (recall \equ(2.17000)):
$$
e^{-\a(\b,\th)}= \frac{\del g_{\b}}{\del m}(\tilde m_{\b,\th},\th). 
\Eq (Palmares)
$$
Since   $\FF$ is invariant by  the $T$-transformation, see \eqv (4.5P), 
interchanging $r\to \infty$ and $r \to -\infty$  in \eqv (AP.5)
there exists  one other family of minimizers
obtained translating $T \bar m$. 
We denote  $$\FF^*= \FF(\bar m)= \FF(T\bar m)>0.  \Eq (min) $$ 
 The functional $\FF$ that enters in the above decomposition into 
a deterministic and a stochastic part,  $\FF+\g \GG$,
is merely a finite volume version of \eqv(AP.1); however  
\eqv(V.1) and  $\FF^*$  will play a crucial role  here. 
\smallskip
The stochastic part of the functional $\cal G$ 
is defined on $\MM_{\d^*}(I)$ (embedded in ${\cal T}$
as previously mentioned) as
$$
\GG(m^{\d^*}_{I})
\equiv
 \sum_{x \in \CC_{\d^*}(I)}
 \GG_{x,m^{\d^*}(x)}(\l(x))
 \Eq(3.14P)
$$
where  for each $x\in\CC_{\d^*}(I)$, $\GG_{x,m^{\d^*}(x)}(\l(x))$ is
the cumulant generating function:
$$
\GG_{x,m^{\d^*}(x)}\left(\l(x)\right)
\equiv -\frac1\b\log  \E_{x,m^{\d^*}(x)}^{\d^*} (e^{ 2\b\th \l(x)\sum_{i\in
D(x)}\s_i}),
\Eq(3.13P)
$$
of the ``canonical" measure on $\{-1,+1\}^{A(x)}$,  defined through
$$
\E_{x,m^{\d^*}(x)  }^{\d^*}(\vf)=
\frac  {\sum_{\s} \vf(\s)\1_{
\{m^{\d^*}(x,\s)=m^{\d^*}(x)\}}} {\sum_{\s}
\1_{ \{m^{\d^*}(x,\s)=m^{\d^*}(x)\}}},
\Eq(3.9P)
$$
the sum being over $\s \in \{-1,+1\}^{A(x)}$.
 Let   $m^{\d^*}_\b$ be one of  the points in
$ \left  \{-1, -1+ \frac {4\g}{\d^*},\dots, 1-
\frac {4 \g}{\d^*},1 \right\}^2$ which is closest to $m_\b$. Given
an interval $I$ we let $m^{\d^*}_{\b,I}$ be the function
which coincides with $m^{\d^*}_\b$ on $I$ and vanishes outside
$I$.    In the analysis of the random fluctuations of our system 
the relevant random quantities will be  
$$
\GG(m_{\b,I}^{\d^*})-\GG(Tm_{\b,I}^{\d^*})=:\sum_{x \in
\CC_{\d^*}(I)}X(x).
\Eq(2.169)
$$
One important property of the random variables $X(x)$ is their symmetry.
The explicit expression of $X(x)$ that one gets  using \eqv(2.169), 
\eqv(3.13P), and \eqv(3.9P) is  almost  useless. One can think about 
making an expansion in $\b\th$ as we basically did in [\rcite{COP1}],
Proposition 3.1  
where $\b\th$ was assumed to be as small as needed. 
 Since here we assume \eqv(2.17), one has to find 
another small quantity. Looking at the term $\sum_{i\in D(x)}\s_i$
in \eqv(3.13P) and  setting   
$$
p(x)\equiv p(x,\o)= |D(x)|/|B^{\l(x)}(x)| = 2\g
|D(x)|/ \d^*,
\Eq(2.1690)
$$
it is easy to see that  
for  $I\subseteq \R$,
if 
 $
\left(\frac{2\g}{\d^*}\right)^{1/2} \log \frac{|I|}{\d^*}\le \frac {1}{32} 
$,
we have
$$
\P\left [ \sup_{x\in \CC_{\d^*}(I)} p(x) > (2\g/\d^*)^{\frac{1}{4}}
\right]
\le e^{-\frac{1}{32}\left(\frac{\d^*}{2\g}\right)^{\frac 12}}.
\Eq(2.17003)
$$
\remark 
Note at this point that the choice of $\d^*$ as
$\g\log\log (1/\g)$ we made  in [\rcite{COP1}], for volume $I$ of order 
$\g^{-1}$ does not satisfied the previous  restriction.
 
Now on the set $\left\{ 
\sup_{x\in \CC_{\d^*}(I)} p(x) \le  (2\g/\d^*)^{\frac{1}{4}}\right\}
$, $p(x)$ is a small parameter (recall $\d^*\g^{-1}\uparrow \infty$). 
It will be proved in Proposition \eqv(PP), see remark \eqv(P10),
that 
on the set   $ \{ p(x) \le (2\g/\d^*)^{1/4}\} $, 
the  quantity $X(x)$ can be written as:
$$
X(x)=
-\l(x)|D(x)|\left[ \log\frac{1+m^{\d^*}_{\b,2}\tanh(2\b\th)}
{1-m^{\d^*}_{\b,1}\tanh(2\b\th)}+ \Xi_1(x,\b\th,p(x))\right]
-\l(x) \Xi_2(x,\b\th,p(x))
\Eq(2.170)
$$
with  
$$
\left|\Xi_1(x,\b\th,p(x))\right|\le
 64  \frac {\b\th (1+ \b \th) }  {(1-m_{\b,1})^2 (1-\tanh(2\b\th))}   (2 \frac \g {\d^*})^{1/4}.      
\Eq(2.1692)
$$
and 
$$
\left|\Xi_2(x,\b\th,p(x)\right|
\le  (2\frac \g {\d^*})^{1/4} \left [ 36+2 c(\b \th) \right ]
\Eq(2.1700)
$$
where $c(\b\th)$ is given in \eqv (3.25PM).

Thus, calling
$$
V(\b,\th)= \log\frac{1+m_{\b,2}\tanh(2\b\th)}
{1-m_{\b,1}\tanh(2\b\th)},
\Eq(2.1701)
$$
on the 
event $ \{ p(x)\le(2\g/\d^*)^{\frac{1}{4}}\}$,
when $ \frac  \g {\d^*} \downarrow 0$
the leading term in \eqv(2.170)  is simply 
$$
 -\l(x)|D(x)| V(\b,\th)
=
 -V(\b,\th)\sum_{i\in A(x)} h_i,
\Eq(2.173)
$$
and,  from \eqv(2.170), we have
$$
\eqalign{
\E[X(x)\1_{\{p(x) \le (2\g/\d^*)^{1/4}\}}]&=0,\cr
\E[X^2(x)\1_{\{p(x) \le (2\g/\d^*)^{1/4}\}} ]&= \frac {\d^*} {\g} 
c(\b,\th,  \g /\d^*)}
\Eq (B.20)
$$
where, if $\g/\d^*<d_0(\b,\th)$ for suitable $0<d_0(\b,\th)$,
$ c(\b,\th, \g /\d^*)$  satisfies:
$$
 V^2(\b,\th)
\left[1-(\g/\d^*)^{1/5}\right]^2 \le  c(\b,\th, \g /\d^*) \le   
V^2(\b,\th) \left[1+
(\g/\d^*)^{1/5}\right]^2.      \Eq(2.21). 
$$

Our final aim is to control the behavior  of the random field over 
intervals of (macroscopic)
length  of order  larger  or equal to  $\frac 1 \g $. To achieve this, 
it is convenient to
consider  blocks of (macroscopic)
length $\e/\g$, with the basic assumption that $\e/\g>\d^*$.
To avoid rounding problems we assume $\e/\g\d^*\in \N$ and we define,
for $\a \in \Z$
$$
\chi^{(\e)}(\a) \equiv \g\sum_{x;\d^*x \in \tilde A_{\e/\g}(\a)}
 X(x)\1_{\{p(x) \le (2\g/\d^*)^{1/4}\} }
\Eq(2.174)
$$
where, according to the previous notation $\tilde
A_{\e/\g}(\a)=((\a-1)\sfrac \e \g,\a\sfrac{\e}{\g} ] \subseteq \R$
and for sake of simplicity the $\g,\d^*$ dependence is not
explicit. To simplify further, and if no confusion arises, we
shall write simply $\chi(\a)$. Note that  $\chi(\a)$ is a
symmetric random variable and assuming that $I\supset
\tilde A_{\e/\g}(\a)$ for all $\a$ under consideration 
$$
\eqalign{ \E(\chi(\a))&=0,\cr
\E(\chi^2(\a))&=\e c(\b,\th, \g
/\d^*),} \Eq(2.20)
$$
as it follows from \equ(B.20) since there are $\e(\g\d^*)^{-1}$
terms in the sum in \eqv(2.174).

\smallskip
As in [\rcite {COP1}], the description of the profiles is based on
the behavior of local averages of
$m^{\d^*}(x,\s)$ over $k$ successive blocks in the block spin
representation, where
$k\geq 2$ is a positive integer. Let $\d=k\d^*$  such that  $1/\d
\in \N$  and let
${\cal C}_\d(\ell)\equiv {\cal C}_\d((\ell-1,\ell])$ defined as before.
Given    $\z \in (0,m_{2,\b}]$ and  $\ell\in \Z$,
we define the random variable
$$
\eta^{\d,\z}(\ell)=\cases{
&1 if $\,\forall_{u \in {\cal C}_{\d}(\ell)}\,
\frac {\d^*}{\d}\sum_{x\in {\cal C}_{\d^*}(((u-1)\d,u\d])}
\|m^{\d^*}(x)-m_\b\|_1 \le\z$;\cr
 &-1 if $\,\forall_{u \in {\cal C}_{\d}(\ell)}\,
\frac {\d^*}{\d}\sum_{x\in {\cal C}_{\d^*}((u-1)\d,u\d])}
\|m^{\d^*}(x)-Tm_\b\|_1 \le\z$;\cr
&0  otherwise.\cr }
\Eq (2.190)
$$
We say that   a magnetization profile $ m^{\d^*}(\cdot)$,
in an interval  $I \subseteq \R$, 
is close to  the equilibrium phase $\t$, $\t=1$ or $\t=-1$,  with  tolerance
$\z$, when 
$$ \{ \eta^{\d,\z}(\ell) = \t,\, \forall  \ell \in I \cap \Z \}  \Eq (equi) $$
In the following we will use always the letter $\ell$
to indicate an element of $\Z$.  This will allow to 
 write
\eqv (equi)  as 
  $ \{ \eta^{\d,\z}(\ell) = \t,\, \forall    \ell \in I  \}   $.

 Given a realization of $h$, we would like to know if ``typically''
with respect to the Gibbs measure we have, as an example, $\eta^{\d,\z}(0)=1$
or $\eta^{\d,\z}(0)=-1$.  The alternative depends on this 
realization of $h$.
Here typically means with an overwhelming  Gibbs measure but having in mind
a exponential convergence.  
First of all, one has to  accept to throw 
away some realizations of $h$ that are not ``typical'' with respect to the 
$\P$--probability. However, depending on  the probabilistic sense
 of ``typical''
one can easily convince himself that the results will be completely different. 
Here we just want that  the $\P$-probability of the realizations of $h$ that 
we throw away goes to zero when $\g\rightarrow 0$.
Some $\P$-almost sure  results can be found in [\rcite{COP1}].
It happens that to  give an answer to such a simple question 
we must know if $\eta^{\d,\z}(0)$ belongs to 
a run of $\eta^{\d,\z}=1$ or to a run of $\eta^{\d,\z}=-1$.
It is rather clear that we have to understand 
the localization of the beginning and the end of consecutive  runs with 
alternating sign.  However to define the beginning and the end of a run,
we have to take into account that  some  messy configurations with 
$\eta^{\d,\z}=0,\pm 1$  could occurs in between two such runs.
So in the first theorem we erase deterministically pieces around 
what we expect to be the endpoints  of the run that countains 
the origin. In the second theorem we consider consecutive 
runs with erased endpoints. In the last theorem we prove that in 
the erased regions between two runs there is just a single run 
of $\eta^{\d,\z}=0$ which is rather short.

The main result of this paper is the following:
\vskip0.5cm 
\noindent{\bf \Theorem (t1)} {\it Given  $(\b,\th)$ that
satisfies \eqv(2.17), $a>0$, 
$\k(\b,\th)>0$  satisfying \eqv(2.19), there exist
$0<\g_0=\g_0(\b,\th)<1$, $ 0<d_0=d_0(\b,\th,a)<1$,
 and $0<\z_0=\z_0(\b,\th)<1$,
such that for all $0<\d^*<1$,    $ 0<\g \le \g_0 $, 
$\g/\d^*\le d_0$,
if $ \z_4$ is  such that    $\z_0 \geq \z_4 >8\g/\d^*$,
$g$ is a positive increasing 
function such that $g(x)\ge 1$,
 $\lim_{x\uparrow \infty}g(x)=+\infty$ and  $\frac {g(x)} x \le 1 $, 
$\lim_{x\uparrow
\infty}x^{-1} g^{38}(x)=0 $,
$$\z_4 > \frac {1}{[\k(\b,\th)]^{1/3}g^{1/6}(\frac{\d^*}{\g})},
\Eq(TE.1)$$
and
$$
\frac{(\d^*)^2}{\g}
g^{3/2}(\frac{\d^*}{\g}) \le 
 \frac 1{\b \k (\b,\th) e^3 2^{13}},
\Eq(TE.2)
$$
then there  exists $ \O_{\g, \d^*}$ with
$$
\P[\O_{\g,\d^*}] \ge 
1-16 \g^2- 
160 \left( g(\frac {\d^*}{\g})\right)^{-\frac{a}{4(2+a)}}
\Eq(5.02)$$
such that   for all realizations of the fields  $ \o\in
\O_{\g,\d^*}$, for
$\e=\left(\frac{5}{g(\d^*/\g)}\right)^4$,
 we can construct explicitly   a random { measurable} pair
$\big(I(\o), \t (\o) )$ where
$$
\t (\o)= {\rm {sgn}}\Big(\sum_
{\a\in \CC_{\e/\g}(I(\o))} \chi(\a) \Big) \in \{-1,+1\} $$
$I(\o)$ is a suitable random macroscopic interval that contains the
origin such that for all  $x>0$
$$
\P ( \o \in \O_{\g, \d^*} \colon\g |I(\o)| >x ) \le 4 
e^{-\frac  x { 8 C_1(\b,\th, \FF^*) } ( 1- \frac {\log3} {\log4})},
   \Eq(2.210a) $$    
$$ \quad
\P(
\o\in \O_{\g, \d^*}\colon \g|I(\o)| <x ) \le 2 
e^{-\frac {(\FF^*)^2}{ 18 x V^2(\b, \th)}} 
 \Eq(2.210b)
$$
where  $C_1(\b,\th, \FF^*)$ is given in \eqv (3P.50), $ \FF^*$ 
 in  \eqv (min)   and $V(\b,\th)$ in \eqv (2.1701). 
The
interval $I(\o)$ is measurable with respect to the $\s$-algebra
$\s\big( \chi(\a), \a \in \CC_{\e/\g}( [ -\frac{ Q}{\g}, \frac {Q}
{\g}])\big)$ where $ Q=\exp\frac{\log g({\d^*}/{\g})}{
\log \log g({\d^*}/{\g})} $, and we have
$$
\mu_{\b,\th,\g}\Big[\forall \ell \in I(\o) \cap \Z,
\eta^{\d,\z_4} (\ell)=\tau (\o)\Big]
\ge 1-  e^{- \sfrac {\b }{\g}
\sfrac{1}{ g(\d^*/\g)}}.
\Eq(2.2100)
$$
here $\d=5^{-1}(g(\d^*/\g))^{-1/2}$.
Moreover the interval $I(\o)$ is maximal, in the following sense:
$\forall J\in \R$, $ I(\o)\subseteq J $,   $|
J \setminus I(\o)| \ge 2\frac \rho \g $, 
with
$\r=\left(\frac{5}{g(\d^*/\g)}\right)^{1/(2+a)}$,

$$
\mu_{\b,\th,\g}\Big[\forall \ell \in J \cap \Z , \eta^{\d,\z_4}
(\ell)=\tau\Big]
\le
e^{- \sfrac {\b }{\g} \sfrac{1}{ g(\d^*/\g)}}.
\Eq(2.2?)
$$
}
\noindent {\bf { Remark.    (Choice of the parameters)}}  The main parameters 
appearing in the problem, besides
$\b, \th $  and $\g$, (we take $\b, \th $, in all the paper, to satisfy 
\eqv (2.17)  and $ \g>0$  small enough to control the
range and the strength of the Kac interaction), are the  smallest coarse 
grained
scale $\d^*$   and the tolerance $\z_4$   around  the ``equilibrium" 
$m_{\b}$ or $Tm_{\b}$.   
We choose a specific  $\d$ for simplification.
There exists an important constraint on how small can $\d^*$ be
taken.
The convergence of the cluster expansion
requires $\frac{ {\d^*}^2}\g \le \frac 1 {6e^3 \b}$, cf. Theorem  
  \eqv (62P).  
 The constraint on $\d^*$ appearing in  
Theorem
\eqv (t1),
\eqv (TE.2), is stronger since to estimate the random field we need to compute
the Lipschitz norm of the multibody term coming from the
cluster expansion and stronger requirements are needed.  

 We decide to write
 the results in term of a rather general function $g$,
verifying the requirements written in  Theorem \eqv (t1). A prototype 
can be  $ g(x)= 1 \vee  \log x$ or any iterated of it.   
The main reason to do this is to have the simplest expression for the 
Gibbs measure estimate \eqv (2.2100). As a consequence, the $ \P $
probability estimate in \eqv (5.02)  is also expressed in term of this 
function  $ g$ as well as all the constraints on the parameters. 
The condition   $ \lim_{x\uparrow
\infty}x^{-1} g^{38}(x)=0 $ comes from an explicit choice of an auxiliary
 parameter $\z_5 $ that will be introduced in Section 5 and the
constraint  \eqv (5.003) that has to be satisfied. 
Notice that taking $g(x)=1\vee \log x$ and $\d^*=\g^{\frac{1}{2} +d^*}$ 
for some
$0<d^*<1/2$ implies that \eqv(TE.2) is satisfied.  

\noindent
Finally the choice of the numerical constants (such as $2^{13}$) is
 never critical and largely irrelevant. We have made no efforts to make
the choices close to optimal. 

\noindent
{\bf { Remark.}}  The endpoints  of the random interval $I(\o)$  
are not stopping times, as it can be seen in Section 3. However, the 
interval
$I(\o)$ is measurable with respect to the
$\s$-algebra
$\s\big( \chi(\a), \a \in \CC_{\e/\g}( [ -\frac{ Q}{\g}, \frac {Q}
{\g}])\big)$,  where $Q$  is given in Theorem \eqv(t1).  Therefore, 
in order to decide if typically $\eta^{\d,\z_4}(0)=+1$ or $-1$,
it suffices to know the realization of the random magnetic fields 
in a volume which, choosing for example,  $ g(x)= 1 \vee  \log x$, 
is  of the order 
$\frac 1 \g \left ( \log \frac 1 \g  \right )^{ \frac 1 
{ \log \log \log \frac 1 \g }}$ 
in macroscopic scale. 

\smallskip
Our next result is a simple
extension  of the
previous theorem.

\noindent{\bf \Theorem (t2)} {\it Under the same hypothesis of
Theorem \eqv(t1),
for all $k \in \N$,
there exists $\O_{\g,\d^*,k}$, with
$$
\P[\O_{\g,\d^*,k}] \ge  
1-32k\g^2- 
320 k \left( g(\frac {\d^*}{\g})\right)^{-\frac{a}{4(2+a)}}
\Eq(2.21000)
$$
 such that for   $\o \in
\O_{\g,\d^*,k}$,
we can construct explicitly a random $(2k+2)$--tuples
$$
\Big(I_{-k}(\o),\dots,I_{k}(\o),
{\rm {sgn}}\Big(\sum_{\a\in \CC_{\e/\g}(I_0(\o))} \chi(\a)\Big) \Big)
\Eq(2.2100000)
$$
where
$I_j(\o), -k\le j\le k$ are suitable disjoint random
intervals, $I_0(\o)$ contains the origin and they satisfy for all $x>0$
$$
\P \left [ \sup_{ -k\le j\le k} \g |I_j(\o)| >x \right ] \le 4(2k+1) 
 e^{-\frac  x { 8 C_1(\b,\th, \FF^*) } ( 1- \frac {\log3} {\log4})},
\Eq(2.2102a) $$    
$$ \qquad \P \left [\inf_{ -k\le j\le k} \g|I_j(\o)| <x \right ] 
 \le (2k+1)
2 e^{-\frac {(\FF^*)^2}{ 18 x V^2(\b, \th)}},  \Eq(2.2102b)
$$
where  $C_1(\b,\th, \FF^*)$ is given in \eqv (3P.50), $ \FF^*$  in 
 \eqv (min)   and $V(\b,\th)$ in \eqv (2.1701). 
The sequence
$(I_{-k}(\o),\dots,I_k(\o))$ is measurable with respect to the
$\s$-algebra
$\s\big( \chi(\a), \a \in \CC_{\e/\g}( [-kQ, kQ])\big)$, 
and
$$
\Big|[\inf (I_{-k}(\o)),\sup(I_k(\o))] \setminus \bigcup_{j=-k}^{k}
I_j(\o)\Big|
\le (2k+1) \frac{\r}{\g}.
\Eq(2.2103)
$$
Moreover for all $-k\le j\le k$,
$$
\mu_{\b,\th,\g}\Big[ \eta^{\d,\z} (\ell)=(-1)^j {\rm sgn}
(\sum_{\a\in \CC_{\e/\g}(I_0(\o))} \chi(\a)),\; 
\forall j\in \{-k,+k\},
\forall \ell \in I_j(\o)
\Big]
\ge
 1-2k e^{- \sfrac {\b }{\g} \sfrac{1}{ g(\d^*/\g)}}.
\Eq(2.2104)
$$
}

In the previous  theorem nothing is said about what  happens in
the region between two consecutive intervals with different signs,
a  region that has a macroscopic length smaller than $ \r/\g$ by
\eqv(2.2103), see before \eqv(2.2?) for $\r$. 
To describe it  we need to introduce the notion of a
single change of phases in a given interval.

\noindent {\bf \Definition  (1)} {\it 
Given an interval $[\ell_1,\ell_2]$ and a  positive integer
 $R_2< | \ell_2-\ell_1|$,
 we   say that a  single change of phases occurs within $[\ell_1,\ell_2]$
on a length
$R_2$
if there exists $ \ell_0\in [\ell_1,\ell_2]$ so that $\eta^{\d,\z}
(\ell)=\eta^{\d,\z}(\ell_1)\in \{-1,+1\}, \forall \ell 
 \in[\ell_1, \ell_0 -R_2]$,
$\eta^{\d,\z}(\ell)=\eta^{\d,\z}(\ell_2)=-\eta(\ell_1), \forall  \ell \in
[\ell_0+R_2,\ell_2]$, and 
$\{\ell \in [\ell_0 -R_2,\ell_0+R_2]: \eta^{\d,\z}(\ell)=0\}$ is a set of
consecutive integers.
We denote by $ \WW_1([\ell_1,\ell_2], R_2,\z)$  the  set 
of all  configurations of $\eta^{\d,\z}$ that satisfies
this properties.}

In other words, there is an unique
run of $\eta^{\d,\z}=0$, with no more than $R_2$ elements,  inside
the interval $[\ell_1,\ell_2]$.  

Our next  result is

\noindent{\bf \Theorem (t3)} {\it
Under the same hypothesis as in Theorem  { \eqv(t2)} 
 and on the same probability space
$\O_{\g,\d^*,k}$, for
$$
R_2=\frac {20(5+ \FF^*)160^3}{\k(\b,\th)} 
\left(g(\frac{\d^*}{\g})\right)^{7/2}
\Eq(2.22)
$$
we have}
$$
  \mu_{\b,\th,\g}\Big[ \bigcap_{-k\le j\le k-1} \WW_1(
[\sup(I_j(\o)), \inf(I_{j+1})],R_2,\z_4)\Big]\ge
 1- 2ke^{-\frac{\b}{\g} \frac{1}{g(\d^*/\g)}}.
\Eq(2.2106)
$$

\smallskip

Note that the regions where the changes  of phases
occur  have  at most  length 
$R_2$ (in macroscopic units)
and we are able to
localize it only within an interval of length $\r/\g>>R_2$.
This means that up to a
small  probability subset,  we are able to give an explicit way of
constructing an interval of length $\r/\g$ where
we have a change of phases that occurs on a scale $R_2$, but
we are not able to   determine where it
occurs within this interval.

\bigskip
\chap{3 Probabilistic estimates}3 \numsec= 3 \numfor= 1
\numtheo=1

In this  section we  construct  a random interval
$J(\o)$, to which the interval $I(\o)$ appearing in Theorem  \eqv
(t1)  is simply related. The construction involves a discrete
random walk obtained from the variables $\chi(\a),\a\in
\Z$, defined by \eqv(2.174) and satisfying \eqv(B.20).
If $\D$ is a finite
interval in $\Z$ we set $\YY(\D)= \sum_{\tilde \a \in \D}
\chi(\tilde \a)$. For convenience we write
$$
\YY_\a \equiv \cases {\;\;\YY(\{1,\dots, \a\}), & if
$\a\ge 1$;\cr \;\;0 & if $\a=0$\cr -\YY(\{\a+1,\dots,0\}),
& if  $\a\le -1$.\cr}\Eq(3.C6)
$$
so that if  $\D=\{\a_1+1,\dots, \a_2\}\equiv(\a_1,\a_2]$, with
$\a_1<\a_2$ integers,  we have
$\YY(\D)=\YY_{\a_2}-\YY_{\a_1}$.
\medskip
As $\g \downarrow 0$, we assume $\e \downarrow 0$ but
$\e/\g\d^*\uparrow+\infty$. In this regime, $\YY_{[\cdot/\e]}$
converges in law to a bilateral Brownian motion (no drift,
diffusion coefficient $V(\b,\th)$).

 Given a real positive number $f$,
$0<f<\FF^*/4$ where $\FF^*$ is defined in \eqv(min), we denote
$$
\DD(f,+) \equiv \DD(f,+, \o)\equiv \left\{ \D \colon\,
 \YY(\D) \ge 2\FF^* +f\, ,\,
 \inf_{\D' \subset \D}\YY(\D') \ge -2 \FF^* +f \right\},
\Eq(3.3)
$$
the set of random (finite) intervals $\D\subseteq \Z$ with  
an (uphill) increment of size at least
$2\FF^* +f$, and such that no interval within $\D$ presents 
a (downhill) increment
smaller  than
$-2 \FF^* +f $. Such an interval $\D\subseteq \Z$ is said to 
give rise to a {\it positive
elongation}, and we set ${\rm sgn}\,\D=+1$.

Similarly,
$$
\DD(f,-) \equiv \DD(f,-,\o)\equiv 
\left\{ \D\colon\, \YY(\D) \le -2\FF^* -f\, ,\,
 \sup_{\D' \subset \D}\YY(\D') \le 2 \FF^* -f \right\},
\Eq(3.4PP)
$$
and such an interval is said to give rise to  a {\it negative elongation}.
If $\D \in \DD(f,-)$, we set ${\rm sgn}\,\D=-1$. We call
 $$\DD(f,\o)\equiv \DD(f,+,\o) \cup\DD(f,-,\o) \Eq (3.C5)$$

\remark $\DD(f,+)\cap \DD(f,-)=\emptyset$ since $f>0$, so that the
above definition of ${\rm sgn}\D$ is  well posed. However, we
may have intervals $\D_1 \in \DD(f,+)$ and $\D_2\in \DD(f,-)$ such
that $\D_1\cap\D_2\neq \emptyset$.

 \vskip0.5cm
 \noindent
Given  $Q>0$ and  writing   $A^c=\O\setminus A$, we let
$$
\PP_0(f,Q)=\left\{
 \exists \D \in \DD(f,\o), \;
\D \subseteq [-Q/\e, Q/\e ] \right\}^c,
\Eq(3.5)
$$
be the set of realizations of the random field that neither  give
rise to a positive   nor to a negative  elongation in the interval
$ [-Q/\e, Q/\e]$. As we will see later, cf. Theorem \equ(T1),
$\P[\PP_0(f,d)]$ is small provided $Q$ is large, {\it uniformly}
on $0 < f \le  \FF^*/4$. (The uniformity is trivial since from the
definitions $\DD(f,\pm) \subseteq \DD(\tilde f,\pm)$ if $0<\tilde
f<f$.)

Deciding if a given interval gives rise to a positive or negative
elongation is a local procedure, in the sense that it depends only
on the values of $\chi(\a)$, with $\a$ in the considered interval.
But, since our goal is to find the beginning and the end of
successive runs of $\eta^{\d,\z}=+1$, and runs of
$\eta^{\d,\z}=-1$, we should    determine contiguous  elongations
with alternating signs. For this we first need (not necessarily
contiguous) elongations with alternating signs. We set, for $k
\in \N$:
$$
\eqalign{
 B_+(f,k, Q) \equiv  \big\{ &\o \in \O\colon \exists\, 0\le a_1 <b_1\le a_2 <
b_2\le\dots\le a_k <b_k \le Q/\e, (a_i,b_i] \in \DD(f),\cr &i=1,..,k;\;
{\rm sgn}(a_i,b_i]= -{\rm sgn}(a_{i+1},b_{i+1}],
 \, i=1,..,k-1  \big\},\cr
}\Eq(3.51)
$$
$$
\eqalign{B_-(f,k, Q) \equiv  \big\{ &\o \in \O\colon \exists\,
 0\ge b_1>a_1\ge b_2>
a_2\ge\dots\ge b_k>a_k\ge -Q/\e, (a_i,b_i] \in \DD(f),\cr &i=1,..,k;\;
{\rm sgn}(a_i,b_i]= -{\rm sgn}(a_{i+1},b_{i+1}],
 \, i=1,..,k-1  \big\},\cr}
\Eq(3.521)
$$
and
$
\PP_1(f,k,Q) \equiv \left ( B_+(f,k,Q)\cap B_-(f,k,Q)\right )^c
\supseteq \PP_0(f,Q).
$
In Theorem \equ(T1) we shall prove that $\P[\PP_1(f,$ 
$ k,kQ)]$ is
small, {\it uniformly} in $0< f \le \FF^*/4$, and $k \ge 1$,
provided $Q$ is taken large enough.

For reasons that will be clear later we set:
$$
\eqalign{ \PP'_2(f,Q)=\{\exists &\a_1<\a_2<\a_3<\a_4 \in
[-Q/\e,Q/\e]\colon |\YY _{\a_1}-\YY_{\a_3}|\vee |\YY
_{\a_2}-\YY_{\a_4}|\le 3f,\cr &||\YY
_{\a_1}-\YY_{\a_2}|-2\FF^*|\le 3f,\cr & \YY_{\a} \in
[\YY_{\a_1}\wedge \YY_{\a_2}-3f,\YY_{\a_1}\vee
\YY_{\a_2}+3f],\forall \a\in [\a_1,\a_4]\}\cr}
$$
and 
$$
\PP''_2(f,Q)=\PP'_2(f,Q)\cup\{\max_{\a \in [-Q/\e,Q/\e]}|\chi(\a)|>
f\}.
\Eq(30.54)
$$

To construct the
  previously described $J(\o)$, with $0 \in J(\o)\subseteq [-Q/\g,Q/\g]$,
  it will suffice to have
$\o \in \big(\PP_1(f,3,Q) \cup \PP''_2(f,Q)\big)^c$. Having fixed $Q$
sufficiently large so that $\P(\PP_1(f,3,Q))$ is suitably small
for any $0< f \le \FF^*/4$, we shall take $f$ small enough and
$\e$ suitably small so that $\P(\PP''_2(f,Q))$ is also suitably
small, as stated in Theorem \equ(T1).
\medskip

Let $\o \in (\PP_1(f,3,Q) \cup \PP''_2(f,Q))^c$. Starting  at
$\a=0$, and going to the right we tag the ``first" interval in
$\Z$ which provides an elongation. We then use an explicit way to
construct {\it contiguous} intervals that provide elongations with
alternating signs. $J(\o)$ will be defined with the help of such
elongations. Having a discrete random walk, different types of
ambiguities appear in this construction  and we need to
estimate the probability of their  occurrence. We discuss a
possible construction.

Let us define for each $a,b \in [-Q/\e,Q/\e]\cap \Z$:
$$
\eqalign{
b_-(a) &\equiv \inf\{ b' >a \colon (a,b'] \in \DD(f,\o) \}\cr
b_+(a) &\equiv \sup \{ b' >a \colon (a,b'] \in \DD(f,\o) \}\cr
a_+(b) & \equiv \sup \{ a'<b \colon (a',b] \in \DD(f,\o)  \}\cr
a_-(b) & \equiv \inf \{ a'<b \colon (a',b] \in \DD(f,\o)\},}\Eq(3.5400)
$$
with the infima and suprema taken on $[-Q/\e,Q/\e]\cap \Z$; thus,
if the corresponding set is non-empty we have a minimum or
maximum; otherwise we make the usual convention: $\inf\emptyset=
+\infty$ and $\sup\emptyset=-\infty$.

We see at once:

$\bullet$ if $b_-(a)<+\infty$ then $ a_-(b_-(a)) \le a \le  a_+(b_-(a))$;

$\bullet$ if $a_+(b)>-\infty$ then $b_-(a_+(b)) \le b \le b_+(a_+(b))  $.

Let us set
$
a_0\equiv\inf\left\{a\ge 0 \colon b_-(a) < +\infty \right\}.
$
Since
 $\o \in  B_+(f,3,Q)\subseteq B_+(f,1,Q)$, we have
$0\le a_0 < b_-(a_0)\equiv b_0\le Q/\e$, and $(a_0,b_0]$ is an
elongation. Also, $(a_-(b_0),b_0]\supseteq (a_0,b_0]$
is an elongation with the same sign. To fix ideas we assume
$+1={\rm sgn}(a_0,b_0]$. This
will
serve as starting point for the
construction.
We now set, for $b<b_0$:
$$\eqalign
{\tilde a_+(b)&= \sup \{a < b\colon (a,b] \in \DD(f,-)\},\cr
b_{-1}&=\sup\{b<b_0 \colon \tilde a_+(b)>-\infty\}, \text{ and }\quad a_{-1}=
\tilde a_+(b_{-1}).\cr
} \Eq(3.7)
$$
Since $\o \in B_-(f,3,Q)\subseteq B_-(f,2,Q)$ we have 
$-Q/\e \le a_{-1} < b_{-1}$, and
from the construction, we easily check  $a_{-1}<0$. 
Observe that in \equ(3.7) we need to consider $b<b_0$ 
(instead of $b \le a_0$) due to the possibility
of non-empty overlap among elongations with different signs.
 We make the following:

\noindent {\bf Claim 1.} If $\o \in (\PP_1(f,3,Q) \cup
\PP''_2(f,Q))^c$  we have $b_{-1}\ge a_-(b_0)$.
\smallskip

\noindent {\bf Proof of Claim.} We prove it by contradiction. For that, we 
suppose that  $b_{-1} < a_-(b_0)$, and  
consider two cases:

\noindent (I) $\YY_{\a} \le \YY _{a_-(b_0)}$ for some $\a \in
[-Q/\e, a_-(b_0))$;

\noindent (II) $\YY_{\a} > \YY _{a_-(b_0)}$ for all $\a \in
[-Q/\e, a_-(b_0))$.
\smallskip
\noindent In case (I), letting $\a_0=\max\{\a< a_-(b_0)\colon
\YY_{\a} \le \YY_{a_-(b_0)}\}$, we take: $\a_3$ any point of
(global) minimum of $\YY_{\cdot}$ in $[a_-(b_0),b_0]$;
$\a_4=\min\{\a\in[\a_3,b_0]\colon \YY_{\a}-\YY_{\a_3}\ge
2\FF^*+f\}$, which exists since ${\rm sgn}(a_-(b_0),b_0]=+1$;
$\a_2=\max\{\a \in [\a_0,\a_-(b_0)]\colon
\YY_{\a_3}-\YY_{\a}<-2\FF^*+f\}$, which exists in this case,
otherwise $(\a_0,b_0]$ would be a { positive}  elongation, 
contradicting the  definition of  $a_-(b_0)$.

We see that starting from $\a_2$ and
moving backwards in time, the process $\YY$ must take  a value
below $\YY_{\a_2} -2\FF^*+3f$  before it reaches a value above
$\YY_{\a_2}+2f$ (otherwise $b_{-1} \ge a_-(b_0)$)); taking $\a_1$
as the ``first" (backwards) such time, we are in the situation
described in $\PP'_2(f,Q)$, contradicting our assumption on $\o$.

\smallskip
\noindent In case (II), let $\a_4$ be any point of minimum of
$\YY(\cdot)$ in $[a_-(b_0),b_0]$. Due to the assumption that $\o
\in B_-(f,3,Q)$, there exists a positive elongation contained in
$[-Q/\e,a_-(b_0)]$. 
Together with the assumption in (II) this
allows to define $\a_1=\max\{\a<a_-(b_0)\colon
\YY_{\a}\ge  \YY_{\a_4}+2\FF^*+f\}$, and $-Q/\e \le \a_1 <a_-(b_0)$.
Taking  $\a_3=\sup\{\a<\a_4\colon \YY_{\a}-\YY_{\a_4}\ge 2\FF^*-f\}$
{ which exists otherwise $[\a_1,\a_4]$ would be a negative 
elongation contradicting $b_{-1}<a_-(b_0)$}.
{ Moreover $\a_3\ge\a_1$.}
We see that starting {from}  $\a_3$ and moving ``backwards"
in time, $\YY_{\cdot}$ has to make a downwards increment of at
least $2\FF^*-3f$ ``before" $\a_1$ [otherwise $b_{-1}\ge a_-(b_0)$].
and we get  $\a_2$ as the
``first" such time, { we are in the situation described in 
$\PP'_2(f,Q)$, contradicting our assumption on $\o$}.

\medskip
Having assumed that $\o \in (\PP_1(f,3,Q) \cup \PP''_2(f,Q))^c$ in
this construction, the previous claim tells us that $b_{-1}\ge
a_-(b_0)$. For ${\rm sgn} (a_0,b_0]=+1$ we  define
$$
\a^*_0= \min\{\a \in [a_-(b_0),b_{-1}]\colon\quad  \YY_{\a}=
\min_{a_-(b_0)\le \tilde \a \le b_{-1}} \YY(\tilde \a) \},\Eq(3.8)
$$
In this situation  $(a_{-1},\a^*_0]$ and $(\a^*_0,b_0]$ are contiguous
elongations, with alternating signs ($-1$ and $+1$ resp.). The same holds for
$(a_-(\a^*_0),\a^*_0]\supseteq(a_{-1},\a^*_0]$ and
$(\a^*_0,b_+(\a^*_0)]\supseteq (\a^*_0, b_0]$.

\noindent {\bf Remark.} Though not needed, one can check that
$\YY_{\a^*_0}= \min_{a_{-1} \le \a \le b_0} \YY_{\a}$.
\medskip

With $\o \in  (\PP_1(f,3,Q) \cup \PP''_2(f,Q))^c$ we may  proceed
one step to the right, where the next ``breaking point" will be a
maximum in a suitable interval. We first set, for $a>\a^*_0$:
$$
\eqalign
{\tilde b_-(a)&= \inf \{b>a\colon (a,b] \in \DD(f,-)\}\cr
a_1=& \inf \{a>\a^*_0 \colon \tilde b_-(a)<+\infty \}, \text 
{ and} \quad b_1=\tilde b_-(a_1)
\cr} \Eq(3.88)
$$
and since $\o \in B_+(f,3,Q)\subseteq B_+(f,2,Q)$ we have $0< a_1<
b_1\le Q/\e$. Moreover, as before we have:

\noindent {\bf Claim 2.} For $\o \in (\PP_1(f,3,Q) \cup
\PP''_2(f,Q))^c$  we must have $a_1\le b_+(\a^*_0)$.

Claim 2 is proven in the same way as the previous one, and we omit details.
It allows to define, for such $\o$:
$$
\a^*_1= \min  \{\a \in [a_1,b_+(\a^*_0)] \colon \YY_{\a}=
\max_{a_1 \le \tilde \a \le b_+(\a^*_0)}\YY_{\tilde \a}\}\Eq(3.89)
$$
so that $(\a^*_0,\a^*_1]$, and $(\a^*_1, b_1]$ are contiguous elongations
with alternating signs ($+1$ and $-1$ resp.). Also  ${\rm sgn }(\a^*_1,
b_+(a_1)]={\rm sgn} (\a^*_1, b_1]$,
and, similarly to previous observation, we see that
$
\YY_{\a^*_1}= \min_{a_{0} \le \a \le b_1} \YY_{\a}.
$

If $\a^*_0< 0$ we set $J(\o)=(\sfrac{\e\a^*_0}{\g},\sfrac{\e\a^*_1}{\g})$.
If instead, $\a^*_0\ge 0$, in order to
determine $J(\o)$
we need to extend the construction one more step to the left. In this case,
we may consider for any $b<\a^*_0$:
$$\eqalign
{\tilde a_+(b)&= \sup \{a < b\colon (a,b] \in \DD(f,+)\},\cr
b_{-2}&=\sup\{b<\a^*_0 \colon \tilde a_+(b)>-\infty\}, 
\text{ and }\quad a_{-2}= \tilde a_+(b_{-2}).\cr
}\Eq(3.90)
$$
Since $\a^*_0\ge 0$, ${\rm sgn}(a_-(\a^*_0),\a^*_0]=-1$, and  $\o
\in B_-(f,3,Q)\subseteq B_-(f,2,Q)$ we have $-Q \le b_{-2} \le
\a^*_0$ and $-Q \le a_{-2}$. Moreover, from the construction
$a_{-2}<a_-(\a^*_0) \le a_{-1}$. As before, we can prove the
following:

\noindent {\bf Claim 3.} For $\o \in (\PP_1(f,3,Q) \cup
\PP''_2(f,Q))^c$  we must have $b_{-2}\ge a_-(\a^*_0)$.

The proof of Claim 3 is omitted, since it follows the same argument of
Claim 1, under the previous assumptions.
Having $b_{-2}\ge a_-(\a^*_0)$ we may split the intervals through
$$
\a^*_{-1}=\inf\{ \a \in [a_-(\a^*_0),b_{-2}]\colon
\YY_{\a}=\max_{a_-(\a^*_0) \le \tilde \a \le b_{-2}} \YY_{\tilde
\a}\} \Eq(3.91)
$$
so that $(a_{-2}, \a^*_{-1}]$ and $(\a^*_{-1},\a^*_0]$ are
elongations with alternating signs. As in the previous steps, we
see that $b_{-2}< a_-(\a^*_0)$  is not possible if $\o \notin
\PP'_2(f,Q)$. Moreover, from the construction it follows that
$\a^*_{-1}<0$, otherwise it would contradict the definition of
$a_0$ and ${\rm sgn} (a_0,b_0]=+1$. Thus, for  $\a^*_0\ge 0$ we
set $J(\o)=(\sfrac{\e\a^*_{-1}}{\g},\sfrac{\e\a^*_0}{\g})$. Though
not used in the sequel, we may again check that,
$
\YY_{\a^*_{-1}}= \min_{a_{-2} \le \a \le b_{-1}} \YY_{\a}.
$

\medskip
Under the assumptions on $\o\in (\PP(f,3,Q)\cup \PP''_2(f,Q))^c$
we have constructed contiguous elongations $(a_{-2},\a^*_{-1}]$,
$(\a^*_{-1},\a^*_0]$, $(\a^*_0,\a^*_1]$, and $(\a^*_1,b_1]$, with
alternating signs. 

Starting from $(a_-(\a^*_{-1}),\a^*_{-1}]$ and
$(\a^*_1,b_+(\a^*_1)]$, the construction may be continued to the
left and right respectively, if $\o \notin \PP_1(f,k,Q) \cup
\PP''_2(f,Q)$ for larger $k$. For Theorem \equ(t2) it
suffices to
have $\o \in (\PP_1(f,3(2k+1),Q) \cup \PP''_2(f,Q))^c$.
\medskip
\noindent {\bf Remark.} We have chosen $\a^*_0$, $\a^*_1$, etc...
as the first minimizer or maximizer, respectively, since the
random walk may have multiple maximizers on the intervals
considered  there. In fact the random walk can oscillate, being
always below or equal to the maximum. 
Since in the limit $\e \downarrow 0$,
the random walk  converges in law to a Brownian motion
where  the local maxima are always distinct, see
[\rcite{RY}] p. 108, we can expect that for a random  walk such a
result holds approximately. 
A way to do it is to 
accept an error 
on the location of the beginning or the end of the runs of
$\eta^{\d,\z}(\ell)$. For this we need to prove that if $\a_1$ and
$\a_2$  are the locations of  two  local maxima of $\YY(\cdot)$
and the distance between $\a_1$ and $\a_2$ is larger than
$\rho/\e$, then  $\P[|\YY_{\a_1}-\YY_{\a_2}|\le \tilde \d]$ goes
to zero in the limit $\e\downarrow 0$, for a suitable choice of
the parameters $\rho=\rho(\e)$, $\tilde \d=\tilde
\d(\rho,\e)=\tilde \d (\e)$ 
both vanishing as $\e \to 0$.

We define, for $\rho$ and $\tilde \d$ positive,
$$
\eqalign{ \PP_2(f,+,Q,a_{-1},b_0,\rho,\tilde \d) \equiv \{ &\o\in
(\PP_1(f,3,Q)\cup \PP''_2(f,Q))^c ;\, \exists\, \tilde \a \in
[a_{-1},b_0],\cr &|\tilde \a -\a^*_0|>\rho/\e, |\YY_{\tilde
\a}-\YY_{\a^*_0}|\le \tilde \d\}},\Eq (3.9)
$$
$$
\eqalign{ \PP_2(f,+,Q,a_{0},b_1,\rho,\tilde \d) \equiv \{ &\o\in
(\PP_1(f,3,Q)\cup \PP''_2(f,Q))^c ;\, \exists\, \tilde \a \in
[a_{0},b_1],\cr &|\tilde \a -\a^*_1|>\rho/\e, |\YY_{\tilde
\a}-\YY_{\a^*_1}|\le \tilde \d\}},\Eq(3.11)
$$
and
$$
\eqalign{ \PP_2(f,+,Q,a_{-2},b_{-1},\rho,\tilde \d) \equiv \{
&\o\in (\PP_1(f,3,Q)\cup \PP''_2(f,Q))^c ;\,\a^*_0>0, \exists\,
\tilde \a \in [a_{-2},b_{-1}],\cr &|\tilde \a -\a^*_{-1}|>\rho/\e,
|\YY_{\tilde \a}-\YY_{\a^*_{-1}}|\le \tilde \d\}}\Eq(3.11')
$$

We will show that  
the previous three sets have $\P$-probability
as small as we want provided we choose the parameters $\e$,
$\rho$, $\tilde \d$ in a suitable way. 

We recall that we have defined the random interval $J(\o)$ as follows:
$$
J(\o)=\cases{\left(\frac{\e\a_0^*}{\g},\frac{\e\a_1^*}{\g}\right),&if
 $\a_0^*<0$;\cr
\left(\frac{\e\a_{-1}^*}{\g},\frac{\e\a_0^*}{\g}\right),& if $\a_0^*\ge 0$.\cr}
\Eq(inter)
$$

There
is some arbitrariness when $\a^*_0=0$, but accepting to make an error
$\rho/\e$ on the location of the maximizers or minimizers, 
we will show that the set 
 $$
\PP_3(f,Q,\rho) 
\equiv \left\{\o\in (\PP(f,3,Q)\cup
\PP''_2(f,Q))^c; \a^*_0\, {\rm or}\, \a^*_{-1} \in
[-2\frac{\rho}{\e},2\frac{\rho}{\e}] \right\} \Eq(3.12)
$$
has a very small probability.

\medskip
\noindent {\bf Remark.} Always assuming $\o\in (\PP(f,3,Q)\cup
\PP''_2(f,Q))^c$, but instead ${\rm sgn}(a_0,b_0]=-1$, we perform
the obvious modifications of the construction.

Recalling that all over this work, $\b>1$ and 
$\th>0$ satisfy \eqv(2.17), 
the control on the various exceptional sets is summarized in the
following:

{\bf \Theorem(T1)}{\it
There exist  positive constants $Q_0=Q_0(\b,\th)$,
$f_0=f_0(\b,\th)$, $\rho_0=\r_0(\b,\th)$ and 
$\g_0=\g_0(\b,\th)$ such that for all 
$ 0 < \g\le \g_0$, $0<\r\le \r_0$, and $0<f\le f_0$, 
for all $\e$ such that 
$$
\d^*\g<\e\le \frac{2}{V^2(\b,\th)\log (1944)}
\left( \rho^{4+2a}\wedge f^2\right)
\Eq(Zorro.1)
$$
for an arbitrary  given $a>0$, 
we have the following: For all
integers $k>1,Q\geq Q_0(\b,\th)$,
$$
\P\left[ \PP_0\big(f,Q\big)\right]
  \le 3e^{-\frac{Q}{2C_1}}
+\frac{1}{\log 2}\frac{ 2f+ 9V(\b,\th)\sqrt{\e \log \frac
{C_1}\e}}{2\FF^*-f} 
\log  \frac {2\FF^*-f} 
{2f+2V(\b,\th)\sqrt{ \e \log \frac {C_1} \e}}
\Eq(3.14)
$$
where $V(\b,\th)$ is given by \equ(2.1701) and $C_1=C_1(\b,\th)$
is given in \eqv(3P.50) with $b=2\FF^*$;
$$
\P \left [  \PP_1\left(f,k, Q \right) \right ] \le
(k+5)e^{-\frac{Q}{2kC_1}}
+\frac{k}{\log 2} 
\frac{ 2f+ 9V(\b,\th)\sqrt{\e \log \frac
{C_1}\e}}{2\FF^*-f} \log  \frac {2\FF^*-f} 
{2f+2V(\b,\th)\sqrt{ \e \log \frac {C_1} \e}}
\Eq(3.141)
$$
$$
\eqalign{ \P  \left[  \PP''_2(f,Q)\right] & \le
8(2Q+1)^2\frac{2\sqrt{2\pi}}{V(\b,\th)} (9f)^{a/(2+a)} +
(2Q+1)\frac{1296}{V(\b,\th)}  \frac{9f+(2+V(\b,\th))
\sqrt{\e\log \frac{C_1}{\e}}}{(9f)^{3/(4+2a)}} \cr
&\quad +\frac{4Q}{\e} e^{-\sfrac{f}{4\e
V^2(\b,\th)}}.}\Eq(3.1411)
$$
 Moreover, 
for   $\tilde\d(\rho)=\rho^{2+a}$ 
we have 
$$
\P\left[ \cup_{i=-k}^{k}\cup_{s_1\in \{\pm 1\}}
\PP_2(f,s_1,Q, a_i,b_{i+1},\rho,\tilde\d(\rho))\right] \le
(4k+2) 3G_1(\b,\th,\tilde \d(\r),\e)\log 
\frac{4}{G_1(\b,\th,\tilde \d(\r),\e)}
\Eq(3.15)
$$
where
$$
G_1(\b,\th,\tilde \d(\r),\e)\equiv
\frac{2^{16}C_1}{\sqrt{V(\b,\th)}}
\left(\r^{a/2}+
\frac{\sqrt{1+V(\b,\th)}(\e\log \frac{C_1}{\e})^{1/4}}
{\r^{3/4}}\right)
\Eq(3PP.16)
$$
with  $C_1$ as in \eqv(3.14), and if $0<\kappa <1/2$}
$$
 \P\left[ \PP_{3}(f,Q,\rho)\right] \le
6 \r^{\sfrac{1}{2}-\kappa}+
\frac{2}{\G(\frac{1}{2}-\k)}
\big(\frac{\e}{\rho}\big)^{\sfrac{1}{2} -\kappa} +\frac
{\e}{\rho^2} \exp\Big(8^{\sfrac{C(\b,\th)}{\kappa^2}} 2 \log
\frac{C(\b,\th) }{\kappa^2}\Big)
\Eq(3.16)
$$
{ \it where $C(\b,\th)$ is a suitable constant that depends on $V(\b,\th)$
and $\G(\cdot)$ is the Euler Gamma function.
}

The proof will be given at the end of this section. 

\remark The quantities $a_i$ and $b_i$ are random variables, but
none is a stopping time. 
As   $\e \downarrow 0$, and then $ \rho
\downarrow 0$ \eqv(3.15) reduces to  the well known fact that with
probability one, the Brownian path does not have two equal local
maximum (or minimum) over any finite interval (see [\rcite{RY}]
pg 108). 

To simplify the writing of the above estimates, we made 
the following choice:
$$
\r=\e^{\frac{1}{4(2+a)}},\quad
f=\e^{\frac{1}{4}},\quad
\kappa=1/4.
\Eq(Zorro.2)
$$
Then, calling
$$
\eqalign{
\PP(k,\e,Q)
=&\PP_1(f=\e^{\frac{1}{4}},k, Q)
\cup
\PP''_2(f=\e^{\frac{1}{4}},Q)
\cup
\PP_{3}(f=\e^{\frac{1}{4}},a_{-2},b_{-1},\r=\e^{\frac{1}{4(2+a)}})\cr
&\cup
\left(\cup_{i=-k}^{k}\cup_{s_1\in \{\pm 1\}}
\PP_2(f=\e^{\frac{1}{4}}
,s_1,Q, a_i,b_{i+1},\rho=\e^{\frac{1}{4(2+a)}}
,\tilde\d(\rho)=\e^{\frac{1}{4}})\right),\cr
}\Eq(Zorro.3)
$$
after simple estimates one gets

{\bf \Corollary(zorro)}{\it 
There exist  positive constants $Q_0=Q_0(\b,\th)$, 
$\g_0=\g_0(\b,\th)$ and $\e_0(\b,\th)$ 
such that for all $ 0 < \g\le \g_0$,   
for all $\e$ that satisfies 
$ \d^*\g<\e\le \e_0$, for all $Q>Q_0$, $k>1$ we have
$$
\P[\PP(k,\e,Q)]\le 
(k+5)e^{-\frac{Q}{2kC_1}}+k\e^{\frac{a}{16(2+a)}}
+Q^2\e^{\frac{a}{8+2a}}+ Q e^{-\frac{1}{2\e^{3/4}V^2(\b,\th)}}
.
\Eq(Zorro.4)
$$
where $a>0$ is a given arbitrary positive number.
}

Recalling \eqv(inter), 
the following Proposition will be used for proving \eqv(2.210a)
and \eqv(2.210b). It will be proved at the end of this section.

{\bf \Proposition(Terminator)}{\it 
For all $0<x<(\FF^*)^2/(V^2(\b,\th)18\log 2)$
we have
$$
\P[\g |J|\le x] \le 2e^{-\frac{(\FF^*)^2}{18 x V^2(\b,\th)}}
\Eq(Terminator1)
$$
while for all $x>0$ we have 
$$
\P[\g|J|\ge x]\le 4e^{-\frac{x}{8C_1(\b,\th,\FF^*)}(1-\frac{\log 3}{\log 4})}.
\Eq(Terminator2)
$$
where $C_1(\b,\th,\FF^*)$ is defined in \eqv(3P.50).
}

\remark Note that for $x\ge (\FF^*)^2/(V^2(\b,\th)18\log 2)$
the right hand side of \eqv(Terminator1) is larger than 1. Therefore
\eqv(Terminator1) is trivially satisfied also in this case.

\medskip
\vskip0.5cm
\noindent{\bf Basic estimates.}

Several  probabilistic  estimates   are needed for Theorem
\eqv(T1) and are summarized in the following  Lemmata and Proposition.
The variables $\chi(\a), \,\a \in \Z$ defined by \equ(2.174), with
$X(x)$ given by \eqv(2.170), constitute the basic objects in the
following analysis. We recall that we always  assume that 
$\b>1$ and $\th>0$ satisfy \eqv(2.17). 
Recalling \equ(2.21) we set 
$$
V^2_-=  V^2(\b,\th)\Big(1- (\sfrac{\g}{\d^*})^{1/5}\Big)^2
\quad{\text{and}}\quad V^2_+=V^2(\b,\th)
\Big(1+(\sfrac{\g}{\d^*})^{1/5}\Big)^2.
 \Eq(4.22001)
$$
\medskip
\remark Throughout this section we shall assume that
$0<\g/\d^*\le d_0(\b,\th)\wedge 2^{-5}$ so that $V(\b,\th)/2\le
V_-\le \sqrt{c(\b,\th,\g /\d^*)}\le V_+ \le 3V(\b,\th)/2$ where
$V(\b,\th)$ is given in \equ(2.1701).

 We need some further simple estimates  concerning the  variables
$\chi(\a)$ that are not difficult  to prove just recalling 
that $\chi(\a)$ is a sum 
over $\e(\g\d^*)^{-1}$ independent symmetric random variables $X(x)$. 
\eqv(4.13) is proved using \eqv(4.9).

\noindent{\bf \Lemma (lem1)} {\it 
There exists a $d_0(\b,\th)>0$, such that if $\g/\d^*
\le d_0(\b,\th)$ then
$$
\E\left[ e^{\l\chi(\a)}\right] \le e^{\sfrac{\l^2}{2}\e V^2_+},
\quad \quad \forall \l \in \R \Eq(4.3)
$$
with $V^2_+$ defined in \equ(4.22001).
If $\,0<\l < [\e V^2_+ ]^{-1}$,
we have
$$
\E\left[ e^{\sfrac {\l}{2} |\chi(\a)|^2}\right] \le  \frac{ 1 }{
1-\e \l V^2_+}.\Eq(4.9)
$$
For all $k\ge 3$ and  $p=1,2,4$:
$$
\E\left[ \max_{\a=1,\dots,k}|\chi(\a)|^p \right] \le  (4 \e V^2_+
\log k)^{p/2} ( 1+\frac {p}{\log k})^{\sfrac{p}{2}\vee 1}.
\Eq(4.13)
$$
}
\medskip

In order to have an elongation, as previously described,
it is necessary  to find suitable  uphill or downhill increments of
height $2\FF^*+{f}$.  

A constructive way to locate
elongations, though it might miss some of them,  is
related to the following stopping
times:

Given $b>0$ ($b=\FF^*+\frac f 2$ later),  we set
$\tau_0=0$, and define, for $k\ge 1$:
$$
\eqalign{
\t_k&=\inf \{ t>\t_{k-1}\colon
|\sum_{\a=\t_{k-1}+1}^{t} \chi(\a)| \ge b \},\cr
\t_{-k}&=
 \sup \{ t<\t_{-(k-1)}\colon
|\sum_{\a=t+1}^{\t_{-(k-1)}}\chi(\a)| \ge b\}.}
\Eq(4.20)
$$

Clearly, the random variables $\D\t_{k+1}:= \t_{k+1}-\t_{k}$, $k\in
\Z$, are
independent and identically distributed.
(Recall that $\D\t_1=\t_1$ from the definitions.)
We define,
 $$ S_k= {\rm sgn} \Big (
\sum_{j=\t_{k-1}+1}^{\t_{k}}
\chi(j)\Big); \qquad
 S_{-k}= {\rm sgn} \Big ( \sum_{j=\t_{-k }+1}^{\t_{-k+1}}
\chi(j)\Big)\qquad \text{for  } k \ge 1
\Eq(4.23)
$$

We need probabilistic estimates for the  variables $\D\t_k$ and $\t_k$,
which are obtained by standard methods. An upper bound on the
tail of their distribution can be given as follows:

\medskip
\noindent{\bf \Lemma (lem41)}
{\it 
There exists a  positive constant
$d_0(\b,\th)$ such
that for all  integer $v$, $\g/\d^*<d_0(\b,\th)$ and  
$0<\e<\e_0(\b,\th, b)$  
where 
$$
\e_0(\b,\th,b):=\frac{1}{3^8}\left(\P\big[ Y\ge \frac{ 4b }
{V(\b,\th)}\big]\right)^{2},
\Eq(4.35003)
$$
we have
$$
\P\left[ \t_1\ge \frac {v}{\e} \right] \le
 \exp\left( -v\P\left[ Y\ge \frac{
4b}{V(\b,\th)}\right] \right), \Eq(4.2401)
$$
where $Y$ is standard Gaussian and $V(\b,\th)$
as in \equ(2.1701).}

\remark
{ For future use, note that
$\e_0(\b,\th,b)$ is a decreasing function of $b$.}

\proof
Since the $\chi(\a)$ are i.i.d. random variables,
for any positive integer $v$, we have:
$$
\P\left[\t_1\ge \frac {v}{\e} \right]
\le \P\left[ \max_{k=0,\dots,v-1}
\left|\sum_{\a=k/\e+1}^{(k+1)/\e}\chi(\a)\right|<2b \right]
=\left( \P\left[
|\YY(1/\e)| \le
2b \right]\right)^v
\Eq(4.35001)
$$

We can use  \eqv(4.3) to get an estimate of the fourth moment of
$\chi(\a)$ and  apply
Berry--Essen Theorem ([\rcite{CT}] p. 304)  to control the
right hand side in \eqv(4.35001).
Consequently,  there exists a
constant $C_{BE}=C_{BE}(\b,\th)$  which, according to Berry-Essen
inequality may be taken as
$$
C_{BE}
=0.8 \sup_{0<\g/\d^* \le d_0(\b,\th),\e>\d^*\g}
\E(|\chi(1)|^3)/\E(|\chi(1)|^2)^{3/2}
\le 3^4
\Eq(4.35002)
$$
assuming at the last step that $\g/\d^*\le d_0(\b,\th)<(1/2)^5$.
Therefore
$$
\P\left[
|\YY(1/\e)| \le
2b\right]
\Big]
\le 1 -2\P\big[ Y\ge \frac{2 b}{\sqrt{c(\b,\th,\g/\d^*)}}\big]  +
  3^4\sqrt{\e} \le 1-\P\big[ Y\ge \frac{ 4b}{V(\b,\th)}\big]
\Eq(4.360)
$$
where $Y$ is a standard Gaussian, using $0<\e<\e_0(\b,\th,b)$
and  \eqv(4.35003)
for the last inequality in \equ(4.360). 
Using $1-x \le e^{-x}$,  we get \eqv(4.2401) \eop

\smallskip

The following lemma gives bounds  for the mean of $\t_1$ and follows
easily from the Wald Identity,  see [{\rcite{Ne}}], pg 83, and \eqv(4.13).

\noindent{\bf \Lemma (lem42)} {\it If
$$
C_1=C_1(\b,\th,b)=\frac{2}{\P[Y>4b/V(\b,\th)]},
\Eq(3P.50)
$$
where $Y$ is standard gaussian  and $0<\e<
\e_0(\b,\th,b)$ cf. \eqv(4.35003), {there exists $d_0(\b,\th)$ 
such that for $\g/\d^*<d_0(\b,\th)$} we have}
$$
\frac{b^2}{\e V^2(\b,\th)}(1-(\g/\d^*)^{1/5})^2 \le \E[\t_1]  \le
\frac {b^2}{\e V^2(\b,\th)} (1+(\g/\d^*)^{1/5})^2 
\left(1 + 9\frac{V(\b,\th)}{b} \sqrt{\e\log \frac {C_1}{\e}} \right)^2.
\Eq(4.371)
$$
\medskip
\remark For future use, note that $C_1(\b,\th,b)$ is increasing with $b$.

We need  exponential estimates for
the probability that  a Ces\`aro average over  $k$ terms of the
previous $\D\t_i$'s is outside an interval that contains
the mean $\E[\t_1]$. The result is:

\noindent {\bf \Lemma (lem4')} {\it For all $0<s<b^2 [4(\log 2)
V^2_+]^{-1}$, 
for all positive
integers $k$
we have
$$
\P\left[ \t_k \le \frac {k s}{\e} \right] \le  e^{-k
\frac{b^2}{4sV^2_+}}, \Eq(4.250000)
$$
where $V^+$ is defined in \eqv(4.22001). Moreover, for
$\e_0=\e_0(\b,\th,b)$ as \eqv(4.35003),  for all $0<\e <\e_0$,
for all positive
integers $k$, and for all $s>0$ we have
$$
\P\Big[ \t_k \ge \frac {k  }{\e}
(s +\log 2)C_1\Big]
\le   e^{- s k}
\Eq(4.3401)
$$
where $C_1=C_1(\b,\th,b)$ is given in \eqv(3P.50).}

\proof \eqv(4.250000) is an immediate consequence of the Markov 
exponential 
inequality together with the exponential Wald identity see [\rcite{Ne}], 
pg 81.
\eqv(4.3401) is an immediate consequence of the Markov exponential
 inequality
together with \eqv(4.2401) to estimate the Laplace transform.\eop

\medskip
\smallskip
As we shall check, the above stopping times with $b=\FF^*+\sfrac
f2$, provide a simple way to catch elongations. 
It will be enough 
to find successive indices $k\ge 1$ ($k \le -2$) such that
$S_k=S_{k+1}$ and eliminating a set of small
probability, see Lemma \eqv(lem111), 
$(\t_{k-1},\t_{k+1}]$ ($(\t_{k},\t_{k+2}]$
respectively) will provide an elongation  which is positive if
$S_k=+1$, or negative otherwise. 
Still, if $S_{-1}=S_1$, then $(\t_{-1},\t_1]$
is an elongation. Not all elongations are of this form, as one
simply verifies, but what matters is that this procedure catches 
enough of them, sufficient to prove Theorem \equ(T1).
The basic ingredient is given  in the next two lemmas.

\smallskip
\noindent{\bf \Lemma (lem6)}
{\it Let  
$\e_0=\e_0(\b,\th,b)$ be given by \eqv(4.35003).
For  all $0<\e <\e_0$,
all integer $k\ge 1$, and all $s>0$ we have}
$$
\P\left[\t_k \le \frac { k (s +\log 2)C_1 }{\e}
; \exists i\in\{1,\dots,k-1\}, S_i=S_{i+1}\right]
\ge \left(1- e^{- sk}\right)(1-\frac 1{2^{k-1}}).
\Eq(4.43)
$$

\proof
It follows at once from the fact, due to the symmetry, 
that  conditionally on $
\D\t_{i}$'s the
variables $ S_i, i \ne 0$'s form a family of
i.i.d.
Bernoulli symmetric random variables (see \equ(4.23)), with the
trivial observation that for i.i.d. symmetric Bernoulli random variables
$$
\P\left[\exists i \in \{1, \dots,k-1\}  
\colon S_i=S_{i+1}\right]=1-\frac{1}{2^{k-1}}. 
\Eq(4.45)
$$
Together with \eqv(4.3401), this entails \eqv(4.43).\eop
\medskip

To deal with  the case where more than one elongation is involved,
 we define to the right of the origin
$$
\eqalign{ i^*_1&\equiv \inf \left\{i\ge 1: S_i=S_{i+1}\right\}\cr
i^*_{j+1}&\equiv\inf\left\{i \ge (i^*_j+2): S_i=S_{i+1}=-S_{i^*_j}
\right\} \quad \qquad  j \ge 1,}
\Eq(4.46)
$$
and  to the left
$$
i^*_{-1} \equiv
\cases{
-1 & if
$S_{-1}=S_1=-S_{i^*_1}$,\cr 
 \sup \left\{i\le -2: S_i=S_{i+1}=-S_{i^*_1}\right\} & if
$S_{-1}\neq S_1$ or $S_1=-S_{i^*_1}$,\cr }
$$
$$
i^*_{-j-1}\equiv\sup\left\{i \le i^*_j-2: S_i=S_{i+1}=-S_{i^*_j} 
 \right\} \quad  \qquad  j \ge 1,
\Eq(4.46')
$$
we then have:

\noindent{\bf \Lemma(lem7)}
{\it
Let  
$\e_0=\e_0(\b,\th,b)$ be given by \eqv(4.35003).  For  all $0<\e <\e_0$,
all $k$ and  $L$ positive integers,  $L$ even,
(just for simplicity of writing) and all $s>0$ we have:}
$$
\P\left[ \t_{kL-1} \le \frac{ (kL-1)(s
+\log2)C_1}
{\e },\forall_{1\le j \le k}\,\, i^*_j < jL\right]
\ge \left(1-e^{-s (kL-1)}\right)\left(1-\sfrac{1}{2^{L-1}}\right)
\left( 1-\left(\sfrac34 \right)^{L/2}\right)^{k-1}    
\Eq(4.48)
$$
{\it and}
$$
\eqalign{
&
\P\left[ \t_{-kL} \ge \frac{ -kL(s
+\log2)C_1}
{\e},\, 
\t_{L-1} \le \frac{ (L-1)(s
+\log2)C_1}
{\e}, \, i^*_1 <L, \,
\forall_{1\le j \le k}\,\, i^*_{-j} > -jL\right] 
\cr&\quad 
\ge \left(1-e^{-s (kL-1)}\right)\left(1-\sfrac{1}{2^{L-1}}\right)
\left( 1-\left(\sfrac34 \right)^{L/2}\right)^{k}.
    \cr
}\Eq(4.480)
$$

\proof We prove \eqv (4.48);  \eqv (4.480) is done similarly.
We again use   that   conditionally on $
\D\t_{i}$'s,
  the variables $ S_i$'s  are i.i.d. Bernoulli symmetric random
variables.
Recalling Lemma
\eqv (lem4'),  it is  then sufficient   to prove   that
$$
 \P[ i^*_1 < L, i^*_2<2L, \dots,i^*_k<kL] \ge
\left(1-\sfrac{1}{2^{L-1}}\right)
\left( 1-\left(\sfrac34 \right)^{L/2}\right)^{k-1}.  \Eq(4.4800)
$$
When $k=1$ this is just \equ(4.45). On the other side, using the
above mentioned properties of the
random variables $S_i$ we easily see that
$$
\P[i^*_{j+1} -i^*_j \le L \mid i^*_1,\dots,i^*_j]
\ge 1-\left(\sfrac34 \right)^{L/2} \quad {\text{a.s.}}
$$
from where \equ(4.48) follows at once.
\eop
\medskip
\smallskip
Next we  verify that the above described method
provides elongations, with overwhelming probability.
Recalling \eqv(4.46) let us assume, to fix ideas, that  
$S_{i^*_1}=S_{i^*_1+1}=1$.
 From the definition of $\t_i$, see \eqv(4.20), 
{with $b=\FF^*+(f/2)$},     we have that
$$
\YY((\t_{\{i^*_1-1\}},\t_{\{i^*_1+1\}}])=
\sum_{\a=\t_{\{i^*_1-1\}}+1}^{\t_{\{i^*_1+1\}}} \chi(\a) \ge 2\FF^*+f.
\Eq(4.97)
$$
Therefore $(\t_{\{i^*_1-1\}},\t_{\{i^*_1+1\}}]$ automatically
 satisfies one of the two
conditions  to give rise to  an elongation, cf.  \eqv (3.3).

Let us
see that, except on a set of small probability,
the  other requirement is fulfilled, i.e.,
$$
\inf_{\t_{\{i^*_1-1\}}< \a_1<\a_2\le \t_{\{i^*_1+1\}}}
\sum_{\a=\a_1}^{\a_2} \chi(\a) \ge -2 \FF^*+f.  \Eq (3.C2)
$$
On the event $\{S_i=1\}$, we readily see that
$$
\inf_{\t_{\{i-1\}}+1 \le \a_\le \t_i}
\sum_{\bar \a=\t_{\{i-1\}}+1}^{\a} \chi(\bar\a)\;\ge\; -\FF^*-f/2,
\,\hbox {and}\,
\inf_{\t_{\{i-1\}}+1 \le \a\le \t_i}
\sum_{\bar \a=\a}^{\t_i} \chi(\bar\a)\;\ge \,0.
 \Eq (3.C1)
$$
Since
$
\sum_{\a=\a_1}^{\a_2} \chi(\a)=
\sum_{\a=\a_1}^{\t_i} \chi(\a) +\sum_{\a=\t_i+1}^{\a_2}\chi(\a),
$
on $\{S_i=S_{i+1}=1\}$  we have
$$
\inf_{\t_{\{i-1\}}+1 \le \a_1\le \t_i<\a_2\le \t_{\{i+1\}}}
\sum_{\a=\a_1}^{\a_2} \chi(\a)\ge\;-\FF^*-f/2 \ge\; -2\FF^*+f.
\Eq(Term5)
$$
In the last inequality we used $ f<\FF^*/4<2\FF^*/3$. 
Therefore, 
it remains to evaluate $\P[\JJ(i^*_1)\cup \JJ(i^*_1+1),S_{i^*_1}=1]$, where
$$
\JJ(i):=\left\{\inf_{\t_{\{i-1\}}+1 \le \a_1<\a_2\le \t_{\{i\}}}
\sum_{\a=\a_1}^{\a_2} \chi(\a) < -2 \FF^*+f\right\}.\Eq(4.102)
$$
Note that on $\{S_{i}=1\}$, we have
$
\inf_{\t_{\{i-1\}}+1\le \a_1<\a_2 \le
\t_{i}}\sum_{\tilde \a=\a_1}^{\a_2} \chi(\tilde \a) \ge -2 \FF^* -f,
$
where we used \equ(3.C1) and 
$
\sup_{\t_{\{i-1\}}+1\le \a_1 \le \t_i}\sum_{\tilde
\a=\t_{i-1}+1}^{\a_1 -1}
\chi(\tilde \a) \le \FF^*
+\frac f2.
$
As a consequence,  for any integer $i$:
$$
\{\JJ(i),S_i=1 \}\subseteq \left\{-2\FF^*-f \le 
\inf_{\t_{\{i-1\}}+1\le \a_1<\a_2 \le
\t_{i}}\sum_{\tilde \a=\a_1}^{\a_2} \chi(\tilde \a)\le -2\FF^*+f\right\}.
$$
An analogous inequality (with a sup instead of an inf) holds in the case 
$S_{i^*_1}=-1$. Therefore we  need to prove the following:

\noindent{\bf \Lemma (lem111)} {\it 
Let  $\e_0=\e_0(\b,\th,2\FF^*)$ be given by \eqv(4.35003) and 
$C_1=C_1(\b,\th,2\FF^*)$ be given by\eqv(3P.50).  For all 
$0<f<\FF^*/4$ and for all $0<\e<\e_0$ we have
$$
\eqalign{ &
\P\left[\cup_{j=i^*_1,i^*_1+1}\left\{ 2\FF^* -f < 
\sup_{ \t_{j-1} <\a_1<\a_2\le
\t_{j}}
\left |   \sum_{\tilde \a=\a_1}^{\a_2} \chi(\tilde \a)\right | < 2\FF^* +f
\right\}\right]
  \cr &   \quad \quad \quad \quad\quad   
\le  \frac{2G (\b,\th, \e, f)}
{\log 2} \log
\frac{1}  {G(\b,\th, \e, f)  } 
}
\Eq(4.103)
$$
where}
$$
G (\b,\th, \e, f) \equiv \frac{ 2f+9V(\b,\th)\sqrt{ \e \log \frac
{C_1}{\e}}}{2\FF^*-f}. \Eq (3.C3)
$$
 \remark
Clearly $i^*_1$ is
anticipating, and
$\t_{i^*_1-1}$ and $\t_{i^*_1}$ are not stopping times.

\proof  Since $\P[i^*_1=i, S_{i^*_1}=1]= 2^{-i+1}$,  we have
$$
\P\big[\JJ(i^*_1), S_{i^*_1}=1\big]\le  \sum_{i=1}^{i_0}
\P\big[\JJ(i),S_i=1\big] + 2^{-i_0}
\Eq(4.105)
$$
where $i_0$  will be suitably chosen.
 To treat the sum, we define the stopping times
$$ T_{\FF^*-\sfrac {3f}{2} }=\inf\big\{\a> \t_{\{i-1\}}; \sum_{\tilde
\a=\t_{i-1}+1}^{\a} \chi(\tilde \a) \ge \FF^* -\sfrac {3f}{2} \big\}
\Eq(4.106) 
$$
$$ T_{\FF^*+\sfrac f2}=\inf\big\{\a> \t_{\{i-1\}}; \sum_{\tilde
\a=\t_{\{i-1\}}+1} ^{\a}
\chi(\tilde \a) \ge \FF^* +\sfrac f2 \big\}
\Eq(4.107)
$$
$$ T^-_{\FF^*-\sfrac {3f}{2}}=\inf\big\{\a> \t_{\{i-1\}}; \sum_{\tilde
\a=\t_{\{i-1\}}+1}
^{\a} \chi(\tilde \a) \le - \FF^* +\sfrac {3f}{2}\big\}
\Eq(4.108)
$$
By inspection we verify that
$
 \left\{\JJ(i),S_i=1\right\} \subseteq \SS (i)\equiv
\{T_{\FF^*-\sfrac {3f}{2}}\le T^-_{\FF^*-\sfrac {3f}{2}}\le
T_{\FF^*+\sfrac f2}\},
$
and by the strong  Markov property, we have
$$
\P\big[ \SS(i)\big] 
\le
\int_{\FF^*-\sfrac {3f}{2}}^{\FF^*+\sfrac f2}\,
\P\big[ \tilde T^-_{\FF^*-\sfrac {3f}{2}+x} < \tilde T_{\FF^*+\sfrac f2
-x}\big]
\P\big[\sum_{\a=\t_{i-1}+1}^{T_{\FF^*-\sfrac {3f}{2}}}\chi(\a)\in dx\big]
\le  \P\big[ \tilde T^-_{2\FF^*-3f} < \tilde T_{2f}\big]
\Eq(4.111)
$$
where, we have written
$
\tilde T_{x}\equiv \inf\big\{\a \ge 1\colon  \YY_{\a}\ge x\big\}$,
$\tilde T^-_x\equiv \inf\big\{\a \ge
1\colon\YY_{\a}\le-x\big\}.
$

At this point we need the estimate \eqv(4.76), in Lemma \equ(lem9)
below, it gives
$$
 \P\big[ \tilde T^-_{2\FF^*-3f} < \tilde T_{2f}\big]
\le \frac{ 2f +9V(\b,\th)\sqrt{ \e \log \frac{C_1}{\e}}}
{2\FF^*-f} \equiv G (\b,\th, \e, f). \Eq(4.113)
$$
{ with $C_1=C_1(\b,\th,2\FF^*)\ge (C_1(\b,\th,(2\FF^*-3f)\vee (2f))$
if $0<\e<\e_0(\b,\th,2\FF^*)\le\e_0(\b,\th,(2\FF^*-3f)\vee (2f))$.
Here  we have used that $\e_0(\b,\th,b)$ is decreasing with $b$ and that 
$C_1(\b,\th,b)$ is increasing with $b$.
}

Consequently, cf. \eqv(4.105),  \eqv(4.111) and \eqv(4.113) we have
$$ 
\P\big[\JJ(i^*_1),\,S_{i^*_1}=1\big] \le  \sum_{i=1}^{i_0}
\P\big[\SS (i)\big] + 2^{-i_0} 
\le  i_0  G (\b,\th, \e, f)
 + 2^{-i_0}
\Eq(4.105P)
$$
Taking $ i_0=  \log \frac 1 { G (\b,\th, \e, f)}  {[\log 2]}^{-1} $  
we obtain \eqv (4.103),
since the same works for $i^*_1+1$.
 \eop
\medskip

\medskip

To show that  \equ(3.15) holds, we need to bound  the probability
of finding two extrema in an interval $[\tau_{i^*_j},\tau_{i^*_{j+1}}]$, 
 at distance larger than $\rho/\e$
and whose values are within
$\tilde \d $.

 We fix the interval $[\tau_{i^*_{-1}},\tau_{i^*_1}]$
(the peculiarity of having fixed the origin will not bother), and 
for any  given
$h , k $  positive integers  we   denote
$$
\EE(k,h,+)= \left\{\o\in \O: i^*_{-1}=-h, i^*_1=k, S_k=-1 \right\},\Eq(4.51)
$$
where for   definiteness we are considering only the case of maxima,
i.e., we have assumed that $S_k=S_{k+1}=-1, 
S_{-h}=S_{-h+1}=+1$ on $\EE(k,h,+)$.
The case of minima is  similar.  Recall that
$\P \left [\EE(k,h,+) \right ] \le 2^{-(k+h)}$.

The positive integers  $h , k $ in \eqv (4.51) determine a random interval
$\{\t_{-h},\dots,\t_{k+1}\}\subseteq \Z$ in which  the 
index $\a$ of the variables
$\chi(\a)$ varies.  Using
 Lemma \eqv(lem4'), on a set of probability larger than $\left ( 1-e^{-s k}
\right ) \left ( 1-e^{-s h} \right )$, we can replace this random
interval  by a larger deterministic one. In particular, {assuming  
$s\ge \log 2$,} 
except for a set of probability at most $4e^{-s}$, {for all $h,k\ge 1$},
$\{\t_{-h},\dots,\t_{k+1}\}\subseteq \{\LL(-h,\e),\dots,
\LL(k+1,\e)\}$  where
$$
\LL(r,\e)\equiv  r\frac{(s+\log2)C_1}{\e} \quad r\in\Z
\Eq(4.512)
$$
 with $C_1=C_1(\b,\th,2\FF^*)\ge 
C_1(\b,\th,\FF^*+(f/2))$ as in \eqv(3P.50). 

We now  
partition the interval $[ \LL(-h,\e),\LL(k+1,\e)]$ into blocks of
length $\rho/\e$, where $\rho$ was already introduced in
\eqv(3.12). Assuming,  as always, that we do not have rounding off
problems, the number of such blocks inside $[
\LL(-h,\e),\LL(k+1,\e)]$ is 
$ \LL(k+1,\rho)-\LL(-h ,
\rho)$, i.e., of order $(k+h+1)\rho^{-1}$, with $\LL(\cdot,\rho)$
defined as in \equ(4.512) with $\e$ replaced by $\r$.

Given $\underline \a\equiv\LL(-h,\e) \le \a_1 \le \a_2 \le  \LL(k+1,\e)$, let:
$$
\YY^*(\underline \a,  \a_1,  \a_2)
\equiv
\max_{ \a_1 \le \tilde \a\le  \a_2}
\sum_{\a=\underline \a}^ {\tilde \a} \chi(\a). \Eq(4.52)
$$
Given $\tilde\d>0$,   $\rho>0$,
and $\ell$ such that $\LL(-h,\rho) \le
\ell\le \LL(k-1,\rho)$,
let us define  the event
$$
\eqalign{
\DD(k,h,\rho,\tilde\d,+,\e) &\equiv
\Big\{\o \in \O:
\exists \ell, \ell',\,\LL(-h,\rho) \le \ell<\ell'\le \LL(k-1,\rho) ;\cr
 &\quad\quad|\YY^*(\underline \a, \sfrac{\rho \ell}{\e}, 
\sfrac{\rho (\ell+1)}{\e})
  -
 \YY^*(\underline \a,  \sfrac{\rho \ell'}{\e},\sfrac{\rho
(\ell'+1)}{\e})| \le
 2\tilde\d
\Big\}.\cr}
\Eq(4.521)
$$
We now prove the following estimate:

\smallskip
\noindent  {\bf\Lemma(lem8)} {\it 
There exist positive  constants
$ \g_0 (\b,\th) $ and $\r_0(\b,\th)$ such that
for  all  $\g \le \g_0 (\b,\th)$, for $0<\r<\r_0(\b,\th)$, for
$\tilde \d=\r^{2+a}$ with $a>0$, 
for $\d^*\g <\e\le \e_0(\b,\th,\r)$, where
$$
\e_0(\b,\th,\r)=\sfrac{4(\r)^{2(2+a)}}{2V^2(\b,\th)\log (1944)},
\Eq(Term19)
$$
and  for all $s>0$
we have}    
$$
\eqalign{ &\P\left[ \cup_{k,h \ge 1}\left(\EE(k,h,+)\cap
\DD(k,h,\rho,\tilde \d,+,\e)\right)\right]
\le\cr  &\quad\quad \
\frac{2^{16}C_1(\b,\th,2\FF^*)}{\sqrt{V(\b,\th)}}(s+\log 2)
\left(\r^{a/2}+
\frac{\sqrt{1+V(\b,\th)}(\e\log \frac{C_1(\b,\th,2\FF^*)}{\e})^{1/4}}
{\r^{3/4}}\right)
.}\Eq(4.54)
$$

\proof
By Schwartz inequality
$$
\P\left[ \cup_{k,h\ge 1}\EE(k,h,+)\cap \DD(k,h,\rho,\tilde\d,+,\e)\right]
\le 
\sum_{h,k\ge 1} \left(\P\left[ \EE(k,h,+)\right]\right)^{1/2}
\left( \P\left[\DD(k,h,\rho,\tilde\d,+,\e)\right]\right)^{1/2}.
\Eq(4.541)
$$
Since
$$
\P\left[ \EE(k,h,+)\right]^{\frac 12}\le 2^{-\frac {(k+h)}  2}
\Eq(4.55)
$$   will be  summable  in $h,k$,
it remains to properly estimate the second term into parenthesis in
\eqv(4.541).
From \eqv(4.521) we just write
$$
\P\left[\DD(k,h,\rho,\tilde\d,+,\e)\right]
\le \sum_{\ell=\LL(-h,\rho)}^{\LL(k-1,\rho)-1} 
\sum_{\ell'=\ell+1}^{\LL(k-1,\rho)}
\P\left[|\YY^*(\underline \a,\sfrac{\rho \ell'}{\e},\sfrac{\rho
(\ell'+1)}{\e})-\YY^*(\underline \a,\sfrac{\rho \ell}{\e},\sfrac{\rho
(\ell+1)}{\e})| \le 2\tilde \d\right]
\Eq(4.56)
$$
and estimate each summand on
 the r.h.s. of \equ(4.56). 
If $\ell+1 <\ell'$ we   write:
$$\eqalign{
&\YY^*(\underline \a,\sfrac{\rho \ell'}{\e},\sfrac{\rho
(\ell'+1)}{\e})-\YY^*(\underline \a,\sfrac{\rho \ell}{\e},\sfrac{\rho
(\ell+1)}{\e})= \cr
&\sum_{\a=\sfrac{\rho(\ell+1)}{\e}+1} ^{\sfrac{\rho\ell'}{\e}}\chi(\a)
+\max_{\sfrac{\rho \ell'}{\e}+1 \le \tilde \a \le \sfrac{\rho
(\ell'+1)}{\e}} \sum_{\a =\sfrac{\rho \ell'}{\e}
+1}^{\tilde \a}\chi(\a)+\min_{\sfrac{\rho \ell}{\e}
\le \tilde \a \le \sfrac{\rho  (\ell+1)}{\e}}
\sum_{\a =\tilde \a+1}^{\sfrac{\rho(\ell+1)}{\e}}\chi(\a),}
$$
and using the independence of the $\chi(\a)$ we easily see that:
$$
\eqalign{
\P\left[|\YY^*(\underline \a,\sfrac{\rho \ell'}{\e},\sfrac{\rho
(\ell'+1)}{\e})-\YY^*(\underline \a,\sfrac{\rho \ell}{\e},\sfrac{\rho
(\ell+1)}{\e})| \le 2\tilde \d \right]
&\le\sup_x
\P\left[\sum_{\a=\sfrac{\rho(\ell+1)}{\e}+1}^{\sfrac{\rho
\ell'}{\e}} \chi(\a) \in [x, x+2\tilde \d]\right] \cr
&\le
\frac{4\tilde \d \sqrt{2\pi}}{V(\b,\th) \sqrt{(\ell'-\ell-1)\rho}}.\cr
} \Eq(4.56''')
$$
In the last inequality we have
used the concentration inequality of Le Cam (e.g. [\rcite{LC}], p.407) 
for the
symmetric random variables $\chi(\alpha)$ 
and {assumed} $0<\e<\e_0(\b,\th,\r)$ see \eqv(Term19).
This condition comes from a lower estimate of  what Le Cam called
$B^2(\t)$.  In our case 
$B^2(2\tilde \d)=(\ell'-\ell-1)\frac{\rho}{\e} 
\E[1\wedge(\chi(1)/2\tilde \d)^2]$.
A short computation gives
$$
\eqalign{
\E[1\wedge(\chi(1)/2\tilde \d)^2]&
\ge \frac{\E[(\chi(1))^2]}{4\tilde \d^2}
\left(1-\frac{\E[(\chi(1))^2\1_{\{|\chi(1)|>4\tilde \d\}}]}{\E[(\chi(1))^2]}
\right). 
}\Eq(AbouAmmar2)
$$
Using \eqv(2.20), \eqv(4.22001), Schwarz inequality, and that  
$\P[|\chi(1)|>4\tilde \d]\le 2 e^{-2\tilde \d^2/(\e V^2_+(\b,\th))}$, 
which 
follows from \eqv(4.3), a short computation shows that for 
$0<\e<\e_0(\b,\th,\r)$  the last term inside  parenthesis
in \eqv(AbouAmmar2) is bounded from below by $1/2$.

When $\ell'=\ell+1$, 
we bound the corresponding term on the r.h.s. of \equ(4.56)
as:
$$
\sup_x \P\left[\YY^*(\sfrac{\r}{\e}) \in [x,x+2\tilde \d]\right] \Eq (4.56'')
$$
where  $\YY^*(\a)\equiv \max_{1\le \tilde \a\le \a} \YY_{\tilde
\a}=\YY^*(1,1,\a)$ if $\a \ge 1$, and $\YY_{\a}$ given in  \eqv(3.C6).
Putting together \eqv(4.512),  \equ(4.56), \equ(4.56''') and \equ(4.56''), we
get
$$
\eqalign { & \P\left[\DD(k,h,\rho,\tilde\d,+,\e)\right] \le  
\big(C_1(\b,\th,2\FF^*){(s+\log2)}\big)^2 2(h+k+1)^2\frac{2
\sqrt{2\pi}}{V(\b,\th)}\sfrac{\tilde \d}{\rho^2} \cr & + \big (
C_1(\b,\th,2\FF^*){(s+\log2)}\big )^2\sfrac{(h+k+1)}{\rho}\sup_x
\P\left[\YY^*(\sfrac{\rho}{\e}) \in [x,x+2\tilde \d]\right]
}\Eq(4.62)
$$
The first term on the r.h.s. of \equ(4.62) suggests to take 
$\tilde \d= \rho^{2+a}$
with $a>0$. The last term
will be  estimated in the next Lemma \eqv(lem30), cf. \equ(3.paris2) below.

Recalling \eqv(4.541), \eqv(4.55),  \eqv(4.56), \eqv(4.62), and 
using \eqv(3.paris2) a short computation entails \eqv(4.54).
\eop
\vskip0.5cm

\noindent{ \bf\Lemma (lem30)} {\it  
There exist positive  constants
$ \g_0 (\b,\th) $ and $\r_0(\b,\th)$ such that
for  all  $\g \le \g_0 (\b,\th)$, for $0<\r<\r_0(\b,\th)$, for
$\tilde \d=\r^{2+a}$ with $a>-1/2$ , 
such that for $\d^*\g <\e\le \e_0(\b,\th,\r)$ with 
$\e_0(\b,\th,\r)$ given in \eqv(Term19), we have
$$
\frac1\rho\sup_x\P\left[\YY^*(\sfrac{\rho}{\e}) \in [x,x+2\tilde
\d]\right]\le 
\frac{1296}{V(\b,\th)}
\left(\frac{\tilde \d+(2+V(\b,\th))
\sqrt{\e\log\frac{ C_1}{\e}}
}{\r^{3/2}}\right)
\Eq(3.paris2)
$$
where  $C_1=C_1(\b,\th,2\FF^*)$ is given by \eqv(3P.50).}

\vskip0.5cm

\proof
 Let
 $\tilde T_x $ be the stopping time  given after
\eqv (4.111). We  write
$$
\P\left[\YY^*(\sfrac{\rho}{\e}) \in [x,x+2\tilde
\d]\right]= 
\P\left[\tilde T_x \le
\frac{\r}{2\e},\tilde T_{x +2\tilde \d} >\frac\r\e\right]+
\P\left[\frac{\r}{2\e} <\tilde T_x < \frac \r \e < \tilde
T_{x+2\tilde \d}\right].
\Eq(m1)
$$
Observe that for any $\tilde \d>0$, we have  
$
\left\{\frac{\r}{2\e} <\tilde T_x < \frac \r \e < \tilde
T_{x+2\tilde \d}\right\}= \left\{\YY^*(\sfrac {\r}{2\e} )<x,
\max_{\sfrac{\r}{2\e} \le \a \le \sfrac \r \e} \YY_{\a} \in
[x,x+2\tilde \d]\right\}
$
therefore if $0<\e<\e_0(\b,\th,\r)$, we obtain 
$$
\P\left[\frac{\r}{2\e} <T_x < \frac \r \e <
T_{x+2\tilde \d}\right] \le \P\left[\max_{\sfrac{\r}{2\e} \le \a
\le \sfrac \r \e} \YY_{\a} \in [x,x+2\tilde \d]\right]
\le
\sup_{u \in \R} \P\left[\YY_{\sfrac{\r}{2\e}} \in [u,u+2\tilde
\d]\right] \le  \frac{4\tilde \d\sqrt{2\pi}}{V(\b,\th)\sqrt{\r/2}}.
\Eq(3PP.1)$$
In the second inequality in \eqv(3PP.1), we  used that 
the law of $\max_{\sfrac{\r}{2\e} \le \a
\le \sfrac \r \e} \YY_{\a}$ is the convolution of the law of
$\YY_{\sfrac {\r}{2\e}}$ with another probability (the law of
$\YY^*(\sfrac {\r}{2\e})$, in this case).

Let us now consider the first summand on the r.h.s.  of \equ(m1).
Decomposing according to the value  of $\YY_{\tilde T_x}$, $\tilde
T_x$  and using the fact the variables $\chi(\cdot)$ are i.i.d. we
get
$$
\P\left[\tilde T_x \le \frac{\r}{2\e},\tilde T_{x +2\tilde \d}
>\frac\r\e\right]=
 \sum_{k=0}^{\rho/2\e} \int_{x}^{x+2\tilde\d} \,
\P\Big[ \tilde T_x=k, \YY_k\in dy \Big]
\P\Big[ \tilde T_{x+2\tilde\d-y} > \frac {\rho}{\e}-k \Big]
$$
Since $x-y \le 0$ we can write:
$$
\P\Big[ \tilde T_{x+2\tilde\d-y} > \frac {\rho}{\e}-k\Big] \le \P\Big[
\tilde T_{2\tilde\d} > \frac {\rho}{\e} -k\Big].
$$
Integrating in $y$ we then have:
$$
\P\left[\tilde T_x \le \frac{\r}{2\e},\tilde T_{x +2\tilde \d}
>\frac\r\e\right]\le \P\left[\tilde T_{2\tilde\d} >\frac{\r}{2\e}\right],
\Eq(4.66)
$$
and  collecting \eqv(m1), \eqv(3PP.1), and \eqv(4.66), we get
$$
\sup_x\P\left[\YY^*(\sfrac{\rho}{\e}) \in [x,x+2\tilde
\d]\right]\le \P\Big[ \tilde T_{2\tilde\d} > \frac {\rho}{2\e}
\Big] +
\frac{4\tilde \d\sqrt{2\pi}}{V(\b,\th)\sqrt{\r/2}}.
\Eq(4.71000)
$$
Now, it is easy to check that
$$
\P\left[ \tilde T_{2\tilde\d} > \frac {\rho}{2\e} \right] \le
\P\left[ \tilde T^-_{c\sqrt{\rho/2}} \le \tilde
T_{2\tilde\d}\right] +\P\left[ \tilde T^-_{c\sqrt{\rho/2}} \wedge
\tilde T_{2\tilde\d}\ge \frac {\rho}{2\e}\right],  \Eq(4.75)
$$
where  $T^-_{c\sqrt{\rho/2}}$ is the stopping time defined after 
\eqv(4.111) 
for a constant $c$ to be chosen soon.
Then we apply inequalities \eqv(4.76) and \eqv(4.77) given in the next lemma,
with  $a=c\sqrt{\r/2}$, 
$d = \rho/2$, and $x = 2\tilde\d$.
Collecting all together the estimates for  
$\P\left[\YY^*(\sfrac{\rho}{\e}) \in [x,x+2\tilde
\d]\right]$,  we have:

$$
\eqalign{ 
& \frac1\rho\,\sup_x\P\left[\YY^*(\sfrac{\rho}{\e}) \in
[x,x+2\tilde \d]\right]\le \frac{ 2\tilde\d + 9V(\b,\th)\sqrt{\e
\log \frac{C_1}{\e}}}
 {\rho(2 \tilde\d+
c\sqrt{\rho/2})}+
 \frac{ 8\sqrt 2 \tilde\d
c}{V^2(\b,\th)\rho^{3/2}} +
\cr &\quad + 
\frac{72}{\r^{3/2}V^2(\b,\th)}
\sqrt{\e\log \frac{C_1}{\e}}
\left(9(2\tilde \d+c\sqrt{\r/2})+V(\b,\th)
\sqrt{\e\log \frac{C_1}{\e}}\right)
}\Eq(4.96)
$$
with $C_1= C_1(\b,\th, (2\tilde \d)\vee c \sqrt{\r/2})$
see \eqv(3P.50). Taking $c=V(\b,\th)$ and assuming  that $\r_0(\b,\th)$
is  small enough, we have $C_1(\b,\th, (2\tilde \d)\vee c \sqrt{\r/2})
\le C_1(\b,\th,2\FF^*)$, and   a short computation entails 
\equ(3.paris2).
\eop

\noindent{ \bf\Lemma (lem9)} {\it
For all $x>0$, $a>0$, $C_1=C_1(\b,\th,x\vee a)$ as 
in \eqv(3P.50),  $\e_0(\b,\th, x\vee a)$ as in \eqv(4.35003), 
and  if  $\d^*\g<\e\le \e_0(\b,\th,x\vee a)$, we have: }
$$
\P\left[ \tilde T^-_a \le \tilde T_{x}\right] \le \frac{  x +
9V(\b,\th)\sqrt{\e \log \frac{C_1}{\e}} } {x+ a}, \Eq(4.76)
$$
$$
\P\left[ \tilde T^-_{a} \ge \tilde T_{x}\right] \le \frac{    a +
9V(\b,\th)\sqrt{\e \log \frac{C_1}{\e}} } {x+a}, \Eq(4.760001)
$$
$$
\P\left[ \tilde T^-_a \wedge \tilde T_{x}\ge \frac
d {\e}\right]
 \le \frac{ 4x a }{V^2(\b,\th) d } +
+
 \frac{36}{V^2(\b,\th) d } \sqrt{\e \log \frac{
C_1}{\e}}\,\,\left(
 9(x +a)+
V(\b,\th)\sqrt{\e \log \frac{C_1}{\e}}\right).
\Eq(4.77)
$$
The proof of the previous lemma is a standard
application of \equ(4.13) and  \equ (4.2401) together with Wald
identity applied to the martingales $\YY_{\a},\a \ge 0$ and
$(\YY_{\a})^2- \e c(\b,\th, \g /\d^*)\a$, and also the bound
\equ(2.21). Details are left out.

\medskip

To prove  \equ(3.16) in Theorem \equ(T1) we
need a classical result on the distribution
of the localization of the mimimum or the maximum of a simple random walk.
Since their distribution is the same, it is enough to consider the case of 
maximum. So, recalling \eqv(4.52),  let us denote
$
L_{\r/\e}=\inf \{\a>0: \YY_\a=\YY^*(0,0,\r/\e)\}.
$
 Such kind of result  was proved  by E. Sparre Andersen [\rcite{SA}]. 
Following step by step  the very nice computations  he did, 
see Theorem 3 of [\rcite{SA}], and using the Berry-Essen theorem
to estimate what is there  denoted by $\P\{S_n>0\}$, 
we can evaluate by the Cauchy integral formula 
the constant called $C$ at pg. 208, 
3 lines before (5.17) of [\rcite{SA}]. After simple,  however lengthy 
computations, we obtain the  following result.

\noindent {\bf \Proposition(prop1)} {\it There exists a constant
$C(\b,\th)$ (related to $V(\b,\th)$) and $\r_0=\r_0(\b,\th)$ such that 
for all $0<\r<\r_0$
there exists $\e_0=\e_0(\b,\th)$ such that for all $0<\e\le \r \e_0$, for all 
$0<\kappa \le 1/2$, for all interval $0<a<a'\le 1$ such that 
$a'-a \ge \sfrac{2\e}{\rho}\,$, 
$$
\eqalign{ &\Big| \P\left[ L_{\rho/\e} \in
[a\rho/\e,a'\rho/\e]\right] -\frac{\cos(\pi\kappa)}{\pi}
\int_{a(\e,\rho)}^{a'(\e,\rho)} \frac{ dx}{x^{\sfrac{1}{2}+\kappa}
(1-x)^{\sfrac{1}{2}-\kappa}} \Big|\cr & \quad \le
\frac{1}{\G(\frac 12 -\kappa)} \Big(\frac{\e}{\rho}\Big)^{\sfrac12
+\kappa} +  \frac{1}{\G(\frac 12 +\kappa)}
\Big(\frac{\e}{\rho}\Big)^{\sfrac12 -\kappa}
+ \frac
{\e}{\rho (a'-a)} \exp\left(  8^{\sfrac {C(\b,\th)}{\kappa^2}} 2\log
\sfrac {C(\b,\th)}{\kappa^2}\right)\cr }\Eq(4.72)
$$ where $ x(\rho,\e)= (\rho x+\e)(\rho +\e)^{-1}$ for $x=a,a'$}

\noindent{\bf Proof of Theorem \eqv(T1)}

We start proving \eqv(3.14).
 For any $Q>Q_0=4\log 2 C_1(\b,\th,\FF^*)$, if  we take $Q/\e$
blocks of length $\e/\g$ on the right of the origin, then using 
Lemma \eqv(lem6) with $s=\log 2$ and 
$k=1+[q/(2C_1(\b,\th,2\FF^*)\log2)]$
where $[\cdot]$ is the integer part, with a $\P$--Probability 
at least $(1-3e^{-Q/2C_1(\b,\th,\FF^*)})$ there is at least one index $i$ 
among $1,\dots, [Q/(2C_1(\b,\th,2\FF^*)\log2)]$ such that $S_i=S_{i+1}$.
From Lemma  \eqv(lem111)  with $\P
\ge 1-     G (\b,\th, \e, f) \log  G (\b,\th, \e, f)  \frac 2 {\log 2}$
with  $G (\b,\th, \e, f)$  defined in \eqv (3.C3) we have an 
elongation there.  Therefore the probability of not having any elongation
on the right of the origin within  $Q/\e$ blocks of length $\e/\g$ is 
less than
$$
3e^{-\frac{Q}{2C_1(\b,\th,\FF^*)}} 
+\frac{2}{\log 2}G (\b,\th, \e, f) \log  G (\b,\th, \e, f)
\Eq(4.114)
$$
which implies \eqv(3.14).

The  proof of  \eqv(3.141) is done in a similar way. We first
apply  Lemma  \eqv(lem7) with $s=\log 2$ and 
$L=1+[Q/(kc(\b,\th,2\FF^*)2\log 2)]$
then Lemma \eqv (lem111).


To prove  \eqv(3.15), we recall  Lemma \eqv(lem8) and the arguments 
that precede it. Taking $\tilde \d(\r)=\r^{2+a}$ and 
recalling \eqv(3PP.16)
we have 
$$
\P\left[\PP_2(f,s_1,Q,a_i,b_{i+1},\r,\tilde \d(\r))\right]
\le 4e^{-s}+(s+\log 2) G_1(\b,\th,\tilde \d(\r),\e).
\Eq(4P.116)
$$
Choosing  $s=\log 4/(G_1(\b,\th,\tilde \d(\r),\e))$ and taking 
$\r_0(\b,\th)$ and $\e_0(\b,\th,\r)$ small enough, we get
\eqv(3.15).

For the proof of \equ(3.1411), { recalling \eqv(30.54)}  we write
$$ \P  \left[  \PP''_2(f,Q)\right]
\le \P\left[\PP'_2(f,Q)\cap \left\{\max_{\a \in
[-Q/\e,Q/\e]}|\chi(\a)|\le f\right\}\right] +\P\left[\max_{\a \in
[-Q/\e,Q/\e]}|\chi(\a)|>f\right],
$$
{ and taking $\rho'=(9f)^{1/(2+a)}$}, we 
consider the event
$$
\tilde\DD(Q,\rho',\e) 
\equiv \Big\{
\exists
\ell, \ell',\,-Q/\rho' \le \ell<\ell'\le (Q-1)/\rho' ;
|\YY^*(\underline \a, \sfrac{\rho' \ell}{\e}, \sfrac{\rho' (\ell+1)}{\e})
  -
 \YY_*(\underline \a,  \sfrac{\rho' \ell'}{\e},\sfrac{\rho'
(\ell'+1)}{\e})-2\FF^*| \le 9f \Big\}.
$$
where $\YY_*$ is defined as in \equ(4.52) replacing $\max$ by
$\min$. 

Simple observations show that $\PP'_2(f,Q)\cap \{\max_{\a
\in [-Q/\e,Q/\e]}|\chi(\a)|\le f\} \subseteq \tilde\DD(Q,\rho',\e)$.
{Following the arguments that lead to \eqv(4.62), 
assuming  $0<\e\le \e_0(\b,\th,f)=(9f)^{2}/(2V^2(\b,\th)\log 1944)$, 
using Lemma \eqv(lem30) with $2\tilde \d$ replaced by $9f$
one gets \eqv(3.1411).}

 The proof of \eqv(3.16) follows from \eqv(4.72) estimating 
the integral in  the left hand side of \eqv(4.72) by
$
8 (a'-a)^{\sfrac 12 -\kappa}
$ 
which can be obtained by cutting the interval 
$[a(\e,\rho),a'(\e,\rho)]$  into two equal
pieces. Using \eqv(4.72) for $a=0,a'=\r$ and a short
computation entails \eqv(3.16).
\eop

\noindent{\bf Proof of Proposition \eqv(Terminator)}

To prove \eqv(Terminator1), notice   that
$\g J(\o)\supset [\e\t_{-1},0]\cup[0,\e\t_{1}]$.
Therefore, using \eqv(4.250000) and a short computation one gets
$$
\P[\g|J|\le x]\le 2e^{-\frac{(\FF^*)^2}{18 x V^2(\b,\th)}}
\Eq(4P.12)
$$
for $0<x<(\FF^*)^2/(V^2(\b,\th) 18\log 2)$.
\eqv(Terminator2) follows at once, due to  \eqv(4.46), \eqv(4.46'),
and the fact that 
$\g J(\o) \subset [\e\t_{i^*_{-2}},\e\t_{i^*_2}]$.
Therefore  \eqv(4.4800) with $k=2$ entails 
$$
\P[\g |J|\ge x] \le 2\P[\e\t_{2L}\ge \frac x2]+ 
2\left(\frac{1}{2^{L-1}} + \left(\frac34\right)^{L/2}\right).
\Eq(4P.13)
$$
Using now \eqv(4.3401) with $k=2$, $s=\log 2$ one gets
$
\P[\e\t_{2L}\ge 4L C_1(\b,\th,\FF^*) \log 2]
\le e^{-2L\log 2}.
$
Taking  $L=x/(8C_1(\b,\th,\FF^*)\log 2)$ one obtains  
after a short computation
\eqv(Terminator2). \eop

The following lemma will be useful in the next section; it is in
fact an immediate consequence of \eqv(3.15)  and the proof is
omitted .

 {\bf \Lemma (lem112)} {\it Under the hypothesis of 
Corollary \eqv(zorro) and with the same notations
with  $\P$--probability larger than
$
1-\e^{\frac{a}{16(2+a)}}
$ we have
$$
\sum_{\a=\a_1}^{\a*_1  } \chi (\a)  \ge \e^{1/4},
\quad \sum_{\a=\a^*_0}^{\a_1  } \chi (\a)  \ge \e^{1/4},
\Eq(4.115)
$$
provided $\a^*_0$ is the beginning and $\a^*_1$ is the end of a
positive elongation, $\a^*_0 +\frac\r\e <\a_1 <
\a^*_1-\frac\r\e$.}

\bigskip

\chap{4 The block spin representation and the $\e$ rigidity}4
\numsec= 4
\numfor= 1
\numtheo=1
\medskip

 We  start by defining the set of profiles  having runs of 
$+$ or of  $-$, with   length at least  $\frac \e \g$.

{\bf \Definition(def1)} {\it Given  $\frac {\e}{\g} >\d^*$,  an
interval $\D_Q\equiv [Q_1,Q_2]\g^{-1}$  of  length in macroscopic units 
$  \frac Q {\g }=\frac {(Q_2-Q_1)}{\g} $,   $Q>0$  such that
$\frac {Q_1} \e$ and $\frac {Q_2} \e $ are integers,
$ \z_4>\z_1>8\g/\d^*
$, $1 >\d >\d^*>0$,
$R_1>0 $, $\eta = \pm 1$,  we define $\AA_1(\D_Q,\eta)=
\AA_1(\D_Q,\d,\z_1,\z_4,\d^*,\g,\e,R_1, \eta )$ as}
$$ \eqalign {
 \AA_1(\D_Q,\eta ) = \Big\{& m^{\d^*}_{\D_Q}\colon \exists
 k \in \N , \exists r_1,\dots,r_k
\in \Big\{\frac
{Q_1}{\e}+1,\frac{Q_1}{\e}+2,\dots,\frac{Q_2}{\e}-2,\frac{Q_2}{\e}-1
\Big \}; \cr &  r_0=\frac{Q_1}{\e}, r_{k+1}= \frac {Q_2}
{\e},r_1<\dots<r_{k}, \exists\, q_{i} \in [r_{i}\frac \e \g
(r_{i}+1 ) \frac \e \g ] \text{s.t.}\cr
 &\eta^{\d,\z_4} (\ell)=\eta (-1)^{i-1}\, \forall  \ell \in
\CC_{1} ([ (r_{i-1}+1 ) \frac \e \g,q_i-R_1]),\cr&
\eta^{\d,\z_1}(q_i-R_1)=(-1)^{i-1}\eta,
\eta^{\d,\z_1}(q_i+R_1)=(-1)^i\eta,\,\cr & \eta^{\d,\z_4}(\ell
)=\eta (-1)^i \,\,\forall  \ell \in [(q_i+R_1)\wedge \frac Q \g,
\frac \e \g (r_{i+1})], \quad\hbox {for} \; i=1,\dots,k \Big\} \cr
} \Eq (6.1)$$ {\it and}
$$
\AA_1(\D_Q) \equiv
 \cup_{\eta \in \{-1,+1\}}\AA_1(\D_Q,\eta).
\Eq(6.10000)
$$
\goodbreak

\noindent{\bf Remark.}

\item {$\bullet$}  The integer $k\ge 0$ represents the number
of blocks of length $R_1$ within $\D_Q$ where there is at least
one change of phases which means that  $
\eta^{\d,\z_1}(q_i-R_1)=(-1)^{i-1}\eta,\,
\eta^{\d,\z_1}(q_i+R_1)=(-1)^i\eta$. There are  no restrictions on
the profiles within the interval $[q_i-R_1+1, q_i+R_1-1]$.

\item {$\bullet$} $ r_{i}$ is the index of the $i$--th block of length
$\e/\g$ in macroscopic units such that in  $ [ q_i-R_1,  
   q_i+R_1] \subset [r_{i}\frac
{\e}
{\g }-   R_1, (r_{i}+1)\frac
{\e}
{\g} + R_1 ] $ we see at least one change of phases.

\item{$\bullet$} $R_1$ will be chosen as
 an upper bound for the  length of the longest interval
where the system can stay out of ``equilibrium", that is to have a run of
$\eta^{\d,\z_1}=0$. This length is related to the parameters
$\z_1,\d$, by $R_1 \approx (\d\z_1^3)^{-1}$,
see \eqv(5.4810).
\medskip
Another definition is needed to describe what happens in the
intervals $[q_i-R_1,q_i+R_1]$.

{\bf \Definition(def2)} {\it Let $\D_L=[\ell_1,\ell_2]$   
be an interval of length $L$ in macroscopic units   and
$\d>0,\;\z_4>\z_1>8\g/\d^*
$ as above.  For $\eta=+1$ or $\eta=-1$ we set
$$
\WW^{\z_1,\z_4 }(\D_L,\eta)\equiv \left\{m^{\d^*}_{\D_L}\colon
\eta^{\d,\z_1}(\ell_1)= \eta^{\d,\z_1}(\ell_2)=\eta,\; \exists
\tilde \ell,\; \ell_1< \tilde \ell < \ell_2\quad
\eta^{\d,\z_4}(\tilde \ell)=-\eta \right\} \Eq(5.31)
$$  and
$
\WW^{\z_1,\z_4}(\D_L)\equiv\WW^{\z_1,\z_4}(\D_L,+1)\cup\WW^{\z_1,\z_4}(\D_L,-1)
$.}

\vskip0.5cm
Given a positive  integer  $L_2 $    we denote by
$
\BB (\D_Q,L_2)=
\cap_{L=3}^{L_2}\cap_{\D_L \subset \D_Q} (\WW^{\z_1,\z_4}(\D_L))^c.
$
The 
profiles in this set  do not have
two changes of phases 
within an interval of length  smaller than
$L_2$,  uniformly along  intervals that are within $\D_Q$.
We set
$$
\AA(\D_Q)=\AA_1(\D_Q) \cap \BB(\D_Q, L_2)
\Eq(ME.3)
$$
If   $L_2>2R_1$ the profiles in $\AA(\D_Q)$ have exactly one
change of phase within each interval $[q_i -R_1,q_i+R_1]$.
The main result of this Section is the following:

\noindent {\bf \Theorem (61)} {  \it Let $\b,\th$ satisfy
\eqv(2.17). We take $\k(\b,\th)>0$ verifying {\rm \eqv(2.19)},
$\FF^*$ is defined in \eqv(min), and $V(\b,\th)$
given by {\rm \eqv(2.1701)}. There exist
$0<\g_0=\g_0(\b,\th)<1$, $ 0<d_0=d_0(\b,\th)<1$, 
 and $0<\z_0=\z_0(\b,\th)<1$,
such that for all $0<\g\le \g_0$,
  for all $\d^*,\d,\z_4,\z_1$ with $\d^*\ge \g$,  $\g/\d^*\le d_0$,
$1>\d  > \d^*>0$,  $\z_0 \geq \z_4 >\z_1>8\g/\d^*
$, and  $Q>3$  
that
satisfy the following conditions
$$
\frac {32}{ \k(\b,\th)} \z_1 \le\d\z_4^3, \Eq(6.03)
$$
$$
\frac{128 (1+\th)}{\k(\b,\th)} 
\frac{2(5+\FF^*)}{\FF^*} \sqrt {\frac{\g}{\d^*}}
<\d\z_1^3, \Eq(6.02)
$$
$$
\z_1\ge \left(5184 (1+c(\b\th))^2 \sqrt{\frac{\g}{\d^*}}\right) \vee
\left(12 \frac{e^3\b }{c (\b,\th)} \frac{(\d^*)^2}{\g}\right)^2
\Eq(6.004)
$$
for constants  $c(\b,\th)$ given in {\rm \eqv(mars.76)}, and
$c(\b\th)$ given in {\rm \eqv(3.25PM)},
$$
\sqrt {\g} \log Q \le 
\frac{\sqrt{6e^3 \b}}{256}
\Eq(6.00410)
$$
if we call
$$
R_1 = \frac{4(5+\FF^*)}{\k(\b,\th)\d\z_1^3}
\Eq(6.04)
$$
and
$$
 L_2= \frac{\FF^*}{32(1+\th)} \sqrt {\frac{\d^*}{\g}},
\Eq(5.042)
$$
then for any  interval $\D_Q $  of length $ \frac {Q} {\g } $ and
any $\e> \g \d^* $, there exists
 $\O_{4} = \O_{4} (\g,\d^*, \D_Q,\e,\d,\z_1,\z_4) $
with
$$
\P[\O_{4}] > 1 -6\g^2 -\frac{6Q}{\e}
\exp\left\{ -\frac{ (\FF^*)^2}{\e2^{6}V^2(\b,\th)}\right\}
\Eq (M.1) $$
and  for all $\o\in \O_{4} $, we have}
$$
\mu_{\b,\th,\g}\left ( \AA (\D_Q)\right ) \ge 1- (\frac
{3Q}{\g^2})^5 e^{-\frac \b \g \left[(\frac{\k(\b,\th)}{4}
\d\z_4^3)\wedge \FF^*\right]}. \Eq (6.2)
 $$

To prove Theorem \eqv(61), we represent  the system in terms  of block 
spins. This  representation  was used  also in [\rcite{COP1}]. 
However, the way to treat some error terms that appear at the very beginning
of the computations is  different, see \eqv(3.6P) and \eqv(C.1P). 
\vskip.3cm
\noindent
{\bf Analysis of the block-spin representation}
\goodbreak

With $C_{\d^*}(V)$ as in Section 2, let $\S_V^{\d^*}$ denote the
sigma--algebra of $\SS$ generated  by $m^{\d^*}_{V}(\s)\equiv
(m^{\d^*}(x,\s)$, $x\in \CC_{\d^*}(V))$, where $ m^{\d^*}(x,\s)=
(m^{\d^*}(+,x,\s),m^{\d^*}(-,x,\s))$, cf. \equ(2.10).

We take  $I= (i^-,i^+]\subset \R$ with $i^{\pm} \in \Z$. The
interval $I$ is assumed to be  $\DD_{\d^*}$--measurable and we set
$\partial^+ I\equiv \{x\in \R\colon i^+< x\le i^+ +1\}$,
$\partial^- I\equiv \{x\in \R\colon i^- -1< x\le  i^-\}$, and
$\partial I=\partial^+I\cup\partial^-I$.

For  $(m^{\d^*}_I,m^{\d^*}_{\partial I})$  in $\MM_{\d^*}(I \cup
\partial  I)$, cf. \eqv(2.14), we set
 $\widetilde
m^{\d^{*}}(x)=({m_{1}^{\d^{*}}(x)+m_{2}^{\d^{*}}(x)})/{2}$,
$$ E(m^{\d^*}_I)\equiv -\frac {\d^*}2 \sum_{(x,y)\in \CC_{\d^*}(I)\times
\CC_{\d^*}(I)} J_{\d^*}(x-y) \tilde m^{\d^*}(x)\tilde m^{\d^*}(y),
\Eq(3.4)
$$
and
$$ E (m^{\d^*}_I, m^{\d^*}_{\partial^{\pm}I}) \equiv
-\d^*\sum_{x\in \CC_{\d^*}(I)}\sum_{y\in
\CC_{\d^*}(\partial^{\pm}I)} J_{\d^*}(x-y) \tilde
m^{\d^*}(x)\tilde m^{\d^*}(y), \Eq(3.5P)
$$
where $J_{\d^*}(x)= \d^*J (\d^*x)$.  It is easy to see that
$$
   H_\g(\s_{\g^{-1}I}) + \th \sum_{i\in \g^{-1}I} h_i\s_i = \frac
1{\g} E(m^{\d^*}_I) + \frac 1 \b  \log \Big [ \prod_{x \in
\CC_{\d^*}(I)}\prod_{y \in \CC_{\d^*}(I)} e^{\b  U
(\s_{A(x)},\s_{A(y)})  } \Big], \Eq(3.6P)
$$
where
$$
U (\s_{A(x)},\s_{A(y)})= -  \sum_{i \in A(x),j  \in A(y)}   \g
  \big[ J(\g|i-j|)-J(\d^*|x-y|) \big ]\s_i\s_j.  \Eq (C.1P)
  $$
Since the interaction is only between  adjacent blocks of
macroscopic length 1, see \eqv(2.2), we see that for all   
intervals $I$, for $s=+$ or $s=-$
$$
\sup_{\s_{\g^{-1}I}\in M^{\d^*}(m^{\d^*}_{I})}
\sup_{\s_{\g^{-1}\partial^{s}I}\in
M^{\d^*}(m^{\d^*}_{\partial^{s}I})} \Big| W_\g(\s_{\g^{-1}I}|
\s_{\g^{-1}\partial^{s}I  } )  -\frac 1{\g} E(m^{\d^*}_I,
m^{\d^*}_{\partial^{s}I})\Big| \le  \d^* \g^{- 1}, \Eq(3.7P)
$$
where
$M^{\d^*}(m^{\d^*}_I)\equiv \{\s\in \g^{-1} I:
m^{\d^*}(x,\s)=m^{\d^*}(x),\,\forall {x \in \CC_{\d^*}(I)}\}$.

 Recalling \eqv(2.13),
and using \eqv (3.6P) and \eqv(3.7P),  if  $F^{\d^*}$ is  a
$\S_I^{\d^*}$-measurable bounded function and $m^{\d^*}_{\partial
I} \in \MM_{\d^*}(\partial I)$, and 
$\mu_{\b,\th,\g}\left(F|\S^{\d^*}_{\partial I}\right)$ denotes the
conditional expectation of $F^{\d^*}$ given the $\s$--algebra
$\S^{\d^*}_{\partial I}$, we have
$$
\eqalign{ &\mu_{\b,  \th,\g} \left(
F^{\d^*}\bigm|\S^{\d^*}_{\partial I}\right) (m^{\d^*}_{\partial
I}) = \frac {  e^{\pm  \frac  \b \g  2\d^* } }
{Z_{\b,\th,\g,I}(m^{\d^*}_{\partial I})} \times
 \cr &\times
\sum_{m^{\d^*}_I\in {\cal M}_{\d^*}(I) }F^{\d^*}(m^{\d^*}_I)
e^{-\frac \b\g \left( E(m^{\d^*}_I)+ E(m^{\d^*}_I,
m^{\d^*}_{\partial^{ }I})- \frac {\th \d^*}2 \sum_{x\in
\CC_{\d^*}(I)} (m^{\d^*}_1(x)-m^{\d^*}_2(x) )\right)}  \cr &
\times
  \sum_{\s_{\g^{-1}I}}
\prod_{x_1 \in \CC_{\d^*}(I)} \1_{
\{m^{\d^*}(x_1,\s)=m^{\d^*}(x_1)\}}\,
 e^{ 2 \b\th \l(x_1)\sum_{i\in D(x_1)}\s_i}\cr & \times \prod_{x_2 \in
\CC_{\d^*}(I)}\prod_{y_2 \in \CC_{\d^*}(I)} e^{-\b U
(\s_{A(x_2)},\s_{A(y_2)}) }
 } \Eq(3.8P)
$$
where
$$
\eqalign{Z_{\b, \g,\th,I}(m^{\d^*}_{\partial I}) &=
\sum_{m^{\d^*}(I)\in {\cal M}_{\d^*}(I) } e^{-\frac \b\g \left(
E(m^{\d^*}_I)+ E(m^{\d^*}_I, m^{\d^*}_{\partial^{ }I})- \frac {\th
\d^*}2 \sum_{x\in \CC_{\d^*}(I)} (m^{\d^*}_1(x)-m^{\d^*}_2(x)
)\right)}  \cr &\times
  \sum_{\s_{\g^{-1}I}}
\prod_{x_1 \in \CC_{\d^*}(I)} \1_{
\{m^{\d^*}(x_1,\s)=m^{\d^*}(x_1)\}}\,
 e^{ 2 \b\th \l(x_1)\sum_{i\in D(x_1)}\s_i}\cr & \times \prod_{x_2 \in
\CC_{\d^*}(I)}\prod_{y_2 \in \CC_{\d^*}(I)} e^{-\b U
(\s_{A(x_2)},\s_{A(y_2)}) }.
 } \Eq(3.8P1)
$$
Equality  \eqv (3.8P) has to be interpreted as an upper bound for
$\pm=+1$ and a lower bound for $\pm=-1$.
 Given $ m^{\d^*}_I$, we define the probability measure
on $\{-1,+1\}^{\g^{-1}I}$ by
$$
\E_{ m^{\d^*}_I}[f]\equiv
 \frac { \sum_{\s_{\g^{-1}I}  } \prod_{x_1
\in \CC_{\d^*}(I) } \1_{ \{m^{\d^*}(x_1,\s)=m^{\d^*}(x_1)\}}e^{ 2
\b \th \l(x_1) \sum_{i\in D(x_1)} \s_i } f(\s) }
 { \sum_{\s_{\g^{-1}I}  }  \prod_{x_1 \in
\CC_{\d^*}(I)} \1_{ \{m^{\d^*}(x_1,\s)=m^{\d^*}(x_1)\}}e^{ 2 \b\th
\l(x_1) \sum_{i\in D(x_1)} \s_i }  }. \Eq (3C.1)
$$

Inside the sum $\sum_{m^{\d^*}_I}$ in \eqv(3.8P), we divide and
multiply by
$$  \sum_{\s_{\g^{-1}I}}
\prod_{x_3 \in \CC_{\d^*}(I) } \1_{
\{m^{\d^*}(x_3,\s)=m^{\d^*}(x_3)\}}e^{ 2 \b \th \l(x_3) \sum_{i\in
D(x_3)} \s_i }
    $$
to get 
 $$
\eqalign{ &\mu_{\b,\th,\g}\left(F^{\d^*}\bigm|\S_{\partial^{ }
I}\right) (m^{\d^*}_{\partial^{ }I}) = \frac {e^{\pm
 \frac{\b}{\g} 2\d^* }}  {Z_{\b,\th,\g,I}(m^{\d^*}_{\partial I})}
\cr & \times\sum_{m^{\d^*}_{I} \in {\cal M}_{\d^*}(I)}
F^{\d^*}(m^{\d^*}) e^{-\frac \b\g \left( E(m^{\d^*}_I)+
E(m^{\d^*}_I, m^{\d^*}_{\partial^{ }I})- \frac {\th \d^*}2
\sum_{x\in \CC_{\d^*}(I)} (m^{\d^*}_1(x)-m^{\d^*}_2(x) )\right)}
\cr & \times e^{\log \E_{ m^{\d^*}_{ I}}[\prod_{x_2 \neq y_2}
e^{\b U (\s_{A(x_2)},\s_{A(y_2)})}]} \cr &\times
 \sum_{\s_{\g^{-1}I}  }
\prod_{x_3 \in \CC_{\d^*}(I) } \1_{
\{m^{\d^*}(x_3,\s)=m^{\d^*}(x_3)\}}e^{ 2 \b \th \l(x_3) \sum_{i\in
D(x_3)} \s_i }.  }\Eq(3.150)
$$

If we notice that the last sum $    \sum_{\s_{\g^{-1}I}  } $
factors out into a product over the intervals
 of length $\d^*\g^{-1}$, indexed by $\CC_{\d^*}(I)$, 
we get that for  each $x\in\CC_{\d^*}(I)$
 $$
 \sum_{\s \in \SS_{\d^*\g^{-1}}}
\1_{ \{m^{\d^*}(x,\s)=m^{\d^*}(x)\}}={{\d^*\g^{-1}/2} \choose
{{{\frac{1+m^{\d^*}_1(x)}{2}\d^*\g^{-1}/2}}}}
  {{\d^*\g^{-1}/2} \choose {{{\frac{1+m^{\d^*}_2(x)}{2}\d^*\g^{-1}/2}}}},
 \Eq(3.10P)
 $$
 and recalling  the probability measure on $\{\-1,+1\}^{A(x)}$ 
defined through \equ(3.9P),  \equ(3.150) becomes 
$$
\mu_{\b,\th,\g}\left(F^{\d^*}\bigm|\S_{\partial^{ }
I}\right) (m^{\d^*}_{\partial^{ }I}) = \frac {e^{\pm \frac \b \g 2
\d^*}}  {Z_{\b,\th,\g,I}(m^{\d^*}_{\partial^{ }I})} 
\sum_{m^{\d^*}_{I} \in {\cal M}_{\d^*}(I)}
F^{\d^*}(m^{\d^*}) e^{-\frac \b \g \left\{{ \widehat
\FF}(m^{\d^*}_I|m^{\d^*}_{\partial^{ }I} )+\g \GG(m^{\d^*}_I)+
\g  V (m^{\d^*}_I)\right\}},
\Eq(3.15P)
$$
where
$$
\eqalign{ {\widehat \FF}(m^{\d^*}_I|m^{\d^*}_{\partial^{ }I} )=&
 E(m^{\d^*}_I)+ E(m^{\d^*}_I,m^{\d^*}_{\partial^{ }I})  -\frac
{\th \d^*}2 \sum_{x\in \CC_{\d^*}(I)} (m^{\d^*}_1(x)-m^{\d^*}_2(
x))\cr  &- \d^* \sum_{x\in \CC_{\d^*}(I)} \frac {\g}{\b\d^*} \log
{{\d^*\g^{-1}/2} \choose
{{{\frac{1+m^{\d^*}_1(x)}{2}\d^*\g^{-1}/2}}}}
  {{\d^*\g^{-1}/2} \choose 
{{{\frac{1+m^{\d^*}_2(x)}{2}\d^*\g^{-1}/2}}}},\cr }\Eq(3.11P)
$$
 $ \GG$ is already defined by \equ(3.14P), \equ(3.13P) and
\equ(3.9P) in Section 2, and
$$
V (m^{\d^*}_I) \equiv V_I (m^{\d^*}_I,h)  = \frac{1}{\b}\log \E_{
m^{\d^*}( I)}[\prod_{\scriptstyle{x\neq y} \atop \scriptstyle {x,y
\in \CC_{\d^*}(I)\times \CC_{\d^*}(I)} } e^{\b U
(\s_{A(x)},\s_{A(y)})}].
  \Eq(3.140)
$$

That is, up to the error  terms ${e^{\pm \frac \b \g 2 \d^*}}$, we
have been able to describe  our system in terms of the  block spin
variables giving  a rather explicit form to the deterministic and
the stochastic part.

The following lemma gives an explicit integral form of the deterministic
part of the block spins system. For $ m\in \cal{T}$, let us call
$$
\eqalign{ \tilde \FF(m_I|m_{\del I})&= \int_I f_{\b,\th}(m(x))\,dx
+ \frac  14 \int_I\int_I J(x-y)[\tilde m(x)-\tilde m(y)]^2\,dxdy\cr &+\frac 12
\int_I\,dx \int_{I^c} J(x-y) [\tilde m(x)-\tilde m(y)]^2\, dy\cr }\Eq(5.Rio2)
$$
which is obviously related to \eqv(AP.1).

\medskip
\noindent{\bf \Lemma(MM)}{\it  
If $m^{\d^*}_{I\cup \del I}\in
\MM_{\d^*}(I\cup \del I)$ and $m(r)=m^{\d^*}(x)$ for $r \in
((x-1)\d^*,x\d^*]$ and $x\in \CC_{\d^*}(I\cup \del I)$, one has
$$
|\widehat \FF(m^{\d^*}_I|m^{\d^*}_{\partial I} )- \tilde
\FF(m_I|m_{\del I})+ 
\frac{\d^*}{2} \sum_{y \in \CC_{\d^*}(\del
I)} \big[\tilde m^{\d^*}(y)\big]^2 \sum_{x \in \CC_{\d^*}(I)}
J_{\d^*}(x-y)
|\le |I|\frac{\gamma}{\d^*} \log\frac{\d^*}{\gamma}.
\Eq(5.Rio3)
$$
} 
\proof
 Since
$$|\1_{\{\g|i-j|\le 1/2\}}-\1_{\{\d^*|x-y|\le 1/2\}}|\le
\1_{\{-\d^* +1/2\le \d^*|x-y|\le \d^* +1/2\}} \Eq (S.9)
$$
we have that
$$ |U (\s_{A(x)},\s_{A(y)})| \le   \g   (\frac {\d^*}
\g)^2 \1_{\{1/2-\d^*\le \d^*|x-y|\le 1/2+\d^*\}}.
 \Eq (H.1) $$
Given $m^{\d^*}_I\in \MM_{\d^*}(I)$, we easily obtain from
\eqv(H.1) that, on $M^{\d^*}(m^{\d^*}_{I})$: 
$$
   \left | H(\s_{\g^{-1}I}) + \th \sum_{i\in \g^{-1}I} h_i\s_i -\frac
1{\g} E(m^{\d^*}_I)  \right| = \frac 1 \b
 \left |  \log \Big [ \prod_{x \in \CC_{\d^*}(I)}\prod_{y \in
\CC_{\d^*}(I)}
e^{\b  U (\s_{A(x)},\s_{A(y)})  } \Big] \right|   
\le |I| \d^*  \g^{- 1}.
\Eq (S.2)
$$
Using Stirling formula, see [\rcite{Ro}],  we get
$$\eqalign{
\Big| \d^* \sum_{x\in \CC_{\d^*}(I)} \frac
1{2\b}\left(\II(m_1^{\d^*})+\II(m_2^{\d^*})\right) - \d^*
\sum_{x\in \CC_{\d^*}(I)} &\frac {\g}{\b\d^*} \log
{{\d^*\g^{-1}/2} \choose
{{{\frac{1+m^{\d^*}_1(x)}{2}\d^*\g^{-1}/2}}}}
  {{\d^*\g^{-1}/2} \choose {{{\frac{1+m^{\d^*}_2(x)}{2}\d^*\g^{-1}/2}}}}
  \Big |\cr & \le \frac
{1}{\b}|I| \frac{\g}{\d^*}\log \frac{\d^*}{\g},} \Eq(3.11P1)
$$
where $\II(\cdot)$ is defined after  \eqv (2.14000).
Recalling the definition of $f_{\beta,\th}(m)$,
 cf. \equ(2.14000) the lemma is proven. \eop

Concerning the stochastic part in \eqv(3.15P),
note that there are two random terms in \eqv(3.15P):
$\GG(m^{\d^*}_I)$ and $V(m^{\d^*}_I)$. To treat them we will use
the following classical deviation  inequality for Lipschitz
function of Bernoulli random variables. See [\rcite{LT}] or 
[\rcite{COP1}] for a short proof.

\medskip
\noindent{\bf \Lemma(41)}{\it Let $N$ be a positive integer and
$F$ be a real function on $\SS_N=\{-1,+1\}^{N}$ and for all $i\in
\{1,\dots,N\}$ let
$$
\|\partial_i F\|_{\infty}\equiv \sup_{(h,\tilde h): h_j=\tilde
h_j,\forall j\neq i} \frac{\left|F(h)-F(\tilde
h)\right|}{|h_i-\tilde h_i|}. \Eq(4.13PP)
$$
 If $\P$ is the symmetric Bernoulli measure and
$\|\partial(F)\|^2_{\infty} =\sum_{i=1}^N
\|\partial_i(F)\|^2_{\infty}$ then, for all $t>0$
$$
\P\left[ F-\E(F)\geq t\right]\le  e^{-\frac{t^2}{4
\|\partial(F)\|^2_{\infty}}} \Eq(4.14)
$$ and also}
$$
\P\left[ F-\E(F)\le -t\right]\le  e^{-\frac{t^2}{4
\|\partial(F)\|^2_{\infty}}}. \Eq(4.140)
$$
\medskip
 For $F(h)=|2\th \l(x_3)\sum_{i\in D(x_3)} \s_i|$, as it  appears  in
\eqv(3.150), Lemma \eqv(41) implies the following rough estimate:

\noindent{\bf \Lemma(60002)}({\bf The rough estimate}) {\it For all 
$\d^*>\g >0$ and   for all positive integer $p$,
that satisfy
$$ 64(2+p){\d^*}\log\frac 1\g\le 1
\Eq(4.00510)
$$
 there exists $\O_{RE}=\O_{RE}(\g,\d^*, p)\subseteq \O$  with $\P[\O_{RE}]\ge 1-\g^2$
such that on $\O_{RE}$ we have:
$$
\sup_{I\subset [-\g^{-p},\g^{-p}]}
 \frac{\sum_{x\in \CC_{\d^*}(I)} |D(x)|-\E[|D(x)|]}
{\sqrt{|I|}} \le \frac {\sqrt{64(2+p)}}{\g} \sqrt{\g \log
\frac{1}{\g}} \Eq(4.0051)
$$
and, uniformly with respect to all intervals
 $I\subseteq [-\g^{-p},\g^{-p}]$, 
$$
\sup_{\s_I \in \{-1,+1\}^{I}}\g \left| \sum_{x\in
\CC_{\d^*}(I)} 2\th \l(x) \sum_{i\in D(x)} \s_i\right| 
\le 2\th
\left( |I| \sqrt{\frac{\g}{\d^*}} + \sqrt {64(p+2)} \sqrt{|I|\g
\log\frac 1\g} \right)
\le  4\th |I| \sqrt
{\frac{\g}{\d^*}}.
\Eq(4.0052)
$$
} \medskip
 This Lemma  is a direct consequence of Lemma
\eqv(41), since  $
|2\th \l(x)\sum_{i\in D(x)} \s_i| \le 2\th (
|D(x)|-\E[|D(x)|])$ $+2\th \E[|D(x)|] 
$
, $|D(x)|=|\sum_{i\in A(x)}h_i|$, and $\E[|D(x)|]\le \sqrt{{\d^*}/{\g}}$ by
Schwarz inequality.

\vskip .5truecm For the function $V(m_I^{\d^*})$ in \eqv(3.140),
the previous rough estimate is useless. In  Theorem
\eqv(62P), with the help of the cluster expansion, we prove the
following

\medskip
\goodbreak
\noindent{\bf \Lemma(62)}{\it For any finite interval $I$, let
$$
\|\del_i V_I\|_\infty \equiv
 \sup_{(h,\tilde h): h_j=\tilde h_j,\forall j\neq i}
\frac{\left |  V_I (m^{\d^*}_I, h) -  V_I (m^{\d^*}_I,\tilde h)
\right|} {|h-\tilde h|}. \Eq (H.220)
 $$
Then,  for all  $\b>0$, for all $\d^*>\g> 0$,
such that
$$
\frac {(\d^*)^2 }{\g}\le  \frac{1}{6e^3\b} \Eq(H.2201)
$$
we have
$$
\sup_{I\subset \Z }\sup_{i\in I} \|\del_i V_I\|_\infty \le \frac 1
\b \frac {S }   {1-  S   }, \Eq (H.22)
 $$
where  $S$ is given in \eqv (paris1), 
  $  0<S \le 6 e^3 \b  \frac {(\d^*)^2} {\g} $. }
\vskip0.5cm

Together with the above estimates for $V_I$, we also need an
explicit expression for $\GG(m_I^{\d^*})$.
Since  $D(x)\subset B^{-\l(x)}(x)$, 
$ \GG_{x,m^{\d^*}(x)}\left(\l(x)\right)$, see \eqv(3.13P),  depends  
only on one component of $ m^{\d^*}(x)$, precisely on 
$m^{\d^*}_{\frac{3+\l(x)}{2}}$.
In fact, we have
$$
\GG_{x,m^{\d^*}(x)}\left(\l(x)\right)
= -\frac{1}{\b}\log
\frac
{\sum_{\s\in\{-1,+1\}^{B^{-\l(x)}(x)}}
\1_{\{m^{\d^*}_{\frac{3+\l(x)}{2}}(x,\s)=m^{\d^*}_{\frac{3+\l(x)}{2}}\}}
e^{2\b\th\l(x)\sum_{i\in D(x)}\s_i}}
{\sum_{\s\in\{-1,+1\}^{B^{-\l(x)}(x)}}
\1_{\{m^{\d^*}_{\frac{3+\l(x)}{2}}(x,\s)=m^{\d^*}_{\frac{3+\l(x)}{2}}\}}},
\Eq(3.13PM)
$$
since the sums over the spin configurations in $\{-1,+1\}^{B^{\l(x)}(x)}$
-- the ones that depend on $m^{\d^*}_{\frac{3-\l(x)}{2}}$-- 
cancel out between the numerator and denominator in \eqv(3.9P).

Depending on the  values of $m^{\d^*}_{\frac{3+\l(x)}{2}}$,
$\GG_{x,m^{\d^*}(x)}\left(\l(x)\right)$ has a behavior that
corresponds  to the classical Gaussian, Poissonian, or  
Binomial regimes, as explained in [\rcite{COP1}]. 
However,  as we shall  see in  Remark \eqv(P), 
we need accurate  estimates only in the  Gaussian regime.

Let
$g_0(n)$ be a positive  increasing real function with
$\lim_{n\uparrow \infty}g_0(n)=\infty$ such that $g_0(n)/n$ 
is decreasing to $0$ when $n\uparrow \infty$.

\noindent {\bf \Proposition(PP)}{\it
For all $\b,\th$ that satisfy \eqv(2.17),
there exist $\g_0=\g_0(\b,\th)$ and
$d_0(\b)>0$ such that for $0<\g\le \g_0$,
$\g/\d^*\le d_0(\b)$,
on the set
 $\{\sup_{x\in\CC_{\d^*}(I)} p(x) \le (2\g/{\d^*})^{1/4}\}$, if
$$
|m^{\d^*}_{\frac{3+\l(x)}{2}}(x)|\le
1-
\left(\frac{g_0(\d^*\g^{-1}/2)}{\d^*\g^{-1}/2}\vee \frac{16p(x)\b\th}
{1-\tanh(2\b\th)}\right),
\Eq(4.PM0)
$$
then
$$
\eqalign{
&\GG_{x,m^{\d^*}(x)}\left(\l(x)\right)=
-\frac{1}{\b}
\log
\frac {\Psi_{\l(x)2\b\th,p(x),m^{\d^*}_{\frac{3+\l(x)}{2}}(x)}}
 {\Psi_{0,0,m^{\d^*}_{\frac{3+\l(x)}{2}}(x)}}
\cr
&
-\frac{1}{\b}|D(x)|\left[ \log\cosh(2\b\th)+
\log\left(1+\l(x)m^{\d^*}_{\frac{3+\l(x)}{2}}(x)\tanh(2\b\th)\right)
+
\hat\varphi(
m^{\d^*}_{\frac{3+\l(x)}{2}}(x),2\l(x)\b\th, p(x) )\right],\cr
}\Eq(4.PM)
$$
where
$$
\left|\hat\varphi(
m^{\d^*}_{\frac{3+\l(x)}{2}}(x),2\l(x)\b\th, p(x) )\right|
\le
\left(\frac{2\g}{\d^*}\right)^{1/4}
 \frac{32\b\th(1+\b\th)}
{(1-|m^{\d^*}_{\frac{3+\l(x)}{2}}(x)|)^2
(1-\tanh(2\b\th))}
\Eq(5.PM)
$$
and
$$
\left|\log
\frac {\Psi_{\l(x)2\b\th,p(x),m^{\d^*}_{\frac{3+\l(x)}{2}}(x)}}
 {\Psi_{0,0,m^{\d^*}_{\frac{3+\l(x)}{2}}(x)}}\right|\le
\frac{18}{g_0(\d^*\g^{-1}/2)} + \left(\frac{2\g}{\d^*}\right)^{1/4}c(\b\th),
\Eq(6.PM)
$$
with  $c(\b\th)$ given in \eqv(3.25PM).}

\noindent{\bf \Remark(P10)}
Recalling \eqv(2.1700), we have 
$$
\Xi_2(x,\b\th,p(x))\equiv
-\l(x)\log\frac{ \Psi_{2\l(x)\b\th,p(x),\l(x)m^{\d^*}_{\b,2}}
\Psi_{0,0,-\l(x)m^{\d^*}_{\b,1}}}
{ \Psi_{2\l(x)\b\th,p(x),-\l(x)m^{\d^*}_{\b,1}} 
\Psi_{0,0,+\l(x)m^{\d^*}_{\b,2}}}
\Eq(nano)
$$
and choosing $g_0(n)=n^{\frac{1}{4}}$, \eqv(2.1700) follows 
from \eqv(6.PM). \eqv(2.1692) follows from \eqv(5.PM).

\noindent{\bf Proof:} The general strategy of the proof is similar 
to that of Proposition 3.1 in [\rcite{COP1}]. However, since there are 
important differences  we give some details.  We introduce the 
``grand canonical'' measure  on
$\{-1,+1\}^{B^{-\l(x)}(x)}$, with chemical potential $\nu\in \R$,
given by
$$
\E_{x,\nu}(f)=\frac {\E_{\s_{B^{-\l(x)}(x)}}\left[f(\s) e^{\nu
\sum_{i\in B^{-\l(x)}(x)}\s_i}\right]}
{\E_{\s_{B^{-\l(x)}(x)}}\left[ e^{\nu \sum_{i\in
B^{-\l(x)}(x)}\s_i}\right]} \Eq(3.14PM)
$$
{ where $\E_{\s_{B^{-\l(x)}(x)}}$ is the Bernoulli uniform on
$\{-1,+1\}^{B^{-\l(x)}(x)}$.}
Then defining
$$
\Psi_{\l(x)2\b\th, p(x), m^{\d^*}_{\frac{3+\l(x)}{2}}(x)}
\equiv
\frac
{\E_{x,\nu_2}\left[e^{\l(x)2\b\th\sum_{i\in D(x)} \s_i}
\1_{\{
\sqrt{\frac{2\g}{\d^*}}\sum_{i\in B^{-\l(x)}(x)}
(\s_i-m^{\d^*}_{\frac{3+\l(x)}{2}}(x))=0
\}}\right]}
{\E_{x,\nu_2}\left[e^{\l(x)2\b\th\sum_{i\in D(x)} \s_i}\right]}
\Eq(3.18PM)
$$
and
$$
\eqalign{
&\phi(m^{\d^*}_{\frac{3+\l(x)}{2}}(x),\l(x)2\b\th,p(x))\cr
&\quad\equiv
\frac{\d^*}{2\g}
\left((\nu_1-\nu_2)m^{\d^*}_{\frac{3+\l(x)}{2}}(x)+
p(x)\log \frac{\cosh(\nu_2+\l(x)2\b\th)}{\cosh(\nu_1)}
+(1-p(x))
\log \frac{\cosh(\nu_2)}{\cosh(\nu_1)}\right),\cr
}
\Eq(3.19PM)
$$
a simple computation gives
$$
\GG_{x,m^{\d^*}(x)}\left(\l(x)\right)
=
-\frac 1\b \log
\frac
{\Psi_{\l(x)2\b\th, p(x), m^{\d^*}_{\frac{3+\l(x)}{2}}(x)}}
{\Psi_{0, 0, m^{\d^*}_{\frac{3+\l(x)}{2}}(x)}}
-\frac 1 \b \phi(m^{\d^*}_{\frac{3+\l(x)}{2}}(x),\l(x)2\b\th,p(x)).
\Eq(3.20PM)
$$

We choose $\nu_1$ such that $m^{\d^*}_{\frac{3+\l(x)}{2}}(x)=\tanh \nu_1$
and $\nu_2$ such that
$$
m^{\d^*}_{\frac{3+\l(x)}{2}}(x)=p(x) \tanh(\nu_2+\l(x)2\b\th)
+(1-p(x))\tanh \nu_2.
\Eq(3.16PM)
$$

By using elementary  formulae on hyperbolic tangents and  cosines,
one can  check the following identity  
$$
\eqalign{
&\phi(m^{\d^*}_{\frac{3+\l(x)}{2}}(x),\l(x)2\b\th,p(x))\cr
&=|D(x)|
\left[\log\cosh 2\b\th+\log\left(1+\l(x)
m^{\d^*}_{\frac{3+\l(x)}{2}}(x)
\tanh(2\b\th)\right)+
\hat\varphi(m^{\d^*}_{\frac{3+\l(x)}{2}}(x),
2\l(x)\b\th, p(x) )\right],\cr
}\Eq(3.21PM)
$$
where
$$
\eqalign{
&
|D(x)| \hat\varphi(
m^{\d^*}_{\frac{3+\l(x)}{2}}(x),2\l(x)\b\th, p(x) )\cr
&\quad\quad =
\frac{\d^*}{2\g}
(\nu_1-\nu_2)
m^{\d^*}_{\frac{3+\l(x)}{2}}(x)
+
\frac{\d^*}{2\g}
\log  \left(1+
m^{\d^*}_{\frac{3+\l(x)}{2}}(x)
\tanh(\nu_2-\nu_1)\right)\cr
&\quad\quad+
\frac{\d^*}{2\g}
\log\cosh(\nu_2-\nu_1)\cr
&\quad\quad+
\frac{\d^*}{2\g}
p(x) \log \left[
1+ \frac{\l(x)\tanh(2\b\th)(1-
(m^{\d^*}_{\frac{3+\l(x)}{2}}(x))
^2)\tanh(\nu_2-\nu_1)}
{\left(1+
\l(x)m^{\d^*}_{\frac{3+\l(x)}{2}}(x)
\tanh(2\b\th)\right)\left(1+
m^{\d^*}_{\frac{3+\l(x)}{2}}(x)
\tanh(\nu_2-\nu_1)\right)}\right].
\cr
}\Eq(mars.68P)
$$

To study \eqv(3.21PM), we need  extensions of results
proved in [\rcite{COP1}].
Defining
$$
\s^2_{\l(x)2\b\th}
\equiv
p(x) \frac 1{\cosh^2(\nu_2+\l(x)2\b\th)}+
(1-p(x)) \frac 1{\cosh^2(\nu_2)},
\Eq(3.22PM)
$$
and  using again elementary formulae on hyperbolic
 tangents and cosines one can check that
$$
\s^2_{\l(x)2\b\th}=\left(1-(m^{\d^*}_{\frac{3+\l(x)}{2}}(x))^2\right)
\left[1-p(x)(1-p(x))S(p(x),m^{\d^*}_{\frac{3+\l(x)}{2}}(x))\right],
\Eq(3.23PM)
$$
where
$$
0\le
S(p(x),m^{\d^*}_{\frac{3+\l(x)}{2}}(x))
\le
\left(1-(m^{\d^*}_{\frac{3+\l(x)}{2}}(x))^2\right)
c(\b\th),
\Eq(3.24PM)
$$
with
$$
c(\b\th)=\frac{\tanh^2(2\b\th)(1+\tanh^2(2\b\th))^2}
{[1-\tanh^2(2\b\th)]^2[1-\tanh(2\b\th)]^6}.
\Eq(3.25PM)
$$

Assuming that $\g/\d^*<d_0(\b)$ for some well chosen $d_0(\b)$, and
following the arguments of  the proof
of Lemma 3.3 in [\rcite{COP1}], we  check that
$$
|\nu_2-\nu_1|\le
\frac{4p(x)\b\th}{
1-(m^{\d^*}_{\frac{3+\l(x)}{2}}(x))^2}.
\Eq(3.27PM)
$$
Using the fact that \eqv(4.PM0) implies that
$
\frac{4p(x)\b\th}{
(1-(m^{\d^*}_{\frac{3+\l(x)}{2}}(x))^2)
(1-\tanh(2\b\th))}\le \frac 14,
$
 recalling  \eqv(3.21PM), and  using Taylor expansion we get
$$
\eqalign{
&\left|
\frac{\phi(m^{\d^*}_{\frac{3+\l(x)}{2}}(x),\l(x)2\b\th,p(x))
}{ 
\frac{\d^*}{2\g}
}
- p(x)\left[\log\cosh 2\b\th+\log\left(1+\l(x)
m^{\d^*}_{\frac{3+\l(x)}{2}}(x)
\tanh(2\b\th)\right)\right]\right|\cr
&
\quad\quad\quad \le \frac{32p^2(x)\b\th(1+\b\th)}
{(1-|m^{\d^*}_{\frac{3+\l(x)}{2}}(x)|)^2
(1-\tanh(2\b\th))}.\cr
}\Eq(3.29PM)
$$
A short computation  concludes  the proof of Proposition \eqv(PP). \eop
\medskip

To prove Theorem \eqv(61), we need results that have been proven 
in [\rcite{COP1}]. We  first define the subsets of 
the complementary 
of $\AA(\D_Q)$ which  will be treated in a  similar 
way to that in  [\rcite{COP1}].

 Let $\D_L=[\ell_1,\ell_2]$ be an interval of length
 $L=\ell_2-\ell_1 \in \N$.
Let $\d>\d^*$, $\z_4>\z_1>8\g/\d^*$ be positive  real numbers. 

\Definition(def3)
{\it We set 
\def\OO{{\cal O}}
$$
\OO_0^{\d,\z_1}(\D_L)
\equiv
\left\{\eta^{\d ,\z_1}(\ell)=0,\quad \forall \ell \in \D_L \cap
\Z\right\}. \Eq (ME.2)
$$
Taking $\tilde L \le L$ a positive  integer, 
 let  $  \D_{\tilde L} =  [\tilde \ell_1,  \tilde \ell_2] $,    
$ \D_ { \tilde L} \subset  \D_L$. Define  for
$\eta=+1$ or $\eta=-1$. 
$$
  \RR_{0,\eta}^{\d,\z_1,\z_4}( \D_L,\tilde L )  \equiv 
\big\{  
\eta^{\d,\z_1}(\ell_1)=\eta^{\d,\z_1}(\ell_2) =\eta;\, \big \} \cap 
\OO_0^{\d,\z_4}([\ell_1+1, \ell_2-1]) \cap 
\OO_0^{\d,\z_1}(\D_{\tilde L})   \Eq(5.41) 
$$
and
$ \RR_{0}^{\d,\z_1,\z_4}(\D_L,\tilde L ) \equiv \RR_{0,+
}^{\d,\z_1,\z_4}(\D_L,\tilde L )\cup \RR_{0,-
}^{\d,\z_1,\z_4}(\D_L,\tilde L ).
$}

\medskip
\noindent Note that $\RR_{0,\eta}^{\d,\z_1,\z_4}(\D_L,\tilde L )$
decreases in $\tilde L$, therefore 
$\cup_{\tilde L:1\le \tilde L\le L}
\RR_{0,\eta}^{\d,\z_1,\z_4}(\D_L,\tilde L)=
\RR_{0,\eta}^{\d,\z_1,\z_4}(\D_L,1 )$.

We set 
$$
\RR_{0}^{\d,\z_1,\z_4}(I)\equiv \bigcup_{L:\;2\le L\le
|I|}\;\bigcup_{\D_L\subset I}\; \RR_{0}^{\d,\z_1,\z_4}(\D_L,1),
\Eq(5.48)
$$
$$
\OO_0^{\d,\z_1}(I, R_1)\equiv \bigcup_{R:\;R_1\le R\le |I|}\;
\bigcup_{\D_R\subset I} \OO_0^{\d,\z_1}(\D_{R}), \Eq(5.46)
$$
and recalling Definition \eqv(def2), 
$$
\WW^{\z_1,\z_4}(I,L_2)\equiv \bigcup_{L:\;2\le L\le L_2}\;
\bigcup_{\D_L\subset I}  \WW^{\z_1,\z_4}(\D_L). \Eq(5.47)
$$ 

\medskip
\noindent{\bf \Theorem(t62)} {\it Given  $\b,\th $ as in \eqv
(2.17), there exist $\g_0=\g_0(\b,\th)>0$, $d_0=d_0(\b,\th)>0$,
and $0<\z_0(\b,\th)<1$ such that if $0<\g\le \g_0$, $\d^*
>\g$, $ \g/\d^* \le d_0$, and  $p$ is a positive integer such that
$$
(p+2) \d^* \log \frac 1 \g \le \frac 1{64}
\Eq(5.448)
$$
there exists  $\O_{RE}=\O_{RE}(\g,\d^*, p)$ with $\P[\O_{RE}] \ge
1- \g^{2}$, such that for all $\d,\z_1 ,\z_4$ with
 $1>\d>\d^*>0$,  $\z_0(\b,\th)>\z_4> \z_1>8\g/\d^*
$,
and
$$
\d \z_1^3 > \frac{128(1+\th)}{\k(\b,\th)} (\d^*\vee
\sqrt{\frac{\g}{\d^*}}), \Eq(5.480)
$$
$$
\d \z_4^3 > \frac{32}{\k(\b,\th)} \z_1, \Eq(5.0480)
$$
where $\k(\b,\th)>0$ satisfies \eqv(2.19), on $\O_{RE}$ we
have
$$
\mu_{\b,\th,\g}\left(\cup_{I\subset[-\g^{-p},\g^{-p}]}\Big(
\OO_0^{\d,\z_1}(I,R_1) \cup \WW^{\z_1,\z_4}(I,L_2)  \cup
\RR_{0}^{\d,\z_1,\z_4}(I)\Big) \right) \le \frac {3^4}{\g^{5p}}
e^{ -\frac \b \g\left[ \big(\frac{
\k(\b,\th)}{4}\d\z_4^3\big)\wedge \big(\FF^*\big)\right]},
\Eq(5.481)
$$
with $\FF^*$ given in  \eqv(min),
$$
R_1=\frac{4(5+\FF^*)}{\k(\b,\th)\d\z_1^3}, \Eq(5.4810)
$$
and}
$$  L_2=\frac{\FF^*}{64(1+\th) }\frac {1}{\d^*\vee \sqrt{
\frac{\g}{\d^*}}}\,. \Eq(5.482)
$$

The proof of Theorem \eqv(t62) is the same as the proof of
Corollary 5.2, Corollary 5.4, and Corollary 5.6 in [\rcite{COP1}],
with $\D\FF$ in 
[\rcite{COP1}] is equal to $2\FF^*$ here. Moreover with a little work, one can 
make explicit the constants depending  on $\b,\th$ that 
appear  in [\rcite{COP1}].
Note  that the condition  \eqv(2.17) on $\b,\th$ is
weaker than the condition used in [\rcite{COP1}], however this will
make no difference at all since we just use the rough estimate, see 
Lemma \eqv(60002) to treat the random field.

Let 
$
\BB_0([-\g^{-p},\g^{-p}],R_1,L_2)\equiv
\cap_{I\subset[-\g^{-p},\g^{-p}]}\big( \OO_0^{\d,\z_1}(I) \cup
\WW^{\z_1,\z_4}(I,L_2)  \cup  \RR_{0}^{\d,\z_1,\z_4}(I)\Big) \big)^c.
$
 On this set we can only have runs of
$\eta^{\d,\z_1}=0$, with length at most $R_1$ and runs of
$\eta^{\d,\z_4}(\ell)=\eta\in \{-1,+1\}$, with length at least
$L_2$. The next step is to prove that the length of the previous
runs of $\eta^{\d,\z_4}=\eta\in \{-1,+1\}$ is indeed bounded from
below by $\e/\g$.

\Definition(def7) {\it For $\eta\in \{+1,-1\}$, $ \ell_1< \tilde
\ell_1< \tilde \ell_2< \ell_2$ with $ 3 \le \tilde \ell_1
-\ell_1\le R_1$ $3\le  \ell_2-\tilde\ell_2 \le R_1$, let}
$$
\eqalign{
 \widetilde \WW_{\eta}^{ \z_1,\z_4}(\ell_1,\tilde
\ell_1,\tilde\ell_2,\ell_2) \equiv \big\{
m^{\d^*}_{[\ell_1,\ell_2]}& \colon \eta^{\d,\z_1}
(\ell_1)=\eta^{\d,\z_1}(\ell_1 +1)= \eta^{\d,\z_1}(\ell_2 -1
)=\eta^{\d,\z_1}(\ell_2)=\eta, \cr &
\quad\quad\eta^{\d,\z_4}(\ell)=-\eta,\,  \forall \ell\in
[\tilde\ell_1-1,\tilde \ell_2+1]
   \big\}.\cr }\Eq(5.49)
$$

\noindent{\bf \Proposition (59)} {\it  Let $\b,\th$  satisfy
\eqv(2.17). We take $\k(\b,\th)>0$ as in \eqv(2.19),
$\FF^*>0$ as in \eqv(min),  $V(\b,\th)$ as  in
\eqv(2.1701), and $c(\b)$ as in \eqv(3.25PM). 
There exist $\g_0=\g_0(\b,\th)>0$,
$d_0=d_0(\b,\th)>0$,
 and $0<\z_0=\z_0(\b,\th)<1$
such that if $ 0<\g\le \g_0$, $\d^*>\g$, $\g/\d^* \le d_0$, and
$0<\z_1< \z_4<\z_0$, $1>\d>\d^*>0$  verify the following conditions
$$
\d\z_1^3 \ge \frac {128(1+\th)(5+\FF^*)}{\k(\b,\th) \FF^*}
\sqrt{\frac{\g}{\d^*}}
\Eq(5.5000100)
$$
$$
\z_1\ge \left(5184 (1+c(\b\th))^2 \sqrt{\frac{\g}{\d^*}}\right)
\vee \left(12\frac{e^3\b}{c (\b,\th)} \frac{(\d^*)^2}{\g}\right)^2
\Eq(5.5000101)
$$
for a  constant $c(\b,\th)$ given in \eqv(mars.76),
if $\D_Q$ is an interval containing the origin, of length $Q/ \g$
in macroscopic units, { with}
$$
\sqrt {\g} \log Q \le \frac{\sqrt{6e^3\b}}{256}, \Eq(5.50001011)
$$
and $\e>\g \d^*$, then there exists  $ \O_4
=\O_4(\b,\th,\g,\z,\d,\D_Q,\e) $ with
$$
\P[\O_4  ]\ge 1-3 \g^{ 2} -  \frac{2 Q
 }{\e}
   e^{-
\sfrac{( \FF^*)^2}{\e 2^{11} \z_4 c^2(\b,\th)}}
 -  \frac{4  Q
 }{ \e}e^{-\sfrac{(\FF^*)^2}{\e2^{6}V^2(\b,\th)} }
\Eq(5.501)
$$
such that on $\O_4$,  we have, for $\eta=\pm 1$
 $$
\mu_{\b,\th,\g}\left( \cup^{*}_{[\ell_1,\ell_2]
\subset I \subset \D_Q}\cup^{**}_{[\tilde
\ell_1,\tilde \ell_2] \subset [\ell_1,\ell_2]}
\widetilde \WW_{\eta}^{\z_1,\z_4}
(\ell_1,\tilde \ell_1,\tilde\ell_2,\ell_2)  \right )
\le \frac {R_1^2 
Q}{\g^3} e^{-\frac \b \g \FF^* }. \Eq (5.51)
$$
In \equ(5.51), the union $\cup^*$ has the constraint $|I|=\e/\g$
while $\cup^{**}$ refers to the extra constraints $ 2 \le \tilde
\ell_1 -\ell_1 \le R_1$, $\ell_2-\tilde\ell_2 \le R_1$, with $R_1$
given by \eqv(5.4810).
 }

\noindent{\bf Remark.}
\item{$\bullet$} The constraint
\eqv(5.50001011) is present since  we use the rough estimate,
Lemma \eqv(60002), to control some terms. Note that taking   $p=2+
[{\log Q}/{\log(1/\g)}]$, \eqv(5.5000101) and   \eqv(5.50001011)
imply   $64 (p+2) \d^* \log (1/\g) \le 1$,  which is the condition
\eqv(4.00510) for the rough estimate. We will see that $\O_4 
\subset \O_{\rm RE}$.

\item{$\bullet$} The constraint $\ell_2-\ell_1 \le \epsilon
\g^{-1}$ enters into play in \equ(5.501), giving  the terms 
proportional to $\e^{-1}$ into the exponential.

\item{$\bullet$}
The uniformity with respect to the  intervals inside $\D_Q$ gives
the {prefactors}  $\frac {Q}{\e}$ in \eqv(5.501) and not
$\frac {Q}{ \g}$, since a maximal inequality is  used. The union
in \eqv(5.51) contains at most $R_1^2\e^2 Q\g^{-3}$ terms.
\medskip

\proof { We split  it in 4 steps.}
\medskip
\noindent{\bf Step 1: reduction to finite volume}

Recalling \eqv(5.49), we define  
$$  \RR(\eta)\equiv\RR^{\d,\z_4}(\eta)\equiv
\RR^{\d,\z_4}(\tilde \ell_1,\tilde\ell_2,\eta) =
\left\{ m^{\d^*}_{[\tilde \ell_1,\tilde \ell_2]}\colon
 \eta^{\d,\z_4}(\ell)=\eta,
\forall \ell \in [\tilde  \ell_1,  \tilde \ell_2]\right\},
\Eq(2.29b)
$$
and
$$
\WW_{\eta}^{\z_1,\z_4} (\ell_1+1,\tilde \ell_1,\tilde\ell_2,\ell_2-1)
\equiv\left\{
\eta^{\d,\z_1} (\ell_1+1)=\eta^{\d,\z_1}(\ell_2-1)=\eta \right\}
\bigcap \RR(-\eta).
\Eq(MRS1)
$$
We can write
$$
\widetilde \WW_{\eta}^{\z_1,\z_4} (\ell_1,\tilde \ell_1,\tilde\ell_2,\ell_2)
=\left\{
\eta^{\d,\z_1} (\ell_1)=\eta^{\d,\z_1}(\ell_2)=\eta \right\}
\bigcap
\WW_{\eta}^{\z_1,\z_4} (\ell_1+1,\tilde \ell_1,\tilde\ell_2,\ell_2-1).
\Eq(MRS2)
$$
Let us first consider a volume $\L$ such that $\g\L \supset
\D_Q$. Recalling \eqv(2.2) and \eqv(2.300), multiplying and
dividing by $Z_{\b,\th,\g,\g^{-1}[\ell_1+1,\ell_2-1]}
^{\s_{\g^{-1}\del[\ell_1+1,\ell_2-1]}}$ we have
$$
\eqalign{
&\mu_{\b,\th,\g,\L}\left( \widetilde \WW_{\eta}^{\z_1,\z_4} (\ell_1,\tilde
\ell_1,\tilde\ell_2,\ell_2)
 \right)
 =\cr
&\quad \frac{1}{Z_{\b,\th,\g,\L}} \sum_{\s_{\L\setminus
\g^{-1}[\ell_1+1,\ell_2-1]}} e^{-\b H(\s_{\L\setminus
\g^{-1}[\ell_1+1,\ell_2-1]})} \1_{\{ \eta^{\d,\z_1}
(\ell_1)=\eta^{\d,\z_1}(\ell_2)=\eta \}}
Z_{\b,\th,\g,\g^{-1}[\ell_1+1,\ell_2-1]}
^{\s_{\g^{-1}\del[\ell_1+1,\ell_2-1]}} \cr
&\quad\sum_{\s_{[\ell_1+1,\ell_2-1]}} \1_{ \WW_{\eta}^{\z_1,\z_4}
(\ell_1+1,\tilde \ell_1,\tilde\ell_2,\ell_2-1)} \frac{ e^{-\b
H_\g(\s_{\g^{-1}[\ell_1+1,\ell_2-1]}) -\b
W_{\g}(\s_{\g^{-1}[\ell_1+1,\ell_2-1]},\s_{\g^{-1}\del
[\ell_1+1,\ell_2-1]} ) } }
{Z_{\b,\th,\g,\g^{-1}[\ell_1+1,\ell_2-1]}
^{\s_{\g^{-1}\del[\ell_1+1,\ell_2-1]}} }. }\Eq(5.51P)
$$
 Since
 $\eta^{\d,\z_1}(\ell_1)=\eta^{\d,\z_1}(\ell_1 +1)=
\eta^{\d,\z_1}(\ell_2 -1)=
\eta^{\d,\z_1}(\ell_2)=\eta$, using \eqv(3.7P) and recalling \eqv(3.8P)
,  we get
$$
\eqalign{ & \sum_{\s_{[\ell_1+1,\ell_2-1]}}
 \1_{\WW_{\eta}^{\z_1,\z_4} (\ell_1+1,\tilde
 \ell_1,\tilde\ell_2,\ell_2-1)}
\frac{ e^{-\b H_\g(\s_{\g^{-1}[\ell_1+1,\ell_2-1]}) -\b
W_{\g}(\s_{\g^{-1}[\ell_1+1,\ell_2-1]},\s_{\g^{-1}\del
[\ell_1+1,\ell_2-1]} ) } }
{Z_{\b,\th,\g,\g^{-1}[\ell_1+1,\ell_2-1]}
^{\s_{\g^{-1}\del[\ell_1+1,\ell_2-1]}} }
\cr &\quad \le
e^{+\frac{\b}{\g} 4(\d^*+\z_1)}
\mu_{\b,\th,\g}\left(
\1_{ \WW_{\eta}^{\z_1,\z_4} (\ell_1+1,\tilde \ell_1,\tilde\ell_2,\ell_2-1)}
\big| \S^{\d^*}_{\del[\ell_1+1,\ell_2-1]}\right)
(m^{\d^*}_{\del[\ell_1+1,\ell_2-1]}=m_{\eta})\cr
} \Eq(MRS3)
$$
where $m_+$ ($m_-$)
is the constant function on $\partial^+I$ or $\partial^-I$ with
value $m_{\b}^{\d^*}$ (resp. $Tm_{\b}^{\d^*}$).

Notice  that for any $\L$ such that $\g\L\supset \D_Q$
$$
\eqalign{ &\frac{1}{Z_{\b,\th,\g,\L}} \sum_{\s_{\L\setminus
\g^{-1}[\ell_1+1,\ell_2-1]}} e^{-\b H(\s_{\L\setminus
\g^{-1}[\ell_1+1,\ell_2-1]})} \1_{\{ \eta^{\d,\z_1}
(\ell_1)=\eta^{\d,\z_1}(\ell_2)=\eta \}}
Z_{\b,\th,\g,\g^{-1}[\ell_1+1,\ell_2-1]}
^{\s_{\g^{-1}\del[\ell_1+1,\ell_2-1]}}\cr &\quad \le
\mu_{\b,\g,\th,\L}( \1_{\{ \eta^{\d,\z_1}
(\ell_1)=\eta^{\d,\z_1}(\ell_2)=\eta \}})\le 1. }\Eq(MRS5)
$$

Therefore, inserting \eqv(MRS3) in
\eqv(5.51P) and  taking the limit $\L\uparrow \Z$ we get
 $$
\eqalign{&
\mu_{\b,\th,\g}\left( \tilde \WW_{\eta}^{\z_1,\z_4} (\ell_1,\tilde
\ell_1,\tilde\ell_2,\ell_2)
 \right)
\cr &\quad\quad
\le e^{\frac \b \g 4(\z_1+\d^*)}
\mu_{\b,\th,\g}\left(
\WW_{\eta}^{\z_1,\z_4} (\ell_1+1,\tilde
\ell_1,\tilde\ell_2,\ell_2-1)
\bigm|\S^{\d^*}_{\partial [\ell_1+1,\ell_2-1]}\right)
(m^{\d^*}\del[\ell_1+1,\ell_2-1]=m_\eta).
}
\Eq(5.52)
$$

To continue, recalling \eqv(3.8P1) and writing   $m^{\d^*}_{\del
I}= (m^{\d^*}_{\del^- I},m^{\d^*}_{\del^+ I})$, we set simply
$$
Z_{\b,\th,\g,I}\left(
m^{\d^*}_{\del^- I}=m_{s_1}
,m^{\d^*}_{\del^+ I}=m_{s_2}\right)\equiv
Z_I^{m_{s_1},m_{s_2}}
\Eq(5.52P1)
$$
when $(m_{s_1},m_{s_2}) \in \{m_-,0,m_+\}^2$ 
where $m_+$  and $m_-$ are as above, and for  $m_{s_1}=0$, we set 
in \eqv(3.8P1)
$E(m^{\d^*}_I,m^{\d^*}_{\del^- I})=0$ while 
for $m_{s_2}=0$ we set $E(m^{\d^*}_I,m^{\d^*}_{\del^+ I})=0$.
In a similar way, recalling \eqv(3.15P), if $F$ is
$\S_I^{\d^*}$--measurable we set 
$$
\frac{Z_I^{m_{s_1},m_{s_2}}(F)}{Z_I^{m_{s_1},m_{s_2}}} \equiv
\frac{ \sum_{m^{\d^*}_{I} \in {\cal M}_{\d^*}(I)} F(m^{\d^*}_I)
e^{-\frac \b \g \left\{ {\widehat
\FF}(m^{\d^*}_I|m^{\d^*}_{\partial^{-}I}=m_{s_1}
m^{\d^*}_{\partial^{ +}I}=m_{s_2}) +\g \GG(m^{\d^*}_I)+ \g  V
(m^{\d^*}_I)\right\}}} {Z_I^{m_{s_1},m_{s_2}}}. \Eq(5.52P4)
$$

 Using the fact that $\eta^{ \d,\z}(\tilde
\ell_1)=\eta^{\d,\z} (\tilde \ell_1-1)$ and $\eta^{\d,\z}(\tilde
\ell_2 +1)=\eta^{\d,\z}(\tilde \ell_2)$ we can decouple the
contribution coming from the interval $[\tilde \ell_1-1, \tilde
\ell_2+1]$ and restrict the configuration in the denominator in a
suitable way to get
$$
\eqalign{
&\mu_{\b,\th,\g}\left(
\WW_{\eta}^{\z_1,\z_4} (\ell_1+1,\tilde
\ell_1,\tilde\ell_2,\ell_2-1)
\bigm|\S^{\d^*}_{\partial [\ell_1+1,\ell_2-1]}\right)
(m^{\d^*}\del[\ell_1+1,\ell_2-1]=m_\eta)\cr
&\quad \le e^{\frac \b\g 8\z_1}
\frac {Z^{m_{\eta},m_{-\eta}}_{[\ell_1+1,\tilde
\ell_1-1]}(\eta^{\d,\z_1}(\ell_1+1)=\eta)}
{Z^{m_{\eta},m_{\eta}}_{[\ell_1+1,\tilde \ell_1
-1]}(\eta^{\d,\z_1}(\ell_1+1)=\eta)} \frac {Z^{0,0}_{[\tilde
\ell_1,\tilde \ell_2]}(\RR(-\eta))} {Z^{0,0}_{[\tilde
\ell_1,\tilde \ell_2]}(\RR(\eta))} \frac {Z^{m_{-\eta},m_{\eta}}_{[\tilde
\ell_2+1, \ell_2-1]}(\eta^{\d,\z_1}(\ell_2-1)=\eta)}
{Z^{m_{\eta},m_{\eta}}_{[\tilde
\ell_2+1,\ell_2-1]}(\eta^{\d,\z_1}(\ell_2-1)=\eta)}.
}\Eq(5.52P6)$$

The first and the third ratio on the right hand side of \eqv(5.52P6)
are easily estimated. Since 
$0<  \ell_1 -\tilde \ell_1 \le R_1$, $0< \ell_2-\tilde\ell_2 \le R_1$
with $R_1$ given by \eqv(5.4810), 
using the rough estimate Lemma \eqv(60002), it can be checked 
that  on $\O_{RE}$,
 uniformly over  all intervals $[\ell_1,\tilde \ell_1] 
\subset [-\g^{-p},\g^{-p}]$, we have
$$
\frac {Z^{m_{\eta},m_{-\eta}}_{[\ell_1+1,\tilde
\ell_1-1]}(\eta^{\d,\z_1}(\ell_1+1)=\eta)}
{Z^{m_{\eta},m_{\eta}}_{[\ell_1+1,\tilde \ell_1
-1]}(\eta^{\d,\z_1}(\ell_1+1)=\eta)}
 \le e^{\sfrac \b \g
(8(1+\th)R_1\sqrt{\frac{\g}{\d^*}})} 
\frac{
e^{-\frac \b \g
\inf_{m^{\d^*}_{[\ell_1+1,\tilde \ell_1-1]}\in
\{\eta^{\d,\z_1}(\ell_1+1)=\eta\}} \tilde
\FF(m^{\d^*}_{[\ell_1+1,\tilde \ell_1-1]}|m_{\eta},m_{-\eta}) }}
{e^{-\frac \b \g
\tilde
\FF(T^{\frac{1-\eta}{2}}
m^{\d^*}_{\b,[\ell_1+1,\tilde \ell_1-1]}|m_{\eta},m_{\eta}) }},
\Eq(5.52P7)
$$
where $\tilde \FF(\cdot)$ is given in \equ(5.Rio2) 
{ and we have used 
the fact 
that since $ \widetilde m^{\d^*}_\b=-\widetilde {Tm^{\d^*}_\b}$
the boundary terms,  see  \eqv(5.Rio3),
$$
\frac{\d^*}{2} \sum_{y \in \CC_{\d^*}(\del[\ell_1+1,\tilde \ell_1-1])}
 \big[\tilde m^{\d^*}(y)\big]^2 
\sum_{x \in \CC_{\d^*}([\ell_1+1,\tilde \ell_1-1])}
J_{\d^*}(x-y) 
\Eq(5.maggio)
$$
cancel between the  numerator and the denominator in \eqv(5.52P7).
}

It can be proved that 
$$
\eqalign{
&\inf_{1 \le \tilde \ell_1 -\ell_1 \le R_1}
\inf_{m^{\d^*}_{[\ell_1+1,\tilde \ell_1-1]}\in
\{\eta^{\d,\z_1}(\ell_1+1)=\eta\}} \tilde
\FF(m^{\d^*}_{[\ell_1+1,\tilde \ell_1-1]}|m_{\eta},m_{-\eta})
-
\tilde
\FF(T^{\frac{1-\eta}{2}}
m^{\d^*}_{\b,[\ell_1+1,\tilde \ell_1-1]}|m_{\eta},m_{\eta}) 
\cr
&\qquad
 \ge \FF^*
-(4L_0+2R_1)(1+\th)\left(\d^*\vee\sqrt{\frac{\g}{\d^*}}\right), \cr
}
\Eq(5.52P8)
$$
where $\FF^*$ is defined in \eqv(min) and 
$L_0=\frac{2}{\a(\b,\th)}\log\frac{\d^*}{\g}$
with $\a(\b,\th)$ as in \eqv(Palmares).
 A similar argument can
be used for the third ratio  in \eqv(5.52P6), { and} we get
$$
\frac {Z^{m_{\eta},m_{-\eta}}_{[\ell_1+1,\tilde
\ell_1-1]}(\eta^{\d,\z_1}(\ell_1+1)=\eta)}
{Z^{m_\eta,m_{\eta}}_{[\ell_1+1,\tilde \ell_1
-1]}(\eta^{\d,\z_1}(\ell_1+1)=\eta)} \frac {Z^{m_{-\eta},m^{\eta}}_{[\tilde
\ell_2+1, \ell_2-1]}(\eta^{\d,\z_1}(\ell_2-1)=\eta)}
{Z^{m_{\eta},m_{\eta}}_{[\tilde
\ell_2+1,\ell_2-1]}(\eta^{\d,\z_1}(\ell_2-1)=\eta)} \le
e^{-\frac{\b}{\g}(2\FF^*-32(1+\th)(R_1+L_0) \sqrt\frac{\g}{\d^*})}.
 \Eq(5.52P9)
$$

It remains to treat the second ratio in \eqv(5.52P6), that is 
$$
\frac
{Z^{0,0}_{[\tilde \ell_1,\tilde \ell_2]}(\RR(-\eta))}
{Z^{0,0}_{[\tilde \ell_1,\tilde \ell_2]}(\RR(\eta))}
\equiv \frac {\sum_{m^{\d^*}_{\tilde I_{12}} \in {\cal
M}_{\d^*}(\tilde I_{12})} \1_{\{\RR(-\eta)\}} e^{-\frac
\b \g\left\{{\widehat \FF}(m^{\d^*}_{\tilde
I_{12}}|0)+{\g}\GG(m^{\d^*}_{\tilde I_{12}}) +  \g  V
(m^{\d^*}_{\tilde I_{12}})\right\}}  } {\sum_{m^{\d^*}_{\tilde
I_{12}} \in {\cal M}_{\d^*}(\tilde I_{12})} \1_{\{\RR(\eta)\}}
e^{-\frac \b\g\left\{{\widehat \FF}(m^{\d^*}_{\tilde
I_{12}}|0) +\g\GG(m^{\d^*}_{\tilde I_{12}})+ \g  V
(m^{\d^*}_{\tilde I_{12}})\right\}}, } \Eq(4.S4) 
$$
where ${\widehat \FF_\g}(m^{\d^*}_{\tilde I_{12}}|0)$ is as
\eqv(3.11P) for $I=\tilde I_{12}=[\tilde \ell_1,\tilde \ell_2]$
but with the term
$E(m^{\d^*}_I,m^{\d^*}_{\del I})\equiv 0$ and, recalling
\equ(4.5P), we have set $T\RR(\eta)=\RR(-\eta)$ and
 $ \1_{\{\RR(-\eta)\}} \equiv \1_{\{\RR(-\eta)\}}(m^{\d^*}_{\tilde
 I_{12}})$.

Notice that if we flip $h_i$ to  $-h_i$, for all $i$,  then
$\l(x)\rightarrow -\l(x)$, $B^+(x)\rightarrow B^-(x)$ while
$|D(x)|$ does not change. Therefore,
$$
\frac {Z^{0,0}_{[\tilde \ell_1,\tilde \ell_2]}(\RR(-\eta))}
{Z^{0,0}_{[\tilde \ell_1,\tilde \ell_2]}(\RR(\eta))}(h) = \frac
{Z^{0,0}_{[\tilde \ell_1,\tilde \ell_2]}(\RR(\eta))}
{Z^{0,0}_{[\tilde \ell_1,\tilde \ell_2]}(\RR(-\eta))}(-h),
\Eq(4.7) 
$$ which implies  that $\log
\frac {Z^{0,0}_{[\tilde \ell_1,\tilde \ell_2]}(\RR(-\eta))}
{Z^{0,0}_{[\tilde \ell_1,\tilde \ell_2]}(\RR(\eta))}(h)$ is  a
symmetric  random variable and in particular has mean zero.
\medskip
\noindent{\bf Step 2: Extraction of the leading stochastic part.}

Recalling
\eqv(2.169),  we introduce
$$
\D^{\eta} \GG(m^{\d^*}_{\b,\tilde I_{12}})
\equiv \eta 
\left[\GG(m^{\d^*}_{\b,
\tilde I_{12}})-\GG(Tm^{\d^*}_{\b, \tilde I_{12}})\right]
\Eq(4.800)
$$
where $m^{\d^*}_{\b,\tilde I_{12}}$  was defined before
\equ(2.169). By definition, $|m^{\d^*}_\b-m_\b| \le 8 \g/\d^*$ and
taking $d_0$ small enough \eqv(5.5000101) implies
$|m^{\d^*}_\b-m_\b|\le  8 \g/\d^*\le \z_1 $. 
Thus, the block spin configuration  constantly equal  to 
$ m^{\d^*}_\b$ (resp. $T m^{\d^*}_\b$) is in $\RR^{\d,\z_4}(+1)$, (resp
$\RR^{\d,\z_4} (-1)$). 
Using the fact that the functional $\widehat \FF$ is
left invariant by $T$, we write
$$
\frac
{Z^{0,0}_{[\tilde \ell_1,\tilde \ell_2]}(\RR(-\eta))}
{Z^{0,0}_{[\tilde \ell_1,\tilde \ell_2]}(\RR(\eta))}(h)
\equiv e^{ \b  \D^\eta\GG(m^{\d^*}_{\b,\tilde I_{12}})}
\frac{Z_{-\eta,0,\d,\z_4}(\tilde
I_{12})}{Z_{\eta,0,\d,\z_4}(\tilde I_{12})}
\Eq(5.61)
$$
where 
$$
\frac{Z_{-\eta,0,\d,\z_4}(\tilde I_{12})}
{Z_{\eta,0,\d,\z_4}(\tilde I_{12})}
\equiv
\frac {\sum_{m^{\d^*}_{\tilde I_{12}} \in {\cal M}^{\d^*}(\tilde I_{12})}
\1_{\{\RR^{\d,\z_4}(\eta)\}} e^{-\frac \b \g\left\{{\widehat
\FF}(m^{\d^*}_{\tilde I_{12}},0)+ \g\D_0^{-\eta}\GG(m^{\d^*}_{\tilde
I_{12}}) + \g  V (Tm^{\d^*}_{\tilde I_{12}})\right\}}
 }
{\sum_{m^{\d^*}_{\tilde I_{12}} \in {\cal M}_{\d^*}( \tilde I_{12})}
\1_{\{\RR^{\d,\z_4}(\eta)\}} e^{-\frac \b\g\left\{{\widehat
\FF}(m^{\d^*}_{\tilde I_{12}},0)+ \g\D_0^\eta \GG(m^{\d^*}_{\tilde
I_{12}}) + \g  V (m^{\d^*}_{\tilde I_{12}})\right\}} }  
\Eq(4.9PP)
$$
with
$$
\D^\eta_0 \GG(m^{\d^*}_{\tilde I_{12}})\equiv
\sum_{x\in \CC_{\d^*}(\tilde I_{12})} \D^\eta_0 \GG^h_{x,m^{\d^*}(x)}
\Eq(4.90P)
$$
and, recalling \eqv(3.13P),
$$
\D^\eta_0 \GG^h_{x,m^{\d^*}(x)}=
\GG_{x,T^{\frac{1-\eta}2}m^{\d^*}(x)}(\l(x))-
 \GG_{x, T^{\frac{1-\eta}2}m^{\d^*}_{\b}(x)}(\l(x))
\Eq(4.10)
$$ with $T^0$ equal to the identity.
\medskip
\noindent{\bf Step 3: Control of the remaining stochastic part.}

To estimate the last term in \eqv (5.61),  we use  Lemma \eqv(41).
A control of  the Lipschitz norm is needed. Since it is rather involved to
do it, we postpone the proof of the next Lemma  to the end of the 
section. 
 
\noindent{\bf \Lemma(401) }
{\it Given $\b, \th $ that satisfy \eqv(2.17),
there exist
$\g_0=\g_0(\b,\th)>0$, $d_0=d_0(\b,\th)>0$, and
$\z_0=\z_0(\b,\th)$ such that for all $0<\g\le \g_0$, for all $\d^*>\g$ with
$\g/\d^* \le d_0$, for all $0<\z_4<\z_0$
that satisfy the following condition
$$
\z_4 \ge \left(5184 (1+c(\b\th))^2 (\frac{\g}{\d^*})^{1/2}\right)
\vee \left(12\frac{e^3\b}{c(\b,\th)} \frac{(\d^*)^2}{\g}\right)^2
\Eq(4.12009)
$$
where $c(\b\th)$ is given in \eqv(3.25PM) and $c(\b,\th)$ is given in
\eqv(mars.76), then  for all $a> 0$, 
$$ \P\left[{\max_{I\subset  \D_Q}}^*
\max_{\tilde I_{12}\subset I}\left|
\log \frac{Z_{-\eta,0,\d,\z_4}(\tilde I_{12})}
{Z_{\eta,0,\d,\z_4}(\tilde I_{12})}\right|
\ge \b\frac{4a+12\z_4}{\g}
\right]
\le
\frac{2  Q }{\e }
\frac{  e^{-\sfrac{u}{\e} }}{1-e^{-\sfrac{u}{\e} }}
\Eq(5.78)
$$ 
where ${\max_{I\subset  \D_Q}}^*$ denote the maximum over the intervals 
$I\subseteq \D_Q$ such that $|I|=\e\g^{-1}$ and 
$
u\equiv
\frac{a^2 \b^2}{8 \z_4 c^2(\b,\th)}.
$
}

\medskip
\noindent{\bf Step 4 Control of the leading stochastic part.}

To estimate the first term in the right hand side of
\eqv(5.61), we recall
$\D^\eta \GG(m^{\d^*}_{\b,\tilde I_{12}})= - \eta
\sum_{x\in \CC_{\d^*}(\tilde I_{12})} X(x)
$ where $X(x)$ is defined in \eqv(2.170). Using Lemma \eqv(lem1), 
exponential Markov inequality, and the Levy inequality
we get
$$
\P\Big[{\max_{I\subset \JJ}}^*\max_{\tilde I_{12}\subset I}
\big|\g\sum_{x\in \CC_{\d^*}(I)}X(x)\big|\ge 2 s\Big]
\le
  \frac{ 4 Q}{ \e} e^{-\sfrac{ s^2 }{2\e V^2_+}}.
 \Eq(5.970)
$$

 Then we  collect  \eqv
(5.78), \eqv(5.970)  and  make  the choice $a=\FF^*/16$, $s=\FF^*/32$.
Using the hypothesis \eqv(5.5000100) and { the definition} \eqv(5.4810), 
choosing $d_0$  small enough, we get 
$ 32(1+\th)(R_1+L_0)  \sqrt{{\g}/{\d^*}} +4\d^*\le \FF^*/2$.
Taking $\z_0$ small enough to have $28 \z_4 \le \FF^*/8$, we get
 $$
\mu_{\b,\th,\g}\big(\widetilde \WW_{\eta}^{ \z_1,\z_4}
 (\ell_1,\tilde \ell_1,\tilde \ell_2,\ell_2)
 \big)
\le e^{ -\frac {\b}{\g}
 ( 2\FF^* -32(1+\th)(R_1+L_0) \sqrt{\frac{\g}{\d^*}} -4(\z_1+\d^*)
- 24 \z- 4 a- 4s )}
\le e^{-\frac \b \g  \FF^*}
\Eq(5.98)
 $$
with $\P$--probability at least
$$  1-3 \g^{ 2} -  \frac{2  Q}{ \e}
\frac { e^{-\sfrac{u}{\e}}}
{1-e^{-\sfrac{u}{\e}}} - \frac{ 4Q}{\e}
e^{-\sfrac{ (\FF^*)2 }{2^6\e V^2_+}},
 \Eq (H.10)$$
where
$$
u\equiv
\frac{(\FF^*)^2}{2^{11} \z_4 c^2(\b,\th)
}.
\Eq(H.11)
$$
The unions in \eqv(5.51) involves at most $R_1^2\e^2 Q\g^{-3}$
terms.
This ends the proof of   Proposition \eqv(59).\eop

\noindent{ \bf Proof of Theorem  \eqv(61)}:

  It  is an immediate consequence of  Theorem  \eqv(t62)  and
 Proposition  \eqv(59) assuming $\z_0$ small enough to have
$u\le{ (\FF^*)^2 }/(2^6 V^2_+)$.
\eop

Lemma \eqv(41) is  the basic ingredient to prove Lemma \equ(401).
An estimate of Lipschitz norms is given in the next lemma. 
Then an Ottaviani type inequality will be used
to take care of the max in \eqv(5.78). 
We state Lemma \eqv(410) for a general $\z$ since it will be used in 
Section 5  with a $\z$ different from $\z_4$.

\noindent{\bf \Lemma (410)} {\it Let $\b>1,\th>0$ that
satisfy \eqv(2.17). We take $c(\b)$ as in \eqv(3.25PM).
There exist 
$\g_0=\g_0(\b,\th)>0$, $d_0(\b,\th)>0$, and
$\z_0(\b,\th)$ such that for all
$0<\g\le \g_0$, for all $\d^*>\g$ with $\g/\d^* <d_0$, and for all
$0<\z\le \z_0$,  
that satisfy
$$
\z> \left(5184 (1+c(\b\th))^2 (\frac{\g}{\d^*})^{1/2}\right)
\vee
\left(\frac{12e^3\b}{c(\b,\th)} \frac{(\d^*)^2}{\g}\right)^2
\Eq(H.4)
$$
where  $c(\b\th)$ is defined in \eqv(3.25PM) and 
 $$
c(\b,\th)=
257 \left(\frac{1}{(1-\tanh(2\b\th))^2}
+\frac{1}{1-m_{\b,1}}\right)+
e^{4\b\th} \frac{1+\tanh(2\b\th)}{1-\tanh(2\b\th)}
e^{
257 \left(\frac{1}{(1-\tanh(2\b\th))^2}
+\frac{1}{1-m_{\b,1}}\right)}
\Eq(mars.76)
$$
then
$$
\left\|\partial_i\log
\frac{Z_{+,0,\d,\z}(\tilde
I_{12})}{Z_{-,0,\d,\z}(\tilde I_{12})}\right\|_{\infty}
\le  \sqrt{\z} c(\b,\th) + 12e^3\b \frac {(\d^*)^2}{\g}
\le 2 \sqrt{\z} c(\b,\th)
\Eq (H.3)$$
where $
\frac{Z_{+,0,\d,\z}(\tilde
I_{12})}{Z_{-,0,\d,\z}(\tilde I_{12})}$ is defined as in \eqv(4.9PP) with
$\z_4$ replaced by $\z$
}
\vskip0.5cm

The proof of Lemma \eqv (410) is done similarly to the
corresponding estimates in Section 4 of [\rcite{COP1}].  
The main differences is that the explicit form of $\D^\eta_0\GG$
in \eqv(4.9PP) is not the same, and
we use the cluster expansion method to  estimate
the  Lipschitz  factors coming from  $ V (m^{\d^*}_{\tilde I_{12}})$.
Since we did not see a   simple way to modify the proof 
given in  [\rcite{COP1}]
we prefer to  start from the very beginning of the computations .

Given $i\in \g^{-1}\tilde I_{12}$, let $x(i)=\left[\g i/\d^*\right]$
be the index of the block of length
$\d^*$ that  contains  $\g i$, and let  $u(i)=\left[x(i)\d^*/\d\right]$
be the index of the block of length $\d$ that contains $x(i)$.

Let us denote
$$ \CC_{\d/\d^*}(u(i)) \equiv \CC_{\d/\d^*}(i)
\equiv\left\{ x \in \Z, \left[\frac{x(i)\d^*}{\d}\right]\frac{\d}{\d^*}<x\le
\left[\frac{x(i)\d^*}{\d}\right]\frac{\d}{\d^*} + \frac{\d}{\d^*} \right\}
\Eq(MRS12)
$$
{\it i.e.,}  the set of indices of those blocks of length $\d^*$ that are
inside  the block of length $\d$ indexed by $u(i)$.

Given a sample of $h$, let us denote $h^{(i)}$ the configuration
$h^{(i)}_j=h_j$ for $j\neq i$, $h^{(i)}_i=-h_i$.
To simplify the notations, we do not write explicitly
the $\d,\z$ dependence of $Z_{\pm,0,\d,\z}$ and  we write 
the Lipschitz factors as
$$
\del_i\log \frac{Z_{+,0,\d,\z}}{Z_{-,0,\d,\z}}=
\log\frac{Z_{+,0}(\tilde I_{12})(h)}{Z_{+,0}(\tilde I_{12})(
h^{(i)})}-
\log\frac{Z_{-,0}(\tilde I_{12})(h)}{Z_{-,0}(\tilde I_{12})(
h^{(i)})}
\Eq(4.260)
$$
To continue we need a simple observation: if
$
 \sum_{x\in \CC_{\d/\d^*}(i)}
\|m^{\d^*}(x)-m_\b\|_1
\le \frac {\d}{\d^*} \z,
$ then, given $g_1(\z)$ decreasing  such that
$ \lim_{\z\downarrow 0}g_1(\z)=0$
but $\frac {\z}{g_1(\z)}< 1$, and if $\z\le 1$, we have
$$
\sum_{x\in \CC_{\d/\d^*}(i)}\1_{\{\|m^{\d^*}(x)-m_\b\|_1\le
g_1(\z)\}}\ge \frac{\d}{\d^*}(1-\frac{\z}{g_1(\z)}).
\Eq(4.2602)
$$ 
This suggests to  make a partition of $\CC_{\d/\d^*}(i)$ into two
sets, 
$$
{\cal K}(m^{\d^*})\equiv \left\{x\in\CC_{\d/\d^*}(i) :
\|m^{\d^*}(x)-m_\b\|_1\le g_1(\z)\right\}.
\Eq(4.26040)
$$ and $\BB(m^{\d^*})=\CC_{\d/\d^*}(i)\setminus {\cal K}(m^{\d^*})$.
Let $\ell(i)=[i\g]$, for all
$m^{\d^*}\equiv m^{\d^*}_{\ell(i)}$ we write
$$
\1_{\{\eta^{\d,\z}(\ell(i))=1\}}(m^{\d^*})=
\sum_{X\subset\CC_{\d/\d^*}(  i )}\1_{\{{\cal K}=X\}}(m^{\d^*})\1_{\{\BB=
X^c\}} (m^{\d^*})
\1_{\{\eta^{\d,\z}(\ell(i))=1\}}(m^{\d^*})
\Eq(4.2605)
$$ where the sum is over all the subsets of $\CC_{\d/\d^*}( i )$ and
$X^c\equiv \CC_{\d/\d^*}( i )\setminus X$. It follows from
\eqv(4.2602) that $\eta_{\d,\z}(\ell(i))=1$ and
$|X| \le \frac{\d}{\d^*}(1-\frac{\z}{g_1(\z)})$ are incompatible.
Therefore we can  impose that $|X| \ge
\frac{\d}{\d^*}(1-\frac{\z}{g_1(\z)})$ in \eqv(4.2605).
Let $$
\NN(\z)=\sum_{X\subset\CC_{\d/\d^*}( i )}\1_{\{|X| \ge
\frac{\d}{\d^*}(1-\frac{\z}{g_1(\z)})\}}
=\sum_{k=\frac{\d}{\d^*}(1-\frac{\z}{g_1(\z)})}^{\frac
{\d}{\d^*}}{{\frac{\d}{\d^*}}\choose{k}},
\Eq(4.2606)
$$
and notice that \eqv(4.260) is the same as  
$$
\log\frac{Z_{+,0}(\tilde I_{12})(h)}{\NN(\z)^{\frac
12} Z_{+,0}(\tilde I_{12})(h^{(i)})}-
\log\frac{Z_{-,0}(\tilde I_{12})(h)}{\NN(\z)^{\frac 12}
Z_{-,0}(\tilde I_{12})(h^{(i)})}.
\Eq(4.2607)
$$ The two terms are estimated in the same way. We consider the
first one. It is easy to see that, with  self--explanatory notation, 
$$
\frac{Z_{+,0}(\tilde I_{12})(h)}{ \NN(\z)^{\frac 12}
Z_{+,0}(\tilde I_{12})(h^{(i)})}
=\frac
1{\NN(\z)^{\frac 12}
}\QQ_+\left[e^{\frac \b \g\big ( \g \D^+_0\GG_{x(i)}^{h}-\g \D^+_0\GG_{x(i)}^{
h^{(i)}} \big)} e^{ \frac \b \g   \big
(  \g V(\tilde I_{12},h)-\g V(\tilde I_{12},h^{(i)})
\big)}\right],
\Eq(4.2608)
$$ where $\QQ$ is the probability measure
$$
\QQ_+[\Psi]
=\frac {\sum_{m^{\d^*}(\tilde I_{12}) \in {\cal
M}_{\d^*}(\tilde I_{12})}\Psi(m^{\d^*})
\1_{\{\RR(+)\}} e^{-\frac \b \g\left\{{\widehat
\FF}(m^{\d^*}_{\tilde I_{12}},0)+ \g\D^+_0\GG^{
h^{(i)}}(m^{\d^*}_{\tilde I_{12}})+\g  V (m^{\d^*}_{\tilde I_{12}},h^{(i)})  \right\}}}
{\sum_{m^{\d^*}(\tilde I_{12}) \in {\cal M}_{\d^*}(I_{12})}
\1_{\{\RR(+)\}} e^{-\frac \b \g\left\{{\widehat
\FF}(m^{\d^*}_{\tilde I_{12}},0)+\g  \D^+_0\GG^{
h^{(i)}}(m^{\d^*}_{\tilde I_{12}})+ \g
 V (m^{\d^*}_{\tilde I_{12}},h^{(i)})   \right\}}
}.
\Eq(4.2609)
$$
Applying Schwartz inequality to \eqv (4.2608)  we obtain
$$
\frac{Z_{+,0}(\tilde I_{12})(h)}{ \NN(\z)^{\frac 12}
Z_{+,0}(\tilde I_{12})(h^{(i)})}
\le
\left ( \frac 1{\NN(\z)
}\QQ_+\left [e^{\frac \b\g 2
\big( \g\D^+_0\GG_{x(i)}^{h}-\g\D^+_0\GG_{x(i)}^{h^{(i)}}\big)
}
 \right] \right
)^{\frac 12} \left( 
\QQ_+\left [  e^{\frac \b \g 2  \left ( \g V(\tilde I_{12},h)-
\g V(\tilde I_{12},h^{(i)})
\right) } \right] \right)^{\frac 12}.  
\Eq (S.3) $$

The last term on  the right hand side of \eqv(S.3), 
can be immediately estimated  through Lemma \eqv(62), and we obtain 
$$
\left| \frac 12 \log \QQ_+\left [  e^{\frac \b \g 2 \left 
(  \g V(\tilde I_{12},h)-
\g V(\tilde I_{12},h^{(i)})
\right) } \right]  \right|  \le 6e^3 \b
\frac {( \d^*)^2} {\g}.
\Eq (mars.1)
$$
The needed  estimates for the first term in the right hand side of \eqv(S.3)
are summarized in the next Lemma

{\bf \Lemma (610) } {\it  Let   $\z$  and $g_{1}(\z)$ be
the quantities defined before \eqv(4.2602).
 For all $\b,\th$ that satisfy \eqv(2.17),
there exist $\z_0(\b\th)$ and $d_0(\b\th)$
such that for all $0<\z\le \z_0(\b\th)$, for all $\g/\d^*\le d_0(\b,\th)$,
 for all increasing $g_0(n)$ such that
$\lim_{n\uparrow\infty}g_0(n)=\infty$ but $g_0(n)/n$ is decreasing with
$\lim_{n\uparrow\infty}
g_0(n)/n=0$  we have that
$$
\left| \frac 12 \log
\frac 1{\NN(\z)  }\QQ_+\left [e^{\frac \b \g 2\big(
\g\D_0\GG_{x(i)}^{h}-\g\D_0\GG_{x(i)}^{h^{(i)}}\big)}
 \right]\right|
\le f_1(\z)+ \frac{\z}{g_1(\z)}e^{|f_2-f_1(\z)|}
 \Eq(4.2623)
$$
 where
$$ f_1(\z)\le \|h-h^{(i)}\|
256 g_{1}(\z)\left(\frac{1}{(1-\tanh(2\b\th))^2}
+\frac{1}{1-m_{\b,1}}\right)
+\frac{72}{g_0(\d^*\g^{-1}/2)} + \left(\frac{2\g}{\d^*}\right)^{1/4}4c(\b\th)
\Eq(4.117)
$$
with $c(\b\th)$ given in \eqv(3.25PM)
and
$$ f_2\equiv f_2(\b,\th) \le \|h-h^{(i)}\|
\Big(\log\frac{1+\tanh(2\b\th)}{1-\tanh(2\b\th)}+4\b\th \Big).
\Eq(4.1170)
$$
}
\goodbreak

 {\bf Proof:} We insert \eqv(4.2605) within the $[.]$ in the
left hand side of \eqv(4.2623). Then, see (4.56) in
[\rcite{COP1}], it can be checked that if we have an estimate of
the form
$$
\left| \D^+_0\GG_{x(i)}^{h}-\D^+_0\GG_{x(i)}^{h^{(i)}}\right|
\le f_1(\z) \1_{\{x(i)\in {\cal K}\}} +f_2\1_{\{x(i)\in \BB\}}.
\Eq(mars.2)
$$
From \eqv(4.2606) we then get
$$
\left|
\log
\frac 1 {\NN(\z)}
{\QQ_+\left [e^{\frac \b\g 2
\big( \g\D^+_0\GG_{x(i)}^{h}-\g\D^+_0\GG_{x(i)}^{h^{(i)}}\big)
}\right]}
\right|
\le
f_1(\z) +\frac{\z}{g_1(\z)} e^{|f_2-f_1(\z)|}.
\Eq(mars.3)
$$

To get \eqv(mars.2) with $f_1(\z)$ that satisfies \eqv(4.117) and
$f_2$ that { satisfies} \eqv(4.1170), we 
recall \eqv(3.13PM) and denote 
$$
\GG_{x,m^{\d^*}(x)}(\l(x))\equiv 
-\frac{1}{\b}
\log  L_{x,m^{\d^*}_{\frac{3+\l(x)}{2}} (x)}^{\d^*}
(\l(x)2\b\th,D(x)),
\Eq(rosa.1)
$$
so that 

$$
\eqalign{
\b\left(\D^+_0\GG_{x(i)}^{h}-\D^+_0\GG_{x(i)}^{h^{(i)}}\right) =&
-\log  
\frac{L_{x(i),m^{\d^*}_{\frac{3+\l(x(i))}{2}} (x(i))}^{\d^*}
(\l(x(i))2\b\th,D(x(i)))}
{L_{x(i),m^{\d^*}_{\frac{3+\l^{(i)}(x(i))}{2}}
(x(i))}^{\d^*}(\l^{(i)}(x(i))2\b\th,D^{(i)}(x(i)))}\cr
& +\log
\frac{L_{x(i),m^{\d^*}_{\b,\frac{3+\l(x(i))}{2}}
(x(i))}^{\d^*}(\l(x(i))2\b\th,D(x(i)))}
{L_{x(i),m^{\d^*}_{\b,\frac{3+\l^{(i)}(x(i))}{2}}
(x(i))}^{\d^*}(\l^{(i)}(x(i))2\b\th,D^{(i)}(x(i)))},\cr
}\Eq(mars.30)
$$
{ where}
$\l^{(i)}(x(i))$ and $D^{(i)}(x(i))$ are the respective images of
$\l(x(i))$ and $D(x(i))$ by the map $h\rightarrow h^{(i)}$.

The first case to consider is when $\l^{(i)}(x(i))=-\l(x(i))$, 
in which case $|D(x(i))|=|D^{(i)}(x(i))|=1$ and, 
using \eqv(3.13PM),  it can be checked that
$$\eqalign{
&\b\left(\D^+_0\GG_{x(i)}^{h}-\D^+_0\GG_{x(i)}^{h^{(i)}}\right)\cr
&\quad=\log\frac
{1+\l(x)m^{\d^*}_{\frac{3+\l(x(i))}{2}} (x(i))\tanh(\l(x(i))2\b\th)}
{1+\l(x)m^{\d^*}_{\b,\frac{3+\l(x(i))}{2}} (x(i))\tanh(\l(x(i))2\b\th)}
\frac
{1-\l(x)m^{\d^*}_{\b,\frac{3-\l(x(i))}{2}} (x(i))\tanh(-\l(x(i))2\b\th)}
{1-\l(x)m^{\d^*}_{\frac{3-\l(x(i))}{2}} (x(i))\tanh(-\l(x(i))2\b\th)}\cr
}\Eq(mars.60)
$$
Now if $\z_0$ is chosen in such a way that
$g_1(\z) \le
 (1-\tanh(2\b\th))/2$, noticing that \eqv(2.17) implies
$0<\tanh(2\b\th)<1$ when $1<\b<\infty$,
a simple computation gives that
$||m^{\d^*}(x(i))-m^{\d^*}_{\b}||_1\le g_1(\z)$ implies
$$
\left|\b(\D^+_0\GG_{x(i)}^{h}-\D^+_0\GG_{x(i)}^{h^{(i)}})\right|
\le \frac{4 ||m^{\d^*}(x(i))-m^{\d^*}_{\b}||_1}{1-\tanh(2\b\th)}
\le \frac{4g_1(\z)}{1-\tanh(2\b\th)}
\Eq(mars.61)
$$
while without condition on $||m^{\d^*}(x(i))-m^{\d^*}_{\b}||_1$
we have
$$
\left|\b(\D^+_0\GG_{x(i)}^{h}-\D^+_0\GG_{x(i)}^{h^{(i)}})\right|
\le \log\frac{1+\tanh(2\b\th)}{1-\tanh(2\b\th)}
\Eq(mars.62)
$$
therefore  \eqv(4.117) and \eqv(4.1170) are satisfied in this particular case.

The other case to study is when $\l^{(i)}(x(i))=\l(x(i))$ and therefore
$\left||D(x(i))|-|D^{(i)}(x(i))|\right|=1$.

If $x(i)\in \BB$, recalling \eqv(mars.3), we do not need a very
accurate  estimate for  the terms in \eqv(mars.30). 
Recalling 
\eqv(3.13PM), it is not difficult to see that each term in term in the right
hand side of \eqv(mars.30) is bounded by $2\b\th$, so we get
$$
\b\left|\D^+_0\GG_{x(i)}^{h}-\D^+_0\GG_{x(i)}^{h^{(i)}}\right|
\le 4\b\th
\Eq(mars.622)
$$
therefore collecting \eqv(mars.62) and \eqv(mars.622),
we have proven \eqv(4.1170).

{ It remains to consider the case where $x(i)\in {\cal K}$.
Recalling \eqv(mars.2) and \eqv(mars.3) this will give us the term
$f_1(\z)$. 
Here  we  want use the explicit form of $\GG_{x,m^{\d^*}}$ 
given in Proposition \eqv(PP).
To check that  \eqv(4.PM0) is satisfied, let us first note that 
since $g_1(x)$ and  $g_0(x)/x$ are  decreasing,
$\lim_{x\downarrow 0}g_1(x)=0$
and   $\lim_{n\uparrow \infty} g_0(n)/n=0$, if we
choose $\z_0=\z_0(\b,\th)$ such that
$$
g_1(\z_0) + \frac{\z_0 g_0(4/\z_0)}{4}\vee
\frac{16(\z_0/4)^{1/4}\b\th}{1-\tanh(2\b\th)}
\le 1-m_{\b,1}
\Eq(4.26041)
$$
and then we choose $d_0$ such that $\g(\d^*)^{-1} <d_0$ and
\eqv(H.4) implies $\z>8\g(\d^*)^{-1}$, we get
$$
g_1(\z) + \frac{g_0(\d^*\g^{-1}/2)}{\d^*\g^{-1}/2}
\vee \frac{16(2\g/\d^*)^{1/4}\b\th}{1-\tanh(2\b\th)}
 \le 1-m_{\b,1}
\Eq(4.260411)
$$
which   implies  that on ${\cal K}(m^{\d^*})$ and on the 
set $\{sup_{x\in\CC_{\d^*}(I)} p(x)\le (2\g/\d^*)^{1/4}\}$, 
 we have \eqv(4.PM0).}

\Remark(P) The fact that it is enough to have accurate  estimates
only in the Gaussian case comes from the previous sentence
together with \eqv(mars.2), \eqv(mars.3) and \eqv(mars.622).

To estimate \eqv(mars.30), 
we first notice   that the contribution to 
$\b\left|\D^+_0\GG_{x(i)}^{h}-\D^+_0\GG_{x(i)}^{h^{(i)}}\right|$ 
 coming from the terms  that correspond to  \eqv(6.PM)
is bounded by
$$
\frac{72}{g_0(\d^*\g^{-1}/2)} + \left(\frac{2\g}{\d^*}\right)^{1/4}4c(\b\th)
\Eq(mars.63))
$$
with  $c(\b\th)$ the  positive constant given in \eqv(3.25PM).
The terms in \eqv(mars.30)
that come from
$$
-|D(x)|\left[ \log\cosh(2\b\th)+
\log\left(1+\l(x)m^{\d^*}_{\frac{3+\l(x)}{2}}(x)\tanh(2\b\th)\right)
\right]
\Eq(mars.64)
$$
in \eqv(4.PM)
give  a contribution that is bounded by
$$
 \frac{8g_1(\z)}{1-\tanh(2\b\th)}
\Eq(mars.65)
$$
when
$||m^{\d^*}(x(i))-m^{\d^*}_{\b}||_1\le g_1(\z).$
It remains to estimate the contribution to \eqv(mars.30)
of the terms that come from
$$
|D| \hat\varphi(
m^{\d^*}_{\frac{3+\l(x)}{2}}(x),2\l(x)\b\th, p(x) )
\Eq(mars.67)
$$
in \eqv(4.PM). Unfortunately the estimate \eqv(5.PM) is useless
and we have to consider the explicit form of $\hat\varphi$,
see \eqv(mars.68P).
The contribution of $\hat\varphi$ in \eqv(mars.30) can be bounded by
$$
\int_{p(x(i))\wedge p^{(i)}}^{p(x(i))\vee p^{(i)}}
\int_
{m^{\d^*}_{\frac{3+\l(x)}{2}}(x)\wedge
m^{\d^*}_{\frac{3+\l(x)}{2},\b}(x)}
^{m^{\d^*}_{\frac{3+\l(x)}{2}}(x)\vee
m^{\d^*}_{\frac{3+\l(x)}{2},\b}(x)}
\left|\frac{\del^2\left[p|B| \hat\varphi(m,2\l(x)\b\th, p )\right]}
{\del m\del p}\right|
\,dp\, dm.
\Eq(mars.69)
$$
It is just a long task to compute the previous partial derivative,
using \eqv(3.16PM), \eqv(3.22PM) and \eqv(3.25PM) and to
check that the following
estimates are valid if $\z$ is such that $g_1(\z) \le (1-\tanh (2\b\th))/2$
$$
\eqalign{
&\frac{\del\nu_2}{\del p}\le \frac{2}{\s^2_{m}},\quad
\frac{\del \nu_2}{\del m} =\frac {1}{\s^2_{m}},\cr
&\left|\frac{\del^2 \nu_2}{\del p\del m} \right|\le \frac{4}{\s^2_{m}},
\quad 0<\frac {1}{\s^2_{m}} -\frac{1}{1-m^2} \le \frac{pc(\b\th)}{\s^2_{m}}.\cr
}\Eq(mars.70)
$$

It is clear that  unpleasant looking  terms like
$(1+m\tanh(\nu_2-\nu_1))^{-1}$  appear in the computations.
Using \eqv(3.27PM), the fact that we can assume that $\z_0=\z_0(\b,\th)$
is small enough to get that if $\z\le \z_0$ then
$||m-m_{\b}||_1\le g_1(\z)$ implies $1-|m| \ge (1-m_{\b,1})/2$.
Then,  assuming $d_0(\b,\th)$ to be  small enough in order to have that 
$\g/\d^*\le d_0(\b,\th)$ implies
$4\b\th (\g/\d^*)^{1/4}/(1-m_{\b,1}) \le 1/2$, we get
$$
1+m\tanh(\nu_2-\nu_1))> 1-\frac{4m_{\b,1}\b\th p(x)}
{1-m_{\b,1}} > \frac 14
\Eq(mars.71)
$$
for all $m$ and $p$ that occur   in the integral in \eqv(mars.69). 
So,  these terms do not present any problem.
We get 
$$
\left|\frac{\del^2\left[p|B| \hat\varphi(m,2\l(x)\b\th, p )\right]}
{\del m\del p}\right| \le |B|
256 \left(\frac{1}{(1-\tanh(2\b\th))^2}
+\frac{1}{1-m_{\b,1}}\right). 
\Eq(mars.72)
$$
Notice  that
$$
\int_{p(x(i))\wedge p^{(i)}}^{p(x(i))\vee p^{(i)}}
\int_
{m^{\d^*}_{\frac{3+\l(x)}{2}}(x)\wedge
m^{\d^*}_{\frac{3+\l(x)}{2},\b}(x)}
^{m^{\d^*}_{\frac{3+\l(x)}{2}}(x)\vee
m^{\d^*}_{\frac{3+\l(x)}{2},\b}(x)} \,dp\,dm
\le ||m^{\d^*}-m^{\d^*}_{\b}||_1 \frac 2{B}
\Eq(mars.73)
$$
Thus, inserting \eqv(mars.72) in \eqv(mars.69),  using \eqv(mars.73)
and then  collecting \eqv(mars.63) and \eqv(mars.65) we get \eqv(4.117).
\eop

\vskip0.5cm

{\bf Proof of Lemma \eqv (410) }
We  recall \eqv (4.260), \eqv (4.2607), and \eqv(S.3)  and apply
 Lemma \eqv (610) and \eqv (mars.1).
The presence of $\z$ in \eqv(4.117) and $\z/g_1(\z)$ in \eqv(4.2623)
suggests to take $g_1(\z)=\sqrt{\z}$.
The presence of $(g_0(\d^*\g^{-1}/2))^{-1}$ and $(2\g/\d^*)^{1/4}$
in \eqv(4.117) suggests  to choose $g_0(n)=n^{1/4}$.
Thus,  calling
$$
c_1\equiv c_1(\b,\th)\equiv 256 \left(\frac{1}{(1-\tanh(2\b\th))^2}
+\frac{1}{1-m_{\b,1}}\right)
\Eq(mars.75)
$$
and
$$
c_2\equiv c_2(\b,\th)=e^{4\b\th} \frac{1+\tanh(2\b\th)}{1-\tanh(2\b\th)}
\Eq(mars.761)
$$
we get that the left hand side of \eqv(4.2623) is bounded by
$$
\sqrt{\z} \left(c_1 +c_2e^{\sqrt{\z} c_1+ 72(1+c(\b\th))(\frac{2\g}{\d^*})^{1/4}}
\right)+
72(1+c(\b\th))(\frac{2\g}{\d^*})^{1/4}
\Eq(mars.77)
$$
from which we get the first term on the right hand side
of\eqv(H.3) with the $c(\b,\th)$ given in  \eqv(mars.76).
\eop

\noindent{\bf Proof of Lemma \eqv(401)}

Using Lemma \eqv(41) and Lemma \eqv(410), we get
after a simple computation,
for all $a> 0$, for all intervals
$\tilde I_{12}=[\tilde \ell_1,\tilde \ell_2]$
$$
\P\Big[
\Big|\log\frac{Z_{-\eta,0,\d,\z_4}(\tilde
I_{12})}{Z_{\eta,0,\d,\z_4}(\tilde I_{12})}\Bigr|\ge
\frac a {\g}\Big]
\le \exp \left(-\frac { a^2}{8  \g|\tilde \ell_1-\tilde \ell_2|
\z c^2(\b,\th)}\right).
\Eq(4.120)
$$

To get \eqv(5.78),  
we need the following  modification of the Ottaviani inequality
done in [\rcite{COP1}], see Lemma (5.8) there.
Given an interval $\tilde I\subset I$, calling 
$ Y(\tilde I) \equiv
\log \frac{Z_{-\eta,0,\d,\z}(\tilde I)}{Z_{\eta,0,\d,\z}(\tilde
I)}$, then for all $a>0$, for all $\z>8\g(\d^*)^{-1}$, we have 
$$
\P\left[\max_{\tilde I_{1,2}\subset I}
\left|
Y(\tilde I_{12})
\right|
\ge \b\frac {4 a+ 12 \z}\g\right]\le
 \frac  {\P\left[\left|Y(I)\right|
\ge \b\frac a \g\right]} {\inf_{\tilde I_{12}\subset I}\P\left[
\left|Y(\tilde I_{12})\right|\le \b\frac a {\g}\right]}.
\Eq(5.65)
$$
Then    for all $a >0$, setting $  \tilde x = 4a+12\z$, we obtain
$$
\P\left[{\max_{I\subset  \D_Q}}^*\max_{\tilde I_{12}\subset I}
|Y(\tilde I_{12})|\ge
 \b\frac {\tilde x } {\g}\right]
\le
\frac {2  Q   }{ \e }
\P\left[\max_{\tilde I_{12}\subset \wh I_{[0,2]}} |Y(\tilde
I_{12})|\ge
\b\frac {\tilde x} {\g} \right],
\Eq(5.750)
$$
where $\wh I_{[0,2]}=[0,2\e\g^{-1}]$. This implies \eqv(5.78) 
after a short computation. \eop

\def\e{\epsilon}
\chap{5 Proof of Theorems }5
\numsec= 5
\numfor= 1
\numtheo=1
In this section we prove Theorems \eqv(t1), \eqv(t2), and \eqv(t3).
They will be  derived from
Proposition \eqv (510P) stated and proved below.
  We   will use the following strictly positive finite
quantities: $\k (\b,\th)$ that satisfies  \eqv(2.19), $\FF^*$
defined  in \eqv (min), $V(\b,\th)$ in  \eqv(2.1701),
$c(\b,\th)$   in \eqv(mars.76)  and $c(\b\th)$ in
\eqv(3.25PM).
We denote
$$
\a(\b,\th,\z_0)
\equiv -\log
\frac{\del g_{\b}}{\del m}(\tilde m_{\b,\th}-\frac{\z_0}{2},\th)>0
\Eq(5.P1)
$$
where $\z_0=\z_0(\b,\th)$ is a small quantity that satisfies
requirements written before \eqv(7.13).
Recalling  \eqv(Palmares), we have
$
\a(\b, \th) \ge \a(\b,\th,\z_0)
$.
 The results  from Sections 3,4, and 6 require relations among
various parameters. For $\g_0,d_0,\z_0$ sufficiently small 
depending on
$\b,\th$
as stated in Theorem \eqv(t1),  $0<\g\le \g_0$, $\g/\d^*<d_0$
$1>\d>\d^*>0$, $\z_0>\z_4>\z_1>\z_5>8\g/\d^*$, $Q>1$, $\e>0$,
we assume that  the following constraints are satisfied:

\noindent {  \bf  The $\CC_0 $ constraints: }
 $$
\frac{128 (1+\th)}{\k(\b,\th)}
\frac{2(5+\FF^*)}{\FF^*}
\sqrt {\frac{\g}{\d^*}}
<\d\z_1^3,
\Eq(5.001P)
$$

$$
\frac {32}{ \k(\b,\th)} \z_1 \le\d \z_4^3
\Eq(5.002)
$$
$$
\left(5184 (1+c(\b\th))^2 \sqrt{\frac{\g}{\d^*}}\right)
\vee \left(12\frac{e^3\b}{c(\b,\th)} \frac{(\d^*)^2}{\g}\right)^2
\le \z_5
\Eq(5.004P)
$$
$$
\frac{512(1+\th)}{\k(\b,\th)\a(\b,\th,\z_0)}
\sqrt {\frac{\g}{\d^*}} \log \frac{\d^*}{\g}
<\d\z_5^3
\Eq(5.003)
$$

$$
\sqrt {\g} \log Q \le 
\frac{\sqrt{6e^3 \b}}{128}
\Eq(5.00410)
$$
$$
\frac{\FF^*}{32(1+\th)} \sqrt{\d^*\g}\le \e
\Eq(5.00411)
$$
\noindent {\bf Remark.} The constraints  \equ(5.001P), \equ(5.002),
\equ(5.004P), and \eqv(5.00410) come from Theorem \equ(61), where
\equ(5.004P) was written for $\zeta_5$ replaced by a larger 
value $\zeta_1$ ; now
we  impose the stronger restriction \eqv(5.004P),
as  it will be needed later. Notice that \eqv(5.00411) and \eqv(5.001P)
imply that $\e\g^{-1}>2R_1$.  \eqv(5.003) comes from \eqv(5.099) in 
Corollary \eqv(6002).

\Remark(P39)
 Note that in \eqv(5.001P)  one can take $\d=\d_1$, in \eqv(5.002)
$\d=\d_4$ and in \eqv(5.003) $\d=\d_5$, with $\d_5=n_5\d^*$,
$\d_1=n_1\d_5$, and $\d_4=n_4\d^*$ for some positive integers
that will diverge since $\d^*\downarrow 0$.
This would  allow $\d_4$ to be
small without imposing as in Theorem \eqv(t1) that it  goes
to zero. Since this would  introduce   new parameters we have decided,
for simplification,  not to do it. 

 \vskip0.5cm \noindent With the choice of parameters that
satisfy the $\CC_0$ constraints, we apply Theorem \eqv(61),
Corollary \eqv(6002) with $p=2 +[(\log Q)/(\log (1/\g))]$,
Lemma \eqv(lem112), and Corollary \eqv(zorro) with $k=5$, to determine
measurable sets $\O_4= \O_{4} ( \g,\d^*, \D_Q,\e,\d,\z_1,\z_4)$,
$\O_{RE}=\O_{RE}(\g,\d^*,$ $ p)\equiv \O_{RE}(\g,\d^*, Q)$, $\O_{\e}$, and
respctively $\PP(5,\e,Q)$ such that, calling $\O_{51}=\O_4\cap 
\O_{RE} \cap \PP(5,\e,Q)^c\cap \O_{\e}$, we have  
$$
\P[\O_{51}] \geq 1-
10e^{-\frac{Q}{10 C_1}}-5\e^{\frac{a}{16(2+a)}}
-Q^2\e^{\frac{a}{8+2a}}- Q e^{-\frac{1}{2\e^{3/4}V^2(\b,\th)}}
-7\g^2,
\Eq(5.P01)
$$
when $\d^*\g<\e\le \e_0(\b,\th)$ and $a>0$.

 For $\o \in  \PP(5,\e,Q)^c$, the origin belongs to an
unique elongation $[\a^*_{j_0},\a^*_{j_0+1}]$ where $j_0=-1\, {\rm
or}\, 0$, see \eqv(3.8) and \eqv(3.89), moreover on this set, recalling
\eqv(3.12), we have,
$$
\left[-\frac{\rho}{\g},\frac{\rho}{\g}\right]
\subset \left[\frac{\e \a^*_{j_0}}{\g}, \frac{\e \a^*_{j_0+1}}{\g}\right]
\subset\left[-\frac{Q}{\g},\frac{Q}{\g}\right]. 
\Eq(zumba1)
$$
We write, for $\eta\in\{-1,+1\}$
$$
\O^{\eta}(\e,Q ) \equiv\left\{\o\in \PP(5,\e,Q)^c,\,
{\rm sgn}\,\left[\frac{\e \a^*_{j_0}}{\g}, \frac{\e \a^*_{j_0+1}}{\g}\right]
=\eta\right\}, 
\Eq(zumbi2)
$$
For concreteness, we take $j_0=0$ and we  assume that this
elongation is positive,  that is,  we are on  $\O_{51}\cap
\O^+(\e,Q)$.  We have the following result:

\noindent {\bf \Proposition (510P)} {\it
If $\CC_0$ holds and
$$
8 f_1 +4 f_2 +4 f_3 + 32  \z_5^{\frac{1-z}{2}}+16 \z_1  \le
\frac{\e^{1/4}}{2}
\Eq(5.0100)
$$
where

$$
f_1=
10(1+\th)
\frac{1}{\a(\b,\th,\z_0)}
\sqrt{\frac{\g}{\d^*}}\log \frac{\d^*}{\g},
\Eq(5.01)
$$
$$
f_2=8 V(\b,\th)
 \sqrt{\g \log\left(\frac{1}{\g}\right)
\left(\frac{1}{\a(\b,\th,\z_0)}\log \left(\frac{\d^*}{\g}\right)+ R_1\right)}
\Eq(5.0100a)
$$
with $R_1= \frac{4(5+\FF^*)}{\k(\b,\th)\d\z_1^3}$,
$$
f_3= 16 (1+\th) R_1\sqrt{\frac{\g}{\d^*}},
\Eq(5.0001)
$$
and  $0<z<1/2$, there  exists $\O_{5}$ such that
$$
\P[ \O_{5} ] \ge 1 -8\g^{2} 
-\frac{
2\exp(-\frac{\b^2}{2^{6}Q
 {\z_5^z}c^2(\b,\th) })
}
{1-\exp(-\frac{\b^2}{2^{6}Q {\z^z_5}c^2(\b,\th) })
}
\Eq (BB.50)
$$
such that  on  $\O_{5}\cap\O_{51}\cap \O^\eta(\e,Q)$,
$$
\eqalign {&
\mu^\o_{\b,\th,\g}\left(
\exists \ell \in
[\frac {\a_0^*\e}{ \g}  + \frac {\r}{ \g} +R_1,
\frac {\a^*_1 \e}{\g} -\frac {\r}{  \g}-R_1 ],
\eta^{\d,\z_4} (\ell)\neq \eta  \right)\le
\cr &\quad\quad \le
\left(\frac {3Q}{\g^2}\right)^5
e^{-\frac \b \g \left\{(\frac{\k(\b,\th)}{4}
\d\z_4^3)\wedge \FF^*\right\}}
   + 28 R_1^2 \left(\frac{2Q}{\g}\right)^5
e^{-\frac  \b \g     \frac{ \e^{1/4}}{5}}
\exp \left\{\frac{4Q}{\g} e^{-\frac  \b\g  \frac {\e^{1/4}}{5}   } \right\}.
}
\Eq (5.2C)
$$
where $\r\equiv \e^{\frac{1}{4(2+a)}}$.
}

\noindent {\bf Remark} Recalling \eqv(inter) and Proposition 
\eqv(Terminator)
the interval $J=[\frac {\a_0^*\e}{ \g},  \frac {\a^*_1 \e}{\g} ]$
is random, its length being a finite and positive random variable, 
of order $\g^{-1}$.
On the other hand   when choosing the parameters
${\r+\g R_1 }$ will tend to zero.

\noindent{\bf  Proof.}
We assume that $\eta=+1$,  the  case $\eta=-1$ being  similar.
To simplify notation we denote
by $N_1= \frac 1 \g \a_0^*
\e $, $N_2= \frac 1 \g \a^*_1 \e$,
  $I= [N_1+R_1+\frac {\r} \g,N_2-R_1-\frac {\r} \g]$,
$\eta^{\d,\z_4} (\ell) =\eta (\ell)$ and
$B(\ell)=  \left \{ \s: \eta(\ell) \neq 1 \right \}
$
{ Recalling \eqv(ME.3)}, we have that
$$ \eqalign { &
\mu_{\b,\th,\g}\left ( \exists \ell \in   I,
\eta(\ell)\neq  1)\right)  \le
\mu_{\b,\th,\g}\left (\MM_{\d^*}(\D_Q)\setminus\AA(\D_Q)
\right) +  \sum_{\ell \in I}
\mu_{\b,\th,\g}\left ( B(\ell) \cap  \AA(\D_Q) \right), 
 }
\Eq (5.3)
 $$
where we denote by $ \AA(\D_Q)^c$ the complement in
$\MM_{\d^*}(\D_Q)$  of $  \AA(\D_Q)$.

According to  Theorem \eqv (61), for $\o \in \O_{51}\subset \O_4$ we have 
$$  \mu_{\b,\th,\g}\left (
\MM_{\d^*}(\D_Q)\setminus\AA(\D_Q)
\right) \le
\left( \frac{3Q}{\g^2}\right)^5
e^{-\frac \b \g \left\{\left(\frac{\k(\b,\th)}{4}
\d\z_4^3\right)\wedge \FF^*\right\}}. 
\Eq (ME.5) $$
To estimate the other term in \eqv (5.3)  we  need to restrict the infinite volume Gibbs measure  to a
finite volume one.
We write
$$
\eqalign { & \mu_{\b,\th,\g}\left ( B(\ell)\cap  \AA(\D_Q) \right
)\cr & \le \sum_{\bar\eta_1,\bar\eta_2 \in \{-1,1\}^2}
\sum_{\ell_1= N_1}^{N_1+ R_1 }\sum_{\ell_2= N_2- R_1 }^{N_2}
\mu_{\b,\th,\g} \left ( \eta^{\d,\z_1} (\ell_1)=\bar
\eta_1,\eta^{\d,\z_1} (\ell_2)=\bar\eta_2,   B(\ell)\cap \AA(\D_Q)
\right ) \cr & +  \mu_{\b,\th,\g}\left ( \eta^{\d,\z_1} (\ell)=0,
\forall \ell \in  [N_1,N_1+R_1] \right )  + \mu_{\b,\th,\g}\left (
\eta^{\d,\z_1} (\ell)=0, \forall \ell \in [N_2- R_1,N_2] \right ).
} \Eq (5.5)
$$
Using Theorem \eqv(t62), with $p=2 + [(\log Q)/(\log
\g^{-1})]$,
on $\O_{RE} \supset \O_{51}$ we have
$$
\eqalign{
&\mu_{\b,\th,\g}\left ( \forall
 \ell \in  [N_1,N_1+R_1],  \eta^{\d,\z_1} (\ell)=0
\right )  +  \mu_{\b,\th,\g}\left ( \forall \ell
\in  [N_2- R_1,N_2] ,  \eta^{\d,\z_1} (\ell)=0   \right )\cr
&\quad \le
\frac{3^4 Q^5}{\g^{10}}
e^{-\frac \b \g
\left\{\left(\frac{\k(\b,\th)}{4}\d\z_4^3\right)\wedge \FF^*\right\}}
}\Eq (5.9) $$
where
$ R_1= \frac { 4(5+\FF^*)}{\k(\b,\th) \d\z_1 ^3} $
and we have used the fact that our choice of $p$ entails $Q\g^{-1}\le \g^{-p}
\le Q \g^{-2}$ to replace $3^4\g^{-5p}$ in \eqv(5.481) by 
$3^4 Q^5 \g^{-10}$ in \eqv(5.9).

Recalling \eqv(3.7P)
and  using that  $\eta^{\d,\z_1}(\ell_1)=\bar\eta_1$
implies that on the left of $\ell_1$
$$
\big|E(m^{\d^*}_{\g^{-1}(\ell_1-2,\ell_1-1]}(\s),
m^{\d^*}_{\g^{-1}(\ell_1-1,\ell_1]}(\s'))-
E(m^{\d^*}_{\g^{-1}(\ell_1-2,\ell_1-1]}(\s),
m^{\d^*}_{T^{\frac{1-\bar\eta}{2}}
\b,\g^{-1}(\ell_1-1,\ell_1]})\big |\le \z_1
\Eq(Yasmina10)
$$
for $\s'$ such that $\eta^{\d,\z_1}(\ell_1)=\eta^{\d,\z_1}(\ell_1)(
\s'_{\g^{-1}(\ell_1-1,\ell_1]})=\bar\eta_1$
and similarly on the right of $\ell_2$, we get
$$
\eqalign {
& \mu_{\b,\th,\g} \left ( \eta^{\d,\z_1} (\ell_1)=\bar \eta_1,
\eta^{\d,\z_1}
(\ell_2)=\bar\eta_2,   B(\ell),  \AA(\D_Q)\right ) \cr
& \le e^{\frac \b \g ( 4 \z_1+\d^*)}
 \frac {  Z_{[\ell_1, \ell_2]}^{0,0} \left (\eta^{\d,\z_1} (\ell_1)=\bar
\eta_1, \eta^{\d,\z_1} (\ell_2)=\bar
\eta_2, B(\ell),\AA([\ell_1,\ell_2])  \right)} 
{  Z_{[\ell_1, \ell_2]}^{0,0} \left
(
\eta^{\d,\z_1} (\ell_1)=\bar \eta_1, \eta^{\d,\z_1} (\ell_2)=\bar \eta_2
\right )  }. } \Eq (5.10)
$$

To get  an upper bound for  \eqv(5.10), we
restrict the denominator to profiles { that  we expect to be
typical for }the Gibbs measure under the constraint
$\eta^{\d,\z_1} (\ell_1)=\bar \eta_1, \eta^{\d,\z_1} (\ell_2)=\bar
\eta_2 $ given that we are inside a positive elongation. Without
the constraints, taking into account only the presence of a
positive elongation, the profiles we expect to be typical are of
course $\eta^{\d,\z_4}=1$ for all $\ell\in [\ell_1,\ell_2]$, this
is also the case for $(\bar \eta_1,\bar \eta_2)= (+1,+1)$.
To take into account the cases  $(\bar \eta_1,\bar \eta_2)\ne
(+1,+1)$, we  leave  intervals $[\ell_1,\ell_1 + L_0 ]$ and/or $[
\ell_2- L_0,\ell_2]$, where $L_0$ is a positive integer to be
chosen later  { to allow the profiles to change from,} say
$\eta^{\d,\z_1} (\ell_1)=\bar  \eta_1=-1$ to $\eta^{\d,\z_1}
(\ell_1 + L_0 ) =+1$.  We actually require the profiles to satisfy
$\eta^{\d,\z_5} (\ell_1 + L_0 ) =+1$, with $\z_5 <\z_1 $  for a
reason that we explain  later.

To proceed on this  it is convenient to define:
 given  $N_1 \le \ell_1 <\ell_2 \le  N_2$ and $\bar \eta \in \{-1,+1\}$,
for $i=1$ and $i=5$
$$
\widetilde  \RR_i(\bar\eta, \ell_1,\ell_2)
= \Big\{ m^{\d^*}_{[\ell_1,\ell_2]}\colon \eta^{\d,\z_i} (\ell_1)=\bar
\eta =\eta^{\d,\z_i} (\ell_2) \Big\}, \Eq(5.L)
$$
$$
\EE (+1,\ell_1,\ell_2,\bar \eta_1 ,\bar \eta_2) \equiv
 \left \{  \eqalign   { &  \widetilde
\RR_1(+1,\ell_1,\ell_2) \cap \{\eta^{\d,\z_5}(\ell_1+L_0)=
 \eta^{\d,\z_5}(\ell_2-L_0)=+1 \}
\,{\rm for}\, \bar \eta_1=-1=\bar \eta_2{\bf ;} \cr &  \widetilde
\RR_1(+1,\ell_1,\ell_2) \cap \{\eta^{\d,\z_5}(\ell_2-L_0)=+1\}
\,{\rm for}\, \bar \eta_1= 1,\bar \eta_2=-1 {\bf ;}\cr &
\widetilde \RR_1(+1,\ell_1,\ell_2) \cap
 \{\eta^{\d,\z_5}(\ell_1+L_0)=+1\} \,{\rm for}\, 
\bar \eta_1=-1,\bar \eta_2=1 {\bf ,}   }  \right.
\Eq(5.12002)
$$
where the $+1$ on  the left hand side  is associated to the
sign of the elongation, chosen  here to be  positive.
We then estimate the expression in \eqv (5.10)  as in Section 4 
(see \eqv(5.52P6)),     to obtain 
$$
\eqalign {
& \mu_{\b,\th,\g} \left ( \eta^{\d,\z_1} (\ell_1)=\bar \eta_1,\eta^{\d,\z_1}
(\ell_2)=\bar\eta_2,   B(\ell), \AA(\D_Q)  \right )\cr
& \le  e^{  \frac { \b}{\g}4(\z_1+\z_5+2 \d^*)}
 \frac {  Z_{[\ell_1, \ell_2]}^{0,0} 
\left (\eta^{\d,\z_1} (\ell_1)=\bar \eta_1,
\eta^{\d,\z_1} (\ell_2)=\bar \eta_2, B(\ell), \AA([\ell_1,\ell_2])
\right) }
{ Z^{0,0}_{[\ell_1,\ell_2]}
\big (\EE(+1,\ell_1,\ell_2,\bar \eta_1 ,\bar \eta_2 )
\big )  }\times \cr
& \quad
\times  \frac     { Z^{0,m_+}_{[\ell_1, \ell_1+L_0-1]}
 \big(
 \eta^{\d,\z_1} (\ell_1)= +1
 \big )} { Z^{0,m_+}_{[\ell_1,
\ell_1+L_0-1]} ( \eta^{\d,\z_1} (\ell_1)= \bar\eta_1) }
\frac   { Z^{m_+,0}_{[  \ell_2-L_0+1,\ell_2]}
 \big ( \eta^{\d,\z_1} (\ell_2)=+1
 \big ) }{
Z^{m_+,0}_{[  \ell_2-L_0+1,\ell_2]} (\eta^{\d,\z_1} (\ell_2)=\bar
\eta_2) } .} \Eq (5.12) $$
 To apply    Lemma \eqv(T1C)  to the last two terms 
 in \eqv (5.12), we take
$$
L_0=\frac{1}{\a(\b,\th,\z_0)} \log \frac{\d^*}{\g}
\ge  \frac{1}{\a(\b,\th)} \log \frac{\d^*}{8\g}.
\Eq(5.12P)
$$
  Replacing
the $f_{11}$ of Lemma \eqv(T1C)  by   $f_1$  defined in \eqv(5.01),
since   here  $\sqrt{\frac{\g}{\d^*}} \ge \d^*$,
       we obtain
$$ \eqalign {
 \mu_{\b,\th,\g} \big (
\eta^{\d,\z_1} (\ell_1)&=\bar \eta_1,\eta^{\d,\z_1}
(\ell_2)=\bar\eta_2,  B(\ell),\AA \big )
\le
e^{\frac {\b}{\g} 4  (\z_1+\z_5+2 \d^*) }
e^{\frac \b \g  \left (  \FF^* +  2 f_1  \right)
\left[\frac 12  (|\bar \eta_1-1|+|\bar \eta_2-1|)\right]}
\times \cr
&
\times \frac {  Z_{[\ell_1, \ell_2]}^{0,0}
\left (\eta^{\d,\z_1} (\ell_1)=\bar \eta_1, \eta^{\d,\z_1}
(\ell_2)=\bar \eta_2, B(\ell),\AA([\ell_1,\ell_2])
\right) }
{  Z_{[\ell_1, \ell_2]}^{0,0} \big 
(\EE (+1,\ell_1,\ell_2,\bar \eta_1 ,\bar \eta_2)
\big ) }
}\Eq (5.16)
$$

To treat the last term in \eqv(5.16), we make a  partition
of the set of profiles in
$\AA([\ell_1,\ell_2]) $  distinguishing   the profiles
according to the number and
the location of the changes of phases in $ [\ell_1,\ell_2]$.
$$\AA([\ell_1,\ell_2]) =   \cup_{n=0}^{\bar N}
\cup_{\{A: |A|=n\}}\AA([\ell_1,\ell_2], A,n)\Eq (M.1E)$$
{ where}
$ \bar N$ is  the  number  of the $\frac \e \g$ blocks  in
$[\ell_1, \ell_2]$,  i.e., 
 $$\bar N= \left [ |\ell_2-\ell_1  | \frac \g \e \right ] = \left [
\left ( \frac \e \g [\a^*_1-\a^*_0]-2R_1 \right )  \frac \g \e
\right ] = \left [ [\a^*_1-\a^*_0]-2\frac \g \e R_1 \right ] \Eq
(D.2)$$
$[ x]$ is the integer part of $x$, and the
first  equality follows from \eqv(5.00411) that entails $\e/\g > 2
R_1$.
  Moreover in \eqv(M.1E),
$
A\subset \left\{\frac {Q_1}{\e}+1, \frac{Q_1}{\e}+2,\dots,
\frac{Q_2}{\e}-2,\frac{Q_2}{\e}-1\right\}.
$
The integer $n$ represents the cardinality of the set $A$ and
therefore the number of $\frac \e \g$ blocks where, in each one of
them,   there is  one and only one  interval of length $2 R_1$ in
which only one change  of phases occurs. Recall that in the
definition of $\AA ([\ell_1,\ell_2])$ cf. \equ(6.1)  the $r_i$, 
 $i=1,..,\bar
N $  indicate that in $[ r_i\frac \e \g ,
(r_i+1) \frac \e \g ] $ there is  $q_i$, such  that in
  $[ q_i-R_1,q_i+R_1] $   there is only  one
change of phases { and there is no change in $[ r_i\frac \e \g
, (r_i+1) \frac \e \g ] \setminus [ q_i-R_1,q_i+R_1]$}.
The notation $\AA([\ell_1,\ell_2], A,n)$ is
self-explanatory. When there is no ambiguity we denote
$\AA([\ell_1,\ell_2], A,n) \equiv \AA(A,n)$.
Going back to \eqv (5.3), taking into  account  \eqv(5.9),  \eqv
(5.16) and \eqv (M.1E)
 on  $\O_{51}$,   we have that
$$
 \mu_{\b,\th,\g}\left ( \exists \ell \in   I,
\eta(\ell) \neq  1 \right )  \le
2\left( \frac{3Q}{\g^2}\right)^5
e^{-\frac \b \g \left\{(\frac{\k(\b,\th)}{4}
\d\z_4^3)\wedge \FF^*\right\}}  + e^{\sfrac {4\b(\z_1+\z_5 +2\d^*) } \g}
\sum_{n=0}^{\bar N}\;  \SS_n,   \Eq (5.18) $$
where
$$  \eqalign {
& \SS_n =   e^{\sfrac \b  \g  \sfrac 12
\left (|\bar \eta_1-1|+|\bar \eta_2-1| \right)
\left (   \FF^*  +  2 f_1   \right) } \times \cr
&\sum_{\ell \in I}
\;\sum_{\bar \eta_1,\bar \eta_2\in \{-1,1\}^2}
\;\sum_{\ell_1=N_1}^{N_1+ R_1 }\;\sum_{\ell_2= N_2- R_1}^{N_2}
\sum_{A,|A|=n}
 \frac { Z_{[\ell_1, \ell_2]}^{0,0}
\left ( \eta^{\d,\z_1} (\ell_1)= \bar \eta_1,\eta^{\d,\z_1}
(\ell_2)=\bar
\eta_2, \AA(A,n),  B(\ell)
\right ) } {  Z_{[\ell_1, \ell_2]}^{0,0} \big (
\EE (+1,\ell_1,\ell_2,\bar \eta_1 ,\bar \eta_2)
\big )}.\cr
}\Eq (5.18a) $$

We must estimate $ \SS_n $ for any   $n$,
taking  care of the probability subspaces on which we are
working.  At first sight one could
have thought that the presence of $n$--changes of phases would
simplify the analysis, at least for $n$ large, due to  the
presence of terms proportional to $ \exp({-n \sfrac \b \g \FF^* })$.
Unfortunately  this is not the right picture since
we must control  the local contributions  of the magnetic field.
{ For $ \D' \subset [  \a_0^*, \a_1^*]$ we only know that
$\sum_{\a \in  \D'} \chi(\a) \ge -2 \left ( \FF^* -f \right
)$. }The analysis is therefore  more delicate, being  summarized in
Lemmas \eqv (BB1) and  \eqv (BB2)   below.

To  complete the  estimate of the expression in  \eqv (5.18) we
 need to   sum  up
the upper bounds of  the $\SS_n$, cf.  Lemmas  \eqv (BB1) and   \eqv
(BB2). For  this we use the
following inequalities that follow  from Taylor formula:
for all $x>0$,
$$
(1+x)^N-\sum_{k=0}^{l} {{N}\choose{k}}x^k \le \frac{(xN)^{l+1}}{(l+1)}
e^{(N-l-1)x} \le (xN)^{l+1} e^{Nx}. 
\Eq(5.2010200)
$$

Recall that $\bar N \frac \e \g =\frac \e \g [(\ell_2-\ell_1) \frac \g \e]
\le (\ell_2-\ell_1)\le \frac {2Q} \g$;
$ |I| \le \frac {2Q} \g$.
To simplify the computations, when  necessary, we take
half of negative part in the exponential
to  compensate the positive part. We also use $\z_5^{\frac {1-z}{2}}
>\z_5$.
  Denote $\O_5=\O_{51}\cap\O_{53}$, with $\O_{53}$ as
in Lemma \eqv(BB1). After some easy however lengthy computations,
using \eqv(5.0100), we see that on $\O_{5}\cap\O^+(\e,Q)$,
$$
\mu_{\b,\th,\g}\left ( \exists \ell \in I,
\eta(\ell)\neq  1\right)
\le
2\left( \frac{3Q}{\g^2}\right)^5
e^{-\frac \b \g \left\{(\frac{\k(\b,\th)}{4}
\d\z_4^3)\wedge \FF^*\right\}}
+
28|R_1|^2\left(\frac {2Q} \g \right)^5
 e^{-\frac \b \g \frac{\e^{1/4}}{5}}
e^ { \left \{ \frac {4Q} \g  e^{-\frac \b \g \frac {\e^{1/4}} 5  } \right \} }
\Eq(5.201021)
$$
which is \eqv(5.2C).
 \eqv(BB.50) follows from  \eqv(5.190) since  $\P[\O_{RE}] \ge 1-\g^2$.
This ends the proof of Proposition \eqv(510P) if we assume 
Lemmas \eqv (BB1) and  \eqv (BB2).
\eop

\vskip .5truecm
\noindent{\bf\Lemma (BB1) (n=0)}
 {\it
For $f_1$ given by \eqv(5.01) and $f_2$ given by \eqv(5.0100a),
for $\frac 12 >z>0$,
there exists $\O_{53}$ with
$$
\P[\O_{53}]\ge 1
-4 \g^2-
\frac{
2e^{-\frac{\b^2}{8Q\z_5^z c^2(\b,\th)}}
}
{1-e^{-\frac{\b^2}{8Q\z_5^z c^2(\b,\th)}}}
\Eq(5.190)
$$
such that on $ \O^+(\e,Q)\cap\O_{53}\cap \O_{51}$,
$$ \SS_0 \le
R_1^2 |I|
e^{\frac \b \g  \big(  4f_1 +f_2\big)}
Ge^{-\frac \b \g   \e^{1/4}
}\Eq (BB.8)$$
where }
$$
G=
e^{\frac \b \g  (4\z_5 +2 f_1+  16 \z_5^{\frac{1-z}{2}}) }
\Big(1+e^{-\frac \b \g \frac{\k(\b,\th)}{4} \d\z_5^3}\Big).
\Eq (BS1)
$$

\noindent{\bf  Proof.}
In this case  the profiles have
no change of phases, therefore  we must have  $\bar \eta_1=\bar \eta_2$.
If $\bar \eta_1=\bar \eta_2=+1$  and we take  $|A|=0 $ in \eqv (M.1E), we have
$$   \left\{  \eta^{\d,\z_1} (\ell_1)= \bar \eta_1,\eta^{\d,\z_1}
(\ell_2)=\bar \eta_2, \AA([\ell_1,\ell_2],A,0),  B(\ell)
\right \}  = \emptyset $$
and there is nothing to prove. So we  consider   the case  $\bar \eta_1=
\bar\eta_2=-1$. With this choice
the set  to estimate in   \eqv (5.18a) is
$$ \eqalign { & \left\{  \eta^{\d,\z_1} (\ell_1)=
\bar\eta_1,\eta^{\d,\z_1} (\ell_2)=\bar\eta_2,
\AA([\ell_1,\ell_2],A,0),  B(\ell)
\right \} \cr & = \left\{  \eta^{\d,\z_1}
(\ell_1)=\eta^{\d,\z_1} (\ell_2)=-1, \forall \tilde \ell \in
[\ell_1+1,\ell_2-1], \eta^{\d,\z_4}(\tilde \ell)=-1 \right\}
 \equiv  R_{1,4}(-1, [\ell_1,\ell_2]) } \Eq(5.matteo)  $$
To estimate the quotient of the two partition functions in \eqv (5.18a),
we need to  extract the contribution of the
magnetic field as we did in the proof of   Proposition \eqv(59),
see \eqv(5.61).
If, however, we proceed   exactly as it was done there,
  we should  get $\z_4$ instead of $\z_5$ on  the right hand
side of \eqv(5.190).  Since $\z_4$ is fixed and $Q$ will be large
at the end, such an estimate would be useless.  Therefore an extra
step is needed.
For  $\bar \eta =\pm 1$, $\ell'_1< \ell'_2$ such that 
$\ell'_1-\ell'_2 > 4\ell_0+8$,    $\ell_0>0$   to be 
chosen later, let us denote  
$$
  R_{5}(\bar \eta,[\ell'_1,\ell'_2])
   \quad=\Big\{m^{\d^*}_{[\ell'_1,\ell'_2]}\colon
   \eta^{\d,\z_5} (\ell)= \bar\eta,
\,  \forall
 \ell \in [\ell'_1 ,\ell'_2 ]\Big\} 
 \Eq(5.1913)
$$
and
 $$   R_{1,4,5}(-1,[\ell_1,\ell_2])  
\equiv R_{1,4,5}(-1,[\ell_1,\ell_2])(\ell_0) =
R_{1,4}(-1, [\ell_1,\ell_2])  \cap
R_{5}(-1,[\ell_1+\ell_0,\ell_2- \ell_0]).  \Eq(5.1912)
$$

Then we write, see \eqv (5.18a) and \eqv (5.matteo)
$$
\frac{ Z_{[\ell_1,\ell_2]}^{0,0}(
  R_{1,4}(-1, [\ell_1,\ell_2]))}
{Z_{[\ell_1, \ell_2]}^{0,0}
\big (\EE (+1,\ell_1,\ell_2,-1 ,-1))
\big ) }
=
\frac{
 Z_{[\ell_1,\ell_2]}^{0,0}(
  R_{1,4,5}(-1, [\ell_1,\ell_2]) )}
{Z_{[\ell_1, \ell_2]}^{0,0}
\big (\EE (+1,\ell_1,\ell_2,-1 ,-1)
\big )  }
\times
\frac{ Z_{[\ell_1,\ell_2]}^{0,0}(
  R_{1,4}(-1, [\ell_1,\ell_2]))}
{ Z_{[\ell_1,\ell_2]}^{0,0}(
  R_{1,4,5}(-1, [\ell_1,\ell_2]) )}. 
\Eq(5.19120)
$$
The choice of $\ell_0$     is related to the   needed  length to
go from   $\eta^{\d,\z_4} (0) =  \eta  $ to
$\eta^{\d,\z_5} (l_0) =\eta  $ knowing that we are within a run of
$\eta^{\d,\z_4}=\eta$. It is determined estimating the last term
in \eqv (5.19120) from which we start.
Since $  R_{1,4,5}(-1, [\ell_1,\ell_2])
\subset   R_{1,4}(-1, [\ell_1,\ell_2])$ we have
$$
1\le
\frac{ Z_{[\ell_1,\ell_2]}^{0,0}(
  R_{1,4}(-1, [\ell_1,\ell_2]))}
{ Z_{[\ell_1,\ell_2]}^{0,0}(
 R_{1,4,5}(-1, [\ell_1,\ell_2]))}
\le
1+
\frac{ Z_{[\ell_1,\ell_2]}^{0,0}(
  R_{1,4}(-1, [\ell_1,\ell_2])
\cap
(  R_{1,4,5}(-1, [\ell_1,\ell_2]))^c
)}
{ Z_{[\ell_1,\ell_2]}^{0,0}(
  R_{1,4,5}(-1, [\ell_1,\ell_2]))}
\Eq(5.191210)
$$
 From Corollary \eqv(6002) it follows  that on $\O_{RE}\supset
\O_{51}$,
 if
$$
\d\z_5^3 >
\frac{512(1+\th)}{\k(\b,\th) \a(\b,\th,\z_0)}
\sqrt{\frac{\g}{\d^*}} \log \frac{\d^*}{\g}
\Eq(5.1912101)
$$  where $\a(\b,\th,\z_0)$ is defined in \eqv(5.P1), and   $ \ell_0$
is chosen\footnote{$^{*}$ }{\eightrm
The  $ L_0$ chosen in \eqv (5.12P)  is obtained setting $d=2$ in 
    Corollary \eqv(6002).
 }   as  $ L_0$  defined in \eqv (5.12P),
then
$$
\frac{ Z_{[\ell_1,\ell_2]}^{0,0}(
  R_{1,4}(-1, [\ell_1,\ell_2])
\cap
(  R_{1,4,5}(-1, [\ell_1,\ell_2]))^c
)}
{ Z_{[\ell_1,\ell_2]}^{0,0}(
  R_{1,4,5}(-1, [\ell_1,\ell_2]))}
\le e^{-\frac \b {\g} \frac{ \k(\b,\th)}{4} \d\z_5^3},
\Eq(5.191220)
$$
uniformly with respect to $[N_1, N_2] \subset [-Q \g^{-1}, Q
\g^{-1}]$, $ \ell_1 \in [N_1,N_1+R_1]$, and $\ell_2 \in
[N_2-R_1,N_2]$.
 To treat the first term in the right hand side of \eqv(5.19120), 
recalling that, see    \eqv (5.12002),
$$ \EE (+1,\ell_1,\ell_2,-1 ,-1) =  
\widetilde \RR_1(+1,[\ell_1,\ell_2]) \cap
\widetilde \RR_5(+1, [\ell_1+L_0,\ell_2-L_0])
  $$
we first split the interval $[\ell_1,\ell_2]$ into three intervals
$[\ell_1,\ell_1+L_0-1]$, $[\ell_1+L_0, \ell_2-L_0]$ and $
[\ell_2-L_0+1,\ell_2]$.  On the first and the last interval,
we use  a block spin representation,   the rough estimate
Lemma \eqv(60002) with $p=2 + [(\log Q)/(\log (1/\g))]$,
and then  the symmetry $m\rightarrow Tm$  of the block spin model.
Thus,   on  $\O_{RE}=\O_{RE}(\g,\d^*, Q)\supset \O_{51}$, we get 
for the first term
$$
\eqalign{
&\frac{ Z_{[\ell_1,\ell_1+L_0-1]}^{0,m_{-}}(\eta^{\d,\z_1}(\ell_1)=-1,
\forall \ell \in [\ell_1+1, \ell_1+L_0-1],\, \eta^{\d,\z_4}(\ell)=-1)}
{Z_{[\ell_1,\ell_1+L_0-1]}^{0,m_{+}}(\eta^{\d,\z_1}(\ell_1)=1)}\cr
& \le
e^{\frac{\b}{\g} 6(1+\th)\frac{1}{\a(\b,\th,\z_0)}
(\log\sfrac{\d^*}{\g})\d^*\vee \sqrt{\sfrac{\g}{\d^*}}}
\le e^{\frac{\b}{\g} f_{1}}
}\Eq(5.ma1)
$$
and in the very same way for the other term.
  Therefore, on $\O_{RE}\supset \O_{51}$, we have
$$
\eqalign { & \frac{
 Z_{[\ell_1,\ell_2]}^{0,0}(
  R_{1,4,5}(-1, [\ell_1,\ell_2])(L_0))}
{Z_{[\ell_1, \ell_2]}^{0,0}
\big (  \EE (+1,\ell_1,\ell_2,-1,-1) )  \big ) }  \le
e^{\frac{\b}{\g} 4\z_5}
e^{\frac{\b}{\g} 2 f_1}
\frac
{Z_{[\ell_1+L_0,\ell_2-L_0]}^{0,0}
(R_{5}(-1, [\ell_1+L_0,\ell_2-L_0]))}
{Z_{[\ell_1+L_0,\ell_2-L_0]}^{0,0}
( R_5 (1, [\ell_1+L_0,\ell_2-L_0] ))}
\cr &  =
e^{- \b\D\GG(m^{\d^*}_{\b,[\ell_1+L_0, \ell_2-L_0]} )}
\frac
{Z_{-1,0}([\ell_1+L_0,
\ell_2-L_0])}{Z_{+1,0,}([\ell_1+L_0, \ell_2-L_0])} }
  \Eq(5.1924)
$$
where
$
\D \GG(m^{\d^*}_{\b,[\ell_1+L_0, \ell_2-L_0]} )= 
\sum_{x\in \CC_{\d^*}([\ell_1+L_0, \ell_2-L_0])} X(x)
$
and  the  remaining term  is defined in \eqv(4.9) with $R(\eta)$ 
replaced by
$R_5(+,[\ell_1+L_0,\ell_2-L_0])$.
The equality in  \eqv (5.1924) is obtained   by  extracting the
 main contribution of the random field as we did in
 \eqv(5.61).

To estimate the last  term in \eqv(5.1924), we  use
Proposition \eqv(401) and  \eqv(5.65)
with $\z=\z_5$, $a=\z_5^{\frac{1-z}{2}}$, for some $0<z<1/2$.
Using \eqv(5.004P), this entails  that on a subset $\O_{54}$, with
$$
\P[\O_{54}] \ge
1-
\frac {2e^{-\frac{\b^2}{8 Q \z_5^zc^2(\b,\th)}}}
{1-e^{-\frac{\b^2}{8 Q \z_5^zc^2(\b,\th)}}}
\Eq(5.1925)
$$
we have
$$
\max_{[\ell_1,\ell_2]\subset[-Q\g^{-1},Q\g^{-1}]}
\frac
{Z_{-1,0}([\ell_1+L_0,
\ell_2-L_0])}{Z_{+1,0,}([\ell_1+L_0, \ell_2-L_0])}
\le e^{\frac{\b}{\g} 16 \z_5^{\frac{1-z}{2}} }
\Eq(5.1926)
$$

Some care is necessary  to estimate the contribution of the first
factor of the r.h.s.  of \eqv(5.1924). By definition, on
$\O^+(\e,Q)$, we have $\D^+ \GG(m^{\d^*}_{\b,
[\a^*_0,\a^*_1]})\ge  2 \FF^* +f\equiv 2\FF^*+\e^{1/4}
$. However the random
contribution we   extracted in \eqv(5.1924) is merely $\D^+
\GG(m^{\d^*}_{\b,[\ell_1+L_0, \ell_2-L_0]} )$, with $\ell_1 \in
[N_1, N_1+R_1]$, $\ell_2 \in [N_2-R_1, N_2]$.
It is easy to check
that there exists a subset $\O_{55}$, { that depends  on
$(\g,\d^*, Q)$ with $\P[\O_{55}] \ge 1-8 \g^2$, such that on
$\O_{55}$,  uniformly with respect to $[N_1, N_2] \subset [-Q
\g^{-1}, Q \g^{-1}]$, and $ \ell_1 \in [N_1,N_1+R_1]$, $\ell_2 \in
[N_2-R_1,N_2]$,  we have}
$$
 e^{-{\b}\D\GG(m^{\d^*}_{\b,[\ell_1+L_0, \ell_2-L_0]} )}
\le e^{-\frac{\b}{\g}\big(2 \FF^* +\e^{1/4}
 - f_2\big)},
 \Eq(5.1927)
$$
where $f_2$ is given in \eqv(5.0100a).
Collecting    \eqv (5.ma1), \eqv (5.1924) and \eqv (5.1926),  on
$\O^+(Q,f)\cap\O_{RE} \cap \O_{54} \cap \O_{55}$,  we have
$$
\frac{
 Z_{[\ell_1,\ell_2]}^{0,0}(
 R_{1,4,5}(-1, [\ell_1,\ell_2]))}
{Z_{[\ell_1, \ell_2]}^{0,0}
\big (\EE (+1,\ell_1,\ell_2,\bar \eta_1 ,\bar \eta_2)
\big ) }
\le
e^{+\frac{\b}{\g} \big(4\z_5 + 2  f_1 + 16 \z_5^{\frac{1-z}{2}}\big) }
e^{-\frac{\b}{\g}\big(2 \FF^* +\e^{1/4} 
- f_2\big)}
\Eq(5.19270)
$$
  Now, collecting \eqv(5.19120), \eqv(5.191210), \eqv(5.191220)
and   \eqv(5.19270), and calling $\O_{53}= \O_{54} \cap
\O_{55}$, on $\O^+(Q,f)\cap \O_{53}\cap \O_{51}$ we have
$$
\frac { Z_{[\ell_1, \ell_2]}^{0,0}
\left ( \eta^{\d,\z_1} (\ell_1)=
\bar\eta_1,\eta^{\d,\z_1} (\ell_2)=\bar\eta_2,
\AA ([\ell_1,\ell_2],A,0),  B(\ell)
\right ) }
{  Z_{[\ell_1, \ell_2]}^{0,0}
\big (   \EE(+1,[\ell_1,\ell_2],
\bar\eta_1,\bar\eta_2)  \big ) }
\le G e^{-\frac{\b}{\g}\big(2 \FF^* +\e^{1/4} 
-f_2 \big)}
\Eq(5.191230)
$$
from which we get easily \eqv(BB.8).\eop

 \smallskip
\noindent{\bf \Lemma (BB2) ($n\ge 1 $) }{\it On $\O(54)\cap
\O_{51}\cap \O^+(\e,Q) \cap \O_\e$, we have 
$$
\SS_1 \le
R_1^2 |I|  [\ell_2-\ell_1]
e^{\frac \b\g
\left (2 f_1
+f_2+ f_3 +4\z_1 \right ) }
G  e^{-  \frac \b \g { \e^{1/4}
}}
 \Eq (BB.14)
$$
$$
\SS_n
\le
R_1^2  |I|
{{ \bar  N}\choose {n}}
\left ( \frac \e \g \right)^n
e^{ n\frac \b \g  \left ( f_3 +4 \z_1 \right )}
G^{\frac n2}
e^{-\frac \b \g \frac n {2} \e^{1/4} } \quad  n
\,\hbox {even}  \quad \bar \eta_1=\bar \eta_2 =1
\Eq (D.40)
$$

$$
 \SS_n
 \le
R_1^{2}
  |I| e^{  \frac \b \g  \left(4f_1
+2f_2  \right)
}
{{ \bar N}\choose {n}}   \left ( \frac \e \g
\right)^n
G^{\frac  n 2 +1 }
e^{n \frac \b  \g
\left ( f_3 +4 \z_1  \right ) }
e^{ -\frac \b\g \left (\e^{1/4}
 [3-n]^+ + \frac {n-2}2 \e^{1/4}
\right ) } \quad  n \ge 2     \quad \bar \eta_1=\bar \eta_2 =-1
\Eq(D.5)
 $$

$$
 \SS_n
\le R_1^2 |I|  e^{\frac \b \g
2 f_1  }
{{ \bar  N}\choose {n}}   \left ( \frac \e \g \right)^n
e^{n \frac \b \g   \left [  f_3 +4 \z_1  \right]}
G^{\frac { n+1}{ 2}}
e^{-\frac \b \g  \left ( \frac {n-1} 2  \e^{1/4}
   \right ) } \quad  n >1 \quad\hbox {odd}
 \Eq (BB.50100)
$$
where $f_1$ is defined in \eqv(5.01) $f_2$ in \eqv (5.0100a),  
$f_3$ in \eqv(5.0001)  and $G$ in \eqv (BS1).}

\smallskip

\noindent {\bf  Proof.} We prove explicitly  the case $n=1$.
The $n>1$ can be   done similarly following the general strategy
outlined later.
When $n=1$ the magnetization profiles have only
one change of phases   and are therefore    compatible only with boundary
conditions $\bar \eta_1 \neq \bar \eta_2$.
Suppose    that $\bar \eta_1=-\bar \eta_2=   1$. The reverse case is done
similarly.
Denote by $r_1$ the index of the $\frac \e \g$
block in which the change of phases occurs. When
$ [r_{1}\frac
\e {\g}-R_1, (r_{1}+1)\frac
\e {\g}+R_1 )] \subset [ N_2-R_1-\frac \r \g,\ell_2]$
we have
$   \left\{  \eta^{\d,\z_1} (\ell_1) = 1,\eta^{\d,\z_1} (\ell_2)=-1,
\AA_{[\ell_1,\ell_2]}(A,1),
B(\ell)
\right \}  = \emptyset $
since $\ell\in [N_1+R_1+\frac{\rho}{\g},N_2-R_1-\frac{\rho}{\g}]$.
Therefore we may assume that
$ [r_{1}\frac
\e {\g}-R_1, (r_{1}+1)\frac
\e {\g}+R_1 )] \subset [\ell_1,  N_2-R_1-\frac \r \g] $.
 We  split  the interval
$[\ell_1,\ell_2]$ into three adjacent intervals
$[\ell_1,q_1-R_1]$, $[q_1-R_1+1,q_1+R_1-1]$ and
$[q_1+R_1,\ell_2]$, assuming that the change of phases happens in
the interval $[q_1-R_1,q_1+R_1]$.    Recalling Definition
\equ(def1)  in Section 4,  one has $\eta^{\d,\z_1}(\tilde
\ell)$ is equal to $+1$ for $\tilde \ell=\ell_1$ and for $\tilde
\ell=q_1-R_1$ while it is equal to $-1$ for $\tilde \ell=q_1+R_1$.
We associate  the interactions between the intervals to the middle
interval. Suitably restricting the denominator  we get

$$
\eqalign{
&\frac { Z_{[\ell_1, \ell_2]}^{0,0} \left ( \eta^{\d,\z_1} (\ell_1)= +1,
\eta^{\d,\z_1} (\ell_2)=-1, \AA_{[\ell_1,\ell_2]}(A,1),  B(\ell)
\right ) }
{Z_{[\ell_1, \ell_2]}^{0,0} \big (
\widetilde \RR_1(+1,\ell_1,\ell_2)\cap\{\eta^{\d,\z_5}
(\ell_2-L_0)=+1 \} \big )}\le
e^{\frac{\b}{\g} 4\z_1}\times \cr
&\quad
\frac { Z^{0,0}_{ [\ell_1,q_1-R_1 ]}
 \big (   R_{1,4}(+1,[\ell_1,q_1-R_1])\big ) }
{ Z^{0,0}_{ [\ell_1,  q_1-R_1]}
\big (\tilde \RR_1(+1,\ell_1,q_1-R_1)\big ) } \times
\frac { Z^{m_+,m_-}_{ [ q_1-R_1+1, q_1+R_1-1 ]}}
{   Z^{m_+,m_+}_{ [q_1-R_1+1,  q_1+R_1-1 ]}
\big ( \widetilde \RR_1(+1,q_1-R_1+1,q_1+R_1-1) \big ) }\times\cr
&\quad
\frac
{ Z_{ [q_1+R_1, \ell_2]}^{0,0}\big( R_{1,4}(-1,[q_1+R_1,\ell_2])\big )}
{ Z_{[q_1+R_1, \ell_2]}^{0,0}  \big (\widetilde
\RR_1(+1,q_1+R_1,\ell_2)\cap
\{\eta^{\d,\z_5}(\ell_2-L_0)=+1\}\big )}
\cr
}\Eq(5.2010)
$$
Since $R_{1,4}(+1,[\ell_1,q_1-R_1]) \subset \tilde
\RR_1(+1,\ell_1,q_1-R_1)$, see \eqv(5.matteo) and \eqv(5.L), the
first ratio on the right
 hand side of \eqv(5.2010) is smaller than $1$.
The second ratio in \eqv(5.2010) is treated in a similar way as in
the proof of Lemma \eqv(T1C). However, since the volume we are
considering is $[q_1-R_1+1,q_1+R_1-1]$, { the error terms} that
come from the block spin approximation and the rough estimates,
see Lemma \eqv(60002), are $ e^{\frac{\b}{\g} f_3}$
with $f_3$ given in \eqv (5.0001).
Therefore,  on $\O_{RE}\supset \O_{51}$, uniformly with
respect to the position of the change of phases in the interval
$[-Q\g^{-1},Q\g^{-1}]$, we have 
$$
\frac { Z^{m_+,m_-}_{ [ q_1-R_1+1, q_1+R_1-1 ]}}
{   Z^{m_+,m_+}_{ [q_1-R_1+1,  q_1+R_1-1 ]}
\big ( \widetilde \RR_1(+1,q_1-R_1+1,q_1+R_1-1) \big ) }
\le  e^{-\frac \b\g (\FF^*  -f_3)  }
\Eq(5.20112)
$$

It remains to treat the last ratio in \eqv(5.2010).
We claim that on $  \O_{54}\cap \O_{RE} \supset \O_{53}$,
  see just before \eqv(5.191230), 
we have
$$
\eqalign{&
\frac
{ Z_{ [q_1+R_1, \ell_2]}^{0,0}
\big(  R_{1,4}(-1,[q_1+R_1,\ell_2])\big )}
{ Z_{[q_1+R_1, \ell_2]}^{0,0}  \big (\widetilde
\RR_1(+1,q_1+R_1,\ell_2)\cap
\{\eta^{\d,\z_5}(\ell_2-L_0)=+1\}\big )}
 \cr & \le
e^{\frac \b \g 4\z_5}
e^{\frac{\b}{\g} 2 f_1}
e^{\frac \b \g 16\z_5^{\frac{1-z}{2}}}
(1+e^{-\frac{\b}{\g}\frac{\k(\b,\th)}{4}
\d\z_5^3})
e^{- \b \D\GG(m^{\d^*}_{\b,[q_1+R_1+L_0+1,\ell_2-L_0-1]})}
=
 Ge^{-\b \D^+\GG(m^{\d^*}_{\b,[q_1+R_1+L_0+1,\ell_2-L_0-1]})}
\cr }
\Eq(5.20113)
$$
where $G$ is defined in \eqv(BS1).

Let us explain where { those terms come from:} We have written {
the ratio} on the left hand side of \eqv(5.20113) as a product of
two ratios in the very same way as  in \eqv(5.19120). The second
ratio gives    the term $(1+e^{-\frac{\b}{\g}\frac{\k(\b,\th)}{4}
\d\z_5^3})$ as in \eqv(5.191210) and \eqv(5.191220),  and
this occurs on $\O_{RE}$. The first ratio was treated by first
splitting  the volume $[q_1+R_1,\ell_2]$ in three  intervals
$[q_1+R_1, q_1+R_1+L_0]$, $[q_1+R_1+L_0+1, \ell_2-L_0-1]$, and
$[\ell_2 -L_0,\ell_2]$. { The first and the last intervals  give
us }the term $\exp ({\frac{\b}{\g} 2 f_1})$ that comes from the
rough estimates, and therefore occurs on $\O_{RE}$. There is also
a term   $\exp({\frac \b \g 4\z_5})$ that comes from the
interactions between the intervals. We remain with a term similar
to the left hand side of \eqv(5.1924) but in the  volume
$[q_1+R_1+L_0+1, \ell_2-L_0-1]$. It give us the term  $\exp({\frac
\b \g 16 \z_5^{\frac{1-z}{2}}})$ and the last term in
\eqv(5.20113) and this occurs on $\O_{54}$.  Collecting
\eqv(5.2010), \eqv(5.20112), and \eqv(5.20113), we have, on
$\O_{51}\cap \O_{53}$
$$
\eqalign {
& \SS_1
  \le
\sum_{\ell_1=N_1}^{N_1+ R_1} \sum_{\ell_2= N_2- R_1}^{N_2}
\sum_{r_1}^*\sum_{q_1=\frac {r_1}\g \e}^{\frac {r_1+1 }\g \e}
   \sum_{\ell= q_1-R_1 }^{\ell_2-\frac \rho \g}
e^{\frac \b\g   \left ( \FF^* +2 f_1 \right) }
e^{-\frac \b\g ( \FF^* -f_3)}
e^{\frac{\b}{\g} 4\z_1}\times\cr
&\quad \times
Ge^{- \b \D^+\GG(m^{\d^*}_{\b,[q_1+R_1+L_0+1,\ell_2-L_0-1]})}\cr &
\le
e^{\frac{\b}{\g} (2f_1+ f_3+4\z_1)}
G
 \sum_{\ell_1=N_1}^{N_1+ R_1} \sum_{\ell_2= N_2- R_1}^{N_2}
\sum_{r_1}^* \sum_{q_1=\frac {r_1}\g \e}^{\frac {r_1+1 }\g \e}
e^{-\b \D^+\GG(m^{\d^*}_{\b,[q_1+R_1+L_0+1,\ell_2-L_0-1]})}. 
\cr } \Eq (E.4)
$$
By  $\sum_{r_1}^*$ we denote  the sum over    blocks  of
length $\frac \e \g$ contained in the interval $ [\ell_1,
N_2-R_1-\frac \r \g ]$, so that  $\sum_{r_1}^*1 \le
[\ell_2-\ell_1] \frac \g \e\le \frac {2Q} \e $. The contribution
of the magnetic field   in \eqv (E.4) is estimated using
Lemma \eqv(lem112) and therefore occurs on $\O_{\e}$. By
definition,  for any value under consideration of 
$q_1\in [ {\frac {r_1}\g \e},{\frac {r_1+1 }\g \e}]$ and  
$r_1$,  we in fact have that 
$$
\eqalign{
  \D\GG(m^{\d^*}_{\b,[q_1+R_1+L_0+1, \ell_2-L_0-1]})& =\sum_{x\in
\CC_{\d^*}([q_1+R_1+L_0 +1, \ell_2-L_0-1])} X(x) \cr
& = \frac 1\g \sum_{\a= \frac {q_1+R_1+L_0+1}\e \g}^{
\a^*_1}
\chi(\a)       - \sum_{x \in
\CC_{\d^*}(\ell_2-L_0,N_2)} X(x)
.} \Eq (SS.8)
$$
The point  is that $\a^*_1 $ is the end of a positive elongation,
that is a maximum,  and by construction
 $ |\a^*_1-\frac {q_1+R_1+L_0+1}\e \g| \ge \frac{\r}{\e}
$, recall $\r=\e^{\frac{1}{4(2+a)}}$.
  Therefore recalling   Lemma \eqv(lem112)  
 on $\O^+(Q,f)\cap \O_{\e}\cap \O_{53}$ we have
$$
 \D\GG(m^{\d^*}_{\b,[q_1+R_1+1, \ell_2]}) \ge
\frac{1}{\g} \left( \e^{1/4}
-f_2\right).
\Eq(SB110)
$$
This entails that on $\O_{51}\cap \O_{53}\cap \O^+(Q,f)$, see \eqv (E.4),
$$
\SS_1 \le R_1^2 |I|  [\ell_2-\ell_1] e^{\frac \b\g \left (2 f_1
+f_2+ f_3 +4\z_1 \right ) } G  e^{-  \frac \b \g { \e^{1/4}
}}.  \Eq (BB.11) $$
  \noindent{\it Remark:}
The fact that $N_2=\a^*_1\e/\g$ is
a maximum is crucial here. In  the case $\bar \eta_1=-1,\bar
\eta_2=+1$, it would be essential that  $\a^*_0$ is a minimum.

\noindent{\it General strategy.} The  estimate of the terms { with
} $n>1$ in \eqv(5.18)  is a simple modification of what we  did in
the cases $n=0$ and $n=1$. Let us summarize the general
strategy: 

\noindent a) Similarly to \eqv(5.2010), if $n$ changes occur 
we bound the ratio of two constrained partition functions by the
product of ratios over the $n$ intervals $[q_i-R_1,q_i+R_1]$ where
the changes occur, a factor $e^{\frac{\b}{\g} 4n\z_1}$ and a
product of ratios over the intervals with no change of phases.

\smallskip
\noindent b) The contribution of a ratio { corresponding to} a
change of phases is estimated by $e^{\frac{\b}{\g}[\FF^*-f_3]}$
where $f_3$ is given in \eqv(5.0001), as we did in \eqv(5.20112).
This holds on $\O_{RE}$ since a rough estimate is used and
therefore on $\O_{51}$.

\smallskip

\noindent c)  The contribution of a ratio over an interval, say $\JJ$,
 where there is no change of phases is bounded by 1 when the
 profile
is  $\z_4$-near  $m_\b$, that is for a run of
$\eta^{\d,\z_4}=+1$. { If, instead, the  profile gives a run of
$\eta^{\d,\z_4}=-1$,
as in \eqv(5.20113), then the corresponding ratio is bounded from
above by
$$
e^{\frac \b \g (4 \z_5+ 2 f_1+  16\z_5^{\frac{1-z}{2}})}(1+e^{-\frac \b \g
\frac{\k(\b,\th)}{4}
\d\z_5^3}) e^{-\frac \b \g \D\GG(m^{\d^*}_{\b,\JJ})}
 = G e^{-\frac \b \g \D\GG(m^{\d^*}_{\b,\JJ})},
\Eq(BBB.111)
$$
 on $\O_{51}\cap
\O_{53}$.}
\smallskip
\noindent d)  The contribution of $\D\GG(m^{\d^*}_{\b,\JJ})$ in
\eqv(BBB.111) depends whether $\JJ$ is between two consecutive
changes of phases or not, with {\bf} $\JJ$ being located  at
an extreme of $I$.
 In the first case we use
$ \sum_{\a\in \D'}\chi (\a)  \ge - \left ( 2 \FF^* -f \right)
\equiv -\left( 2\FF^* -\e^{1/4})\right)
 $
which holds  on $\O^+(\e,Q)$. In the second case,  if the
length of   $\JJ$ is larger than
$\frac{\rho}{\e}\equiv\e^{-\frac{5+4a}{8+4a}}$, we apply  Lemma
\eqv (lem112) as in \eqv(SB110), on $\O_\e$.  This gives
$\D^+\GG(m^{\d^*}_{\b},\JJ) \ge \g[
\e^{1/4}-f_2]$. Otherwise,   we use the fact that
$$  \inf_{\ell_1\le t \le \ell_2 } \sum_{\a=t}^{\a^*_1}
\chi(\a) \ge 0  \qquad \qquad  \inf_{\ell_1\le t \le \ell_2 }
\sum_{\a=\a_0^*}^{t} \chi(\a) \ge 0,  \Eq (BB.40) $$ since
$\a^*_1$ is the location of a  maximum   and $\a_0^*$ the location
of  a minimum.

\noindent e) At least there are two  factors in \eqv(BB.14),
\eqv(D.5), and \eqv(BB.50100)  that come from
$$ \sum_{r_1,...r_n} 1 \le {{ \bar N}\choose {n}} \qquad \qquad
\sum_{q_1,...q_n} 1 \le  \left ( \frac \e \g\right )^n.\Eq (M.12)$$
\eop

\noindent{ \bf Proof of Theorem \eqv (t1)  } The proof of Theorem
\eqv (t1) is a consequence of  Proposition \eqv (510P) and of  the
next  choice of parameters.  Take  $g(\cdot)$  satisfying the
hypothesis of Theorem  \eqv(t1),  i.e., 
 $g(x)$ is increasing, $ g(x) \ge 1$ diverges
as  $x\uparrow \infty$, $ x^{-1} g^{{38}}(x) \le 1$ and
$$
 x^{-1} g^{{38}}(x)\downarrow 0,
 \Eq(Sancho.1)
 $$
$$
\e^{1/4}=\frac 5 {g(\frac {\d^*}{\g})},
\Eq(5.7029)
$$
 $$
Q=\exp \left( \frac{\log g(\frac{\d^*}{\g})}{\log\log g(\frac{\d^*}{\g})}
\right)
,\Eq(5.7030)
$$
$$
\z_5= 
\frac 1{2^{18}c^{6}(\b,\th)}\frac {1}{ g^3(\frac{\d^*}{\g})},\quad z=\frac13,
\Eq(5.97033)
$$
$$
\z_1=\frac 1{160 g(\frac {\d^*}{\g})}\quad{\rm and}\quad
\d= \frac 1{5(g(\frac {\d^*}{\g}))^{1/2}}.
\Eq(5.97037)
$$
 First we have to  check that the $\CC_0$ constraints are  satisfied if the
 parameters are chosen as above.
 \eqv(5.001P) is immediate from \eqv(5.97037) and
\eqv(Sancho.1).  \eqv(5.002) is just \eqv(TE.1) with the
choice in  \eqv(5.97037).  \eqv(5.004P) is just
\eqv(TE.2) with \eqv(5.97033) and \eqv(Sancho.1).
 \eqv(5.003) is immediate from \eqv(5.97033) and
\eqv(Sancho.1).
  \eqv(5.00410) is immediate from \eqv(5.7030), \eqv(Sancho.1),
$\d^*<1$  and $\g/\d<d_0$,  by taking $d_0$ small enough.
 \eqv(5.00411) follows from \eqv(Sancho.1) by taking $\g_0$
and $d_0$ small enough.
 It is immediate  to check that \eqv(5.0100) holds and also that
\eqv(5.2C) implies \eqv(2.2100) after easy simplifications.
 It remains to check  \eqv(5.02).
Notice  that \eqv(5.97033) gives $2^6 Q \z_5^z c^2(\b,\th)=
Q/g$,  and taking $d_0$ small enough we have $e^{-\b^2/(2^6 Q
\z_5^z c^2(\b,\th))}\le e^{-\b^2 \sqrt{g}}$. It is then  easy
to check that with our choice of $Q$ the leading  term in
\eqv(5.P01) is $5\e^{\frac{a}{16(2+a)}}$ from which we 
easily get  \eqv(5.02).

We then set 
$$
I(\o)=\left[ \frac{\e\a^*_0}{\g}+\frac{\rho}{\g}+R_1,
\frac{\e\a^*_1}{\g}-\frac{\rho}{\g}-R_1\right]
$$
and
$\t(\o)=+1$ if $\o\in \O^+(\e,Q)\cap \O_5$ and
$\t(\o)=-1$ if $\o\in \O^-(\e,Q)\cap \O_5$.
The estimates \eqv(2.210a) and \eqv(2.210b) are immediate consequences
of Proposition \eqv(Terminator). \eop

\noindent{\bf Proof of Theorem \eqv(t3)}
Since the proof follows from arguments similar to  the ones we already
used, we  will sketch it.
It is enough to consider two consecutive elongations
$$
\eqalign{
I_0&=\left[\frac{\a_0^*\e}{\g}+R_1+\frac\r\g,
\frac{\a_1^*\e}{\g}-R_1-\frac\r\g\right]\cr
I_1&=\left[\frac{\a_1^*\e}{\g}+R_1+\frac\r\g,
\frac{\a_2^*\e}{\g}-R_1-\frac\r\g\right]\cr
}\Eq(Sherazade.1)
$$
with $\hbox {sgn} I_0=+1$ and  $\hbox {sgn} I_1=-1$.
The main point is to estimate $\mu_{\b,\th,\g}(\CC_{0,1})$
where
$$
\CC_{01}
\equiv
\WW^c_1\left([\frac{\a_1^*\e}{\g}-R_1-\frac\r\g,
\frac{\a_1^*\e}{\g}+R_1+\frac\r\g],R_2,\z_4\right)
\cap \AA(\D_{2Q})
\Eq(5.990)
$$
where $\WW_1$ is defined in Definition \eqv(1).
Using Theorem \eqv(t62), we get
$$
\eqalign{
\mu_{\b,\th,\g}[\CC_{01}]&\le
\sum_{\bar \eta_1,\bar\eta_2 \in\{-1,+1\}}
\sum_{\ell_1=\frac{\e\a^*_1}{\g}-R_1}^{\frac{\e\a^*_1}{\g}}
\sum_{\ell_2=\frac{\e\a^*_1}{\g}}^{\frac{\e\a^*_1}{\g}+R_1}
\mu_{\b,\th,\g}[\CC_{01}\cap\{\eta^{\d,\z_1}(\ell_1)=\bar \eta_1,
\eta^{\d,\z_1}(\ell_2)=\bar \eta_2\}]\,\,+\cr
&\quad\quad\quad +3^4\left(\frac{2 Q}{\g^2}\right)^{5}
e^{-\frac \b \g \left((\frac{\k}{4}\d\z_4^3)\wedge \FF^*\right)}\cr
}\Eq(5.999)
$$
where $Q$ is defined in \eqv(5.7030). To study
$\mu_{\b,\th,\g}[\CC_{01}\cap\{\eta^{\d,\z_1}(\ell_1)=\bar \eta_1,
\eta^{\d,\z_1}(\ell_2)=\bar \eta_2\}]$, we decompose the event in
a way similar to \eqv(M.1E). Consider first the case $\bar
\eta_1=+1$. To be able to use   \eqv(D.40) where there is a positive
elongation, we need to have another $\eta^{\d,\z_1}(\ell)=+1$ for
$\ell$ on the left of $\frac{\a_1^*\e}{\g}-R_1-\frac\r\g$ instead
of the $\eta^{\d,\z_4}(\ell)=1$ that is present by Theorem
\eqv(t1). Using Theorem \eqv(t62), we will find such an $\ell$ in
the interval $[\frac{\a_1^*\e}{\g}-2R_1-\frac\r\g,
\frac{\a_1^*\e}{\g}-R_1-\frac\r\g]$, and  we   apply
\eqv(D.40) in the interval $[\ell,\ell_1]\subset
[\frac{\a_1^*\e}{\g}-2R_1-\frac\r\g,\frac{\a_1^*\e}{\g}]$. As a
consequence, on $\O_{51}\cap\O_{53}$, the Gibbs--probability to
have an even number of changes of phases $n\ge 2$ within $
[\frac{\a_1^*\e}{\g}-2R_1-\frac\r\g,\frac{\a_1^*\e}{\g}-R_1]$ is
bounded from above by
$$
56R_1^2(\frac{2Q}{\g})^5 e^{-\frac \b\g
\frac{\e^{1/4}}{4}}e^{\left\{\frac{Q}{\g} e^{-\frac \b\g\frac
{\e^{1/4}}{4}}\right\}}.  \Eq(5.9990)
$$

Consider now the case $\bar \eta_1=-1$.  Thus,  within
the interval $[\frac{\e\a^*_1}{\g}-R_1-\frac{\r}{\g},
\frac{\e\a^*_1}{\g}]$ the profile makes   an odd number of
changes of phases. When  $n>1$, we can apply
\eqv(BB.50100) and we get that the contribution of these terms is
also bounded from above by \eqv(5.9990).

So,  on the left of $\a^*_1$, there are two cases left from
the previous analysis: no change of phases when $\bar \eta_1=+1$ or
a single change of phases when $\bar \eta_1=-1$.

The same arguments apply on the right of $\a^*_1$ and therefore
 we can have at most
one change of phases on the  left of $\a^*_1$ and at most one
change of phases on its  right. Now we show that to have
simultaneously one change of phases on the right of $\a_1^*$ and
one on its  left has a very small  Gibbs--probability. It only
remains to consider the case 
 $\bar \eta_1=-1,\bar \eta_2=+1$.  Since
$\eta^{\d,\z_4}(\frac{\a^*_1\e}{\g}+R_1+\frac{\r}{\g})=-1$
the profile in $\CC_{01}$ makes two changes of phases on
the right of $\ell_1$
but since we are on $\AA(\D_{2Q})$
this means that there exists an $\ell \in
[\ell_1,\frac{\a^*_1\e}{\g}+R_1+\frac{\r}{\g}]$
 with $\ell-\ell_1 \ge \e/\g$
such that $\eta^{\d,\z_1}(\ell)=+1$.
That is within the  negative elongation that
occurs on the left of $\a^*_1$, we have
$\eta^{\d,\z_1}(\ell_2)=+1, \eta^{\d,\z_1}(\ell)=+1$.
 By using the very same argument
as in   \eqv(D.40), taking  care that here with the same
notations as in \eqv(5.1924), we will merely use
$$
\frac{Z^{0,0}_{[\ell_2+L_0,\ell-L_0]}(R_5(+1,[\ell_2+L_0,\ell-L_0]))}
{Z^{0,0}_{[\ell_2+L_0,\ell-L_0]}(R_5(-1,[\ell_2+L_0,\ell-L_0]))}
=e^{+ \b \D\GG(m^{\d^*}_{\b,[\ell_2+L_0,\ell-L_0]})} \frac{
Z_{+1,0}([\ell_2+L_0,\ell-L_0])}{
Z_{-1,0}([\ell_2+L_0,\ell-L_0])},  \Eq(5.9991)
$$
and since we are within a negative effective elongation we have
$$
\g \D\GG(m^{\d^*}_{\b,[\ell_2+L_0,\ell-L_0]}) \le 2\FF^* -\e^{1/4}.
\Eq(5.9992)
$$
As in   \eqv(D.40), the $2\FF^*$ cancels with the contributions of the two
changes of phases and we get a contribution which is bounded
from above by \eqv(5.9990).

Therefore we are left with  the three cases
$\bar \eta_1=-1, \bar\eta_2=-1$,
$\bar \eta_1=+1, \bar\eta_2=+1$, and
$\bar \eta_1=+1, \bar\eta_2=-1$ that belong to 
$\WW_1\left([\frac{\a_1^*\e}{\g}-R_1-\frac\r\g
\frac{\a_1^*\e}{\g}+R_1+\frac\r\g],R_2,\z_4\right)$.
 This ends the proof of Theorem \eqv(t3). \eop

\chap{6 Functional }6
\numsec= 6
\numfor= 1
\numtheo=1

In this section we  prove some estimates needed in Section  5,  based on
results on a finite volume version of  the excess free energy functional,
$\FF(\cdot)$,  see \eqv(AP.1). They are adaptation to
our case from results in [\rcite{CP}] and [\rcite{BZ}]. More care 
is needed here,
since the profiles  belong to
$ {\cal T} \subset L^{\infty}(\R,[-1,+1])\times L^{\infty}(\R,[-1,+1])$ 
instead of
$L^{\infty}(\R,[-1,+1])$ and the norm involved, see \eqv (7.2a), is stronger
 than the $L^\infty$ norm
used in  [\rcite{CP}] and [\rcite{BZ}].

\noindent{$\bullet$} {\bf I: Minimizers  in finite volume}

As in Section 2,  $\DD_\d$  denotes the partition of $\R$ into the
intervals $((\ell-1)\d,\ell\d ]$, $\ell \in \Z$, for $\d>0$
rational. In particular, if $\d =n\d'$, $n\in \N$,  then $\DD_\d$
 is coarser than $\DD_{\d'}$. For $r \in \R$,
we denote by $D^{\d}(r)$ the interval of $\DD_\d$  that contains
$r$. A function $f(\cdot)$ is $\DD_\d$--measurable if it is
constant on each interval of $ \DD_\d$.
In terms of the notation of Section 2, we
have
 $
D^\d(r)=\tilde A_\d([r/\d]+1),
$
where $[x]$  denotes the integer part of $x$.
We define for  $ m= (m_1,m_2) \in    {\cal T}$, see \eqv (space),
$$ m ^\d_i (r) =  \frac {1}{\d} \int_{D^{\d}(r)} m_i(s) ds
 \qquad \qquad i=1,2. \Eq (7.300) $$
By definition, the functions  $ m_i^\d (\cdot), i=1,2$,  are  constant on
each  $D^{\d}(r)$.
Definition \eqv (2.190)  is extended to functions in $ {\cal T}$,
 and, with an abuse of notation, we denote
$\eta^{\d,\z} (\ell)$,  $\ell \in \N$,
$$
\eta^{\d,\z}(\ell)=\cases{
 +1 &if $\,\forall_{u \in (\ell-1, \ell] }, \,
 \frac {1}{\d} \int_{D^{\d}(u)} ds \|m^{\d^*} (s)-m_\b\|_1 \le\z$;\cr
- 1 &if $\,\forall_{u \in (\ell-1, \ell] }, \,
 \frac {1}{\d} \int_{D^{\d}(u)} ds \|m^{\d^*}(s)- Tm_\b\|_1 \le\z$;\cr
0 & otherwise.\cr }
\Eq (7.2a)
$$
 If $ m^{\d^*}(x) = m^{\d^*} (x,\s)$  for $ x \in \CC_{\d^*}(I)$,
see  Section 2  before \eqv (2.14), and  we identify it with an element 
of ${\cal T}$,
piecewise constant on each
$ ((x-1)\d^*, x\d^* ] $,
and take  $\d=k \d^*$, then \eqv (7.2a) coincides with \eqv (2.190).
Given $L_0 \in \N$, $\d>\d^*>0, \z>0$
 and  $\eta\in \{-1,+1\}$
 we  set
 $$
 \VV_{\d, \z, L_0} (\eta)  = \left \{   m=( m_1,m_2)\colon (\eta
 m_1,\eta m_2) \in  \MM_\infty, \; \eta^{\d,\z}(0)= -\eta,
\eta^{\d,\z}(L_0)= \eta \right \}, \Eq (SS.4)
$$
where $\MM_\infty$ was defined in \eqv(AP.5).

\smallskip
\noindent{\bf \Lemma (5.30)  } {\it
 Let  $\b>1$ and  $\th>0$  satisfy \eqv (2.17). There exist
$\d_0=\d_0(\b,\th)>0$,  $\z_0=\z_0(\b,\th)>0$  such that
for all
$0<\d \le  \d_0$ and  $0<\z\le \z_0$,
  for all integers  $
L_0  \ge  \frac 2 { \a(\b,\th) } \log 1/\z  $, with  $ \l(\b,\th)$ 
 given in \eqv (Palmares)     we have
$$ \inf_{m \in \VV_{\d, \z, L_0}(+1) }
\FF(m) =  \FF^*=  \inf_{m \in \VV_{\d, \z, L_0}  (-1) } \FF(m),
\Eq (P.2)
$$
where  $\FF ^* $   is defined in \eqv (min).   The infimum in the 
first (last) term of \eqv
(P.2) is a minimum, attained at a suitable translate of $\bar m$
($T\bar m$, respectively).}

Lemma  \eqv (5.30)  follows  from the variational result
proven in [\rcite{COP2}] and recalled in  Section 2, once we show
that the set $ \VV_{\d, \z, L_0} (+1)$ ($\VV_{\d, \z, L_0} (-1)$)
contains a suitable translate of $\bar m$ ($T\bar m$,
respectively). Due to the $T$-invariance of the functional $\FF$
it suffices to check the first. This is easily obtained.   Namely, 
from  the exponential decay
properties of $\bar m$, see \eqv (V.1),  $\|\bar
m(r)- m_\b\|_1 \le \z$ for $ r \ge \frac 1 { \a(\b,\th) } \log c/\z  $ and
$\|\bar m(r)- Tm_\b\|_1 \le \z$ for $r \le -\frac 1 { \a(\b,\th) } \log
c/\z  $.
Taking into  account the definition
\eqv (SS.4) we can take
  $L_0 \ge   \frac 2 { \a(\b,\th) } \log c/\z  $  and find  
a translate  of  $\bar
m$ in the set $ \VV_{\d, \z, L_0} (+1)$.


\vskip0.5cm
For  any  interval $I \subset \R $ and $m=(m_1,m_2) \in {\cal T} $, we 
denote by
$m_{I}\equiv m \1_{I }$ the function that coincides with $m$ on
$I$ and vanishes outside $I$. We define
$$
  \FF^0 (m_I) \equiv
  \int_I  \left(f_{\b,\th}(m(r))-f_{\b,\th}(m_\b)\right)\, dr
 +
\frac {1}{4}
\int_{I}\, dr \int_{I}\, dr'
J(r-r') \big[ \tilde m(r)-  \tilde m(r')\big]^2,
\Eq(5.110)
$$
where  $f_{\b,\th}$ is defined in \eqv(2.14000) and
$\tilde m= \frac {m_1+m_2}2$.
For a
 given $\underline m \in  {\cal T}$,
we denote
$$
  \FF (m_I|\underline m_{\del I}) \equiv
\FF^0 (m_I ) +
\frac{1}{2}
\int_{I}\,dr\int_{I^c}\,dr'
J(r-r') \big[ \tilde m(r)-  {\tilde {\underline m }}(r')\big]^2.
\Eq (5.111)
$$
Both   functionals are positive and well defined for all $I\subset
\R$, however they  could be infinite if $I$ is unbounded. Observe
that  when $m_I \equiv m_\b$ (or  $ m_{ I} \equiv Tm_\b$)  then
$\FF^0 (m_I) $ reaches its minimum value $\FF^0(m_\b)=
\FF^0(Tm_\b)=0$  in $I$. The same holds for   $ \FF
(m_I|\underline m_{\del I})  $ when
 $ \underline m_{\del I} \equiv m_\b$
(or
$\underline m_{\del I} \equiv Tm_\b$).  When  the
boundary conditions $ m_{\del I}$ are different from $ m_\b$
(or  $ Tm_\b$) but are suitably
close to them  we will prove that the minimizer
exists and it decays exponentially fast to $ m_\b$  (or  $ Tm_\b$)
with the distance from the boundaries of $I$.
 The value of the functional  at the
minimizer will be, therefore,  close to the null value.
For all $\eta\in \{-1,+1\}$, we denote
$$\MM (\z,\d,\eta ) = \left \{ m=(m_1,m_2)
\in   {\cal T};
  \, \eta^{\d,\z}(\ell)=\eta,  \,\forall \ell \in  \Z \right \}, \Eq (roma.1)
$$
$$\AA (\z,\d,\eta ) = \left \{ m=(m_1,m_2)
\in   {\cal T};
 \, \bar \eta^{\d,\z}(\ell)=\eta,  \,\forall \ell \in  \Z \right \}, 
\Eq (roma.2)$$
where $\eta^{\d,\z}(\cdot)$ was defined in \eqv (7.2a) and
$$
\bar \eta^{\d,\z}(\ell)=\cases{
 +1 &if $\,\forall_{u \in (\ell-1, \ell] }\, \|m^{\d}(u) -m_\b\|_1 \le\z
$;\cr - 1 &if $\,\forall_{u \in (\ell-1 , \ell] }\,
\|m^{\d}(u) -T m_\b\|_1
\le\z$;\cr
0 & otherwise.\cr }
\Eq (7.5a)
$$

\noindent Using $\|m^{\d}(u)-m_\b\|_1\le
\d^{-1} \int_{D^{\d}(u)} ds\|m^{\d^*}(s)-m_\b\|_1$, it is easy to
see that $  \MM (\z,\d,\eta)  \subset \AA (\z,\d,\eta)  $.
 We denote by $\MM_I (\z,\d,\eta) = \{m\1_I \,
\hbox {for }\, m\in \MM(\z,\d,\eta)\}$
and in a similar way  $\AA_I(\z,\d,\eta)$.

\vskip0.5cm
\noindent{\bf \Theorem (7)} {\it   For $(\b,\th)$   that satisfies 
 \eqv (2.17), there
exists
$0<\z_0=\z_0(\b,\th)<1$
 and,    for  $0<\z \le \z_0$, there exists  $\d_0 =\d_0 (\z) >0$,  
such that for
 any $0<\d \le \d_0$,
given a  $\DD_\d$--measurable  interval $I$ and boundary
 conditions $\underline m_{\del I}\in
 \MM_{\del I}  (\z,\d,+1 )$
there exists    an unique   $\psi = (\psi_1,\psi_2)  $ in
 $\MM_I (\z,\d,+1 ) $ such that
$$
\inf_{m_I\in \MM_I (\z,\d,+ )}
\FF (m_I|\underline m_{\del I} )= \FF (\psi |\underline m_{\del I} ).
\Eq(Alrachid.1)
$$
    The
minimizer  $\psi$ is a  continuous function  with
 uniformly bounded first derivative in the interior of $I$,
$\lim_{r \uparrow \partial^+ I} \psi(r)$ and 
$\lim_{r \downarrow \partial^- I} \psi(r)$
exist,
with the further property that
$$
|\psi_1(r)- m_{\b,1}|+ |\psi_2(r)- m_{\b,2}| 
\le  \z   \qquad \forall r \in I
\Eq (7.3a)
$$
$$
|\psi_1(r)- m_{\b,1}|+ |\psi_2(r)- m_{\b,2}| \le 
 \z e^{-\a(\b,\th,\z_0)[2d(r,\del I)]}
\qquad \forall r \in I\quad
\hbox {such that}\quad d(r,\del I)
\ge \frac 12,
\Eq (7.3)
$$
where  $d(r,\del I)$ denotes the distance from  $r$, to  the 
closure of $\del I$,
$[\cdot]$ refers to the integer part,  and $\a(\b,\th,\z_0)$ is 
defined in \eqv (5.P1).
}

\remark An analogous  result, changing
$m_{\b}$ to $T m_{\b}$, holds for $\eta=-1$. 
\vskip0.5cm
\noindent \proof
Since $\MM_I (\z,\d,1)  \subset \AA_I (\z,\d,1)  $, we first
prove that the infimum of $\FF (\cdot|\underline m_{\del I} )$
over $ \AA_I (\z,\d,1)$, a priori smaller than the one
in \eqv(Alrachid.1), is reached at a unique $\psi \in \AA_I
(\z,\d,1)$. Then we prove that $\psi$ can be taken continuous and that
it verifies \eqv(7.3a). This implies that   $\psi \in  \MM_I(\z,\d,1)$, and
therefore \eqv(Alrachid.1) holds.
The proof that the  minimizer of $
\FF (\cdot |\underline m_{\del I})$ over  $\AA_I (\z,\d,1)$ exists is
obtained  dynamically. We study  a system of integral  differential equations
for which  $\FF (\cdot|\underline m_{\del I} )$   is decreasing along 
its solutions:
$$ \eqalign { & { \partial m_1 \over\partial t} =  - m_1 +
\tanh\{\beta
\left (J\star \tilde  m+\th + J \star
 \tilde {\underline m}_{\del I} \right ) \};\cr &
{\partial m_2 \over\partial t} =  - m_2 + \tanh\{\beta
\left (J\star \tilde m-\th + J \star 
\tilde {\underline m}_{\del I}   \right ) \}.}
\Eq (7.6)   $$
Therefore the minimizers of   $
\FF (\cdot |\underline m_{\del I})$  correspond to     
stationary solutions of\eqv (7.6),
i.e:
$$ \eqalign { &   \psi_1 =
\tanh\left\{\beta
\left (J\star \tilde  \psi +\th +   J \star \tilde
 {\underline m}_{\del I} \right ) \right\}; \cr &
  \psi_2 = \tanh\left\{\beta
\left (J\star \tilde \psi -\th + J \star \tilde 
{\underline m}_{\del I}  \right ) \right\}.}
\Eq (7.8)   $$
This method has been already applied  to characterize 
the minimum of  the infinite volume functional  \eqv (AP.1),
see [\rcite{COP2}] and reference therein.
To show \eqv (7.3a) set   $\tilde \psi = \frac 12
(\psi_1 +\psi_2)$ so that, from  \eqv (7.8),
$$\tilde \psi = \frac 12 \tanh\left\{\beta
\left (J\star \tilde \psi  +\th +
J \star \tilde  {\underline m} _{\del I}  \right ) \right\}
+\frac 12 \tanh\left\{\beta
\left (J\star \tilde  \psi  -\th +
 J \star\tilde  {\underline m}_{\del I} \right ) \right\}.
\Eq (7.30) $$
 Since, see \eqv (2.1800),  $g_{\b}(s,\th)<s$ when $s > \tilde m_\b$ and
$g_{\b}(s,\th)>s$ when $0 \le s < \tilde
m_\b$,
it is easy to see  that for
 $ 0< \z \le \tilde  m_{\b}$
 there exists $\d_0(\z)$ such that for $\d \le  \d_0(\z)$,
$| \tilde \psi (r) -\tilde m_{\b}|
\le \frac \z 2 $ for $r \in I$.
 \eqv (7.3a) is then easily derived, once we observe that
$$
\eqalign { & |\psi_1(r)- m_{\b,1}| = |\tanh \b [J \star (\tilde
\psi + \tilde {\underline m}_{\del I}) (r)  +\theta ] - \tanh \b
[\tilde m_\b
 +\theta ] |\cr
&\quad=
   \left   | \int_0^1 ds   \b (1- \tanh ^2 \b [  s J \star 
(\tilde \psi + \tilde {\underline
m}_{\del I}) (r)+ (1-s) \tilde m_\b
   + \theta ]     \left [ J  \star (\tilde \psi + \tilde {\underline
m}_{\del I})
(r) - \tilde m_\b \right ]   \right | . }\Eq (again.4) $$
 Replacing  $ \tilde {\underline
m}_{\del I} $ by $\tilde {\underline m^\d}_{\del I} (r)$, we obtain

$$
\eqalign { & |\psi_1(r)- m_{\b,1}|  \le \b \left [ 1- \tanh^2 \b \{\tilde
m_\b -  \frac \z2 -\d +\th \} \right ]  \left | (J \star (\tilde
\psi + \tilde {\underline m^\d}_{\del I}) (r) + \d -\tilde m_\b (J
\star \1_{I\cup \d I }) (r)   \right |\cr & \le \b \left [ 1-
\tanh^2 \b \{\tilde m_\b -  \frac \z2 -\d +\th \} \right ]   \left
(\frac \z 2 + \d \right ). }     $$
Doing something similar for the other component we obtain
$$ |\psi_1(r)- m_{\b,1}| + |\psi_2(r)- m_{\b,2}|  \le
e^{-\a(\b,\th,\z+2\d)}[\z+2\d], $$
where  we set
$ \a(\b,\th,\z)
= -\log  \frac {\partial
g_\b}{\partial m} (\tilde m_{\b,\th}-\frac {\z} 2,\th)  $,
  $g_\b$  being  defined in \eqv (2.1800).
By the  smoothness of  $g_\b$, since
  \eqv  (2.17000),   there exists $ \z_0= \z_0(\b,\th)$ so that for
$\z\le \z_0(\b,\th)$   and  $\d$
 small enough (depending on $\z$) we have  $ e^{-\a (\z+2\d)}
 [\z +2\d ]\le \z $.
To get    \eqv (7.3) we first show  that  $\tilde \psi$ solution
of \eqv (7.30) has the following
property
$$
|\tilde \psi(r)- \tilde m_\b| \le \frac \z 2 e^{-\a(\b,\th,\z_0)
[2d(r,\del I)]}  \quad \text{if } \quad d(r, I^c) \ge \frac 12,
  \Eq (7.13)
$$
where $[x]$ is the integer  part of $x$.
Since $\tilde m_\b$ is a
solution of \eqv (2.1800), we have:
$$
\left|\tilde \psi_I(r)- \tilde m_\b\right|
\le e^{-\a(\b,\th,\z)} \left| J \star \tilde \psi_I(r)
-\tilde m_\b \right| +e^{-\a(\b,\th,\z)}
J\star|\tilde {\underline m}_{\del I}|(r),
\quad \forall r\in I.
\Eq(alibaba)
$$
 Notice that $( J \star | \tilde  {\underline m}_{\del I}|) (r)=0$  for
$  r \in I,\,d(r, \del I) \ge \frac12$ and,
since $J(r)=\1_{\{|r|\le 1/2\}}$, if $r$  is such that
$d(r,\del I)>N_0/2$
for some $N_0\in \N$, we  have $(J^{\star N_0}\star
|{\underline m}_{\del I}|)(r)=0$.  Therefore, iterating \eqv(alibaba)
$N_0$--times, 
for $r$ such that $(N_0+1)/2\ge d(r,\del I)>N_0/2$, we see that
$$
\left|\tilde \psi_I(r)- \tilde m_\b\right|\le e^{-{N_0}\a(\b,\th,\z)}
\left|\tilde J\star\psi_I(r)-\tilde m_\b\right|\le
e^{-{N_0}\a(\b,\th,\z)}\frac{\z}{2}.
\Eq(exp)
$$
Since $e^{-\a(\b,\th,\z)} <1$ for $\z \le \z_0$,  we
obtain  \eqv (7.13).
Since  $ d(r,\del I) \ge
\frac 12 $
  implies  that $
 (J \star \tilde {\underline m}_{\del I}) (r) =0 $, from \eqv (again.4)
and \eqv (7.13), and doing similarly for the other component, we obtain  that
$$ |\psi_1(r)- m_{\b,1}| + |\psi_2(r)- m_{\b,2}| \le e^{-\a(\b,\th,\z)} \z
e^{- \a(\b,\th,\z_0)[2 d(r,\del I)]} \le \z e^{-\a(\b,\th,\z)
[2d(r,\del I)]}.   \Eq (7.14ab)  $$
\qed
\smallskip
\noindent $\bullet$ {\bf II:  Surface tension.}
\smallskip

\noindent {\bf \Lemma (T1C)} {\it  Given  $(\b,\th)$ that satisfies \eqv(2.17),
there exist   $\g_0=\g_0(\b,\th)>0$, $d_0=d_0(\b,\th)>0$,
$1>\z_0=\z_0(\b,\th) >0$
 such that  for all $0<\g\le \g_0$, all $\d^*>0$ with $\g/\d^* \le d_0$, and
 all positive integer $p$ satisfying
$$
(p+2) \d^* \log \frac 1 \g \le \frac 1{64}
\Eq(5.12PP)
$$
there exists $\O_{RE}= \O_{RE}(\g,\d^*, p)$ with $\P[\O_{RE}]\ge 1-\g^2$ such
that for any $\o \in \O_{RE}$,
any  $1>\d  >\d^*>0$, and any
$\z_0 >\z_1> 8 \g/\d^*$, if  $L_0= 
\frac  d {\a(\b,\th)}  \log(\frac{\d^*}{8\g})$ for some
 $d\ge   2  $
and $\a(\b,\th)$ defined in \eqv(Palmares), we then
we have, uniformly with respect to  the choice  of $[\ell_1,
\ell_1+L_0-1]$ and $[\ell_2-L_0+1,\ell_2]$ inside
$[-\g^{-p},\g^{-p}]$:
$$
\frac   { Z^{m_+,0}_{[  \ell_2-L_0+1,\ell_2]}
 \big ( \eta^{\d,\z_1} (\ell_2)=+1
 \big ) }{
Z^{m_+,0}_{[  \ell_2-L_0+1,\ell_2]} (\eta^{\d,\z_1} (\ell_2)=\bar
\eta_2) }
 \frac
{ Z^{0,m_+}_{[\ell_1, \ell_1+L_0-1]} \big (
 \eta^{\d,\z_1}(\ell_1)=+1
  \big )}
{ Z^{0,m_+}_{[\ell_1, \ell_1+L_0-1]}
(\eta^{\d,\z_1}(\ell_1)=\bar \eta_1) }
\le e^{ \frac \b
\g \big ( \FF^*  +   f_{11}  \big)  \big [
\frac 12 (|\bar \eta_1-1|+|\bar \eta_2 -1|)
\big] },
\Eq(5.1500)
$$ where $\FF^*$ is defined in  \eqv (min) and
$$
  f_{11}  \equiv 10  (1+\th)
(\d^*\vee\sqrt{\frac{\g}{\d^*}}) d \log \frac{\d^*}{8\g}.
\Eq(5.1501)
$$
}
\noindent \proof        We start estimating
 $  \frac { Z^{0,m_+}_{[\ell_1, L_0-1]}
(\eta^{\d,\z_1}(\ell_1)=\bar \eta_1) }
{ Z^{0,m_+}_{[\ell_1, L_0-1]} \big (
 \eta^{\d,\z_1}(\ell_1)=+1
  \big )} $ from below.
  When $ \bar \eta_1=+1$, the previous quantity is
 equal to  1 and there is nothing to prove. We then suppose that  $
\bar \eta_1=-1$ and to simplify notation we set
$\ell_1=0$. We perform a block spin
transformation as in Section 4 and   use Lemma \eqv(MM).  
 For  the random terms we use
the rough estimate,
Lemma \eqv (60002), obtaining  for $\o \in  \O_{RE}$,
$$
\eqalign{
Z^{0,m_+}_{[0, L_0-1]}
(\eta^{\d,\z_1}(0)=\bar \eta_1)
&\ge
e^{-\frac \b\g    L_0 \left ( \d^*+ c \frac{\g}{\d^*}\log \frac{\g}{\d^*} +
4 \th  \sqrt{\sfrac {\g}{\d^*}}\right )  }
 e^{-\frac \b\g \big[\tilde  \FF
(\widehat m^{\d^*}_{[0,
L_0-1]}|\underline m^{\d^*}_{\partial[0,
L_0-1]})\big]}\times \cr
&\quad\times
e^{+\frac{\b}{\g} \big[\frac{\d^*}{2} \sum_{y \in \CC_{\d^*}(\del
[0,L_0-1])} \big[\tilde {\underline m}^{\d^*}(y)\big]^2
\sum_{x \in \CC_{\d^*}([0,L_0-1])}
J_{\d^*}(x-y)\big]}\cr
}\Eq(5.12005)
$$
  where $\underline  m^{\d^*}_{\partial[0,
L_0-1]}$  is  the profile associated to the chosen boundary conditions, i.e.,
$\underline m^{\d^*}_{\partial^-[0, L_0-1]}=0 $,
$\underline  m^{\d^*}_{\partial^+[0,
L_0-1]}$ $=m^{\d^*}_\b  $ and   $ \widehat  m^{\d^*}_{ [0,
L_0-1]} \in \MM^{ -} \equiv \MM_{\d^*}(
[0,L_0-1]) \cap \{  \eta^{\d,\z_1} (0)=-1\} $ will
 be suitable chosen in the following.
In a similar way,   we
estimate the denominator by
$$
\eqalign{ { Z^{0,m_+}_{[0, L_0-1]}
\big(\eta^{\d,\z_1}(0)= 1 \big )} &\le e^{\frac \b \g   L_0
\left ( \d^*  + 4 \th  \sqrt{\frac \g{\d^*}}   \right )  } \times
e^{ \frac {1}{\g} [    L_0 \frac \g{\d^*} \log \frac {\d^*}\g]}
\times e^{ \frac {\b}{\g} [    L_0 c\frac \g{\d^*} \log \frac
{\d^*}\g]}\times \cr &\quad \times  e^{ -\frac \b \g \big[\inf_{\{
m^{\d^*} \in \MM^+\}} \tilde \FF(m^{\d^*}| \underline
m^{\d^*}_{\partial[0, L_0-1]}) \big ] }\times\cr
&\quad\times e^{+\frac{\b}{\g} \big[\frac{\d^*}{2} \sum_{y \in
\CC_{\d^*}(\del [0,L_0-1])} \big[\tilde {\underline
m}^{\d^*}(y)\big]^2 \sum_{x \in \CC_{\d^*}([0,L_0-1])}
J_{\d^*}(x-y)\big]}.\cr }\Eq(5.12009)
$$
The term
$e^{ \frac {1}{\g} [   L_0 \frac \g{\d^*} \log \frac {\d^*}\g]}$
comes  from  counting  the number of configurations of
$m^{\d^*} \in \MM_{\d^*}(
[0, L_0-1])$.
The infimum  in \eqv (5.12009) is over the set
 $
 \MM^+ \equiv
   \MM_{\d^*}(
[0,L_0-1]) \cap
\big \{ \eta^{\d,\z_1}(0)   =
1 \big \}
 $
 and it  is attained on the configuration  
 $ \{ m^{\d^*}(x)=m_\b^{\d^*}, \forall
x \in \CC_{\d^*}([0,L_0-1])\} $,  since  the boundary
conditions are at  one side zero and
at  the other side already  equal to $m_\b^{\d^*}$.
We need only that
$\z_1> 8\g/\d^*$ to be sure that
$\|m^{\d^*}_\b-m_\b\|_1 \le  \z_1 $ entails that the configuration constantly
equal  to $m^{\d^*}_\b$ belongs to $ \MM^+$.
 Taking in account \eqv (5.12005), \eqv (5.12009)   we
obtain
$$
\eqalign {
&  \frac { Z^{0,m_+}_{[0, L_0-1]}
(\eta^{\d,\z_1}(0)=\bar \eta_1) }
{ Z^{0,m_+}_{[0, L_0-1]} \big (
 \eta^{\d,\z_1}(0)=+1
  \big )}  \ge  e^{-\frac \b \g  \frac 12 |\bar \eta_1-1|  \left [  2  L_0
\left (\d^*+ (1+c) \frac \g{\d^*} \log \frac {\d^*}\g +
 4 \th  \sqrt{\frac \g{\d^*}} \right )
\right ]  }  \cr
 &\quad \times  e^{-\frac \b \g  \frac 12 |\bar \eta_1-1| 
\left [\tilde  \FF ( \widehat
m^{\d^*}_{[0,
L_0-1]}|\underline m^{\d^*}_{\partial[0,
L_0-1]})  -\tilde \FF(m_\b^{\d^*}| \underline  m^{\d^*}_{\partial[0,
L_0-1]}) \right ] }. } \Eq (P.100)
$$
The
exponent in the last line of \eqv (P.100)   can be written as
$$ \eqalign { & \left [  \tilde \FF ( \widehat m^{\d^*}_{[0,
L_0-1]}|\underline m^{\d^*}_{\partial [0,
L_0-1]})  -\tilde \FF (m_\b^{\d^*}| \underline  m^{\d^*}_{\partial[0,
 L_0-1]})  \right ]
=
\FF^0 ( \widehat m^{\d^*}_{[0,
 L_0-1]})  +[f(m_\b)-f(m^{\d^*}_\b)] [L_0-1] + \cr
&\qquad +
 \frac{\d^*}{2}
\sum_{x \in \CC_{\d^*}([0,
 L_0-1])}\sum_{y \in \CC_{\d^*}(  \del [0,
 L_0-1] )}
J_{\d^*}(x-y) \big[\tilde { \widehat m} ^{\d^*}(x)- \tilde { \underline
m}^{\d^*}_{\partial[0,
 L_0-1]}(y)\big]^2 \cr &
\qquad -
\frac{\d^*}{2}
\sum_{x \in \CC_{\d^*}([0,
 L_0-1])}\sum_{y \in \CC_{\d^*}(   \del [0,
 L_0-1] )}
J_{\d^*}(x-y) \big[ \tilde {  m }_\b^{\d^*}(x)-  \tilde {\underline
m}^{\d^*}_{\partial[0, L_0-1]} (y)\big]^2 \cr
} \Eq (PM.1)
$$
where    $\FF^0$ is  the functional defined  in \eqv (5.110).
We  take  $\z= 8 \frac \g {\d^*}$ in Lemma \eqv(5.30), assuming
that $\g/\d^*$ is smaller that the  $\z_0$ there,   $L_0 = \frac d {\a(\b,\th)}
\log \frac{\d^*}{8\g}$ with $d \ge   2  $, and
$\a(\b,\th) $  defined in \eqv (Palmares). Then,
Lemma \eqv (5.30) says that  a suitable translate of $\bar m$
belongs to $\VV_{\d,\z,L_0}$, see \eqv(SS.4), provided
   and $0<\d<\d_0$.
 By an abuse of
notation we always denote such translate by $\bar m$. 
Since $\MM^- \subset \VV_{\d,\z,L_0}$,
 we  can  choose $\widehat m^{\d^*} \in \MM^- $  such that  $\|\widehat
m^{\d^*}(r)  -\bar m (r) \|_1\le 8\g/\d^* $  for all $ r\in [0,
L_0-1]$, where $\bar m$ is the previous chosen minimizer.
An easy computation gives
$$
|f(m_\b)-f(m^{\d^*}_\b)]| [L_0-1] + | \FF^0 ( \widehat m^{\d^*} _{[0,
 L_0-1]})- \FF^0 ( \bar  m_{[0,
 L_0-1]})| \le  8L_0(1+\th)\sqrt{\frac{\g}{\d^*}}.
\Eq (PP.1)
$$
Since  $\widehat m^{\d^*} \in \VV_{\d,\z,L_0}$  and 
$\z= 8 \frac \g {\d^*}$,  the difference of the last two sums
in \eqv(PM.1) is bounded from above by
$ 64 \frac{\g}{\d^*}<\sqrt{\frac{\g}{\d^*}}$ and 
$ \frac \g {\d^*}$ is small enough.
Since
$  \FF^0 (\bar m_I) \le  \FF^* $ we obtain
$$
\frac { Z^{0,m_+}_{[\ell_1, \ell_1+L_0-1]}
(\eta^{\d,\z_1}(\ell_1)=-1) }
{ Z^{0,m_+}_{[\ell_1, \ell_1+L_0-1]} \big (
\eta^{\d,\z_1}(\ell_1)=1
  \big )}
\ge  e^{-\frac \b \g   
[ \FF^* + 10(1+\th)(\d^*\vee\sqrt{\frac{\g}{\d^*}})L_0)]}.
\Eq (P.01)
$$

Repeating similar arguments for  the term with $\bar \eta_2$ we 
end the proof. \qed

\smallskip
\noindent $\bullet$ {\bf III: Shrinking of the typical profiles.}
\smallskip

\noindent {\bf \Theorem (6111)}
 {  \it Given $(\b,\th)$ that satisfies \eqv(2.17),
there exist
 $0<\g_0=\g_0(\b,\th)<1$, $ 0<d_0=d_0(\b,\th)<1$ and $0<\z_0=\z_0(\b,\th)<1$,
such that for all $0<\g\le \g_0$,   $\g/\d^*\le d_0$,
for all $p\in\N$ verifying the condition
$$ (p+2)\d^*\log \frac 1 \g \le \frac
1{64},
\Eq(6.3800)
$$
there exists $\O_{RE} = \O_{RE} (\g,\d^*,p) $
with $\P\big[\O_{RE} \big] \ge 1-\g^2$
such that for any $\o \in \O_{RE}$,  $\bar \eta\in \{-1,+1\}$,
  $\ell_0 \in \N $,
 $\d ,\z_4,\z_5$ with $1>\d > \d^*>0$, and any
$\z_0\ge \z_4>\z_1 >\z_5 \ge 8\g/\d^*$, we
have
$$
\eqalign{
& 
\sup_{\D_L \subset [-\g^{-p},\g^{-p}]}
\mu_{\b,\th,\g}\Big(
R_{1,4}(\bar \eta, [\ell_1,\ell_2])
\cap
(  R_{1,4,5}(\bar \eta, [\ell_1,\ell_2])(\ell_0))^c
\Big)
\cr   & \le
\frac{2}{\g^p} e^{-\frac\b \g \Big\{ \frac{ \k (\b,\th)}{2}
\d \z_5^3 -2\z_4 e^{-\a(\b,\th,\z_0)[2 \ell_0]} -
12(1+\th)(4\ell_0+10)[
 \d^* \vee \sqrt {\frac{\g}{\d^*}}] \Big\}},\cr
}\Eq(6.1002)
$$
where
 $R_{1,4,5}(\bar \eta,[\ell_1,\ell_2])(\ell_0)$ is defined
in \eqv(5.1912), and  $R_{1,4}(\bar\eta,[\ell_1,\ell_2])$
in \eqv(5.matteo), $ \k (\b,\th)>0$ satisfies   \eqv (2.19),
$\a(\b,\th,\z_0)$ is defined in \eqv (5.P1)
 and   $\D_L=[\ell_1,\ell_2]$ is an interval
 of  length $L\ge 4\ell_0+10$.
Moreover
$$
\sup_{\D_L \subset [-\g^{-p},\g^{-p}]}
\frac{ Z_{[\ell_1,\ell_2]}^{0,0}(
  R_{1,4}(\bar \eta, [\ell_1,\ell_2])
\cap
(R_{1,4,5}(\bar \eta, [\ell_1,\ell_2])(\ell_0))^c
)}
{ Z_{[\ell_1,\ell_2]}^{0,0}(
  R_{1,4,5}(\bar \eta, [\ell_1,\ell_2]))}
\Eq(laila.10)
$$
satisfies the same estimates as \eqv(6.1002).}

\noindent \remark
 Note the crucial fact that the last  term in
the exponent on  the right hand side of
\eqv(6.1002)  is  proportional to  $4\ell_0+10$ and not to  $L$.

\noindent The following corollary is an immediate consequence 
of Theorem \eqv(6111).
Its proof  consists essentially in choosing an appropriate $\ell_0$
in \eqv(6.1002), see \eqv(6.0491), 
and  taking in account   that, under \eqv(6.0049) and  $\d^*>\g$, we have
$\d^*\vee \sqrt{\frac{\g}{\d^*}}=
\sqrt{\frac{\g}{\d^*}}$.

\noindent {\bf \Corollary (6002)}
{\it
 Under the same hypothesis of Theorem \eqv(6111) with the further requirements
$$
\d \z_5^3 >  \frac{512(1+\th)}{\k (\b,\th)\a(\b,\th,\z_0)}
\sqrt {\frac{\g}{\d^*}}
\log \frac{\d^*}{\g},
\Eq(5.099)
$$

$$
\frac{(\d^*)^2}{\g} \le \frac 1{6e^3\b}
\Eq(6.0049)
$$
where
$\k (\b,\th)>0$ satisfies  \eqv (2.19) and
 $\a(\b,\th,\z_0)$ is defined in \eqv (5.P1).
$$
\ell_0=  \frac {d} {2\a(\b,\th,\z_0)}  \log
\frac{\d^*}{\g} \qquad d>1,
\Eq (6.0491)
$$
then for $\o\in \O_{RE} $ and $\bar \eta\in \{-1,+1\}$, we have
$$
\sup_{\D_L \subset [-\g^{-p},\g^{-p}]}
\mu_{\b,\th,\g}\Big(
\big(
R_{1,4}(\bar \eta, [\ell_1,\ell_2])
\cap
(  R_{1,4,5}(\bar \eta, [\ell_1,\ell_2])(\ell_0))^c
\Big)
\le  e^{-\frac \b \g
{\frac{\k (\b,\th)}{4} \d_5\z_5^3 }}
\Eq(6.000030)
$$
where $\D_L$  is an interval of length $L \ge 4 \ell_0+10$.
Moreover \eqv (laila.10)
satisfies the same estimates as
\eqv(6.000030).
}
\vskip.5cm
{\bf Proof of Theorem \eqv (6111)}
Given an interval $\D_L\equiv [\ell_1,\ell_2]$, with
$\ell_2-\ell_1=L >4\ell_0+10$ for some
$\ell_0$ to be chosen later,
$\ell \in [\ell_1+2\ell_0,\ell_2  -2\ell_0]$, $ \bar \eta=\pm 1$,
we denote
$$
\EE_{\bar \eta}(\ell) \equiv
\Big\{ m^{\d^*}(x), x \in \CC_{\d^*}(\D_L):
\eta^{\d,\z_5}(\ell)=0, \eta^{\d,\z_4}(\ell')= \bar 
\eta\quad  \forall \ell' \in
[\ell-2\ell_0-5,\ell+2\ell_0+5]
\Big\}.
\Eq(5.101)
$$
Since
$$
R_{1,4}(\bar \eta, [\ell_1,\ell_2])
\cap
(R_{1,4,5}(\bar \eta, [\ell_1,\ell_2]))^c
\subset
\cup_{\ell=\ell_1+2\ell_0}^{\ell_2-2\ell_0}
 \EE_{\bar \eta}(\ell)
\Eq(5.103)
$$
 it is enough to estimate
$\mu_{\b,\th,\g}\big(\EE_{\bar \eta}(\ell)\big)$ and we   
assume   $\bar \eta=+1$.
After an easy  computation, calling $I=[\ell-2\ell_0-5,\ell+2\ell_0+5]$,
for   $\o\in \O_{RE}$,  introduced in Lemma \eqv (60002),
for all $\ell \in [-\g^{-p},\g^{-p}]$, we
obtain
$$
\eqalign{
&\mu_{\b,\th,\g}\big(
\EE_{1}(\ell)
\big)\le \frac 1{Z_{\b,\th,\g,\L}}\cr
&\quad
\sum_{\s_{\L\setminus \g^{-1}I}}
e^{-\b H(\s_{\L\setminus  \g^{-1}I})}
\1_{\{\eta^{\d,\z_4}(\ell-2\ell_0 -5)=1 \}}(\s_{\g^{-1}\partial I})
\1_{\{\eta^{\d,\z_4}(\ell+2\ell_0 +5)=1\}}(\s_{\g^{-1}\partial I})
Z_{\b,\th,\g,\g^{-1}I}^{\s_{\g^{-1}\partial I}}
\cr &
\quad\quad\quad\quad\quad\quad \times\frac
{e^{-\sfrac \b\g\big\{
\inf_{\EE_{1}(\ell)}
\widetilde \FF\big (m^{\d^*}_I|m^{\d^*}_{\partial I}(\s)\big)
-8(1+\th) (4\ell_0+10) [\d^* \vee  \sqrt{\sfrac{\g}{\d^*}}]
\big\}}}
{e^{-\sfrac \b\g\big \{
\widetilde \FF\big (\overline m^{\d^*}_I|m^{\d^*}_{\partial I }(\s)
\big)\big\}}},\cr
}\Eq(5.106)
$$
where $ \widetilde \FF $   is given in \eqv(5.Rio2)
and  $\overline m^{\d^*}_I$ is a fixed profile.
This inequality is obtained as follows: writing $\mu_{\b,\th,\g}(
\EE_{1}(\ell))$ as a sum of the expression in \equ(2.300)
over the configurations in $\sigma_{\Lambda} \in \EE_{1}(\ell)$ we 
multiply and divide by
$Z_{\b,\th,\g,I}^{\s_{\L\setminus \g^{-1}I}}$, inside the  sum over
 $\s_{\g^{-1}I}$, perform
a block spin transformation in the volume $\g^{-1}I$
and  roughly  estimate the magnetic field applying Lemma \eqv (60002).
This {last two steps  are}  done  in the numerator and the denominator and 
they produce an  error term
$8(1+\th) (4\ell_0+10) [\d^* \vee  \sqrt{\sfrac{\g}{\d^*}}]$.
We  get an upper bound  restricting  in the denominator the sum over all
profiles to the single one  $\overline m^{\d^*}_I$.
Notice the important fact that the  term
$$
\frac{\d^*}{2}
\sum_{y \in \CC_{\d^*}(\del I)}
\big[\tilde m^{\d^*}(y,\s)\big]^2 \sum_{x \in \CC_{\d^*}(I)} J_{\d^*}(x-y)
\Eq(5.107)
$$
in \eqv(5.Rio3) cancels out in the formula \eqv(5.106), since it is present
both  in the numerator and in the denominator.
We can subtract from the two $\widetilde \FF$ in \eqv(5.106) the
quantity $  f(m^\b) |I| $
obtaining $\FF\big (\cdot |m^{\d^*}_{\partial I}(\s) \big) $
instead of $ \widetilde \FF\big (\cdot |m^{\d^*}_{\partial I}(\s)
\big)$. Therefore to prove  Theorem \eqv(6111), it remains to 
prove that  we can choose $\overline
m^{\d^*}_I$    in such a way that
$$
\inf_{m^{\d^*}_I \in \EE_{1}(\ell) }
\FF  \big (m^{\d^*}_I|\,m^{\d^*}_{\partial I }\big) \ge
\frac {\k (\b,\th)}{2} \d \z_5^3  -
2\z_4 e^{-\a(\b,\th,\z_0) 2\ell_0} -4(4\ell_0 +10)(1+\th) 
\sqrt{\frac{\g}{\d^*}}+
\FF \big (\overline m^{\d^*}_I|\,m^{\d^*}_{\partial I }\big)
\Eq(5.1091)
$$
uniformly with respect to $m^{\d^*}_{\del I} \in
R_{1,4}(+1, [\ell_1,\ell_2])$.
In fact the  terms in  the second line of \eqv(5.106) will
be bounded by $Z_{\b,\th,\g,\L}$ uniformly in $\L$ and we get \eqv (6.1002).
It is rather delicate  to prove \eqv(5.1091).

Using \eqv(5.110) and \eqv(5.111), and splitting
$I= I^- \cup (\ell-1,\ell] \cup I^+$ where 
 $I^-\equiv (\ell -2\ell_0-5 ,\ell -1]$ and $I^+\equiv (\ell,\ell
+2\ell_0+5]$, we get that
for all $ m^{\d^*}_I \in \EE_{1}(\ell)$
$$
\eqalign { \FF (m^{\d^*}_I \,| m^{\d^*}_{\partial I })& \ge
\inf_{m^{\d^*}_{I^-}\in \MM_{I^-}(\z_4,\d,+1)} \FF
(m^{\d^*}_{I^-} |\, m^{\d^*}_{\del^- I},  m^{\d^*}_{
(\ell-1,\ell]} ) +\FF ^0(m^{\d^*}_{ (\ell-1,\ell]}) \cr
&\quad\quad +\inf_{m^{\d^*}_{I^+} \in \MM_{I^+}(\z_4,\d,+1)} \FF
(m^{\d^*}_{I^+} |\, m^{\d^*}_{ (\ell-1,\ell]}, m^{\d^*} _{\del^+
I}),  }\Eq(5.116)
$$
where   $m^{\d^*}_{ (\ell-1,\ell]}\equiv  \{ m^{\d^*}(x),
x\in \CC_{\d^*}((\ell-1,\ell])\} $.
Since $m^{\d^*}_{\del^{\pm} I}$   belongs to  $\MM_{\del^{\pm}I}
(\z_4,\d,+1)$ and $ m^{\d^*}_{ (\ell-1,\ell]} $  to
$\MM_{(\ell-1,\ell]} (\z_4,\d,+1 )$, using  Theorem
\eqv(7), there exist unique  minimizers
$\psi^1_{I^+}\in \MM_{I^+}(\z_4,\d,+1)$ and
  $\psi^2_{I^-}  \in \MM_{I^-}(\z_4,\d,+1)$ such that 
$$
\FF (m^{\d^*}_I|\, m^{\d^*}_{\del I})\ge
\FF (\psi^1_{I^-}|\, m^{\d^*}_{\del^- I}(\s),  m^{\d^*}_{ (\ell-1,\ell]} )
+
\FF^0(m^{\d^*}_{ (\ell-1,\ell]}) +
\FF (\psi^2_{I^+}|\, m^{\d^*}_{\del^+ I}(\s),  m^{\d^*}_{ (\ell-1,\ell]} ),
\Eq(5.122)
$$
for any fixed boundary condition and any $ {m^{\d^*}_I} \in \EE_{1}(\ell) $.
By \eqv (2.19)
$$ \FF ^0(m^{\d^*}_{ (\ell-1,\ell]}) \ge  
 \frac{\k (\b,\th)}{2} \z_5^3 \d_5. \Eq (luminy1) $$
 Denote by  $I^-_1 = (\ell -2\ell_0-5 ,\ell - \ell_0-3]$, 
 $I^-_1 \subset I^-$.
By   the positivity  property  of the functional, see
\eqv (5.111),
$$ \FF (\psi^1_{I^-}|\, m^{\d^*}_{\del^- I}(\s),  
m^{\d^*}_{ (\ell-1,\ell]} )  \ge
\FF (\psi^1_{I^-_1}|\,  m^{\d^*}_{\del^- I}(\s),   
 \psi^1_{ (\ell - \ell_0-3, \ell -
\ell_0-2] } ).
$$
Applying  \eqv (7.3)  of   Theorem \eqv (7)    we  have that
 $$\FF (\psi^1_{_{I^-_1}}|\,
m^{\d^*}_{\del^- I}(\s),    \psi^1_{ (\ell - \ell_0-3, \ell -
\ell_0-2] } ) \ge  \FF (\psi^1_{_{I^-_1}}|\,  m^{\d^*}_{\del^- I}(\s),
 m_\b \1_{ (\ell-\ell_0-3, \ell - \ell_0-2] } ) - 
 \z_4 e^{-\a(\b,\th,\z_0)[2\ell_0]}.   $$
Doing the same computations for  
$\FF (\psi^2_{I^+}|\, m^{\d^*}_{\del^+ I}(\s),  m^{\d^*}_{ (\ell-1,\ell]} )$ 
and setting $I^+_2= (\ell +\ell_0+3 ,\ell + 2 \ell_0+5 ] $,  we obtain
$$  \eqalign { &
\FF (\psi^1_{I^-}|\, m^{\d^*}_{\del^- I}(\s),  m^{\d^*}_{ (\ell-1,\ell]} )
+
\FF (\psi^2_{I^+}|\, m^{\d^*}_{\del^+ I}(\s), 
 m^{\d^*}_{ (\ell-1,\ell]} ) \cr & \ge
\FF (\psi^1_{_{I^-_1}}|\,  m^{\d^*}_{\del^- I}(\s),
 m_\b \1_{ (\ell-\ell_0-3, \ell - \ell_0-1] } )
   + \FF (\psi^2_{{I^+}_2}|\,  m^{\d^*}_{\del^+ I}(\s),
 m_\b \1_{ ( \ell + \ell_0+1, \ell+\ell_0+3]   } )  
-  2\z_4 e^{-\a(\b,\th,\z_0)[2\ell_0]} \cr & =
 \FF (\psi^3_{I}|\, m^{\d^*}_{\del^- I}(\s),  m^{\d^*}_{\del^+ I}(\s) )     
-  2\z_4 e^{-\a(\b,\th,\z_0)[2\ell_0]} }
 \Eq
(luminy2)
$$
where we  set   $ \psi^3_I = \psi^1_{I^-_1} + m_\b \1_{
(\ell - \ell_0-3, \ell +\ell_0+3]}  +
\psi^2_{I^+_2} $.  In the last equality in \eqv (luminy2) we use 
that  $ \FF^0(m_\b)=0$.

By Theorem \eqv (7), there exists an unique 
$\psi^*_I \in  \MM_{I} (\z_4,\d,+1 )$ such that
$$
 \inf_{\psi_I \in  \MM_{I} (\z_4,\d,+1 ) } 
\FF (\psi_I|m^{\d^*}_{\partial I })\equiv
\FF (\psi^*_I|m^{\d^*}_{\partial I }).
\Eq(5.132)
$$
Therefore, since $\psi^3_{I} \in \MM_{I}
(\z_4,\d,+1)$, we have
$$ \FF (\psi^3_{I}|\, m^{\d^*}_{\del^- I}(\s),  m^{\d^*}_{\del^+ I}(\s) )   
\ge \FF (\psi^*_{I}|\, m^{\d^*}_{\del^-
I}(\s),  m^{\d^*}_{\del^+ I}(\s) ). \Eq (fine)  $$

Then, from \eqv (5.122), \eqv (luminy1), \eqv (luminy2), \eqv (fine)  
we obtain
$$\inf_{m^{\d^*}_I \in \EE_{1}(\ell) }
\FF  \big (m^{\d^*}_I|\,m^{\d^*}_{\partial I }\big) \ge  
 \FF (\psi^*_{I}|\, m^{\d^*}_{\del^-
I}(\s),  m^{\d^*}_{\del^+ I}(\s) ) +  
\frac{\k (\b,\th)}{2} \z_5^3 \d_5  -  2\z_4 e^{-\a(\b,\th,\z_0)[2\ell_0]}. 
\Eq(boredom)
$$

Choosing  for $\overline m^{\d^*}_I$ a  $\DD^{\d^*}$--measurable
approximation of  $\psi^*_I$  with values in $\MM_{\d^*}(I)$,
see \eqv(2.14), we get
$$
\FF (\psi^*_I|m^{\d^*}_{\partial I }) \ge
\FF^{\d^*}\big (\overline
m^{\d^*}_I|\,m^{\d^*}_{\partial I }\big) -
4(4\ell_0+10)(1+\th) \left(\d^*\vee \sqrt{\frac{\g}{\d^*}}\right).
\Eq(5.134)
$$
Collecting \eqv(boredom) and \eqv(5.134) we get \eqv(5.1091).
\eop

\bigskip
\chap{   7 Appendix: The cluster expansion} 7
\numsec= 7
\numfor= 1
\numtheo=1
\bigskip

In this  section we   prove Lemma \eqv (62) of Section 4. We  will
write   $  V (m^{\d^*}_I,h) $, 
 defined in \eqv(3.140),
as an absolute  convergent series   and then estimate its 
Lipschitz norm. 

To state the result   we need some preliminary definitions.
 Let $I\subset \R$ be  a bounded,  $\DD_{\d^*}-$ measurable  interval,  
$\AA (I)$ the set of blocks  $ A(x)
$, $ x \in \CC_{\d^*}(I) $.  We denote by $\l= (A,A')$ a pair of
different blocks belonging to $ \AA (I) $ and by $\bar \l =A \cup
A'$ its support. We define a {\it graph} $g$ in $\AA (I) $ as
any  collection of pairs of different blocks $g=
\{\l_1,\l_2,..,\l_m\}$, with  $ 0\le m\le \frac {|\AA(I)| }
2(|\AA(I)|-1) $, such that    $\l_s\neq \l_t$ for all $s\neq t$. A
graph  $g$ will be said to be  connected if, for any pair $B$
and $C$ of disjoint  subsets of $\AA (I) $ such that $B\cup
C= \cup_{s=1}^{m} \bar \l_s $, there is a $\l_s \in g$ such that
$\bar\l_s \cap B \neq \emptyset $  and $\bar\l_s \cap C \neq
\emptyset $.
Given
a graph $g = \{\l_1,\l_2,..,\l_m\}$, $ \l_1,\l_2,..,\l_m $ are
called {\it links } of the graph $g$ and the blocks  $A(x) $
belonging to $\cup_{s=1}^{m} \bar \l_s $ are called vertices of

$g$. We denote $\GG_{\AA (I)}$ the set of all connected graphs of
$\AA(I)$. A connected {\it tree graph} $\t$ (or simply a tree
graph) is a connected graph with $m$ vertices and $m-1$ links. We
denote by $\TT_{\AA (I)} $ the set of all tree graphs in $\AA(I)$.
 Given a tree graph $\t$  the incidence number of the vertex $A(x)  $,
 denoted by $d_{A(x)}$, is the number of links
$\l_s$ in $\t$ such that $ A(x)\cap \bar \l_s
\neq  \emptyset$.
 In the following we denote by a {\it polymer }
  $R$  a subset  of  blocks  of  $ \AA (I) $,  by
$ \CC_{\d^*}(  R) = \{  x\in \CC_{\d^*}(I) \, \hbox {such that}
 \, A(x) \in   R \} $ and
$ m^{\d^*}_R = \{
m^{\d^*}(x);  x\in \CC_{\d^*}(R)  \} $.
 We have the following Theorem.

\noindent{\bf \Theorem (62P)}
{\it For all $\b>0$,  $h\in \O$, for any  bounded    
interval $I \subset \R$,  for  $\d^*>0$,  $\frac{ (\d^*)^2} \g <
\frac 1 {6 e^3 \b}  $,   $  V (m^{\d^*}_I,h) $ can be written as
an absolutely convergent series:
$$  V (m^{\d^*}_I,h)   =
\frac 1 \b \sum_{n=1}^{\infty} \frac 1 {n!}  \sum_{R_1,R_2,\dots,R_n,
 |R_{\ell} | \ge 2 }  \Phi^T(R_1,R_2,\dots,R_n)
 \prod_{{\ell}=1}^n
 \r(R_{\ell}),   \Eq (C.3) $$ where $
\Phi^T(R_1,R_2,\dots,R_n)$ are the Ursell coefficients, see \eqv
(C.1d),  and   $\rho(R_{\ell})$  is   given by
$$
  \r(R_{\ell}) =\r(R_{\ell},h) = \E_{ m^{\d^*}_{( R_{\ell})}} \left 
 [ \sum_{g\in \GG_{R_{\ell}}}
\prod_{(x,y) \in g, x \neq y } \left [ e^{ \b U (\s_{A(x)},\s_{A(y)})}
-1 \right ]
 \right ].     \Eq (C.2)
 $$
 $ \GG_R $ is the set of the  connected graphs in {\bf} $R$ and $x$ 
is a short notation for  $ A(x)$.   (So $
(x,y) \in g$ is a short notation for $(A(x),A(y))\in g $.) 
Moreover
 $$ \left |  V (m^{\d^*}_I,h)  \right| \le |\CC_{\d^*}(I)|
\frac 1 \b \frac  {S }   {1-  S  },     \Eq (H.2)
 $$
where $$  S =\sup_{h} \sup_{x\in \CC_{\d^*}(I)} 
\sum _{R:x \in R}e^{|R|}\rho(R) < 6e^3 \b  \frac {(
\d^*)^2}  {\g} <1 \Eq (paris1)$$
   and
 $$
 \sup_{I\subset \Z }\sup_{i\in I} \|\del_i V_I\|_\infty
\le \frac  { S }   {1-  S   } \frac 1 \b.    \Eq (H.22P)
 $$ }
\bigskip

\proof  The proof is obtained via a standard tool of Statistical
Mechanics, the so called cluster expansion, see [\rcite{Br}]  and
bibliography therein.  This  expansion is done in three
steps:

\noindent  (1)  express the $ \log V$ as a formal series,

\noindent  (2) establish sufficient conditions for the series to
 converge  absolutely,

\noindent  (3) control that under the hypothesis  of Theorem 
\eqv (62P) these conditions  are indeed satisfied.

\noindent We start with the following  identity
$$
\eqalign { &  \E_{ m^{\d^*}_I}[\prod_{x \neq y} e^{ \b U
(\s_{A(x)},\s_{A(y)})}] = \E_{ m^{\d^*}_I} \prod_{x \neq y} [ e^{ \b U
(\s_{A(x)},\s_{A(y)})}-1+1]
  \cr & =1+  \sum_{n=1}^{\infty}\frac 1 {n!} 
\sum_{R_1,R_2,..,R_n, |R_{\ell}| \ge
2 } e^{-\tilde U(R_1,..R_n)}
\prod_{ \ell=1}^n
 \r(R_{\ell}), 
} \Eq (C.1a)
$$
where
$$\tilde  U(R_1,\dots,R_n) = \sum_{1\le {\ell},s\le n} 
\tilde U(R_{\ell}, R_s),  \Eq (C.1b)  $$

$$ \tilde U(R_{\ell}, R_s)= \left \{  \eqalign
 { & 0\qquad  \hbox {if}  \qquad R_{\ell} \cap
R_s=
\emptyset
\cr &
\infty   \qquad
\hbox {if} \qquad  R_{\ell} \cap R_s \neq \emptyset} \right.   \Eq (C.1c)  $$
and
 $
 \r(R_{\ell})$ is given in \eqv (C.2).
Since  $|\AA (I)| < \infty $ the number  of terms contributing 
to \eqv (C.1a) is finite.
We have that  the $\log $ of the right hand side of \eqv (C.1a) can be 
written as  a formal expansion
$$ \eqalign {\b  V(m^{\d^*}_I) &=  \log \left 
[ 1+  \sum_{n=1}^{\infty}\frac 1 {n!} 
\sum_{R_1,R_2,\dots ,R_n, |R_{\ell}| \ge
2 } e^{-\tilde U(R_1,\dots ,R_n)} \prod_{{\ell}=1}^n
 \r(R_{\ell})  \right] \cr &=
\sum_{n=1}^{\infty} \frac 1 {n!}\sum_{R_1,R_2,\dots ,R_n,
|R_{\ell}| \ge 2 } \Phi^T(R_1,R_2,\dots ,R_n)
\prod_{{\ell}=1}^n
 \r(R_{\ell}) ,   }\Eq (C.3P) $$
where $ \Phi^T(R_1,R_2,\dots ,R_n)$ are the Ursell
coefficients
$$\Phi^T(R_1,R_2,\dots ,R_n)=
\left \{  \eqalign { & \sum_{g\in \GG_{R_1,\dots, R_n}}
\prod_{({\ell},s) \in g, {\ell} \neq s}\left [ e^{-\tilde U
(R_{\ell}, R_s) }-1 \right ]
 \qquad  \hbox {if}  \qquad
n \ge 2 \cr & 1 \qquad \hbox {if}\qquad  n=1 } \right.   \Eq
(C.1d)  $$ Observe  that $ \Phi^T(R_1,R_2,\dots ,R_n)=0$
if $g\in \GG_{R_1,...R_n}$ is not connected.

\noindent   We must now prove that the formal series 
\eqv (C.3P) actually converges.
Fix  $x \in C_{\d^*} (I)$ and a  polymer $R$,
such that $ A(x) \in R$.  Recall that
$\Phi^T(R)=1$,
 when $n=1$.  Then, 
   \eqv (C.3P) can be written as  
$$   \b V(m^{\d^*}_I,h)  =
\sum_{ x\in C_{\d^*}(I)}\sum_{R, x\in R,  |R| \ge 2 } \r (R) \left
[1+ \sum_{n\ge 2} \frac 1 {n!} B_n(R) \right ],   \Eq (C.100)
$$ where
$$ B_n(R)=  \sum_{R_2,\dots ,R_n, |R_{\ell}| \ge 2 }
\prod_{{\ell}=2}^n
 \r(R_{\ell}) \Phi^T(R,R_2,\dots ,R_n). \Eq (C.11) $$
From the definition of $ \Phi^T(R,R_2,\dots ,R_n)$  we
see  that $ B_n(R)$ can be written as
$$ B_n(R)= \sum_{g \in \GG_{R,R_2,\dots, R_n} }  \sum_{   f \subset g  }
(-1)^{|f|}  \sum_{R_2,\dots ,R_n, |R_{\ell}| \ge 2 }
 \prod_{{\ell}=2}^n
 \r(R_{\ell}),   \Eq (C.12) $$
where $f \subset g$ means that every link of $f$ is also a link of $g$.
Recall that, from Rota inequality,  see  [\rcite{ROTA}],
$$\left | \sum_{f \subset g  } (-1)^{|f|} \right | \le N(g),  $$
where $N(g) $ is the number of connected tree graphs in $g$.
Setting $\TT_{R,R_2,\dots,R_n}  \equiv \TT_n$,  we have that
$$  \sum_{ g \in \GG_{R_1,\dots,R_n }}  =  \sum_{
\t  \in \TT_n } \sum_{g: \t \in g } \frac 1 {N(g)} $$
and then we can express
$$ B_n(R)= \sum_{
\t  \in \TT_n }  w(\t) \Eq (C.13) $$
where
$$  w(\t) =  \sum_{R_2,\dots ,R_n, |R_{\ell}| \ge 2, 
\t \in g(R,R_2,\dots ,R_n) }
 \prod_{{\ell}=2}^n
 \r(R_i)  \Eq (C.12a ) $$
For any fixed set $R'$ we have the bound
$$  \sum_{R, R \cap R'\neq \emptyset} \le |R'| \sup_{x \in R'} 
\sum_{R: x\in
R} $$
then
 $$  w(\t)  \le |R|^{d_1} \prod_{{\ell}=2}^n
\left [ \sup_{x \in  \CC_{\d^*}(I) }
\sum_{R_{\ell}: x \in  R_{\ell}} |R_{\ell}|^{d_i-1} \left |
\r(R_{\ell}) \right |   \right ],  \Eq (C.13P) $$ where
$d_{\ell}$ is the incidence number of the vertex ${\ell}$ in the
tree $\t$.
Using Caley formula  [\rcite{PLS}],  we get
   $$ \eqalign { &  B_n(R) = \sum_{
\t  \in \TT_n }  w(\t) \cr &
\le \sum_{d_1,...,d_n} |R|^{d_1}  \frac {(n-2)!} {\prod_{{\ell}=1}^n
(d_{\ell}-1)}
\prod_{{\ell}=2}^n \left [
\sup_{x
\in
\CC_{\d^*}(I) }
\sum_{R_i: x \in  R_{\ell}} | R_{\ell}|^{d_{\ell}-1}
 \left | \r(R_{\ell}) \right |    \right ]  \cr &
\le (n-1)! \left [ \sum_{d_1=1}^\infty \frac {|R|^{d_1}} {d_1!} \right ]
 \prod_{{\ell}=2}^n
\left [
\sup_{x\in C_{\d^*}(I)} \sum_{R_{\ell}: x \in R_{\ell}} 
\sum_{d_{\ell}=1}^\infty
|\r(R_{\ell})|  \frac {|R|^{d_{\ell} -1}} {(d_{\ell}-1)!}
\right ]  \cr &
\le (n-1)! \left (e^{|R|} -1 \right ) \prod_{{\ell}=2}^n
\left [
\sup_{x\in C_{\d^*}(I)} \sum_{R_{\ell}: x \in R_{\ell}} |\r(R_{\ell})|
 e^{ |R_{\ell}|}  \right ]
 \le  (n-1)!  e^{ |R|} S^{n-1} ,  }\Eq (C.16) $$
where
in the second inequality  we used  that $ n-1 \ge d_1 $ 
to obtain the factor $ \frac 1
{d_1!} $ and  in the last inequality we set
 $$  S =\sup_{h} \sup_{x\in \CC_{\d^*}(I)}
 \sum _{R:x \in R}e^{|R|}\rho(R).   \Eq (paris10)$$
Thus,  under the condition that $S<1$ we obtain 
$$  \eqalign { &
\sum_{R, |R| \ge 2, x \in R } |\r (R)| \left [1+
\sum_{n\ge 2} \frac 1 {n!} B_n(R) \right ]
\le
\sum_{R, |R| \ge 2, x \in R } |\r (R)| \left [1+  e^{ |R|} 
 \sum_{n\ge 2} \frac 1
n  S^{n-1}
\right ] \cr & \le
\sum_{R, |R| \ge 2, x \in R } |\r(R) | e^{|R|}  \left [1 +
\frac  {S }   {1-  S  } \right ]= \frac S {1-S}.} \Eq (C.22)
$$
  Therefore, recalling \eqv (C.100), we obtain  \eqv (H.2).
 The  important  remark   to prove \eqv (H.22P) is that  to
 obtain the Lipschitz norm we
make the difference of two absolute convergent series  having
the only difference in one site $i$.   We then obtain
$$
\eqalign{ &   V (m^{\d^*}_{\tilde I_{12}},h) -  V
(m^{\d^*}_{\tilde I_{12}},h^i )
  = \cr
& \frac 1 \b  \quad\quad
\sum_{n=1}^{\infty} \frac 1 {n!}\sum_{R_1,R_2,\dots,R_n,
 |R_{\ell}| \ge 2
}
  \Phi^T(R_1,R_2,\dots,R_n) \left [ \prod_{\ell=1}^n
 \r(R_\ell,h) - \prod_{\ell=1}^{n} \r(R_\ell,h^i)  \right ]
 \cr
&  \le  \frac 2 \b   \sum_{n=1}^{\infty} \frac 1 {n!}
\sum_{R_1,R_2,\dots,R_n, |R_l|
\ge 2
\exists l : i \in R_l}
\left |\Phi^T(R_1,R_2,\dots,R_n)\right|  \prod_{\ell=1}^n  \sup_{h}
\left |\r(R_\ell,h) \right|.
}\Eq (S.5P)  $$
  Following the same strategy used above we obtain  \eqv (H.22P).
Next we show that    $S $,  see \eqv (paris10), satisfies \eqv
(paris1).    Taking
into  account  \eqv (H.1) and setting
 $ \Phi(x,y) = \1_{\{\frac 12 -
\d^* \le \d^*|x-y| \le \frac 12 + \d^*\}}   \left (\frac {\b (
\d^*)^2} {\g} \right ) $ we obtain that  if $g$ is a  
connected graph   with support $R$, then:  
$$   \sup_{h} \E_{ m^{\d^*}_R}\left  [
\prod_{(x,y) \in g, x \neq y}\left [ e^{ \b U (\s_{A(x)},
\s_{A(y)}} -1 \right ]
\right]
       \le   \prod_{(x,y) \in g, x \neq y}
\Phi(x,y). \Eq (C.8)$$ In the last estimate we used    \eqv(S.9).
From  \eqv (paris10)    we have that
$$ \eqalign { & S \equiv   \sup_{h} \sup_{x\in C_{\d^*}(I)} 
\sum_{ R: x \in R} |\r(R)|  e^{ |R|}
\cr & \le \sup_{x\in C_{\d^*}(I)} \sum_{ R: x \in R} \sum_{g \in\GG_R}
  e^{|R|} \prod_{(z,y)
\in g, z \neq y} \Phi (z,y). }  \Eq (C.20P)
$$
 An essential fact to prove  \eqv (paris1) is that  $ \Phi
(z,y) \neq 0 $ only   when $\frac 12 - \d^* \le \d^*|z-y| \le
\frac 12 + \d^* $, i.e.,    the block $A(z)$ interacts only
with three  blocks, the $ A(y)$ block which is at   distance
$\frac 1 {2\d^*}$ from it and the two blocks,   to the left   and
to the right of $ A(y)$\footnote{$^{*}$ }{\eightrm This depends on
the particular choice of the  potential, $\1_{|x|\le \frac 12}$.
For general potential, always with support $\{ x: |x|\le \frac
12\}$ this will be not true. In that case $ \Phi (z,y) \neq 0 $
when $
 \d^*|z-y| \le \frac 12 $, therefore the block $A(z)$ will interact
with $\frac 1 {\d^*}$ blocks.  Nevertheless this will not cause
problems to  get \eqv (paris1).  Namely  in this case the  function
$\Phi$,   using  Taylor formula to estimate  the potential,  becomes
  $ \Phi(x,y) = \1_{\{ \d^*|x-y| \le
\frac 12 \}}   \left (\frac {\b ( \d^*)^2} {\g}c \d^* \right )$, where
$c$ is a positive constant depending on the potential. Performing the
sums in \eqv (P.P1) we should replace $3$ with $\frac 1 {\d^*}$.
The result will be  similar. The only difference is given by the
presence of the constant $c$.  }.
\hskip.2cm
 Therefore for any fixed
polymer $R$, $x \in R$,  $|R|=\ell$,  the number of  graphs
that  contribute  to the sum in \eqv (C.20P) is at
most  $ 3^{(\ell-1)}$. Namely,   $\ell-1$ is the number
of links connecting the $\ell$ vertices of the graph and $3$ is
the maximum number of links that a vertex can have with the
others,  since  $ \Phi (z,y) \neq 0 $ only   when $\frac 12 -
\d^* \le \d^*|z-y| \le \frac 12 + \d^* $. Since  $\Phi $ is
translational invariant  we can
assume $x=0$. Then from    \eqv (C.20P)  
we obtain that
 $$ \eqalign { &
S \le  \sum_{ R: 0 \in R} \sum_{g \in\GG_R}
  e^{|R|} \prod_{(z,y)
\in g, z\neq y} \Phi (z,y)  = \sum_{\ell\ge 2}\sum_{R,0 \in R,
|R|=\ell}\sum_{g \in \GG_R}e^{|R|} \prod_{(z,y)\in g, z\neq
y}\Phi(z,y) \cr & \le  \sum_{\ell\ge 2}\ell 3^{(\ell-1)} 
 e^{\ell} [ \frac  \b  \g  (\d^*)^2 ]^{\ell-1}
 < 3\left [ \frac {e^3}
{1- 3 e^2 \frac \b \g  (\d^*)^2 }  \right]  \frac \b \g 
 (\d^*)^2 \le 6 e^3 \frac \b \g  (\d^*)^2.   } \Eq (P.P1) $$
\qed
\vskip0.5cm
\centerline{\bf References}
\vskip.3truecm
\item{[\rtag{Ah}]} A. Aharony. {\it
Tricritical points in systems with random fields. }
 Phys. Rev. B {\bf 18}, 3318--3327 (1978).
\item{[\rtag{AW}]}M. Aizenman, and  J. Wehr. 
{\it Rounding of first or\-der pha\-se tran\-si\-tions 
in sys\-tems with quenched disorder.}
 Com. Math. Phys. {\bf 130}, 489--528 (1990).
\item{[\rtag{AP}]}J.M. G.  Amaro de Matos and J. F. Perez.
{\it Fluctuations in the Curie-Weiss version of the random field Ising model.}
J. Stat. Phys. {\bf 62}, 587--608 (1991).
\item{[\rtag{APZ}]}J M G Amaro de Matos, A E Patrick, and  V A Zagrebnov.
{\it Random infinite volume Gibbs states for the Curie-Weiss 
random field Ising model.}  J. Stat. Phys. {\bf 66},
139--164 (1992).
\item{[\rtag{B}]}A. Berretti. {\it Some properties of random Ising 
models.} J. Stat. Phys, {\bf 38}, 483-496 (1985).
 \item{[\rtag{BRZ}]} P. Bleher, J. Ruiz, and V. Zagrebnov.
{\it One-dimensional random-field Ising model: Gibbs states and structure
of ground states.}
J. Stat. Phys. {\bf 84}, 1077--1093 (1996).
\item{[\rtag{Bo}]} T. Bodineau. {\it Interface in a one-dimensional 
Ising spin
system.} Stoch. Proc. Appl. {\bf 61}, 1--23 (1996).
\item{[\rtag{BGP4}]} A. Bovier, V. Gayrard and P. Picco. 
{\it Distribution of profiles for the Kac-Hopfield model.}
 Comm. Math. Phys. {\bf 186} 323--379 (1997).
\item{[\rtag{BZ}]} A. Bovier, and M. Zahradnik. 
{\it The Law temperature phase of Kac-Ising models.}
J. Stat. Phys. {\bf 87}, 311--332 (1997).
\item{[\rtag{BK}]} J. Bricmont, and A. Kupiainen.
{\it Phase transition in the three-dimensional 
random field Ising model.}
 Com. Math. Phys.,{\bf 116}, 539--572 (1988).
\item{[\rtag{Br}]} D.C. Brydges. {\it A short course on cluster expansions.}
in Critical phenomena, random systems, gauge theories, ed
K. Osterwalder, R. Stora, Les Houches XLlll, North Holland (1986)
\item{ }
R. Kotecky, D. Preiss. 
{\it  Cluster expansion for abstract polymer models.}  
Comm. Math. Phys. {\bf 103} 491--498 (1996). \item{  }
B. Simon. {\it The statistical mechanics of lattice gases.}
 Vol I, Princeton University Press (1993)
\item{[\rtag{LC}]} L. Le Cam. 
{\it Asymptotic methods in statistical decision theory.} 
Springer-Verlag  (1986).
\item{[\rtag{COP1}]} M. Cassandro, E. Orlandi, and P.Picco.
{\it Typical configurations for one-dimensional random field Kac model. } 
Ann. Prob.  {\bf 27}, No 3,
1414-1467 (1999).
\item{[\rtag{COP2}]}   M.Cassandro, E.Orlandi,and  P. Picco. 
{\it  Uniqueness and global stability of the interface in a model 
of phase separation.}
Nonlinearity   {\bf 27} No. 3, 1414-1467 (2002).   
  \item{[\rtag{COP}]} M. Cassandro, E. Orlandi, and E. Presutti.
{\it Interfaces and typical Gibbs configurations for one-dimensional
Kac potentials.} Prob. Theor. Rel. Fields {\bf 96}, 57-96 (1993).
 \item{[\rtag{CP}]} M. Cassandro, and  E. Presutti.
{\it Phase transitions in Ising systems with long but finite range 
interactions.} Markov
Process, Rel. Fields,  {\bf 2}, 241-262  (1996).
 \item{[\rtag{CT}]}Y.S. Chow  and H. Teicher. 
{\it Probability Theory: Independence, Interchangeability, Martingales.}
Springer-Verlag, Berlin and New York, (1978).  
\item{[\rtag{FFS}]} D.S. Fisher, J. Fr{\"o}hlich, and T. Spencer.
{\it The Ising model in a random magnetic field.} 
 J. Stat. Phys.  {\bf 34}, 863--870 (1984).
\item{[\rtag{FP}]}S. Friedli, C.-E. Pfister, 
{\it Non-analiticity and the van der Waals limit.}, J. Stat. Phys. 
(to appear).
\item{[\rtag{I}]}J. Imbrie.
{\it The ground states of the three-dimensional 
random field Ising model.}
 Com. Math. Phys.  {\bf 98}, 145--176 (1985). 
\item{[\rtag{IM}]}Y. Imry, and  S.K. Ma. 
{\it Random-field instability of the 
ordered state of continuous symmetry.}
 Phys. Rev. Lett. {\bf 35}, 1399--1401 (1975).
\item{[\rtag{KUH}]} M. Kac, G. Uhlenbeck, and P.C. Hemmer. {\it  On the van
der Waals theory of vapour-liquid equilibrium. I. Discussion of a
one-dimensional model.} J.  Math. Phys. {\bf 4}, 216--228 (1963);
{\it II. Discussion of the distribution functions.} 
J. Math. Phys. {\bf 4}, 229--247 (1963);
{\it III. Discussion of the critical region.} 
J. Math. Phys. {\bf 5}, 60--74 (1964).
\item{[\rtag{Ku}]} C. K\"ulske.
{\it On the Gibbsian nature of the random field Kac model under 
block-averaging.}
J. Stat. Phys. no 5/6, 991-1012  (2001).
\item{[\rtag{LP}]} J. Lebowitz, and O. Penrose. {\it Rigorous treatment of
the Van der Waals Maxwell theory of the liquid-vapour
transition.} J. Math. Phys. {\bf 7}, 98--113 (1966).
\item{[\rtag{LMP}]} J. Lebowitz, A Mazel and  E. Presutti. {\it
 Liquid-vapor phase transition for systems with finite-range interactions
system.} J. Stat. Phys. {\bf 94}, no 5/6, 955-1025  (1999).
\item{[\rtag{LT}]} M. Ledoux, and M. Talagrand. {\it Probability in 
Banach Spaces.} Springer, Berlin-Heidel\-berg-New York, (1991).
 \item{[\rtag{Ne}]} J. Neveu. {\it Martingales a temps discrete.} 
Masson and C  (1972).
\item{[\rtag{PL}]} O. Penrose, and J.L. Lebowitz.
{\it Towards a rigorous molecular theory of metastability.} 
in Fluctuation Phenomena ( e.W. Montroll and J.L. Lebowitz
ed) North-Holland Physics Publishing (1987).
\item{[\rtag{RY}]} D. Revuz,  and M. Yor.
{\it Continuous Martingales and Brownian Motion.}
Springer Verlag (1991)
\item{[\rtag{Ro}]} H. Robbins. {\it A remark on Stirling's formula.}
Amer. Math. Montlhy {\bf 62}, 26--29 (1955).
\item{[\rtag{ROTA}]} G.C. Rota. {\it On the foundations of the
 combinatorial theory.I Theory
of Moebius functions, Z.}. Wahrscheinlichkeitstheorie Verw. 
Gebiete {\bf 2}, 340--368 (1964).
\item{[\rtag{SW}]} S.R. Salinas, and W.F. Wreszinski.
{\it On the mean field Ising model in a random external field.} 
J. Stat. Phys. {\bf 41}, 299--313 (1985).
\item{[\rtag{SA}]} E. Sparre Andersen.
{\it On the  fluctuations of sums of Random variables II.} 
Math.  Scand. {\bf 2}, 195--233 (1954).

\end